\documentclass[a4paper]{amsart}
\usepackage{mathbbol}

\usepackage{bm}
\usepackage[a4paper,hmargin=36mm,vmargin=48mm]{geometry}
\usepackage[svgnames,table]{xcolor}
\usepackage{booktabs}
\arrayrulecolor{MidnightBlue}
\newcolumntype{C}{>{\columncolor{yellow!10}}c}
\usepackage{xparse,mathtools,etoolbox}
\usepackage{mathrsfs}
\usepackage{amsfonts,amssymb,mathrsfs}

\DeclareFontEncoding{LS2}{}{\noaccents@}
\DeclareFontSubstitution{LS2}{stix}{m}{n}

\DeclareSymbolFont{arrows3} {LS2}{stixtt} {m} {n}
\DeclareMathSymbol{\longrightsquigarrow}{\mathrel}{arrows3}{"0C}

\DeclareFontEncoding{LS1}{}{}
\DeclareFontSubstitution{LS1}{stix}{m}{n}
\DeclareSymbolFont{symbols2}{LS1}{stixfrak} {m} {n}

\DeclareFontFamily{U}  {MnSymbolD}{}
\DeclareSymbolFont{MnSyD}         {U}  {MnSymbolD}{m}{n}
\DeclareFontShape{U}{MnSymbolD}{m}{n}{
    <-6>  MnSymbolD5
   <6-7>  MnSymbolD6
   <7-8>  MnSymbolD7
   <8-9>  MnSymbolD8
   <9-10> MnSymbolD9
  <10-12> MnSymbolD10
  <12->   MnSymbolD12}{}
\DeclareMathSymbol{\sqfrowneqsmile}{\mathrel}{MnSyD}{57}

\DeclareFontFamily{U}  {MnSymbolC}{}
\DeclareSymbolFont{MnSyC}         {U}  {MnSymbolC}{m}{n}
\DeclareFontShape{U}{MnSymbolC}{m}{n}{
    <-6>  MnSymbolC5
   <6-7>  MnSymbolC6
   <7-8>  MnSymbolC7
   <8-9>  MnSymbolC8
   <9-10> MnSymbolC9
  <10-12> MnSymbolC10
  <12->   MnSymbolC12}{}
\DeclareMathSymbol{\triangledown}{\mathbin}{MnSyC}{81}
\DeclareMathSymbol{\triangleup}{\mathbin}{MnSyC}{83}

\usepackage{enumitem}
\setlist[enumerate]{label=\upshape\alph*), ref=\alph*}
\let\realItem\item
\NewDocumentCommand\centeredItem{so}{
   \IfNoValueTF{#2}{\realItem}{\realItem[#2]}%
   \IfBooleanF{#1}{\hfil}%
}
\newlist{Enumerate}{enumerate}{1}
\setlist[Enumerate]{
  label=\textup{\alph*)},
  ref=\theequation\alph*,
}
\newlist{relations}{enumerate}{1}
\setlist[relations]{
  label=$(\text{KLR}_{\arabic*})$,
  before=\let\item\centeredItem
}

\usepackage{tikz}
\usetikzlibrary{arrows.meta,matrix,decorations.markings,decorations.pathreplacing,calc}

\newcommand\doverline[1]{%
\tikz[baseline=(nodeAnchor.base)]{
    \node[inner sep=0] (nodeAnchor) {$#1$};
    \draw[line width=0.1ex,line cap=round]
        ($(nodeAnchor.north west)+(0.0em,0.2ex)$)
            --
        ($(nodeAnchor.north east)+(0.0em,0.2ex)$)
        ($(nodeAnchor.north west)+(0.0em,0.5ex)$)
            --
        ($(nodeAnchor.north east)+(0.0em,0.5ex)$)
    ;
}}

\usepackage{aliascnt}
\usepackage{cite}

\title{Content systems and deformations of cyclotomic {KLR} algebras of type $A$ and $C$}


\def\Email#1{\email{\href{mailto:#1}{#1}}}
\author{Anton Evseev}
\address{School of Mathematics, University of Birmingham,
         Edgbaston, Birmingham B15 2TT, UK}
\author{Andrew Mathas}
\address{School of Mathematics and Statistics,
         University of Sydney \\ NSW 2006 \\ Australia}
\Email{andrew.mathas@sydney.edu.au}

\usepackage{hyperref}
\hypersetup{%
  anchorcolor = red,
  citecolor = blue,
  colorlinks = true,
  linkcolor =blue,
  naturalnames,
  unicode=true,
  urlcolor = blue,
  raiselinks,
  hypertexnames
}

\RenewDocumentCommand\And{s}{%
  \IfBooleanTF{#1}{\quad\text{and}\quad}{\qquad\text{and}\qquad}%
}

\def\map#1#2{\colon #1\!\longrightarrow\!#2}
\newcommand\bijection[1][\sim]{\overset{#1}{\longrightarrow}}

\DeclarePairedDelimiterX\cart[2]{\langle}{\rangle}{#1,#2}
\let\kill\relax
\DeclarePairedDelimiterX\kill[2]{(}{)}{#1|#2}

\newcommand\N{\mathbb{N}}
\renewcommand\k{\mathbb{k}}
\newcommand\K{\mathbb{K}}
\newcommand\kx{\k[\ux]}
\newcommand\Kx{\K[\ux]}
\newcommand\Kxx{\K[\ux^\pm]}

\newcommand\CC{\mathbb{C}}
\newcommand\Q{\mathbb{Q}}
\newcommand\Sym{\mathfrak S}
\newcommand\Z{\mathbb{Z}}
\newcommand\Jell{{J_\ell}}

\newcommand\A{\mathcal{A}}
\renewcommand\AA{\mathbb{A}}

\newcommand\Ncal{\mathcal{N}}
\newcommand\Pcal{\mathcal{P}}

\newcommand\TT{\mathtt{t}}

\newcommand\B[1][0]{\mathcal{B}_{#1}}
\renewcommand\L[1][0]{\mathcal{L}_{#1}}

\newcommand\mull{\mathsf{m}}
\newcommand\Nodes[1][n]{\Ncal^\ell_{#1}}

\newcommand\charge{{\bm{\rho}}}
\newcommand\blam{{\bm{\lambda}}}
\newcommand\bmu{{\bm{\mu}}}
\newcommand\bnu{{\bm{\nu}}}
\newcommand\bsig{{\bm{\sigma}}}
\newcommand\btau{{\bm{\tau}}}
\newcommand\balp{{\bm{\alpha}}}
\newcommand\bbet{{\bm{\beta}}}

\let\uaccent=\u

\def\s{\mathsf{s}}
\def\t{\mathsf{t}}
\def\u{\mathsf{u}}
\def\v{\mathsf{v}}
\def\w{\mathsf{w}}
\def\x{\mathsf{x}}
\def\y{\mathsf{y}}
\def\z{\mathsf{z}}

\newcommand\ux{{\underline{x}}}

\newcommand\bc{\mathsf{c}}
\newcommand\br{\mathsf{r}}

\newcommand\floor[2]{\lfloor\tfrac{#1}{#2}\rfloor}

\let\Dom\triangledown 
\newcommand\doM{\triangleup}
\newcommand\Domeq{\underline{\Dom}}
\usepackage{stackengine}
\newcommand\Domneq{\mathrel{\stackinset{c}{}{c}{}{$/$}{$\Dom$}}}

\newcommand\Domin{{\Dom} \in\set{\Ldom,\Gdom}}
\newcommand\Doming{\set{\Dom,\doM}=\set{\Ldom,\Gdom}}

\newcommand\GDom\blacktriangleright
\newcommand\LDom\blacktriangleleft
\newcommand\GeDom{\mathrel{\underline{\GDom}}}

\newcommand\Gedom{\trianglerighteq}
\newcommand\Gnedom{\ntrianglerighteq}

\newcommand\Gdom{\triangleright}
\newcommand\Ledom{\trianglelefteq}
\newcommand\Ldom{\triangleleft}

\newcommand\DefPol[1][\blam]{\gamma^{{\Ldom}\kern-0.9pt{\Gdom}}_{#1}}

\newcommand\Llex{\mathop{<_{\text{lex}}}}
\newcommand\Glex{\mathop{>_{\text{lex}}}}

\DeclarePairedDelimiterX\CPair[2]{\langle}{\rangle}{#1,#2}
\DeclarePairedDelimiterX\GPair[2]{(}{)^{\Gdom}}{#1,#2}
\DeclarePairedDelimiterX\LPair[2]{(}{)^{\Ldom}}{#1,#2}
\DeclarePairedDelimiterX\DPair[2]{(}{)^{\Dom}}{#1,#2}

\newcommand\Gpsis[1][\s]{\psi^{\Gdom}_{#1}}
\newcommand\Gpsist[1][\s\t]{\psi^{\Gdom}_{#1}}
\newcommand\Gfst[1][\s\t]{f^{\Gdom}_{#1}}
\newcommand\Gfs[1][\s]{f^{\Gdom}_{#1}}
\newcommand\Gylam[1][\blam]{y^{\Gdom}_{#1}}
\newcommand\Gilam[1][\blam]{\bi^{\Gdom}_{#1}}
\newcommand\Gtlam[1][\blam]{\t^{\Gdom}_{#1}}
\newcommand\Ggt[1][\t]{\gamma^{\Gdom}_{#1}}
\newcommand\Gwlam[1][\blam]{d^{\Gdom}_{#1}}
\newcommand\Gdt[1][\t]{d^{\Gdom}_{#1}}
\newcommand\Gbeta{\beta^{\Gdom}}
\newcommand\GAdd{\mathop{\mathrm{Add}}\nolimits^{\Gdom}}
\newcommand\GRem{\mathop{\mathrm{Rem}}\nolimits^{\Gdom}}
\NewDocumentCommand\GSlam{ sO{\blam}}{S^{\Gdom\IfBooleanT{#1}\sgn}_{#2}}
\newcommand\GGnu[1][\bnu]{\mathbb{Y}^{\Gdom}_{#1}}
\newcommand\GHnu[1][\bnu]{\mathbb{D}^{\Gdom}_{#1}}
\newcommand\Gynu[1][\bnu]{\mathsf{y}^{\Gdom}_{#1}}
\newcommand\GXnu[1][\bnu]{\mathsf{X}^{\Gdom}_{#1}}
\NewDocumentCommand\GEnu{O{\bnu}D(){\K}}{E^{\Gdom}_{#1}(#2)}
\newcommand\Gdeg{\mathop{\deg^{\Gdom}}}

\NewDocumentCommand\GeRlam{sD<>{n} O{\blam}d()}%
{\bigl(\IfBooleanTF{#1}{\sRx}{\Rx[#2]\IfNoValueF{#4}{(#4)}}\bigr)^{\Gedom#3}}
\newcommand\Gow[1][n]{{{\bm\omega}_{#1}^{\Gedom}}}
\NewDocumentCommand\Gdilam{ s O{A} D(){\blam} }{d^{\Gdom\IfBooleanT{#1}{\sgn}}_{#2}(#3)}
\NewDocumentCommand\Gglamnu{O{\blam\bnu}D(){\K}}{\mathbb{a}^{#2\Gdom}_{#1}(q)}
\NewDocumentCommand\Ghlamnu{O{\blam\bnu}D(){\K}}{\mathbb{b}^{#2\Gdom}_{#1}(q)}
\NewDocumentCommand\GGlamnu{O{\blam\bnu}}{\mathbb{d}^{\Gdom}_{#1}(q)}
\NewDocumentCommand\GHlamnu{O{\bnu\blam}}{\mathbb{e}^{\Gdom}_{#1}(-q)}
\NewDocumentCommand\Gdlamnu{O{\blam\bnu}D(){\K}}{\mathsf{d}^{#2\Gdom}_{#1}(q)}
\NewDocumentCommand\Gelamnu{O{\blam\bnu}D(){\K}}{\mathsf{e}^{#2\Gdom}_{#1}(-q)}
\NewDocumentCommand\Gclamnu{O{\blam\bnu}D(){\K}}{\mathsf{c}^{#2\Gdom}_{#1}(q)}
\NewDocumentCommand\GCar{O{n}D(){\K}}{\mathsf{C}^{#2\Gdom}_{#1}}
\NewDocumentCommand\GDec{O{n}D(){\K}}{\mathsf{D}^{#2\Gdom}_{#1}}
\newcommand\Gdec{\mathsf{d}^{\Gdom}}

\NewDocumentCommand\GAdj{O{n}D(){\K}}{\mathsf{A}^{#2\Gdom}_{#1}}
\NewDocumentCommand\Gfock{sO{\blam}}
   {\mathsf{s}^{\Gdom\IfBooleanT{#1}{\sgn}}_{#2}}
\NewDocumentCommand\GFock{sO{\A}}{%
    \mathscr{F}^{\Lambda\IfBooleanT{#1}{{}^\sgn}\Gdom}_{#2}}

\newcommand\Lpsis[1][\s]{\psi^{\Ldom}_{#1}}
\newcommand\Lpsist[1][\s\t]{\psi^{\Ldom}_{#1}}
\newcommand\Lfst[1][\s\t]{f^{\Ldom}_{#1}}
\newcommand\Lfs[1][\s]{f^{\Ldom}_{#1}}
\newcommand\Lylam[1][\blam]{y^{\Ldom}_{#1}}
\newcommand\Lilam[1][\blam]{\bi^{\Ldom}_{#1}}
\newcommand\Ltlam[1][\blam]{\t^{\Ldom}_{#1}}
\newcommand\Lgt[1][\t]{\gamma^{\Ldom}_{#1}}
\newcommand\Lwlam[1][\blam]{d^{\Ldom}_{#1}}
\newcommand\Ldt[1][\t]{d^{\Ldom}_{#1}}
\newcommand\Lbeta{\beta^{\Ldom}}
\newcommand\LAdd{\mathop{\mathrm{Add}}\nolimits^{\Ldom}}
\newcommand\LRem{\mathop{\mathrm{Rem}}\nolimits^{\Ldom}}
\NewDocumentCommand\LSlam{ sO{\blam}}{S^{\Ldom\IfBooleanT{#1}\sgn}_{#2}}
\newcommand\LGmu[1][\bmu]{\mathbb{Y}^{\Ldom}_{#1}}
\newcommand\LHmu[1][\bmu]{\mathbb{D}^{\Ldom}_{#1}}
\newcommand\Lymu[1][\bmu]{\mathsf{y}^{\Ldom}_{#1}}
\newcommand\LXmu[1][\bmu]{\mathsf{X}^{\Ldom}_{#1}}
\NewDocumentCommand\LEmu{O{\bmu}D(){\K}}{E^{\Ldom}_{#1}(#2)}
\newcommand\Ldeg{\mathop{\deg^{\Ldom}}}
\newcommand\LRlam[1][\blam]{\bigl(\Rx(\kx)\bigr)^{\Ldom#1}}
\NewDocumentCommand\LeRlam{sD<>{n}O{\blam}d()}%
{\bigl(\IfBooleanTF{#1}{\sRx}{\Rx[#2]\IfNoValueF{#4}{(#4)}}\bigr)^{\Ledom#3}}
\newcommand\Low[1][n]{{{\bm\omega}_{#1}^{\Ledom}}}
\NewDocumentCommand\Ldilam{ sO{A}D(){\blam} }{d^{\Ldom\IfBooleanT{#1}{\sgn}}_{#2}(#3)}
\NewDocumentCommand\Lglammu{O{\blam\bmu}D(){\K}}{\mathbb{a}^{#2\Ldom}_{#1}(q)}
\NewDocumentCommand\Lhlammu{O{\blam\bmu}D(){\K}}{\mathbb{b}^{#2\Ldom}_{#1}(q)}
\NewDocumentCommand\LGlammu{O{\blam\bmu}}{\mathbb{d}^{\Ldom}_{#1}(q)}
\NewDocumentCommand\LHlammu{O{\bmu\blam}}{\mathbb{e}^{\Ldom}_{#1}(-q)}
\NewDocumentCommand\Ldlammu{O{\blam\bmu}D(){\K}}{\mathsf{d}^{#2\Ldom}_{#1}(q)}
\NewDocumentCommand\Lelammu{O{\blam\bmu}D(){\K}}{\mathsf{e}^{#2\Ldom}_{#1}(-q)}
\NewDocumentCommand\Lclammu{O{\blam\bmu}D(){\K}}{\mathsf{c}^{#2\Ldom}_{#1}(q)}
\NewDocumentCommand\LCar{O{n}D(){\K}}{\mathsf{C}^{#2\Ldom}_{#1}}
\NewDocumentCommand\LDec{O{n}D(){\K}}{\mathsf{D}^{#2\Ldom}_{#1}}
\newcommand\Ldec{\mathsf{d}^{\Ldom}}

\NewDocumentCommand\LAdj{O{n}D(){\K}}{\mathsf{A}^{#2\Ldom}_{#1}}
\newcommand\Lfock[1][\blam]{\mathsf{s}^{\Ldom}_{#1}}
\NewDocumentCommand\LFock{sO{\A}}{%
    \mathscr{F}^{\Lambda\IfBooleanT{#1}{{}^\sgn}\Ldom}_{#2}}

\NewDocumentCommand\GDnu{ sO{\bnu} } {D^{\Gdom\IfBooleanT{#1}\sgn}_{#2}}
\NewDocumentCommand\LDmu{ sO{\bmu} } {D^{\Ldom\IfBooleanT{#1}\sgn}_{#2}}
\NewDocumentCommand\DDmu{ sO{\bmu} } {D^{\Dom\IfBooleanT{#1}\sgn}_{#2}}

\newcommand\Dpsis[1][\s]{\psi^{\Dom}_{#1}}
\newcommand\Dpsist[1][\s\t]{\psi^{\Dom}_{#1}}
\newcommand\Dfst[1][\s\t]{f^{\Dom}_{#1}}
\newcommand\Dfs[1][\s]{f^{\Dom}_{#1}}
\newcommand\Dylam[1][\blam]{y^{\Dom}_{#1}}
\newcommand\Dilam[1][\blam]{\bi^{\Dom}_\blam}
\newcommand\Dtlam[1][\blam]{\t^{\Dom}_{#1}}
\newcommand\Dgt[1][\t]{\gamma^{\Dom}_{#1}}

\newcommand\Ddt[1][\t]{d^{\Dom}_{#1}}
\newcommand\Dbeta{\beta^{\Dom}}
\newcommand\DAdd{\mathop{\mathrm{Add}}\nolimits^{\Dom}}

\NewDocumentCommand\DSlam{ sO{\blam}}{S^{\Dom\IfBooleanT{#1}\sgn}_{#2}}

\newcommand\DGmu[1][\bmu]{\mathbb{Y}^{\Dom}_{#1}}
\newcommand\DHmu[1][\bmu]{\mathbb{D}^{\Dom}_{#1}}

\newcommand\DXmu[1][\bmu]{\mathsf{X}^{\Dom}_{#1}}
\NewDocumentCommand\DEmu{O{\bmu}D(){\K}}{E^{\Dom}_{#1}(#2)}
\newcommand\Ddeg{\mathop{\deg^{\Dom}}}
\newcommand\DRlam[1][\blam]{(\Rx)^{\Dom#1}}
\NewDocumentCommand\DeRlam{sD<>{n}O{\blam}d()}%
{\bigl(\IfBooleanTF{#1}{\sRx}{\Rx[#2]\IfNoValueF{#4}{(#4)}}\bigr)^{\Domeq#3}}
\newcommand\Dow[1][n]{{{\bm\omega}_{#1}^{\Domeq}}}
\NewDocumentCommand\Ddilam{ sO{A}D(){\blam} }{d^{\Dom\IfBooleanT{#1}{\sgn}}_{#2}(#3)}
\NewDocumentCommand\Dodilam{ O{A} D(){\blam} }{d^{\doM}_{#1}(#2)}
\NewDocumentCommand\Dclammu{O{\blam\bmu}D(){\K}}{\mathsf{c}^{#2\Dom}_{#1}(q)}
\NewDocumentCommand\Ddlammu{O{\blam\bmu}D(){\K}}{\mathsf{d}^{#2\Dom}_{#1}(q)}
\NewDocumentCommand\Dglammu{O{\blam\bmu}D(){\K}}{\mathbb{a}^{#2\Dom}_{#1}(q)}
\NewDocumentCommand\Dhlammu{O{\blam\bmu}D(){\K}}{\mathbb{b}^{#2\Dom}_{#1}(q)}
\NewDocumentCommand\DGlammu{O{\blam\bmu}}{\mathbb{d}^{\Dom}_{#1}(q)}
\NewDocumentCommand\DHlammu{O{\bmu\blam}}{\mathbb{e}^{\Dom}_{#1}(-q)}
\NewDocumentCommand\Delammu{O{\blam\bmu}D(){\K}}{\mathsf{e}^{#2\Dom}_{#1}(-q)}
\NewDocumentCommand\DCar{O{n}D(){\K}}{\mathsf{C}^{#2\Dom}_{#1}}
\NewDocumentCommand\DDec{O{n}D(){\K}}{\mathsf{D}^{#2\Dom}_{#1}}
\newcommand\Ddec{\mathsf{d}^{\Dom}}
\newcommand\Dcar{\mathsf{c}^{\Dom}}
\NewDocumentCommand\DAdj{O{n}D(){\K}}{\mathsf{A}^{#2\Dom}_{#1}}
\newcommand\Dfock[1][\blam]{\mathsf{s}^{\Dom}_{#1}}
\NewDocumentCommand\DFock{sO{\A}}{%
    \mathscr{F}^{\Lambda\IfBooleanT{#1}{{}^\sgn}\Dom}_{#2}}

\NewDocumentCommand\DCrystal{t! s}
{\mathscr{B}^{\Dom\IfBooleanT{#1}\sgn}_{\IfBooleanTF{#2}{\infty}{0}}%
   (\Lambda\IfBooleanT{#1}{{}^\sgn})}
\NewDocumentCommand\LCrystal{t! s}
{\mathscr{B}^{\Ldom\IfBooleanT{#1}\sgn}_{\IfBooleanTF{#2}{\infty}{0}}%
   (\Lambda\IfBooleanT{#1}{{}^\sgn})}
\NewDocumentCommand\GCrystal{t! s}
{\mathscr{B}^{\Gdom\IfBooleanT{#1}\sgn}_{\IfBooleanTF{#2}{\infty}{0}}%
  (\Lambda\IfBooleanT{#1}{{}^\sgn})}

\NewDocumentCommand\DCry{s}{\mathsf{B}^\Dom(\Lambda)}
\NewDocumentCommand\LCry{s}{\mathsf{B}^\Ldom(\Lambda)}
\NewDocumentCommand\GCry{s}{\mathsf{B}^\Gdom(\Lambda)}

\NewDocumentCommand\Parts{sO{n}}{%
  \mathcal{P}^{\ell}_{\IfBooleanTF{#1}{\bullet}{#2}}}
\NewDocumentCommand\DKlesh{sO{n}}{%
  \mathcal{K}^{\Dom}_{\IfBooleanTF{#1}{\bullet}{#2}}}
\NewDocumentCommand\GKlesh{sO{n}}{%
  \mathcal{K}^{\Gdom}_{\IfBooleanTF{#1}{\bullet}{#2}}}
\NewDocumentCommand\LKlesh{sO{n}}{%
  \mathcal{K}^{\Ldom}_{\IfBooleanTF{#1}{\bullet}{#2}}}

\NewDocumentCommand\DYmu{ O{\bmu} }{Y^{\Dom}_{#1}}
\NewDocumentCommand\GYnu{ O{\bnu} }{Y^{\Gdom}_{#1}}
\NewDocumentCommand\LYmu{ O{\bmu} }{Y^{\Ldom}_{#1}}

\newcommand\Dzlam[1][\blam]{z^{\Dom}_{#1}}
\newcommand\Dzlams[1][\s]{z^{\Dom}_{#1}}
\newcommand\Gzlam[1][\blam]{z^{\Gdom}_{#1}}
\newcommand\Gzlams[1][\s]{z^{\Gdom}_{#1}}
\newcommand\LZlamu[1][\blam]{Z^{\Ldom}_{#1\uparrow}}
\newcommand\Lzlam[1][\blam]{z^{\Ldom}_{#1}}
\newcommand\Lzlams[1][\s]{z^{\Ldom}_{#1}}
\newcommand\Lzlamus[1][\s]{z^{\Ldom}_{#1\uparrow}}

\let\eps\varepsilon
\newcommand\tei[1][i]{\tilde{e}_{#1}}
\newcommand\tfi[1][i]{\tilde{f}_{#1}}

\newcommand\Deps{\varepsilon^\Dom}
\newcommand\Geps{\varepsilon^\Gdom}
\newcommand\Leps{\varepsilon^\Ldom}

\newcommand\Dphi{\varphi^\Dom}
\newcommand\Gphi{\varphi^\Gdom}
\newcommand\Lphi{\varphi^\Ldom}

\NewDocumentCommand\DCan{ sO{\bmu} }{%
  \mathsf{G}^{\Dom}_{\IfBooleanTF{#1}{\infty}{0},#2}}
\NewDocumentCommand\GCan{ sO{\bnu} }{%
  \mathsf{G}^{\Gdom}_{\IfBooleanTF{#1}{\infty}{0},#2}}
\NewDocumentCommand\LCan{ sO{\bmu} }{%
  \mathsf{G}^{\Ldom}_{\IfBooleanTF{#1}{\infty}{0},#2}}

\newcommand\GLa[1][\A]{\mathscr{L}^{\Gdom}_{#1}(\Lambda)}
\newcommand\LLa[1][\A]{\mathscr{L}^{\Ldom}_{#1}(\Lambda)}
\newcommand\DLa[1][\A]{\mathscr{L}^{\Dom}_{#1}(\Lambda)}

\NewDocumentCommand\Ganumu{O{\bnu\bmu}D(){\K}}{\mathsf{a}^{#2\Gdom}_{#1}(q)}
\NewDocumentCommand\Lanumu{O{\bnu\bmu}D(){\K}}{\mathsf{a}^{#2\Ldom}_{#1}(q)}
\NewDocumentCommand\Danumu{O{\bnu\bmu}D(){\K}}{\mathsf{a}^{#2\Dom}_{#1}(q)}

\NewDocumentCommand\DNi{ O{i} }{\Ncal^\Dom_{#1}}
\NewDocumentCommand\GNi{ O{j} }{\Ncal^\Gdom_{#1}}
\NewDocumentCommand\LNi{ O{i} }{\Ncal^\Ldom_{#1}}

\newcommand\qint[1]{[#1]}

\NewDocumentCommand{\setargs}{>{\SplitArgument{1}{|}}m}{\setargsaux#1}
\NewDocumentCommand{\setargsaux}{mm}
{\IfNoValueTF{#2}{#1} {#1\,\delimsize|\,\mathopen{}#2}}

\DeclarePairedDelimiterX{\set}[1]{\{}{\}}{\setargs{#1}}
\DeclarePairedDelimiterX{\gen}[1]{\langle}{\rangle}{\setargs{#1}}

\def\pmod#1{\text{ }(\textrm{mod } #1)\,}

\def\NewTheorem#1{%
  \newaliascnt{#1}{equation}%
  \newtheorem{#1}[#1]{#1}%
  \aliascntresetthe{#1}%
  \expandafter\def\csname #1autorefname\endcsname{#1}%
}
\def\equationautorefname~#1\null{(#1)\null}
\def\relationsiautorefname~#1\null{#1\null}
\def\itemautorefname~#1\null{(#1)\null}

\numberwithin{equation}{subsection}

\newcounter{maintheorem}

\newtheorem{MainTheorem}[maintheorem]{Theorem}

\swapnumbers

\NewTheorem{Proposition}
\NewTheorem{Theorem}
\NewTheorem{Corollary}
\NewTheorem{Lemma}
\NewTheorem{Definition}
\NewTheorem{Notation}
\theoremstyle{definition}
\NewTheorem{Example}
\AtEndEnvironment{Example}{\null\hfill$\Diamond$}%
\NewTheorem{Examples}
\AtEndEnvironment{Examples}{\null\hfill$\Diamond$}%
\theoremstyle{remark}
\NewTheorem{Remark}

\newcounter{case}[equation]

\newcommand\Case[1]{\refstepcounter{case}%
\medskip\noindent\textbf{Case \arabic{case}.} #1:\space}

\newcommand\ei[1][\bi]{\mathsf{1}_{#1}}

\newcommand\bi{\mathbf{i}}
\newcommand\bj{\mathbf{j}}

\newcommand{\DeclareMyOperator}[1]{%
  \expandafter\DeclareMathOperator\csname #1\endcsname{#1}
}
\forcsvlist{\DeclareMyOperator}{%
  Add,
  ch,
  cont,
  diag,
  End,
  END,
  head,
  Hom,
  HOM,
  im,
  inv,
  Irr,
  Deg,
  Ind,
  Rad,
  Rem,
  rad,
  Res,
  res,
  soc,
  supp,
  Shape,
  Std,
  Top,
}

\DeclareMathOperator\Rep{Rep}
\DeclareMathOperator\Proj{Proj}

\NewDocumentCommand\RepN{ sD<>{\K}}{%
  \bigl[\Rep_{#2}\Rx[\bullet](\K[x])\bigr]\IfBooleanT{#1}{{}_{\Z[q]}}}
\NewDocumentCommand\ProjN{ D<>{\K} o}{%
  \bigl[\Proj_{#1}\Rx[\bullet](\K[x])\bigr]\IfNoValueF{#2}{{}_{#2}}}

\newcommand\RepK[1][\K]{[\Rep_{#1}\Rx(#1[x])]}
\newcommand\ProjK[1][\K]{[\Proj_{#1}\Rx(#1[x])]}

\newcommand\gdim{\mathop{\rm dim}\nolimits_q}
\newcommand\defect{\mathop{\rm def}\nolimits}
\newcommand\StI[1][m]{I^{#1}_{\text{Std}}}

\DeclareMathOperator\noarrow{\:\rlap{\hspace*{0.25em}/}\text{---}\:}
\DeclareMathOperator\arrow{\text{---}}

\newcommand\Rxaff[1][n]{\mathsf{R}_{#1}}
\newcommand\Rnaff[1][n]{\mathscr{R}_{#1}}
\newcommand\Raff{\mathscr{R}}
\newcommand{\Saff}{\mathcal{S}}
\newcommand\Snl{\Saff^\ell_n}

\newcommand\Qx{Q^{\ux}}
\NewDocumentCommand\Rx{ D<>{\Lambda} O{n} }{\mathsf{R}^{#1}_{#2}}
\NewDocumentCommand\sRx{ D<>{\Lambda} O{\alpha} }{{}^{\sgn}\mathsf{R}^{#1}_{#2}}
\NewDocumentCommand\sRn{ D<>{\Lambda} O{\alpha} }{{}^{\sgn}\mathscr{R}^{#1}_{#2}}
\NewDocumentCommand\Hn{  D<>{\Lambda} O{n} }{\mathscr{H}^{#1}_{#2}}
\NewDocumentCommand\Rn{  D<>{\Lambda} O{n} }{\mathscr{R}^{#1}_{#2}}
\NewDocumentCommand\Rnaf{ D<>{n} }{\Rn[]<#1>}

\newcommand\h{\mathfrak{h}}
\newcommand\cij[1][ij]{\mathsf{c}_{#1}}
\newcommand\di[1][i]{\mathsf{d}_{#1}}
\newcommand\Uq{U_q(\mathfrak{g}_\Gamma)}
\newcommand\UA[1][\A]{U_{#1}(\mathfrak{g}_\Gamma)}
\newcommand\adq{\operatorname{ad}_{q^i}}

\NewDocumentCommand\Dgood{sO{\bi}}
    {\stackrel{#2\Dom\IfBooleanT{#1}{\sgn}}\longrightsquigarrow}
\NewDocumentCommand\Ggood{sO{\bj}}
    {\stackrel{#2\Gdom\IfBooleanT{#1}{\sgn}}\longrightsquigarrow}
\NewDocumentCommand\Lgood{sO{\bi}}
    {\stackrel{#2\Ldom\IfBooleanT{#1}{\sgn}}\longrightsquigarrow}

\newcommand\iarrow[1][i]{\stackrel{#1}\longrightarrow}

\newcommand\zero{\underline{\mathbf{0}}_\ell}
\DeclareMathOperator\wt{wt}
\newcommand\iInd[1][i]{F_{#1}^\Lambda}
\newcommand\iRes[1][i]{E_{#1}^\Lambda}

\let\sgn\varepsilon

\newcommand\Aone[1][e-1]{A^{(1)}_{#1}}
\newcommand\Cone[1][e-1]{C^{(1)}_{#1}}

\newcommand{\Qbx}[1][I]{\mathbf{Q}^{\ux}_{#1}}
\newcommand{\Qbxs}[1][I]{\mathbf{Q}^{\ux,\sgn}_{#1}}
\newcommand{\Qb}[1][I]{\mathbf{Q}_{#1}}
\newcommand{\Wbx}[1][I]{\mathbf{W}^{\ux}_{#1}}
\newcommand{\Wbxs}[1][I]{\mathbf{W}^{\ux,\sgn}_{#1}}
\newcommand{\Wb}[1][I]{\mathbf{W}_{#1}}

\newcommand{\CycIdeal}[1][\alpha]{\mathscr{W}^\Lambda_{#1}(\Wb)}
\newcommand{\Jstd}{J^n_{\text{Std}}}
\renewcommand{\phi}{\varphi}
\def\iso{\stackrel{\sim}{\longrightarrow}}

\newcount\ncells
\tikzset{
  centered/.style = {
     baseline = {([yshift=#1]current bounding box.center)}
  },
  centered/.default={-0.5ex},
  ->-/.style={
    decoration={
      markings,
      mark=at position 0.6 with {\arrow[thin,scale=2]{>}}
    },
    postaction={decorate}
  },
  -<-/.style={
    decoration={
      markings,
      mark=at position 0.6 with {\arrow[thin,scale=2]{<}}
    },
    postaction={decorate}
  },
  circled/.style = {fill=Tan},
  domstyle/.style = {
    thick,
    scale=#1,
  },
  vertex/.style = {
    circle,
    ball color=MidnightBlue,
    font=\small,
    inner sep=0pt,
    minimum size=2mm
  },
  pics/domm/.style = {
    code = {
      \draw[domstyle=#1](0.8,0)--++(0.6,0.4)--++(0.6,-0.4)--++(-0.6,-0.4)--++(-0.6,0.4);
    }
  },
  pics/dom/.style = {
    code = {
      \draw[domstyle=#1](0.2,0)--(0.8,0)
        --++(0.6,0.4)--++(0.6,-0.4)--++(-0.6,-0.4)--++(-0.6,0.4);
    }
  },
  pics/domeq/.style = {
    code = {
      \pic at (0,0){dom=#1};
      \draw[scale=#1](0.8,-0.5)--++(1.2,0);
    }
  },
  pics/domneq/.style = {
    code = {
      \pic at (0,0){domeq=#1};
      \draw[scale=#1](0.9,-0.75)--++(1.0,0.5);
    }
  },
  pics/mod/.style = {
    code = {
      \draw[domstyle=#1](0,0)--++(-0.6,0)--++(-0.6,-0.4)--++(-0.6,0.4)--++(0.6,0.4)--++(0.6,-0.4);
    }
  },
  pics/Sdom/.style = {
    code = {
      \draw[domstyle=#1](0.2,0)--(0.8,0);
      \draw[scale=#1,fill=black](0.8,0)--++(0.6,0.4)--++(0.6,-0.4)--++(-0.6,-0.4)--++(-0.6,0.4);
    }
  },
  pics/Sdomeq/.style = {
    code = {
      \pic at (0,0){Sdom=#1};
      \draw[scale=#1](0.8,-0.5)--++(1.2,0);
    }
  },
  pics/Sdomneq/.style = {
    code = {
      \pic at (0,0){Sdomeq=#1};
      \draw[scale=#1](0.9,-0.75)--++(1.0,0.5);
    }
  },
}

\usepackage{atableaux}
\usepackage{adecomp}

\usepackage{tabularx}
\usepackage{anindex}
\anindexsetup{
    heading,
    lines,
}

\let\emph\textbf
\let\<=\langle
\let\>=\rangle

\begin{document}

\begin{abstract}
  This paper initiates a systematic study of the cyclotomic KLR algebras
  of affine types~$A$ and~$C$. We start by introducing a graded
  deformation of these algebras and then constructing all of the
  irreducible representations of the deformed cyclotomic KLR algebras
  using \textit{content systems}, which are recursively defined using
  Rouquier's $Q$-polynomials. This leads to a generalisation of the
  Young's seminormal forms for the symmetric groups in the KLR setting.
  Quite amazingly, the same theory captures the
  representation theory of the cyclotomic KLR algebras of affine
  types~$A$ and~$C$, with the main difference being that the definition
  of the residue sequence of a tableau depends on the Cartan type. We use
  our semisimple deformations to construct two ``dual'' cellular bases
  for the non-semisimple KLR algebras of affine types~$A$ and~$C$. As
  applications we recover many of the main features from the
  representation theory in type~$A$, simultaneously proving them for the
  cyclotomic KLR algebras of types~$A$ and~$C$. These results are
  completely new in type~$C$ and we, usually, give more direct proofs in
  type~$A$. In particular, we show that these algebras categorify the
  irreducible integrable highest weight modules of the corresponding
  Kac-Moody algebras, we construct and classify their simple modules, we
  investigate links with canonical bases and we generalise Kleshchev's
  modular branching rules to these algebras.
\end{abstract}

\vspace*{40mm}

\maketitle

\begin{center}\slshape
  We record with deep sadness the passing of Anton Evseev on February 21, 2017.
\end{center}


\section{Introduction}

  The \textit{KLR algebras} are a remarkable family of graded algebras
  that were independently introduced by
  Khovanov-Lauda~\cite{KhovLaud:diagI} and
  Rouquier~\cite{Rouquier:QuiverHecke2Lie,Rouq:2KM}. These algebras are
  now central to many of the recent developments in representation
  theory, not least because these algebras categorify the positive part
  of quantised Kac-Moody
  algebras~\cite{VaragnoloVasserot:CatAffineKLR}.

  \newpage
  The \textit{cyclotomic KLR algebras} are natural finite dimensional
  quotients of the KLR algebras that categorify the irreducible highest
  weight representations of the corresponding quantum
  groups~\cite{KangKashiwara:CatCycKLR,BK:GradedKL,Webster:HigherRep,
  BrundanStroppel:KhovanovII}. These algebras are only well understood
  for quivers of type $\Aone$ and $A_\infty$, where it has been possible
  to bootstrap results using the Brundan-Kleshchev isomorphism
  theorem~\cite{BK:GradedKL}, which shows that the cyclotomic KLR
  algebras of type $A$ are isomorphic to the (ungraded) Ariki-Koike
  algebras.  Using the Brundan-Kleshchev isomorphism, there is now an
  extensive literature in type~$A$ including a categorification
  theorem~\cite{BK:GradedDecomp}, cellular
  bases~\cite{HuMathas:GradedCellular,Bowman:ManyCellular}, and results
  on Specht
  modules~\cite{BKW:GradedSpecht,KMR:UniversalSpecht,HuMathas:GradedInduction}.

  Very little explicit information is known about the cyclotomic KLR
  algebras for other Cartan types and even in type~$A$ our understanding
  is imperfect because it is seen through the lens of the
  Brundan-Kleshchev isomorphism theorem, which does not keep track of
  the grading. Hu and Shi have proved an amazing general formula that
  gives the graded dimensions of the weight spaces of the cyclotomic KLR
  algebras of symmetrisable Cartan type~\cite{HuShi:GradedDimensions}.
  Recent work of the second author and
  Tubbenhauer~\cite{MathasTubbenhauer:Subdivision,MathasTubbenhauer:BAD}
  shows that the cyclotomic KLR algebras of types $A^{(2)}_{2e}$,
  $B_\infty$, $\Cone$ and $D^{(1)}_{e-1}$ are graded cellular algebras,
  in the sense of \cite{GL,HuMathas:GradedCellular},
  using the weighted KLRW algebras pioneered by
  Webster~\cite{Webster:HigherRep,Webster:RouquierConjecture,Webster:WeightedKLR}
  and Bowman~\cite{Bowman:ManyCellular}, who mainly consider type~$A$.
  The combinatorics in this paper is influenced by a beautiful series of
  papers by Ariki and Park~\cite{ArikiPark:A2,ArikiPark:C,ArikiPark:D},
  which determine the representation type of the cyclotomic KLR
  algebras in certain types, and by the attempts of Ariki, Park and
  Speyer~\cite{ArikiParkSpeyer:C} to construct Specht modules for the
  cyclotomic KLR algebras of affine type~$C$. The semisimplicity of the
  cyclotomic KLR algebras of types~$A$ and~$C$ is determined in the
  papers~\cite{Mathas:Singapore,Speyer:SemisimplicityC}.

  The cyclotomic KLR algebras are defined by generators and relations
  with the most important relations being encoded in Rouquier's
  $Q$-polynomials. Modulo a choice of signs, which do not affect the
  algebras up to isomorphism, the ``standard'' $Q$-polynomials in
  literature take the form
  \[
        Q_{i,j}(u,v) = \begin{cases*}
           u-v        & if $i\to j$,\\
           (u-v)(v-u) & if $i\leftrightarrows j$,\\
           u-v^2      & if $i\Rightarrow j$,
        \end{cases*}
  \]
  where $i$ and $j$ are vertices of the underlying quiver and  $u$ and
  $v$ are indeterminates of degree~$2$ (see \autoref{S:Quiver} for more
  detailed definitions.) Our starting point is to consider
  ``deformations'' of these polynomials, such as
  \[
        Q^x_{i,j}(u,v) = \begin{cases*}
           u-v-x^2             & if $i\to j$,\\
           (u-v+x^2)(v-u+x^2) & if $i\leftrightarrows j$,\\
           u-(v-x^2)^2         & if $i\Rightarrow j$,
        \end{cases*}
  \]
  where $x$ is an indeterminate over $\K$ of degree~$1$. (We allow more general
  deformations.) Using the standard $Q$-polynomials $Q_{i,j}(u,v)$, and a dominant
  weight~$\Lambda$, we define the ``standard'' (cyclotomic) KLR
  algebras~$\Rn$ via \autoref{D:KLR}. Using the deformed
  $Q$-polynomials $Q^x_{i,j}(u,v)$, the same definition gives us the
  deformed (cyclotomic) KLR algebras $\Rx$, for $n\ge0$. For quivers of
  types $\Aone$ and $\Cone$ we show that the deformed cyclotomic KLR
  algebras $\Rx$ are split semisimple graded algebras over
  $\Kxx=\K[x,x^{-1}]$. We prove this by introducing \textit{content
  systems}, which are a generalisation of the classical content
  functions from the semisimple representation theory of the symmetric
  groups.  Unlike the classical situation, a content system consists of
  two functions that determine ``contents'' and ``residues'', where the
  content function is determined by the $Q$-polynomials. We use content
  systems to construct irreducible representations of the deformed
  cyclotomic KLR algebras of types~$A$ and~$C$ over $\Kxx$, giving a
  generalisation of Young's seminormal forms in the KLR setting.  The
  appearance of seminormal forms in the representation theory of the KLR
  algebras of type~$A$ is not surprising but, at least for us, this was
  unexpected for the algebras of type~$C$.

  The graded semisimple deformations of the cyclotomic KLR algebras
  gives a new way of approaching the non-semisimple representation
  theory of the cyclotomic KLR algebras, even though these algebras are
  rarely semisimple. The deformed cyclotomic KLR algebras are semisimple
  over $\Kxx$ but they stop being semisimple when $x$ is not invertible,
  which allows us to recover the standard cyclotomic KLR algebras from
  the deformed algebras by specialising $x=0$. In this way, we can use
  the semisimple representation theory of $\Rx$ over $\K[x,x^{-1}]$ to
  understand the non-semisimple representation theory of~$\Rn$
  over~$\K$. In fact, throughout the paper we work mainly with the
  deformed KLR algebra~$\Rx$, both because $\Rx$ is easier to work with
  and because it has a richer representation theory that encodes
  everything about~$\Rn$.

  The first main result of this paper, \autoref{T:kcellular}, is the following.

  \begin{MainTheorem}\label{MT:ACCellular}
    Let $\Rn$ be a cyclotomic KLR algebra of type $\Aone$ or
    $\Cone$. Then $\Rn$ is a graded cellular algebra.
  \end{MainTheorem}

  Knowing that an algebra is cellular gives a framework for
  understanding its representation theory, including a construction of
  the irreducible representations of the algebra. We actually prove
  several enhanced versions of \autoref{MT:ACCellular}.  First, over a
  positively graded ring $K$, such as $\K[x]$, we show that the deformed
  KLR algebra $\Rx$ over $K$ is a graded $K$-cellular algebra, where
  $K$-cellularity further generalises cellular algebras to the category
  of finite dimensional graded algebras that are defined over graded
  rings. Secondly, we construct four different cellular bases of~$\Rx$,
  two of which specialise to give cellular bases of~$\Rn$, and two of
  which give bases for the split semisimple algebra $\Rx$ when we extend
  scalars to~$\K[x^\pm]$.

  The proof of \autoref{MT:ACCellular} starts by using our
  generalisation of Young's seminormal forms to show that $\Rx$ has two
  seminormal cellular bases, $\set{\Lfst}$ and $\set{\Gfst}$, over
  $\K[x^{\pm}]$. The seminormal bases are characterised as bases of
  simultaneous eigenvectors for the generators $y_1,\dots,y_n$ of~$\Rx$,
  where the eigenvalues are given by our content systems.
  The seminormal bases are then used to show that $\Rx$ has two
  ``integral'' cellular bases, $\set{\Lpsist}$ and $\set{\Gpsist}$
  (\autoref{D:psist}), that are defined over~$\K[x]$ and which specialise
  to give cellular bases of~$\Rn$. In type~$A$, the $\psi$-bases
  of~$\Rx$ generalise the $\psi$-bases constructed
  in~\cite{HuMathas:GradedCellular}.  The transition matrix from the
  $\Lfst[]$-basis to the $\Lpsist[]$-basis is unitriangular, as is the
  transition matrix from the $\Gfst[]$-basis to the $\Gpsist[]$-basis,
  so it is very easy to deduce properties of $\psi$-bases from the
  seminormal bases.

  The key difference between the $\Lpsist[]$-basis and the
  $\Gpsis[]$-basis, and between the $\Lfst[]$-basis and the
  $\Gfst[]$-basis, is that one is defined using the \textit{reverse
  dominance order} on the set of $\ell$-partitions and the other is
  defined using the \textit{dominance order}. (Here $\ell$ is the
  \text{level} of the dominant weight~$\Lambda$; see
  \autoref{S:Tableaux}.) That is, by reversing the choice of partial order
  in our definitions we can switch between these two families of
  cellular bases. In turn, this leads to the construction of two closely
  related families of \textit{cell modules}, or \textit{Specht modules},
  $\set{\LSlam[\bmu]}$ and $\set{\GSlam[\bnu]}$, and two families of
  simple $\Rx(F[x])$-modules $\set{\LDmu}$ and $\set{\GDnu}$.
  Throughout the paper we keep track of these two families of modules
  because, aside from the notation, doing this requires almost no extra
  work, with the only real difference being whether we work with the
  dominance or reverse dominance order. In fact, we need to work with
  these two ``dual'' families of modules because some of our main
  results are proved by exploiting the close connections between these
  two families of modules.

  Once we have proved that $\Rx$ and $\Rn$ are cellular algebras, we
  next turn to understanding their representation theory. We first use
  the semisimple representation theory to show that~$\Rx$ (and~$\Rn$),
  is a graded symmetric algebra. There is a natural symmetrising form
  that is defined using \textit{defect polynomials}
  (\autoref{D:DefectPolynomial}), which are graded analogues of the
  \textit{generic degrees} from the representation theory of
  cyclotomic Hecke algebras~\cite{Mathas:gendeg}. In particular, this
  allows us to show that $\LSlam$ is isomorphic to the dual of $\GSlam$,
  up to shift. The \textit{defect} of a Specht module is equal to the
  degree of its defect polynomial. Defect is a key invariant of the
  blocks of the cyclotomic KLR algebras, which generalises the
  $p$-weight of a partition in the modular representation theory of the
  symmetric groups.

  As a second application of the semisimple representation theory, we
  give explicit Specht filtrations of the modules obtained by inducing
  and restricting the Specht modules of~$\Rx$ over an arbitrary ring.
  Together with the combinatorics based on the defect polynomials, the
  graded branching rules for the Specht modules translate into our next
  main result, which is a categorification theorem. To state this we
  need to introduce some notation.

  Let $\K$ be a field and $x$ an indeterminate over $\K$. We consider
  $\K[x]$ as a positively graded ring, with $x$ in degree~$1$, and set
  $\A=\Z[q,q^{-1}]$. Let $\Rep_\K\Rx(\K[x])$ be the category of graded
  $\Rx(\K[x])$-modules that are finite dimensional as $\K$-vector spaces
  and let $\Proj_\K\Rx(\K[x])$ be the full subcategory of projective
  $\Rx(\K[x])$-modules. Let $[\Rep_\K\Rx[\bullet](\K[x])]$ and
  $[\Proj_\K\Rx[\bullet](\K[x])]$ be the direct sum of the Grothendieck
  groups of  these categories for~$n\ge0$, which we consider as free
  $\A$-modules by letting $q$ act as the grading shift functor.

  Suppose that $\Gamma$ is a quiver of type $\Aone$ or type $\Cone$.
  Let $\Uq$ be the corresponding quantised Kac-Moody algebra and let
  $\UA$ be Lusztig's $\A$-form of~$\Uq$. For a dominant
  weight~$\Lambda$, let $L(\Lambda)_\A$ be the $\A$-form of the
  corresponding irreducible integrable highest weight module for $\UA$
  and let $L(\Lambda)^*$ be is dual, with respect to the Cartan pairing.

  \begin{MainTheorem}[Cyclotomic categorification theorem]
    \label{MT:Categorification}
    Suppose that $\Gamma$ is a quiver of type $\Aone$ or $\Cone$ and let
    $\Lambda$ be a dominant weight.  Then, as $\UA$-modules,
    \[
        L(\Lambda)_\A\cong[\Proj_\K\Rx[\bullet](\K[x])]
        \qquad\text{and}\qquad
        L(\Lambda)_\A^*\cong[\Rep_\K(\Rx[\bullet](\K[x])].
    \]
  \end{MainTheorem}

  This result, which is \autoref{T:CyclotomicCat}, is not new. In
  type~$\Aone$ it is one of the main results of \cite{BK:GradedDecomp}.
  More generally,~\cite{KangKashiwara:CatCycKLR} establishes this result
  whenever~$\Gamma$ is a quiver of symmetrisable Cartan type. What is
  new about this result is that it is deduced almost directly from the
  graded branching rules for the Specht modules of~$\Rx(\K[x])$, which
  directly encode the action of $\Uq$ on the Grothendieck groups. This
  explicit link with the representation theory of~$\Rx(\K[x])$ makes it
  much easier to apply this result to the representation theory
  of~$\Rx(\K[x])$. In fact, the information flow is stronger in both
  directions, so we also use the representation theory of~$\Rx(\K[x])$
  to better understand~$L(\Lambda)$. In particular, we are able to give
  detailed information about the canonical bases of~$L(\Lambda)_\A$
  and~$L(\Lambda)_\A^*$ and their role in this categorification theorem.

  Our approach to \autoref{MT:Categorification} is partly based on
  \cite{BK:GradedDecomp}, although our perspective is fundamentally
  different because we work almost exclusively inside the Grothendieck
  groups of the cyclotomic KLR algebras whereas \cite{BK:GradedDecomp}
  works mainly inside a combinatorial Fock space, which we also use. In
  particular, we use \autoref{MT:ACCellular}, and the triangularity of
  the decomposition matrices of~$\Rx(\K[x])$, to show that Lusztig's bar
  involution is triangular on the basis of Specht modules.  Our
  arguments work simultaneously for the algebras of type~$\Aone$ and
  $\Cone$ and, via \autoref{MT:ACCellular}, we obtain two versions of
  \autoref{MT:Categorification} corresponding to the $\Lpsist[]$ and
  $\Gpsist[]$ cellular bases. This gives two explicit realisations of
  the irreducible integrable highest weight $\UA$-modules
  $L(\Lambda)_\A$ and $L(\Lambda)_\A^*$.

  Our next main goal is to classify the irreducible graded
  $\Rx(\K[x])$-modules. Our parallel theories, using the $\Lpsist[]$ and
  $\Gpsist[]$ cellular bases, leads to two combinatorial descriptions of
  the crystal graph of $L(\Lambda)$, which we call the $\Ldom$-crystal
  graph and the $\Gdom$-crystal graphs in this introduction.  To
  describe these, let~$I$ be the vertex set of the quiver~$\Gamma$. The
  paths in the crystal graphs of $L(\Lambda)$ are labelled by $n$-tuples
  $\bi\in I^n$, corresponding to generalisations of Kleshchev's
  \textit{good node sequences} (\autoref{D:NormalGood}).  Each good node
  sequence $\bi$ determines two paths: one path $\zero\Lgood\bmu$ in the
  $\Ldom$-crystal graph  and a second path $\zero\Ggood[\bi]\bnu$ path in the
  $\Gdom$-crystal graph. Here, $\zero$ is the empty $\ell$-partition and
  $\bmu,\bnu$ are $\ell$-partitions of~$n$. Let
  $\LKlesh=\set{\bmu\in\Parts|0\Lgood\bmu\text{ for some }\bi\in I^n}$
  and $\GKlesh=\set{\bnu\in\Parts|0\Ggood\bmu\text{ for some }\bj\in
  I^n}$ be the vertex sets of the two crystal graphs. Calculations
  with the canonical bases in the Grothendieck
  groups~$\Rep_\K\Rx[\bullet](\K[x])$ and $\Proj_\K\Rx[\bullet](\K[x])$ allows
  us to classify the irreducible $\Rx(\K[x])$-modules over a field, for
  $n\ge0$. As \autoref{T:KleshchevSimples}, we prove.

  \begin{MainTheorem}\label{MT:Simples}
    Let $\K$ be a field. Up to shift,
    $\set[\big]{\LDmu|\bmu\in\LKlesh}$ and
    $\set[\big]{\GDnu|\bnu\in\GKlesh}$ are both complete sets of pairwise
    non-isomorphic irreducible $\Rx(\K[x])$-modules.
  \end{MainTheorem}

  In particular, over any field, this result classifies the irreducible
  modules of the cyclotomic KLR algebras of type $\Aone$ and $\Cone$.

  \autoref{MT:Simples} implies that there is a bijection
  $\mull\map{\LKlesh}\GKlesh$ such that $\LDmu\cong\GDnu[\mull(\bmu)]$.
  In \autoref{C:MullinuxMap} we show that if $\bmu\in\LKlesh$ and
  $\zero\Lgood\bmu$ is a path in the $\Ldom$-crystal graph
  of~$L(\Lambda)$ then there is a unique $\ell$-partition
  $\bnu=\mull(\bmu)$ such that $\zero\Ggood[\bi]\bnu$ is a path in the
  $\Gdom$-crystal graph. This gives a way to compute the
  $\ell$-partition $\mull(\bmu)$. In the special case of the symmetric
  groups, this gives another description of the \textit{Mullineux map},
  which describes what happens to the simple modules of the symmetric
  group when they are tensored with the sign representation. We
  introduce a sign representation for the algebras $\Rx(\K[x])$ and show
  in our setting, which generalises that of the symmetric groups, the
  Mullineux map is the function $\bmu\mapsto\mull(\bmu)'$, where $\bmu'$
  is the $\ell$-partition conjugate to~$\bmu$; see
  \autoref{S:Bases}.

  Finally, we show that Kleshchev's modular branching
  rules~\cite{Klesh:III} extend to give branching rules for the simple
  $\Rx(\K[x])$-modules. For $i\in I$, let $\iRes$ and $\iInd$ be the
  corresponding $i$-restriction and $i$-induction functors and let $e_i$
  and $f_i$ be Kashiwara's operators on the crystal graph of
  $L(\Lambda)$. We refer the reader to \autoref{S:ModularBranching} for
  the precise definitions and statements, but the main results take the
  form:

  \begin{MainTheorem}\label{MT:ModularBranching}
    Suppose that $\bmu\in\LKlesh$, $\bnu\in\GKlesh$ and $i,j\in I$. Then,
    up to grading shift,
    \[
      \LDmu[e_i\bmu]=\soc(\iRes\LDmu),\quad
      \LDmu[f_i\bmu]=\head(\iInd\LDmu),\quad
      \GDnu[e_j\bnu]=\soc(\iRes[j]\GDnu)
      \quad\text{and}\quad
      \GDnu[f_j\bnu]=\head(\iInd[j]\GDnu).
    \]
  \end{MainTheorem}

  In type $\Aone$, Brundan and Kleshchev~\cite[Theorem]{BK:GradedDecomp}
  prove this result by lifting Ariki's~\cite{Ariki:can,Ariki:branching}
  and Grojnowski's work~\cite{Groj:control}, from the ungraded
  representation theory, into the KLR world. More generally, for any
  symmetrisable Cartan type, Lauda and
  Vazirani~\cite{LaudaVazirani:CatCrystals} show that analogues of these
  modular branching rules categorify the crystal graph of~$L(\Lambda)$
  by lifting parts of Grojnowski's approach to the KLR setting. Lauda
  and Vazirani's result does not imply \autoref{MT:ModularBranching}
  because it is not clear how their crystal graph is related to the
  labelling of the simple modules given in \autoref{MT:Simples}. Our
  proof of \autoref{MT:ModularBranching} is almost axiomatic in that it
  uses \autoref{MT:Categorification} to lift the result from
  \autoref{MT:Categorification} and properties of the canonical basis.

  Throughout the paper we work almost exclusively with a deformed
  cyclotomic KLR algebra~$\Rx$ that has a content system to prove our
  results, after which the results for~$\Rn$ are obtained by specialising
  the deformation parameters to~$0$. We show by example that every
  cyclotomic KLR algebra of types $\Aone$ and $\Cone$ has a graded
  content system over $\Z[x]$, so our results apply to all cyclotomic KLR
  algebras of affine types~$A$ and~$C$ over any ring. In type~$\Aone$,
  the results we obtain for $\Rn$ were known but those for $\Rx$ are
  new. In type~$\Cone$, all of these results are completely new. As we
  note in \autoref{S:Quiver}, the results in this paper also extend to
  quivers of type $A_\infty$ and $C_\infty$. It likely that the general
  framework that we develop can be modified to work in other types.

  It is quite striking that we are able to prove all of these results
  using a common framework for the cyclotomic KLR algebras of
  type~$\Aone$ and $\Cone$.  Ultimately, the reason why this works is
  that our deformation arguments show that the algebra $\Rx$ over
  $\K[x^{\pm}]$ is isomorphic to a direct sum of matrix algebras that
  depend only on~$n$ and~$\ell$, and not on the choice of dominant
  weight~$\Lambda$ or even on the quiver~$\Gamma$. In fact,
  \autoref{T:SSemisimple} shows that if~$n$ and~$\ell$ are fixed then,
  for any choice of content system, the deformed cyclotomic KLR algebras
  over $\K[x^{\pm}]$ are canonically isomorphic as ungraded algebras.

  An \hyperref[ANotationalIndex]{index of notation} is included at the
  end of the document, before the list of references.


  \subsection*{Acknowledgements}
  First and foremost, I thank Anton for his important contributions to
  this paper. He was an insightful mathematician and it was a pleasure
  to work with him. I thank Jun Hu, Huansheng Li, Huang Lin, Tao Qin,
  Liron Speyer, Catharina Stroppel, Daniel Tubbenhauer and Ben Webster
  for helpful discussions and feedback. I thank the referee for
  their comprehensive report, which significantly improved our
  exposition. Finally, I thank the Max Plank Institute in Bonn, where
  part of this research was completed.

  This research was partially supported by the Australian Research Council.

\section{KLR algebras}\label{S:KLRAlgebras}

\subsection{Graded rings, algebras and modules}\label{S:Rings}
Throughout this paper we work with $\Z$-graded rings, algebras and
modules.  For convenience, we refer to each of these structures as being
\textit{graded}. This section recalls the basic definitions that we need
for modules over graded rings.

All rings in the paper will be commutative integral domains with~$1$.  A
\emph{graded ring} is a ring~$K$ that has a decomposition
$K=\bigoplus_{d\in\Z} K_d$ as an additive abelian group such that
$K_dK_e\subseteq K_{e+d}$. In particular, note that $K_0$ is a subring
of~$K$.
\footnote{We apologise to the readers who instantly think that $K$ is a
field. In the body of the paper we mostly work with a field $\K$, which
is a $\k$-algebra (often the field of fractions of the ring~$\k$), and
we consider modules over the graded rings $\kx$, $\Kx$ and $\Kxx$.
}

Let $K$ be a graded commutative domain. Then:
\begin{itemize}[nosep]
  \item A \emph{graded $K$-module} is a $K$-module $M$ that admits a
  decomposition $M=\bigoplus_{d\in\Z}M_d$ as a $K_0$-module such that
  $K_dM_e\subseteq M_{d+e}$.

  \item A \emph{graded $K$-algebra} is a $K$-algebra $A$ that admits an
  decomposition $A=\bigoplus_{d\in\Z}A_d$ as a graded $K$-module such
  that $K_dA_e\subseteq A_{d+e}$.

  \item A \emph{graded $A$-module} is an $A$-module $M$ that admits a
  decomposition $M=\bigoplus_{d\in\Z}M_d$ as a graded $K$-module such
  that $A_dM_e\subseteq M_{d+e}$.
\end{itemize}
If $R=\bigoplus_dR_d$ is a graded ring, algebra or module let
$\underline{R}$ be the corresponding structure obtained by forgetting
the grading.  An element $m\in R$ is \emph{homogeneous} of degree~$d$ if
$0\ne m\in R_d$, in which case we set $\deg(m)=d$. By definition, $0$ is
not homogeneous. In particular, note that if $r\in R$ and $m\in M$ are
homogeneous then $\deg(rm)=\deg(r)+\deg(m)$. Further, $R$ is
\emph{positively graded} if there are no elements of negative degree
(that is, these are non-negatively graded structures) and $R$ is
\emph{concentrated in degree~$d$} if $R=R_d$.

In this paper the three types of graded rings~$K$ that we consider are:
\begin{itemize}[nosep]
  \item commutative domains~$\k$ with~$1$,
  \item polynomial rings~$\kx=\k[\ux]$, where~$\ux$ is a (possibly
  empty) family of indeterminates over~$\k$ with each indeterminate
  having degree~$1$,
  \item and Laurent polynomial rings~$\Kxx=\K[\ux,\ux^{-1}]$, where $\K$
  is a field that is a $\k$-algebra, such as the field of fractions of~$\k$.
\end{itemize}
In these rings, the elements of~$\k$ and~$\K$ are in
degree~$0$.

A \emph{graded field} is a graded ring in which every nonzero
homogeneous element has a multiplicative inverse.  In
particular, $\K$ and ~$\K[x^{\pm}]$ are graded fields. By
\cite[Theorem~4.1]{vanGeelVanOystaeyen:GradedFields} all graded fields
are of this form.

\notation{$\k$}{A commutative integral domain with $1$, concentrated in degree $0$}
\notation{$\K$}{A field that is a $\k$-algebra, again in degree~$0$}
\notation{$\ux$}{A family of indeterminates over the ground ring, which is normally $\k$}
\notation{$\kx$}{The positively graded polynomial ring $\k[\ux]$, with $x\in\ux$ in degree~$1$}
\notation{$\Kxx$}{The $\Z$-graded Laurent polynomial ring $\K[\ux,\ux^{-1}]$}

If $A$ is a graded $K$-algebra and $M$ is a graded $A$-module then
graded submodules, quotient modules, projective modules, $\ldots$ are
defined in the obvious ways. If $K$ is a graded field and $A$ is a
graded $K$-algebra then an \emph{irreducible graded $A$-module} is a
graded $A$-module that has no non-trivial proper graded $A$-submodules.
We emphasise that irreducible graded modules make sense when the ground
ring is a graded field that is not a field.

\begin{Remark}\label{R:GradedIrreducible}
  Let $K$ be a field and $A$ a graded $K$-algebra. Then a
  graded $A$-module $D$ is an irreducible graded $A$-module if and only
  if $\underline{D}$ is an irreducible $\underline{A}$-module by
  \cite[Theorem 4.4.4 and Theorem 9.6.8]{NV:graded}. In contrast, if
  $A$ is a graded $\Kxx$-algebra then an irreducible graded $A$-module
  is not necessarily irreducible when we forget the grading. For example, if
  $A=\K[x^\pm]$ and $D=\K[x^\pm]$ then $D$ is an irreducible graded
  $A$-module but $D$ is not irreducible as an $\underline{A}$-module
  because, for example, it contains the (non-homogeneous) ideal
  $(x+1)\K[x^\pm]$.
\end{Remark}

If $M$ and $N$ are graded $A$-modules then a \emph{homogeneous
$A$-module homomorphism} of degree~$d$ is an $A$-module
homomorphism $f\map MN$ such that $\deg f(m)=\deg(m)+d$ whenever $m\in
M$ is homogeneous. In this case we write $\deg f=d$. The map $f$ is an
\emph{$A$-module isomorphism} if it is bijective and it is homogeneous
of degree~$0$.

\notation{$\A$}{The ring $\A=\Z[q,q^{-1}]$, where $q$ is an indeterminate}
\notation{$\AA$}{The ring $\Q(q)$ of rational functions in~$q$}
\notation{$q^dM$}{The graded module obtained by shifting the grading on $M$ by~$d$}

Let $q$ be an indeterminate and set $\A=\Z[q,q^{-1}]$ and $\AA=\Q(q)$.
If $M=\bigoplus_d M_d$ is a graded $A$-module and $s\in\Z$ let
$q^sM=\bigoplus_d(q^sM)_d$ be the graded $A$-module that is equal
to~$\underline{M}$ as an ungraded module, has $(q^sM)_d=M_{d-s}$ and
with $A$-action inherited from the action on~$M$.

If $M$ and $N$ are graded $A$-modules let $\Hom_A(M,N)$ be the
homogeneous $A$-module homomorphisms of degree~$0$. Then
$\Hom_A(q^{d}M,N)\cong\Hom_A(M,q^{-d}N)$ is naturally identified with
the set of homogeneous maps $M\to N$ of degree~$d$, for $d\in\Z$. Set
$\HOM_A(M,N) = \bigoplus_{d\in\Z}\Hom_A(q^dM,N)$. Define $\End_A(M)$ and
$\END_A(M)$ similarly.

\notation{$\Hom_A(M,N)$}{The homogeneous $A$-module maps $M\to N$ of degree $0$}
\notation{$\HOM_A(M,N)$}{\textit{All} homogeneous $A$-module maps $M\to N$}
\notation{$\End_A(M)$}
   {The homogeneous $A$-module endomorphisms of $M$ of degree $0$}
\notation{$\END_A(M)$}{\textit{All} homogeneous $A$-module endomorphisms of $M$}

\begin{Remark}\label{R:GradedSimples}
  For geometric reasons, indeterminates are usually put in degree~$2$.
  It is more convenient for us to put the indeterminates in~$\ux$ in
  degree~$1$ because then the graded field $\K[x^{\pm1}]$ has a unique
  irreducible graded representation, namely itself; see
  \autoref{R:GradedIrreducible}. (In contrast, if  we set $\deg(x)=2$
  then $\K[x^{\pm1}]$ and $q\K[x^{\pm1}]$ are non-isomorphic irreducible
  graded $\K[x^{\pm1}]$-modules.) On the other hand, $\set{q^d\K|d\in\Z}$
  is a complete set of pairwise non-isomorphic  unique irreducible
  graded $\K[x]$-modules, where the $\K[x]$-module $q^d\K$ is
  concentrated in degree~$d$ and~$x$ acts as multiplication by $0$.
\end{Remark}

If $A$ is a graded $K$-algebra then we will usually work in the category
$\Rep A$ of finitely generated $A$-modules with homogeneous maps of
degree~$0$. If $K=\bigoplus_dK_d$ and $\K=K_0$ is a field let
$\Rep_{\K}A$ be the full subcategory of $\Rep A$ consisting of
$A$-modules that are \textit{finite dimensional} as $\K$-vector spaces.
Similarly, let $\Proj A$ be the additive subcategory of $\Rep A$
consisting of \emph{projective graded $A$-modules} and let $\Proj_\K A$
be the corresponding subcategory of $\Rep_\K A$. The proofs of our Main
Theorems take place in the categories $\Rep_\K\Rx(\Kx)$ and
$\Proj_\K\Rx(\Kx)$.

Let $[\Rep_\K A]$ and $[\Proj_\K A]$ be the \emph{Grothendieck groups} of
the categories $\Rep_\K A$ and $\Proj_\K A$, respectively. Given a module
$M$ in $\Rep_\K A$, or in $\Proj_\K A$, let $[M]$ be its image in $[\Rep_\K
A]$ or $[\Proj_\K A]$, respectively.  Both $[\Rep_\K A]$ and $[\Proj_\K A]$
are free $\A$-modules where $q$ acts by grading shift.  That is,
$[qM]=q[M]$.

\subsection{Quivers and \texorpdfstring{$Q$}{Q}-polynomials}\label{S:Quiver}

  In this section we fix the Lie theoretic data that will be used
  throughout the paper.  Let $\N=\Z_{\ge0}$ be the set of non-negative
  integers.

  \notation{$\N$}{The set of non-negative integers $\Z_{\ge0}$}
  \notation{$\Gamma$}{A symmetrisable quiver, usually of type $\Aone$ or $\Cone$}
  \notation{$I$}{The vertex set $\set{0,1,\dots,e-1}$ of $\Gamma$}
  \notation{$C=(\cij)$}{Cartan matrix of $\Gamma$}
  \notation{$\di$}{$D=\diag(\di[0],\dots,\di[e-1])$ is the symmetriser of $C$}
  \notation{$\alpha_i$}{Simple root, for $i\in I$}
  \notation{$\Lambda_i$}{Fundamental weight, for $i\in I$}
  \notation{$P^+$}{Dominant weight lattice}
  \notation{$Q^+$}{Positive root lattice}

  Let $\Gamma$ be a symmetrisable \emph{quiver} $\Gamma$ with vertex
  set~$I$. Let $\bigl(C, P, P^\vee, \Pi, \Pi^\vee\bigr)$ be the
  \emph{Cartan data} attached to~$\Gamma$, consisting of:
  \begin{itemize}
    \item A \emph{symmetrisable Cartan matrix}, $C=(\cij)_{i,j\in I}$
    satisfies $\cij[ii]=2$, $\cij\le0$ for $i\ne j$, $\cij=0$ whenever
    $\cij[ji]=0$. Since $C$ is symmetrisable, there exists a diagonal
    matrix $D=\diag(\di|i\in I)$ such that $DC$ is symmetric
    \item The \emph{weight lattice} $P$ is a free abelian group with
    basis the \emph{simple roots} $\Pi=\set{\alpha_i|i\in I}$.
    \item The \emph{dual weight lattice} is $P^\vee=\Hom(P,\Z)$ has
    basis the \emph{simple coroots} $\Pi^\vee=\set{\alpha_i^\vee|i\in I}$.
  \end{itemize}
  The \emph{Cartan pairing} $\cart{\ }{\ }\map{P^\vee\times P}\Z$
  and \emph{fundamental weights} $\set{\Lambda_i|i\in I}\subset P$ are
  given by
  \[
      \cart{\alpha^\vee_i}{\alpha_j} =\cij\quad\text{and}\quad
      \cart{\alpha^\vee_i}{\Lambda_j}=\delta_{ij},
      \qquad\text{for }i,j\in I.
  \]
  The \emph{positive root lattice} is $Q^+=\sum_{i\in I} \N \alpha_i$,
  and $P^+=\bigoplus_{i\in I}\N\Lambda_i$ is the set of \emph{dominant
  weights} of~$\Gamma$. The \emph{height} of $\alpha=\sum_{i\in I} h_i
  \alpha_i\in Q^+$ is the non-negative integer
  $\operatorname{ht}(\alpha)=\sum_{i\in I} h_i$.  Let~$Q^+_n$ be the set
  of all elements of $Q^+$ of height $n$.  Set $\h=\Q\otimes_\Z P^\vee$.
  As $C$ is symmetrisable, there exists a symmetric bilinear form
  $\kill{\ }{\ }$ on~$\h^*$ such that
  \[
    \kill{\alpha_i}{\alpha_j}=\di \cij=\cij\di[j]\qquad\text{and}\qquad
    \cart{\alpha_i^\vee}{\lambda}
         =\frac{2\kill{\alpha_i}{\lambda}}{\kill{\alpha_i}{\alpha_i}},
    \qquad\text{for $\lambda\in\h^*$ and $i\in I$.}
  \]

  Fix $n\in\N$ and let $\Sym_n$ be the \emph{symmetric group} on~$n$
  letters. As a Coxeter group, $\Sym_n$ is generated by the simple
  transpositions $\sigma_1,\dots,\sigma_{n-1}$, where $\sigma_k=(k,k+1)$ for
  $1\le k<n$. Let $L\colon\Sym\to\N$ be the \emph{length function}
  on~$\Sym_n$, so if $w\in\Sym_n$ then $L(w)=l$ if $l$ is minimal such
  that $w=\sigma_{a_1}\dots \sigma_{a_l}$, for some $1\le a_j<n$. A \emph{reduced
  expression} for $w\in \Sym_n$ is any expression
  $w=\sigma_{a_1}\dots \sigma_{a_l}$ with $l=L(w)$.
  \notation{$\Sym_n$}{Symmetric group on $\set{1,2,\dots,n}$}
  \notation{$\sigma_k$}{Simple reflection $\sigma_k=(k,k+1)\in\Sym_n$, for $1\le k<n$}
  \notation{$L(w)$}{Coxeter length of $w\in\Sym_n$}

  The group $\Sym_n$ acts from the left on the set
  $I^n=I\times\dots \times I$ by place permutations: if $w\in\Sym^n$ and
  $\bi=(i_1,\dots,i_n)\in I^n$ write $w\bi=(i_{w(1)},\dots,i_{w(n)})$.

  In this paper we will mainly consider the quivers of type $\Aone$
  ($e\ge2$) and $\Cone$ ($e\ge3)$, for which we use the following
  quivers:

  \notation{$\Aone$}{Affine quiver of type $A$ with vertex set $I$}
  \notation{$\Cone$}{Affine quiver of type $C$ with vertex set $I$}
  \begin{center}
    \begin{tabular}{@{\quad}cccc@{\quad}}\toprule
      Type & Dynkin diagram & $\delta$ & $(\di[0],\dots,\di[e-1])$
      \\\midrule
          $\Aone$ &
          \begin{tikzpicture}[centered]
            \draw(2.5,0)--(4,0)--(2,0.5)--(0,0)--(1.5,0);
            \draw[dashed](1.5,0)--(2.5,0);
            \foreach \x/\y/\z [count=\bc] in {2/0.5/0,4/0/1,3/0/2,0/0/e-1,1/0/e-2} {
              \ifnum\bc=1
                \node[vertex, label=above:$\z$] at (\x,\y){};
              \else
                \node[vertex, label=below:$\z$] at (\x,\y){};
              \fi
            }
          \end{tikzpicture}
          & $\alpha_0+\alpha_1+\dots+\alpha_{e-2}+\alpha_{e-1}$
          & $(1,1,\dots,1,1)$
      \\[8mm]
          $\Cone$ &
          \begin{tikzpicture}[centered]
            \draw[->-, double, double distance=2pt](0,0)--(1,0);
            \draw(1,0)--(1.5,0);
            \draw[dashed](1.5,0)--(2.5,0);
            \draw(2.5,0)--(3,0);
            \draw[-<-, double, double distance=2pt](3,0)--(4,0);
            \foreach \x/\z in {0/0,1/1,3/e-2,4/e-1} {
                \node[vertex, label=below:$\z$] at (\x,0){};
            }
          \end{tikzpicture}
          & $\alpha_0+2\alpha_1+\dots+2\alpha_{e-2}+\alpha_{e-1}$
          & $(2,1,\dots,1,2)$
       \\[5mm]
      \bottomrule
    \end{tabular}
  \end{center}
  Here, $\delta$ is the \emph{null root}, which satisfies
  $\cart{\delta}{\alpha_i^\vee}=0$, for $i\in I$. Notice that for both of these
  quivers we have $I=\set{0,1,\dots,e-1}$. Our arguments apply equally
  well
  to the infinite quivers $A_\infty$ and $C_\infty$ but there is no real
  gain in considering these because the cyclotomic KLR algebras for
  these quivers are isomorphic to cyclotomic KLR algebras for a suitably
  large finite quiver.

  Fix a (graded) commutative domain $K=\bigoplus_{d\in\Z} K_d$ and let
  $u,v$ be indeterminates over~$K$. Following
  Rouquier~\cite[Definition~3.2.2]{Rouquier:QuiverHecke2Lie} and
  Kashiwara-Kang~\cite{KangKashiwara:CatCycKLR}, a family
  of \emph{$Q$-polynomials} for~$\Gamma$ is a collection of polynomials
  $Q_{ij}(u,v)\in K[u,v]$, for $i,j\in I$, such that
  $Q_{i,j}(u,v)=Q_{j,i}(v,u)$, $Q_{i,i}(u,v)=0$ and if $i\ne j$ then
 \begin{equation}\label{E:Qpolynomials}
        Q_{i,j}(u,v)=\sum_{p,q\ge0} t_{i,j;p,q}u^pv^q,
        \qquad\text{where $t_{i,j,-\cij,0}\in K_0^\times$ and $t_{i,j;p,q}\in K_d$, }
 \end{equation}
  where $d=-2\kill{\alpha_i}{\alpha_j}-p\kill{\alpha_i}{\alpha_i}
            -q\kill{\alpha_j}{\alpha_j}$.
  That is, $Q_{i,j}(u,v)$ is homogeneous of degree $d$.  By assumption,
  $Q_{i,j}(u,v)=Q_{j,i}(v,u)$, so $t_{i,j;p,q}=t_{j,i;q,p}$.  One
  standard choice for these polynomials is
  \begin{equation}\label{E:StandardQ}
        Q_{i,j}(u,v) = \begin{cases*}
           u-v        & if $i\to j$,\\
           (u-v)(v-u) & if $i\leftrightarrows j$,\\
           u-v^2      & if $i\Rightarrow j$.
        \end{cases*}
  \end{equation}
  As discussed in the introduction, this paper uses ``deformed
  analogues'' of these standard $Q$-polynomials. More examples can be
  found in \autoref{Ex:ContentSystem} below.
  \notation{$\Qb$}{Family $\Qb=(Q_{i,j}(u,v))_{i,j\in I}$ of Rouquier's
                   $Q$-polynomials}[E:Qpolynomials]

  For $i,j,k\in I$ and indeterminates $u,v$ and $w$ over $K$, define the
  three variable $Q$-polynomials
  \begin{equation}\label{E:Qijk}
    Q_{i,j,k}(u,v,w) = \delta_{ik}\frac{Q_{ij}(u,v)-Q_{jk}(v,w)}{u-w},
  \end{equation}
  where $\delta_{ik}$ is the Kronecker delta.
  It is straightforward to check that $Q_{i,j,k}(u,v,w)\in K[u,v,w]$.

  \subsection{KLR algebras}\label{S:KLR}
  This section defines the (cyclotomic) KLR algebras, which are one of
  the main objects of interest in this paper.  Unless otherwise stated,
  all of our algebras are $K$-algebras, where~$K$ is a (graded)
  commutative integral domain with one.

  As in the last section, let $K=\bigoplus_dK_d$ be a graded commutative
  ring with one and fix algebraically independent indeterminates
  $u_1,\dots,u_n$ over~$K$. The symmetric group $\Sym_n$ acts on
  $K[u_1,\dots,u_n]$ by permuting indeterminates $f\mapsto
  {}^wf=f(u_{w(1)},\dots,u_{w(n)})$, for $w\in\Sym_n$ and
  $f\in K[u_1,\dots,u_n]$.

  Recall from \autoref{S:Quiver}, that $I=\set{0,1,\dots,e-1}$ is the
  (finite) vertex set of the quiver~$\Gamma$ and that we have fixed a
  family $\Qb=\bigl(Q_{ij}(u,v)\bigr)_{i,j\in I}$ of Rouquier's
  $Q$-polynomials. In addition, fix a family of homogeneous
  \emph{weight polynomials} $\Wb = \bigl( W_i(u)\bigr)_{i\in I}$ such that
  \begin{equation}\label{E:Kpolys}
    W_i(u) = \sum_{k=0}^{\kill{\alpha^\vee_i}{\Lambda}}
                      a_{i;k}u^{\kill{\alpha^\vee_i}{\Lambda}-k},
             \qquad\text{ for $a_{i;k}\in K_{\di k}$ and $a_{i;0}=1$}.
  \end{equation}
  The weight polynomials $\Wb$ determine a dominant weight
  $\Lambda=\Lambda_{\Wb}=\sum_{i\in I}l_i\Lambda_i\in P^+$,
  where $l_i=\deg W_i(u)$ for $i\in I$. The \textbf{level}
  of~$\Lambda$ is~$\ell=\sum_il_i$. We assume $\ell\ge1$.
  \notation{$\Wb$}{Family $\Wb=(W_i(u))_{i\in I}$ of weight polynomials, for $i\in I$}[E:Kpolys]
  \notation{$\Lambda$}{The dominant weight in $P^+$ determined by $\Wb$}[E:Kpolys]

  A \emph{cyclotomic KLR datum} is a triple $(\Gamma,\Qb,\Wb)$,
  where~$\Gamma$ is a quiver and $\Qb$ and $\Wb$ are families of
  $Q$-polynomials and weight polynomials for~$\Gamma$, respectively. The
  quiver~$\Gamma$ has vertex set~$I$ and comes equipped with a Cartan
  datum as in \autoref{S:Quiver}.

  If $\alpha\in Q^+_n$ let
  $I^\alpha=\set{\bi\in I^n|\alpha=\alpha_{i_1}+\dots+\alpha_{i_n}}$.

  \notation{$I^\alpha$}{The orbit $\set{\bi\in I^n|\alpha=\alpha_{i_1}+\dots+\alpha_{i_n}}$
  for $\alpha\in Q^+$}
  \notation{$\Rn, \Rn[\alpha]$}{A (standard) cyclotomic KLR algebra}[D:KLR]
  \notation{$\Rnaff, \Rnaff[\alpha]$}{A (standard) KLR algebra}[D:KLR]
  \notation{$\ei$}{An idempotent in, and generator of, $\Rx$ or $\Rn$, for $i\in I$}
                  [D:KLR]
  \notation{$y_1,\dots,y_n$}{Generators of $\Rx$ or $\Rn$}[D:KLR]
  \notation{$\psi_1,\dots,\psi_{n-1}$}{Generators of $\Rx$ or $\Rn$}[D:KLR]

  \begin{Definition}\label{D:KLR}
    Let $(\Gamma,\Qb,\Wb)$ be a cyclotomic KLR datum
    and suppose that  $\alpha\in Q^+_n$.
    The \textbf{KLR algebra} $\Rnaff[\alpha]=\Rnaff[\alpha] (\Qb)$
    is the unital associative $K$-algebra generated by
    \[
      \set{ \ei | \bi \in I^\alpha } \cup \set{ \psi_k | 1\le k<n }
             \cup \set{ y_m | 1\le m\le n}
    \]
    subject to the relations:
    \begin{relations}
     \item $\ei \ei[\bj] = \delta_{\bi,\bj} \ei$
           \quad and\quad
           $\sum_{\bi\in I^\alpha} \ei =1$ \label{R:Ids}
     \item $y_k \ei = \ei y_k$
           \quad and\quad
           $y_k y_m = y_m y_k$ \label{R:y}
     \item $\psi_k y_m = y_m \psi_k$ if $m\ne k,k+1$
           \label{R:psiycomm}
     \item $\psi_k \psi_m =\psi_m \psi_k$ if $|m-k|>1$
           \label{R:psicomm}
     \item $\psi_k \ei = \ei[\sigma_k\bi] \psi_k$, \label{R:psi}
     \item $(\psi_k y_{k+1} - y_k \psi_k)\ei
              = \delta_{i_k,i_{k+1}} \ei
              =(y_{k+1} \psi_k - \psi_k y_k ) \ei$
           \label{R:psiy}
     \item $\psi_k^2 \ei = Q_{i_k,i_{k+1}} (y_k, y_{k+1}) \ei$
           \label{R:quad}
     \item $(\psi_{k+1}\psi_k \psi_{k+1}-\psi_k \psi_{k+1} \psi_k)\ei
                = Q_{i_k,i_{k+1},i_{k+2}}(y_k,y_{k+1},y_{k+2})\ei$
       \label{R:braid}
    \end{relations}
    for all $\bi \in I^\alpha$ and all admissible $k$ and $m$.
    The~\textbf{cyclotomic KLR algebra} is the quotient algebra
    \begin{equation}\label{R:cyc}
                \Rnaff[\alpha]^\Lambda=\Rnaff[\alpha]^\Lambda(\Qb,\Wb)
                     =\Rnaff[\alpha]/\CycIdeal,
    \end{equation}
    where $\CycIdeal$ is the two-sided ideal of~$\Rnaff[\alpha]$
    generated by $\set{W_{i_1} (y_1) \ei |\bi\in I^\alpha}$.

    Set $\Rnaff= \bigoplus_{\alpha\in Q_n^+} \Rnaff[\alpha]$ and
    $\Rn=\bigoplus_{\alpha\in Q_n^+} \Rnaff[\alpha]^\Lambda$.
  \end{Definition}

  We abuse notation and use $\ei$, $y_r$ and $\psi_r$ for both the
  generators of $\Rnaff[\alpha]$ and $\Rnaff[n]$ and for their images
  in~$\Rn[\alpha]$ and~$\Rn$. When we want to emphasise the base
  ring~$K$ we write $\Rnaff(K)=\Rnaff(\Qb,\Wb,K)$ and
  $\Rn(K)=\Rnaff(\Qb,\Wb,K)$.

  Importantly, the algebras $\Rnaff$ and $\Rn$ are graded $K$-algebras with degree
  function
  \[
      \deg\ei=0, \quad
      \deg y_m\ei=\kill{\alpha_{i_m}}{\alpha_{i_{m}}}=2\di,
      \quad\text{and}\quad
      \deg\psi_k\ei=-\kill{\alpha_{i_k}}{\alpha_{i_{k+1}}},
  \]
  for $\bi\in I^n$, $1\le k< n$ and $1\le m\le n$.
  \notation{$\deg$}{Degree function on $\Rn$, $\Rx$, graded rings, and tableaux}

  Inspecting the relations, there is a unique anti-isomorphism~$*$ of
  $\Rnaff$, and of $\Rn$, that fixes each of the generators.
  If $M$ is a graded $\Rn$-module then the \emph{graded dual} of $M$ is
  \begin{equation}\label{E:Dual}
    M^\circledast = \HOM_{\Rn}(M,K),
  \end{equation}
  where the $\Rn$-action on $M^\circledast$ is given by $(af)(m)=f(a^*m)$,
  for $a\in\Rn$, $f\in M^\circledast$ and $m\in M$.
  \notation{$*$}{The unique anti-isomorphism of $\Rx$, or $\Rn$, that fixes each generator}
  \notation{$M^\circledast$}{Graded dual $M^\circledast=\HOM_A(M,K)$ of $M$}[E:Dual]

  We reserve the notation $\Rn$ for the cyclotomic KLR algebras that
  are defined using $Q$-polynomials such that
  $Q_{i,j}(u,v)\in K_0[u,v]$, such as the standard $Q$-polynomials
  given in \autoref{E:StandardQ}. For most of this paper we work
  with cyclotomic KLR algebras $\Rx$ that are defined using
  ``deformations'' of the standard $Q$-polynomials, such as those in
  \autoref{Ex:ContentSystem} below.

  \begin{Remark}
    There is an extensive literature for the cyclotomic KLR algebras of
    affine type~$A$. Almost all of these papers work with the quiver
    $\Aone[e-1]$. In particular, in characteristic~$p>0$ the group algebra
    of the symmetric group is isomorphic to a cyclotomic KLR algebra of
    type~$\Aone[p-1]$. As this paper simultaneously treats affine
    types~$A$ and~$C$, we have chosen our notation to be consistent with
    the literature in affine type $A$ and so that both quivers
    have the same vertex set $\set{0,1,\dots,e-1}$. This is why we work
    with quivers of types $\Aone$ and $\Cone$ even though a more natural
    notation would be to work with quivers of types $\Aone[e]$ and
    $\Cone[e]$.
  \end{Remark}

  When $K$ is positively graded the algebras in this paper fit into the
  general framework developed by Kang and Kashiwara in
  \cite{KangKashiwara:CatCycKLR}. In particular,
  \cite{KangKashiwara:CatCycKLR} proves the following result using an
  intricate induction on~$n$.

  \begin{Proposition}[{%
    Kang-Kashiwara~\cite[Theorem~4.5]{KangKashiwara:CatCycKLR}}]
    \label{P:KLRfree}
    Suppose that $K$ is a positively graded ring. Then $\Rx(K)$ is free as
    a $K$-module.
  \end{Proposition}

  \begin{proof}
    By \cite[Theorem~4.5]{KangKashiwara:CatCycKLR}, $\Rx(K)$ is
    projective as an $\Rx[n-1](K)$-module, which implies that $\Rx(K)$ is
    projective as an $\Rx[0](K)$-module. This gives the result since
    $\Rx[0](K)\cong K$.
  \end{proof}


  A cyclotomic KLR datum $(\Gamma, \Qb, \Wb)$ is \emph{standard} if
  $Q_{i,j}(u,v), W_i(u)\in K_0[u,v]$, for all $i,j\in I$. A (cyclotomic) KLR
  algebra is \emph{standard} if its cyclotomic KLR datum is standard.  Many
  papers in the literature define KLR algebras over positively graded
  rings $K=\bigoplus_{d\ge0}K_d$ but in almost all cases they only
  consider standard $Q$-polynomials, like those in
  \autoref{E:StandardQ}. Non-standard $Q$-polynomials, such as those
  in \autoref{Ex:ContentSystem} below, play an important role in this paper.

  Let $\k$ be a commutative integral domain with~$1$. Let $\K$ be a
  field that is a $\k$-algebra. (Often, $\K$ will be the field of
  fractions of~$\k$.) Let~$\ux$ be a (possibly empty) tuple of
  indeterminates over~$\k$. In this and later sections, we work over the
  polynomial ring~$\kx=\k[\ux]$ and the Laurent polynomial
  ring~$\Kxx=\K[\ux,\ux^{-1}]$ with indeterminates~$\ux$. We
  consider~$\kx$ as a positively graded ring, and~$\Kxx$ as a $\Z$-graded
  ring, with the indeterminates in~$\ux$ all having degree~$1$; compare
  \autoref{R:GradedIrreducible}.

  \notation{$\Qbx$}{Family $(Q^{\ux}_{i,j}(u,v))_{i,j\in I}$ of deformed $Q$-polynomials
                    defining $\Rx$}
  \notation{$\Wbx$}{Family $(W^\ux_i(u))_{i\in I}$ of deformed weight polynomials
                    defining $\Rx$}

  Fix a standard family of standard $Q$-polynomials $\Qb$ together with
  a family of standard weight polynomials $\Wb$, both with coefficients
  in~$\k$.  Let $\Rn(\k)=\Rn(\Qb,\Wb,\k)$ be the corresponding
  cyclotomic KLR algebra over~$\k$. An \emph{$\kx$-deformation} of
  $(\Gamma,\Qb,\Wb)$ is a cyclotomic KLR datum $(\Gamma,\Qbx,\Wbx)$ such
  that $\Qbx=\bigl(Q^\ux_{i,j}(u,v)\bigr)_{i,j\in I}$ is a family of
  $Q$-polynomials with coefficients in~$\kx$ and
  $\Wbx=\bigl(W^\ux_i(u)\bigr)_{i\in I}$ is a family of weight
  polynomials such that the polynomials in $\Qb$ and $\Wb$ are the
  degree zero terms of the polynomials in~$\Qbx$ and $\Wbx$,
  respectively. That is, $\Qb={\Qbx}_{|\ux=0}$ and
  $\Wb={\Wbx}_{|\ux=0}$. (Here, and below, if $f(\ux)\in\kx$ then
  $f(\ux)_{|\ux=0}$ is the constant term of $f(\ux)$.)

  \begin{Notation}\label{N:WbxQbx}
    Suppose that $(\Gamma,\Qbx,\Wbx)$ is a $\kx$-deformation of $(\Gamma,\Qb,\Wb)$.
    Let
    \[
      \Rx(\kx)=\Rn(\Qbx,\Wbx,\kx) \qquad\text{and}\qquad
      \Rx(\Kxx)=\Rn(\Qbx,\Wbx,\Kxx)
    \]
    be the corresponding cyclotomic KLR algebras over $\kx$ and $\Kxx$,
    respectively.
  \end{Notation}
  \notation{$\Rx$}{Deformed cyclotomic KLR algebra determined by $(\Gamma,\Qbx,\Wbx)$}[N:WbxQbx]
  \notation{$\Rx[\alpha]$}{Block of cyclotomic KLR algebra $\Rx$}[N:WbxQbx]

  The $\kx$-deformations $(\Gamma,\Qbx,\Wbx)$ used in this paper are
  part of the data of a \textit{content system}, which is the subject of
  the next section. Non-trivial examples of the polynomials $\Qbx$ and
  $\Wbx$ are given in \autoref{Ex:ContentSystem} below. We will
  sometimes use the deformed KLR algebras
  $\Rxaff(\kx)=\Rxaff(\Qbx,\kx)$ and $\Rxaff(\Kxx)=\Rxaff(\Wbx,\Kxx)$
  determined by the polynomials $\Qbx$.  Let $\Qx_{ijk}(u,v,w)$ be the
  analogue of the three variable $Q$-polynomials in \autoref{E:Qijk}
  determined by $(\Gamma,\Qbx,\Wbx)$.

  As before, let $\Rn(\k)=\Rn(\Qb,\Wb,\k)$ be the standard cyclotomic
  KLR algebra determined by $(\Gamma,\Qb,\Wb)$. By specialising the
  indeterminates in $\ux$ to zero, the relations of $\Rx(\kx)_{|\ux=0}$
  coincide with those of the algebra $\Rn(\k)$, so we have the following
  trivial but useful observation

  \begin{Proposition}\label{P:Specialisation}
    Suppose that $(\Gamma,\Qbx,\Wbx)$ is a $\kx$-deformation of
    $(\Gamma,\Qb,\Wb)$. Consider~$\k$ as a graded $\kx$-module by
    letting $\ux$ act as zero. Then
    $\Rn(\k)\cong\Rx(\k)=\k\otimes\Rx(\kx)$ as graded algebras.
  \end{Proposition}

  That is, the standard cyclotomic KLR algebra $\Rn(\k)$ is isomorphic,
  as a graded algebra, to the specialisation of $\Rx(\kx)$ at $\ux=0$.
  Equivalently, $\Rn(\k)$ is the degree zero component, with respect
  to the~$\ux$-grading, of the algebra $\Rx(\kx)$. Note also that
  $\Rx(\kx)$ is free as a $\kx$-module by \autoref{P:KLRfree}.

  It turns out that the representation theories of the algebras
  $\Rn(\k)$ and~$\Rx(\kx)$ are very similar, with the theory
  for~$\Rx(\kx)$ being slightly richer. In contrast, under the
  assumptions introduced below, the algebra~$\Rx(\Kxx)$ is semisimple,
  which makes it a useful tool for studying the algebras~$\Rx(\kx)$
  and~$\Rx(\k)\cong\Rn(\k)$. Note that $\Rx(\kx)$ embeds into~$\Rx(\Kxx)$
  by \autoref{P:KLRfree}.

  \subsection{Bases of KLR algebras}\label{S:KLRproperties}
  For each $w\in \Sym_n$, fix a \emph{preferred reduced expression}
  $w=\sigma_{a_1} \dots \sigma_{a_l}$ and define $\psi_w= \psi_{a_1} \dots
  \psi_{a_l}$.  In general, $\psi_w$ depends on the choice of the
  preferred reduced expression for $w$.
  \notation{$\psi_w$}{Element of $\Rx$ or $\Rn$ defined by a fixed reduced expression for $w\in\Sym_n$}

  \begin{Theorem}[{Khovanov-Lauda~\cite[Theorem 2.5]{KhovLaud:diagI},
                   Rouquier~\cite[Theorem 3.7]{Rouq:2KM}}]\label{T:KLRBasis}
    The algebra $\Rxaff$ is free as a $K$-algebra with basis
    $\set{ \psi_w y_1^{m_1} \dots y_n^{m_n} \ei |w\in \Sym_n,
                          \, m_1,\dots,m_n\in \N, \, \bi \in I^n }$.
  \end{Theorem}

  Given $1\le k<n$,  define the \emph{divided difference operator}
  \[
    \partial_{k}\map{K[u_1,\dots,u_n]}{K[u_1,\dots,u_n]};
    f\mapsto \frac{f-{}^{\sigma_k}f}{u_k-u_{k+1}}.
  \]
  The next result follows easily from the relations in \autoref{D:KLR}.

  \begin{Lemma}[{Kang-Kashiwara~\cite[Lemma 4.2]{KangKashiwara:CatCycKLR}}]
    \label{L:KKDividedDifference}
  Let $V$ be an $\Rxaff[n]$-module and $f\in K[u_1,\dots,u_n]$ such that
  $f(y_1,\dots,y_n) \ei V=0$, for  $\bi \in I^n$. Suppose that
  $i_{k} = i_{k+1}$, for some $1\le k<n$. Then
  \[
    (\sigma_{k} f)(y_1,\dots,y_n) \ei V=0
    \And*
    (\partial_{k} f)(y_1,\dots,y_n) \ei V=0.
  \]
  \end{Lemma}

  \begin{Lemma}\label{L:DividedDiff}
  Let $f=(u_1-a_1)\dots (u_1-a_t) \in K[u_1,u_2]$, for
  $a_1,\dots,a_t\in K$. Then
  \[
      (\partial_{1} f)(a_1,u) = (u-a_2) \dots (u-a_t).
  \]
  \end{Lemma}

  \begin{proof} This follows easily by induction on $t$ using the
    general identity
    $\partial_{k} (fg) = (^{\sigma_k}f)\partial_{k}g +(\partial_{k}f) g$.
  \end{proof}

  Following~\cite[(1.6)]{KangKashiwaraKim:SymmI}, if $1\le r<n$, define
  $\phi_r=\sum_{\bi\in I^n}\phi_r\ei\in\Raff_n$ by
  \begin{equation}\label{E:phi}
    \phi_r \ei = \begin{cases*}
          \bigl(\psi_r (y_r-y_{r+1})+1\bigr) \ei& if $i_r = i_{r+1}$, \\
          \psi_r \ei& if $i_r \ne i_{r+1}$.
          \end{cases*}
  \end{equation}
  By definition, $\phi_r\ei$ is homogeneous and $\deg\phi_r\ei\ge0$.  If
  $w=\sigma_{a_1} \dots \sigma_{a_m}$ is a reduced expression for $w\in\Sym_d$
  define $\phi_w=\phi_{a_1} \cdots \phi_{a_m}$.  Parts (b) and (c) of
  the next lemma show that $\phi_w$ does not depend on the choice
  of the reduced expression.
  \notation{$\phi_w$}{Element of $\Rx$ or $\Rn$ indexed by $w\in\Sym_n$}

  \begin{Lemma}[\protect{Kang, Kashiwara and
       Kim~\cite[Lemma 1.5]{KangKashiwaraKim:SymmI}}]\label{L:KaKaKi}
    The following identities hold:
    \begin{enumerate}
      \item If $1\le r<n$, then $\phi_r^2 \ei =
      \bigl(Q_{i_r, i_{r+1}} (y_r, y_{r+1}) + \delta_{i_r,i_{r+1}}\bigr) \ei$.
      \item\label{KaKaKi2}
        If $1\le r<n-1$, then
        $\phi_r \phi_{r+1} \phi_r = \phi_{r+1} \phi_r \phi_{r+1}$.
      \item If $|r-s|>1$, then $\phi_r \phi_s = \phi_s \phi_r$.
      \item\label{KaKaKi3}
      If $w\in\Sym_n$ and $1\le t\le n$, then $\phi_w y_t = y_{w(t)} \phi_w$.
      \item\label{KaKaKi4} If $1\le k<n$ and $w(k+1) = w(k)+1$,
      then $\phi_w \psi_k = \psi_{w(k)} \phi_w$.
      \item\label{KaKaKi5}
        If $w\in\Sym_n$, then
        $\phi_{w^{-1}} \phi_w \ei
           =\displaystyle \prod_{\substack{1\le a<b\le n\\ w(a)>w(b)}}
             \bigl(Q_{i_a, i_b} (y_a, y_b) + \delta_{i_a,i_b}\bigr) \ei$.
    \end{enumerate}
  \end{Lemma}

  \section{Content systems for KLR algebras}\label{S:Contents}

  This chapter introduces \textit{content systems}, which are the basic
  combinatorial tool underpinning this paper. Using content systems, we
  will give analogues of Young's seminormal forms for cyclotomic KLR
  algebras of types $\Aone$ and $\Cone$, which are then used to prove
  the main results of this paper.

  \subsection{Content systems}\label{SS:Contents}

  As in \autoref{S:KLR}, in this chapter we let $\k$ be a commutative
  ring with~$1$ and fix a family of indeterminates $\ux$ and work over
  the rings $\kx$. In this chapter, $\K$ is the field of fractions
  of~$\k$ and we will mainly work over $\Kxx$. Let $(\Gamma,\Qbx,\Wbx)$
  be a $\kx$-deformation of the standard cyclotomic KLR datum
  $(\Gamma,\Qb,\Wb)$. This chapter studies the algebras $\Rx(\kx)$ and
  $\Rx(\Kxx)$ under the additional assumption that they come equipped
  with a \textit{content system}, which is the subject of this section.

  As in \autoref{S:KLR}, the cyclotomic KLR datum $(\Gamma,\Qbx,\Wbx)$
  determines  a dominant weight~$\Lambda=\Lambda_{\Wbx}\in P^+$ of
  level~$\ell$.  Fix an $\ell$-tuple
  $\charge=(\rho_1,\dots,\rho_\ell)\in I^\ell$, the
  \textbf{$\ell$-charge}, such that $\Lambda=\sum_{l=1}^{\ell}
  \Lambda_{\rho_l}$.

  Let $\Gamma_\ell$ be the quiver of type
  $A_\infty^{\times\ell}=A_\infty\times\dots\times A_\infty$,
  with~$\ell$ factors. More explicitly, $\Gamma_\ell$ has vertex set
  $\Jell=\set{1,2,\dots,\ell}\times\Z$ and edges $(l,a)\longrightarrow
  (l,a+1)$, for all $(l,a)\in\Jell$.  Given $(k,a),(l,b)\in \Jell$,
  write $(k,a)\arrow(l,b)$ if $(k,a)\ne(l,b)$ and there is an arrow between
  $(k,a)$ and $(l,b)$, in either direction. Similarly, write
  $(k,a)\noarrow(l,b)$ if $(k,a)\ne(l,b)$ and there are no arrows
  between $(k,a)$ and $(l,b)$. By definition, if $k\ne l$ then
  $(k,a)\noarrow(l,b)$.

  \begin{Definition}\label{D:ContentSystem}
  A \textbf{content system} for $\Rx(\kx)$ with values in $\kx$
  is a pair of maps $(\bc,\br)$, with
  \[
      \bc\map{\Jell}\kx \quad\text{and\quad} \br\map{\Jell}I,
  \]
  such that:
  \begin{enumerate}
    \item\label{CA:Initial}
         If $1\le l\le \ell$ then $\br(l,0)=\rho_l$. Moreover,
         if $i\in I$ then
         $W^\ux_i(u)=\displaystyle\prod_{l\in[1,\ell], \rho_l=i}
                            \bigl(u-\bc(l,0)\bigr)$.
    \item \label{CA:Qpoly} If $(k,a)\in\Jell$ and
    $j\in\set{\br(k,a-1),\br(k,a+1)}$ then there exists a unit
    $\epsilon=\epsilon_{k,a,j}\in\k^\times$ such that
    \[
      Q^\ux_{\br(k,a), j}\bigl(\bc(k,a), v\bigr)
          =\epsilon\prod_{\substack{b\in \set{a-1,a+1}\\\br(k,b)=j}}
              \bigl(\bc(k,b)-v\bigr).
    \]
    \item\label{CA:Unique}
    If $(k,a),(l,b)\in\Jell$ with $-n<a,b<n$ then $\br(k,a)=\br(l,b)$ and
    $\bc(k,a)=\bc(l,b)$ if and only if~$(k,a)=(l,b)$.
  \end{enumerate}
  The function~$\bc$ is the \emph{content function} of the content system
  and~$\br$ is the \emph{residue function}.  A content system $(\bc,\br)$
  is \textbf{graded} if $\bc(k,a)$ is homogeneous of degree
  $\kill{\alpha_i}{\alpha_i}=2\di$, where $i=\br(k,a)\in I$ for
  $(k,a)\in\Jell$.
  \end{Definition}
  \notation{$(\bc,\br)$}{A content system for $\Rx$}[D:ContentSystem]

  Almost all of the content systems that we consider will be graded.
  Even though content systems are defined using a quiver of type
  $\Gamma_\ell$, the quiver~$\Gamma$ is not assumed to be of this type.
  Notice that the roots of the polynomials~$W^\ux_i(u)$ are pairwise
  distinct by condition~(a) and~(c) of \autoref{D:ContentSystem}.

  By definition, a content system $(\bc,\br)$ depends on the choices of
  $K=\k[\ux]$, $\Gamma$, $\Qbx$, $\Wbx$, $\charge$ and~$n$. To define a
  content system we need to specify all of this data. As we will see,
  content systems are closely related to semisimple representations. In
  particular, the theory below implies that content systems do not exist
  for most choices of (standard) $Q$-polynomials or over fields of
  positive characteristic. As we explain in \autoref{T:SSemisimple}
  below, if a content system exists then the algebra $\Rx(\Kxx)$ is
  uniquely determined up to non-homogeneous isomorphism.  On the other
  hand, the examples below show that by deforming the standard
  $Q$-polynomials we can always find content systems for any standard
  cyclotomic KLR algebra $\Rn$ of type $\Aone$ or type $\Cone$.

  In the examples below, we give the minimum information necessary to
  specify the $Q$-polynomials. Recall from \autoref{E:Qpolynomials} that
  $Q^\ux_{i,j}(u,v)=Q^\ux_{j,i}(v,u)$, $Q^\ux_{i,i}(u,v)=0$ and that
  $Q^\ux_{i,j}(u,v)=1$ if~$i$ and~$j$ are not connected in~$\Gamma$, so
  we only need to specify one of the polynomials $Q^\ux_{i,j}(u,v)$ and
  $Q^\ux_{j,i}(v,u)$ whenever~$i$ and~$j$ are connected in~$\Gamma$.

  \begin{Example}\label{Ex:ContentSystem}
    The content systems below are completely new, so the use of the
    adjectives \textit{classical} and \textit{reduced} is purely
    descriptive.  For parts (a)--(e), we allow $n\ge0$ to be arbitrary
    and we take $K=\Z[\ux]=\Z[x]$, where $\ux=(x)$ and $x$ is an
    indeterminate of degree~$1$ over~$\Z$. For the examples of
    level~$\ell=1$ we identify~$\Jell$ with $\Z$ via the obvious map
    $(1,a)\mapsto a$ and set $\charge=(0)$. Throughout we use the weight
    polynomials $\Wbx=\bigl(W_i(u)\bigr)$, where
    $W^\ux_i(u)=\prod_{l\in[1,\ell], \rho_l=i}
                            \bigl(u-\bc(l,0)\bigr)$
    in accordance with \autoref{D:ContentSystem}(a).
    If $a,b\in\Z$ with $b\ne0$ let $\floor ab$ be the
    integer part of $\frac ab$ and set $\overline a=a\pmod{e}\in I$.
    \begin{enumerate}
    \item (The quiver $\Gamma_\ell$) Let $\Gamma=\Gamma_\ell$, the
    quiver of type $A_\infty^{\times\ell}$, and let
    $\rho=\bigl((1,0),\dots,(\ell,0)\bigr)$.  Let $\Qbx=\Qb$ be the
    standard $Q$-polynomials for $\Gamma_\ell$ given by
    \autoref{E:StandardQ}. Let
    $\br^\Jell$ be the identity map on $\Jell$ and define $\bc^\Jell$
    to be identically zero. Then $(\br^\Jell,\bc^\Jell)$ is a content
    system for $\Rn=\Rx$, where
    $\Lambda=\Lambda_{(1,0)}+\dots+\Lambda_{(\ell,0)}$.

    \item (Classical contents) Let $\Gamma$ be a quiver a type $\Aone$.
    Define
    \[
          Q^\ux_{i,j}(u,v) = \begin{cases*}
                 (v-u+x^2)(u+x^2-v)   & if $i\leftrightarrows j$,\\
                 (u+x^2-v)   & if $i\rightarrow j$,\\
           \end{cases*}
    \]
    for $i,j\in I=\set{0,1,\dots,e-1}$. Then $\Lambda=\Lambda_0$ and
    $\ell=1$. Then a content system for $\Rx$ is given by the functions
    $\bc(a) = ax^2$ and $\br(a) =\overline a$, for $a\in\Z$. More explicitly,
    $(\bc,\br)$ is given by the table:
    \[
      \begin{array}{l|cCCCCccccCCCCcc}
         a& -1 & 0    & 1       &\ldots
           & e{-}1 & e &  \ldots
           & 2e{-}1 & 2e & \ldots
           & 3e{-}1 & \ldots
           \\\hline
        \br(a) & e{-}1 & 0 & 1 & \ldots& e{-}1 & 0 & \ldots
            & e{-}1 & 0 & \ldots& e{-}1 & \ldots\\
        \bc(a) & -x^2 & 0 & x^2 & \ldots&(e{-}1)x^2 & ex^2 &
           \ldots & 2e x^2 & (2e{+}1)x^2& \ldots&(3e{-}1)x^2 & \ldots
      \end{array}
    \]
    Here, and below, the shading in the table highlights how the content
    function depends on $e=|I|$. The residue function~$\br$ is the
    standard residue function for type~$\Aone$. We call this a
    \textit{classical} content system because we recover the content
    function used in the classical semisimple representation theory of
    the symmetric groups by setting $x=1$. For more details, see
    \autoref{Ex:TableauContents}.

    To verify this example, and the examples that follow, observe that
    if $e>2$ and $\br(a)=i$ and $\bc(a)=cx$ then
    $(c+1)x-v=Q^\ux_{i,i+1}(\bc(a),v)=\epsilon\bigl(\bc(a+1)-v\bigr)$ by
    \autoref{D:ContentSystem}(c), so we require $\bc(a+1)=(c+1)x$ (and
    $\epsilon=+1$). The calculation when $e=2$ is similar except that we
    also need to inductively assume that $\bc(a-1)=(c-1)x$. In this way,
    the content function $\bc$ is completely determined by the
    $Q^x$-polynomials and the ``initial condition'' given by the weight
    polynomial~$W^x_0(u)=u-\bc(0)=u$.

    There is a related content system $(\bc',\br')$ that is, in a
    certain sense, dual to $(\bc,\br)$, which is given by
    $\bc'(a)=\bc(-a)$ and $\br'(a)=\br(-a)$, for $a\in\Z$.  This is a
    special case of a general construction given in
    \autoref{S:SignAutomorphism}, so similar remarks apply to every
    example below.

    \item (Reduced contents) Let $\Gamma$ be a quiver a type $\Aone$. Define
     \[
                Q^\ux_{i,j}(u,v) = \begin{cases*}
                   (u-v)(v+x^2-u) & if $e=2$ and $(i,j)=(0,1)$,\\
                   (u-v-x^2)      & if $e>2$ and $(i,j)=(0,e)$,\\
                   (u-v)        & if $i\rightarrow j\ne e$,\\
             \end{cases*}
      \]
      for $i,j\in I$. As in the last example, $\Lambda=\Lambda_0$
      and $\ell=1$.  Then a content system $(\bc,\br)$ for $\Rx$ is
      given by the functions $ \br(a) =\overline a$ and
      $\bc(a) =\floor{a}{e}x^2$, for all $a\in\Z$.  More explicitly, $(\bc,\br)$
      is given by the table:
      \[
        \begin{array}{l|cCCCCccccCCCCcc}
           a& -1 & 0    & 1       &\ldots
             & e-1 & e & e+1  & \ldots
             & 2e-1 & 2e & 2e+1  & \ldots
             & 3e-1 & 3e &  \ldots
             \\\hline
          \br(a) & e-1 & 0 & 1 & \ldots& e-1 & 0 & 1 & \ldots
               & e-1 & 0 & 1 & \ldots& e-1 & 0 & \ldots\\
          \bc(a) & -x^2 & 0 & 0 & \ldots& 0 & x^2 & x^2&  \ldots & x^2
                & 2x^2 & 2x^2& \ldots& 2x^2 & 3x^2& \ldots
        \end{array}
      \]

    \item (Classical contents) Let $\Gamma$ be a quiver a type $\Cone$.
    Define
    \[
      Q^\ux_{i,j}(u,v) = \begin{cases*}
         u-(v-x^2)^2& if $i=0\Rightarrow1=j$,\\
         (u+x^2)^2-v& if $i=e-1\Leftarrow e=j$,\\
         (u-v+x^2)   & if $i\rightarrow j$,\\
       \end{cases*}
    \]
    for $i,j\in I$.  As in the last example, $\Lambda=\Lambda_0$ and
    $\ell=1$. For an integer $a$ set $a'=\floor{a}{e-1}$ and
    let~$\doverline{a}$ be the unique integer such that
    $a\equiv\doverline{a}\pmod{2(e-1)}$ and $0\le \doverline a<2e-1$. A
    content system $(\bc,\br)$ for~$\Rx$ is given by the functions
    \[
       \bc(a) = \begin{cases*}
         (a+1)^2x^4& if $\overline a=0$,\\
         (-1)^{a'}(a+1)x^2 & if $\overline a>0$\\
       \end{cases*}
       \And*
       \br(a) = \begin{cases*}
         \overline a   & if $\doverline{a}<e$,\\
         \overline{-\doverline{a}-2} & otherwise,
       \end{cases*}
    \]
    for $a\in\Z$.  More explicitly, $(\bc,\br)$ is given by the table:
    \[
      \begin{array}{l|CcCCCcCCCcCCc}
         a& -1 & 0    & 1       &\ldots
           & e-2 & e-1 & e  & \ldots
           & 2e-3 & 2e-2 &  2e-1 &\ldots
           \\\hline
        \br(a) & 1 & 0 & 1 & \ldots& e-2 & e-1 & e-2 & \ldots
            & 1 & 0 & 1 &\ldots \\
            \bc(a) & 0x^2 & 1^2x^4 & 2x^2 & \ldots&(e{-}1)x^2 &e^2 x^4
           &-(e{+}1)x^2&
            \ldots &-(2e{-}2)x^2 & (2e{-}1)^2x^4 &2ex^2&\ldots
      \end{array}
    \]
    Notice that we cannot set $\bc(0)=0$ because this would force
    $c(-1)=x^2=\bc(1)$, which would violate \autoref{D:ContentSystem}(c).
    As we will see, the residue function~$\br$ is the type~$\Cone$
    residue function used by Ariki, Park and
    Speyer~\cite{ArikiParkSpeyer:C}. (Again, compare with
    \autoref{Ex:TableauContents}.)

    \item (Reduced contents) Let $\Gamma$ be a quiver a type $\Cone$.
    Define
    \[
      Q^\ux_{i,j}(u,v) = \begin{cases*}
         u-(v-x^2)^2& if $i=0\Rightarrow1=j$,\\
         (u+x^2)^2-v& if $i=e-2\Leftarrow e-1=j$,\\
         (u-v)   & if $i\rightarrow j$,\\
       \end{cases*}
    \]
    for $i,j\in I$.  As in the last example, $\Lambda=\Lambda_0$
    and $\ell=1$. A content system $(\bc,\br)$ for $\Rx$ is
    given by the functions
    \[
         \bc(a) = \begin{cases*}
           (2a'+1)^2x^4& if $\overline a=0$,\\
           (-1)^{a'}(2a'+2)x^2 & if $\overline a>0$\\
         \end{cases*}
         \And*
         \br(a) = \begin{cases*}
           \overline a & if $\doverline{a}<e$,\\
           \overline{-\doverline{a}-2} & otherwise,
         \end{cases*}
    \]
    for $a\in\Z$. More explicitly, $(\bc,\br)$ is given by the table:
    \[
      \begin{array}{l|CcCCCcCCCcCCCc}
         a& -1 & 0    & 1       &\ldots
           & e-2 & e-1 & e  & \ldots
           & 2e-3 & 2e-2 & 2e-1  & \ldots
           \\\hline
        \br(a) & 1 & 0 & 1 & \ldots& e-2 & e-1 & e-2 & \ldots
            & 1 & 0 &  1 &\ldots\\
        \bc(a) & 0x^2 & 1^2x^4 & 2x^2 & \ldots& 2x^2 & 3^2x^4 & -4x^2
           & \ldots & -4x^2 &5^2x^4 & 6x^2& \ldots
      \end{array}
    \]

    \item(Higher levels, many parameters) We extend the examples of
    content systems for level one algebras given in Examples (b)--(e) to
    algebras of level $\ell>1$.  Let $\Gamma$ be a quiver of type $\Aone$
    or $\Cone$, as above, and let $\Lambda\in P^+$ be a dominant weight
    with $\ell$-charge $\charge\in I^\ell$.  Fix a family of
    indeterminates $\ux=(x, x_1,\dots,x_\ell)$ over $\Z$ and set
    $K=\Z[\ux]$. Let $\Qbx$ be one of the families of $Q$-polynomials
    given in Examples (b)--(e) and let $(\br_0,\bc_0)$ be the
    corresponding level one content system for~$\Lambda=\Lambda_0$. A content
    system for the algebra~$\Rx$ is then given by setting
    $\br(k,a)=i=\br_0(\rho_k+a)\in I$ and $\bc(k,a)=\bc_0(\rho_k+a)+x_k^{2\di}$,
    for $(k,a)\in\Jell$.

    \item(Higher levels, one parameter) We can tweak the last example to
    give a content system that is defined over $\Z[x]$ for any
    $\ell\ge1$.  For example, in type~$\Aone$ to satisfy
    \autoref{D:ContentSystem}(c) we can fix integers
    $c_1>c_2+2n>\dots>c_\ell+2n$, and then specialise $x_k$ to $c_kx^2$ in
    example~(f), for $1\le k\le\ell$. For type~$\Cone$, we need
    $c_1>c_2+2n^2>\dots>c_\ell+2n^2$. More generally, if~$\k$ is a ``large enough''
    ring such that $2n\cdot1_\k\ne0$ then a higher level content system
    with values in $\k[x]$ is given by defining $\bc(k,a)=(c_k+a)x$, for
    suitable choices $c_1,\dots,c_\ell\in\k$ such that $c_k+a=c_l+b$
    only if $(k,a)=(l,b)$ for $-n<a,b<n$ and $1\le k,l\le\ell$. The
    content system in \autoref{Ex:ContentSystem}(d)--(f) extend to higher
    levels in essentially the same way except that extra care is
    required in choosing the ``initial contents'' $\bc(k,0)$, for $1\le
    k\le\ell$, to ensure that \autoref{D:ContentSystem}(c) is satisfied.
    We leave the details to the reader.

    \item(Non-graded content systems) In characteristic zero, the
    content systems given in Examples (a)--(f) are all graded content
    systems for any $n\ge0$.  By \autoref{P:Specialisation}, the
    standard cyclotomic KLR algebra $\Rn$ is isomorphic to the algebra
    $\Rx/\ux\Rx$ obtained by specialising all of the indeterminates
    at~$0$. We can obtain ungraded content systems for $\Rx$ over~$\Z$
    by specialising the indeterminates to a fixed prime $p$. Reducing
    modulo~$p$, it follows that the algebra $\Rx/p\Rx$ is isomorphic to
    the corresponding standard cyclotomic KLR algebra $\Rn(\Z/p\Z)$,
    defined over the finite field $\Z/p\Z$.

    \item(Finite type) It is possible to construct content systems for
    some quivers of finite type, such as type $A_e$, but we do not
    consider these here. The main difference is that in finite type the
    irreducible modules defined in \autoref{P:SeminormalForm} below
    exist only for certain $\ell$-partitions.
    \end{enumerate}
    \vspace*{-2.3\baselineskip}
  \end{Example}

  In particular, (b)--(e) and (g) of \autoref{Ex:ContentSystem} show the following:

  \begin{Lemma}\label{L:IntegralContentSystems}
    Let $\Gamma$ be a quiver of type $\Aone$ or $\Cone$ and
    suppose that $(\Gamma,\Qb,\Wb)$ is a standard cyclotomic KLR datum for
    $\Rn(\Z)$. Then there exists a $\Z[x]$-deformation
    $(\Gamma,\Qbx,\Wbx)$ of $(\Gamma,\Qb,\Wb)$ such that the algebra
    $\Rx=\Rx(\Qbx,\Wbx,\Z[x])$ has a content system $(\bc,\br)$ with values in
    $\Z[x]$.
  \end{Lemma}

  If $\k$ is a field of characteristic $p>0$ then the functions
  $(\bc,\br)$ from \autoref{Ex:ContentSystem}(b)--(h) define content
  systems only for ``small'' values of~$n$ because the uniqueness
  requirement of \autoref{D:ContentSystem}(c) fails whenever $n$ is too
  large.  For example, in characteristic~$2$ examples~(c) and~(d) define
  contents systems in type~$\Cone$ only when $n=1$. However, since
  content systems for cyclotomic KLR algebras of types $\Aone$ and
  $\Cone$ always exist over~$\Z[x]$ we can use content systems to
  construct cellular bases for these algebras by base change
  from~$\Z[x]$.

  \begin{Lemma}\label{L:AdjacentCS}
    Suppose that $(\bc,\br)$ is a content system and  $i=\br(l,a)$ and
    $j=\br(l,a+1)$, for $(l,a)\in\Jell$. Then $j\arrow i$ and, in
    particular, $i\ne j$. Moreover,
    $j=\br(l,a-1)$ if and only if $i\Longrightarrow j$ or
    $j\leftrightarrows i$.
  \end{Lemma}

  \begin{proof}
    By \autoref{D:ContentSystem}(b), $Q^\ux_{i,j}(\bc(k,a),v)$ is a nonzero
    polynomial in~$v$, so $i\ne j$ and $\kill{\alpha_i}{\alpha_j}\ne0$
    by \autoref{E:Qpolynomials}.
    Hence, $j\arrow i$. If, in addition, $\br(l,a-1)=j$ then
    $Q^\ux_{i,j}(\bc(k,a),v)$ is a polynomial of degree~$2$ in~$v$.
  \end{proof}

  \autoref{L:AdjacentCS} implies that if $(\bc,\br)$ is a content system
  for $\Rx$ and  $\Gamma$ is a quiver of type $\Aone$ and $1\le l\le
  \ell$ then either $\br(l,a)=\overline{\rho_l+a}$ or
  $\br(l,a)=\overline{\rho_l-a}$, for all $a\in\Z$. Similarly, if
  $\Gamma$ is of type $\Cone$ then $\br(l,a)=\br(\rho_l+a)$ or
  $\br(l,a)=\br(\rho_l-a)$, where $\br$ is the level one residue function
  used in~(c) and~(d) of \autoref{Ex:ContentSystem}. As sketched in
  example~(b) above, the content function is almost uniquely determined
  by the cyclotomic KLR datum $(\Gamma,\Qbx,\Wbx)$ because $\bc(l,0)$ is
  a root of the polynomial $W^\ux_{\br(l,0)}(u)$ and $\bc(l,a+1)$ is a
  root of the polynomial $Q^\ux_{i,j}(\bc(l,a),v)$, where $i=\br(l,a)$ and
  $j=\br(l,a+1)$. So, defining a content system $(\bc,\br)$ amounts to finding
  a $\kx$-deformation $(\Gamma,\Qbx,\Wbx)$ of the cyclotomic KLR
  datum.

\subsection{Tableau combinatorics}\label{S:Tableaux}
By \autoref{D:ContentSystem}, a content system $(\bc,\br)$ with values
in $\kx$, is just a pair of functions. This section extends these
functions to maps on $\ell$-partitions and standard tableaux, and the
next section uses this combinatorics to construct irreducible graded
representations of the deformed KLR algebra $\Rx$ over $\Kxx$. These
representations, which are modelled on Young's seminormal forms, are the
foundations that this paper are built on. We start by setting up the
required combinatorics.

A \textbf{partition} is a weakly decreasing sequence of positive
integers. If $\lambda=(\lambda_1,\dots,\lambda_r)$ is a partition, then
the \textbf{size} of $\lambda$ is $|\lambda|=\sum_{t=1}^r \lambda_t$,
and we set $\lambda_t=0$ for $t>r$.  An \textbf{$\ell$-partition} is an
ordered tuple $\blam= (\lambda^{(1)}|\dots|\lambda^{(\ell)})$ of
partitions. The \textbf{size} of $\blam$ is $|\blam|=\sum_{c=1}^{\ell}
|\lambda^{(c)}|$.  Let~$\Parts$ be the set of $\ell$-partitions of
size~$n$. We identify partitions and $1$-partitions in the obvious way.

If $\blam,\bmu\in\Parts$ then $\blam$ \emph{dominates} $\bmu$, written
$\blam\Gedom\bmu$, if
\[
      \sum_{c=1}^{k-1}|\lambda^{(c)}| + \sum_{r=1}^s\lambda^{(k)}_r
        \ge\sum_{c=1}^{k-1}|\mu^{(c)}| + \sum_{r=1}^s\mu^{(k)}_r,
        \qquad\text{ for $1\le k\le \ell$ and $s\ge1$.}
\]
Similarly, the \emph{reverse dominance order} $\Ledom$ is defined by
$\blam\Ledom\bmu$ if $\bmu\Gedom\blam$.  Write $\blam\Gdom\bmu$  and
$\bmu\Ldom\blam$ if $\blam\Gedom\bmu$ and $\blam\ne\bmu$.
\notation{$\Parts$}{The poset of $\ell$-partitions of $n$}
\notation{$\Ldom,\Gdom$}{Reverse dominance and dominance orders on $\Parts$}
\notation{$\Dom,\doM$}{Throughout, $\Dom\in\set{\Ldom,\Gdom}$ and $\set{\Dom,\doM}=\set{\Ldom,\Gdom}$}

In this paper, we consider the set of $\ell$-partitions $\Parts$ both as
the poset $(\Parts,\Gedom)$, under dominance, and as the poset
$(\Parts,\Ledom)$, under reverse dominance. As we will see, the interplay
between the dominance and reverse dominance partial orders corresponds
to a duality in the representation theory.

Let $\Nodes=\set{(k,r,c)|1\le k\le\ell\text{ and } r,c\in\Z_{>0}}$ be
the set of \emph{nodes}, which we consider as a totally ordered set
under the \emph{lexicographic order}~$\ge$. We also use the reverse
lexicographic order~$\le$. (We emphasize that our use of, and notation
for, the lexicographic and reverse lexicographic orders coincides with
how we use the dominance and reverse dominance orders.)
Identify an $\ell$-partition $\blam\in\Parts$ with its \emph{Young
diagram}, which is the set of nodes:
\[
   \blam=\set[\big]{(k,r,c)|1\le k\le \ell\text{ and } 1\le c\le \lambda^{(k)}_r}.
\]

\begin{Remark}
  In this paper the node $(k,r,c)\in\Nodes$ sits in component $k$,
  row~$r$ and column~$c$ of an $\ell$-partition. This is different to
  the conventions of \cite{DJM:cyc}, where the components of the nodes
  are indexed in order $(r,c,k)$. The convention used in this paper is
  preferable because many places in this paper order the nodes
  lexicographically, or reverse lexicographically, looking first at the
  component index and then at the row and column indices.
\end{Remark}

A \textbf{$\blam$-tableau} is a bijection
$\t\map{\blam}\set{1,2,\dots,n}$.  The group $\Sym_n$ naturally acts
from the left on the set of all $\blam$-tableaux. A $\blam$-tableau $\t$
is \textbf{standard} if $\t(k,r,c)<\t(k,r+1,c)$, and
$\t(k,r,c)<\t(k,r,c+1)$, whenever these nodes are in~$\blam$. That is,
the entries in each component of a standard tableau increase along rows and
down columns. Let $\Std(\blam)$ be the set of standard $\blam$-tableaux.
For $\Pcal\subseteq\bigcup_{n\ge0}\Parts$, set
\[
   \Std(\Pcal)= \set[\big]{\s|\s\in\Std(\blam)\text{ for }\blam\in\Pcal}
    \qquad\text{and}\qquad
   \Std^2(\Pcal)= \set[\big]{(\s,\t)|\s,\t\in\Std(\blam)\text{ for }\blam\in\Pcal}.
\]
Write $\Shape(\t)=\blam$ if $\t\in \Std(\blam)$.  Given
$\t\in \Std(\Parts)$ and $1\le m\le n$ let $\t_{\downarrow m}$ be the
subtableau of~$\t$ containing the numbers in~$\set{1,\dots,m}$. That is,
$\t_{\downarrow m}$ is the restriction of $\t$ to $\t^{-1} (\set{1,\dots,m})$.
\notation{$(k,r,c)$}{The node in component $k$, row $r$ and column $c$}
\notation{$\le,\ge$}{Lexicographic orders on the set of nodes $\set{(k,r,c)}$}
\notation{$\Std(\blam)$}{Standard tableau of shape $\blam\in\Parts$}
\notation{$\Std^2(\Pcal)$}{Pairs of standard tableaux $\bigcup_{\blam\in\Pcal}\Std(\blam)\times\Std(\blam)$, for $\Pcal\subseteq\Parts$}

Armed with this notation, we can now extend $(\bc,\br)$ to functions on
$\ell$-partitions and tableaux.

\begin{Definition}\label{D:Contents}
  Let $A=(k,r,c)\in\Nodes$ be a node. The
  \emph{content} of $A$ is $c(A)=\bc(k,c-r)\in\k[\ux]$ and the
  \emph{residue} of $A$ is $\br(A)=\br(k,c-r)\in I$. If $i\in I$, then
  $A$ is an \emph{$i$-node} if $\br(A)=i$.
\end{Definition}

Let $\t\in \Std(\blam)$ a standard $\blam$-tableau, for
$\blam\in\Parts$. Fix $1\le m\le n$. Define
\[
  \bc_m(\t)=\bc(\t^{-1}(m))\qquad\text{and}\qquad
  \br_m(\t)=\br(\t^{-1}(m)),
\]
which are the \emph{content} and \emph{residue} of $m$ in $\t$, respectively.
Similarly, the \emph{content sequence} and the \emph{residue sequence}
of~$\t$ are
\[
  \bc(\t)=\bigl(\bc_1(\t),\dots,\bc_n(\t)\bigr) \in \kx^n
  \quad\text{and}\quad
  \br(\t)=\bigl(\br_1(\t),\dots,\br_n(\t)\bigr)\in I^n,
\]
respectively. Let $\Std(\bi)=\set{\t\in\Std(\Parts)|\br(\t)=\bi}$ be the
set of standard tableaux with residue sequence~$\bi$.
\notation{$\Std(\bi)$}{Set of standard tableaux with residue sequence~$\bi$}

\begin{Example}\label{Ex:TableauContents}
    Suppose that $\ell=1$ and let $\blam=(5,3,2)$. Using the content
    systems from parts~(b)--(e) of \autoref{Ex:ContentSystem} for the quivers $\Aone[2]$
    and $\Cone[2]$, the different residues and contents in $\blam$ are:
    \[
    \begin{array}{cccc}\toprule
        \text{Quiver} & \text{\autoref{Ex:ContentSystem}} & \text{Contents} & \text{Residues}
        \\\midrule
        \Aone[2] & \text{(b)} & \Tableau[scale=0.77]{{0,x,2x,3x,4x},{-x,0,x},{-2x,-x}}
                 & \Tableau[scale=0.77]{{0,1,2,0,1},{2,0,1},{1,2}}\\[11mm]
        \Aone[2] & \text{(c)} & \Tableau[scale=0.77]{{0,0,0,x,x},{-x,0,0},{-x,-x}}
                 & \Tableau[scale=0.77]{{0,1,2,0,1},{2,0,1},{1,2}}\\[11mm]
        \Cone[2] & \text{(d) and (e) }& \Tableau[scale=0.77]{{x^2,2x,3^2x^2,4x,5^2x^2},{0,x^2,2x},{-2x^2,0}}
                 & \Tableau[scale=0.77]{{0,1,2,1,0},{1,0,1},{2,1}}\\[10mm]
        \bottomrule
      \end{array}
    \]
\end{Example}

The symmetric group $\Sym_n$ acts on~$I^n$ and~$\kx^n$ by
place permutations. Write $w\bc(\t)$ and $w\br(\t)$ for
the content and residue sequences obtained by acting with~$w$, for $w\in\Sym_n$.

From \autoref{S:Quiver}, recall that $\sigma_j=(j,j+1)\in\Sym_n$, for $1\le j<n$.

\notation{$\bc(k,r,c)$}{Content $\bc(k,c-r)$ of the node $(k,r,c)$}[D:Contents]
\notation{$\br(k,r,c)$}{Residue $\br(k,c-r)$ of the node $(k,r,c)$}[D:Contents]
\notation{$\bc(\t)$}{Content sequence $\bc(\t)=(\bc_1(\t),\dots,\bc_n(\t))$ of the tableau $\t$}
\notation{$\br(\t)$}{Residue sequence $\br(\t)=(\br_1(\t),\dots,\br_n(\t))$ of the tableau $\t$}

\begin{Lemma}\label{L:Separation}
   Suppose that $\s\in\Std(\blam)$ and $\t\in\Std(\bmu)$, for
   $\blam, \bmu\in\Parts$.
   \begin{enumerate}
      \item We have $\s=\t$ if and only if $\bc(\s) = \bc(\t)$ and
      $\br(\s)=\br(\t)$.
      \item\label{Swap} Suppose $\blam=\bmu$,  $\bc(\s)=\sigma_m\bc(\t)$ and
      $\br(\s)=\sigma_m\br(\t)$, for some $1\le m<n$. Then $\s=\sigma_m\t$.
   \end{enumerate}
\end{Lemma}

  \begin{proof}
  (a) If $\s\ne\t$ then let $m$ be minimal such that $\s_{\downarrow
  m}\ne\t_{\downarrow m}$. Set $\bmu=\Shape(\s_{\downarrow(m-1)})$ and
  let $A=(k,r,c)=\s^{-1}(m)$ and $B=(l,s,d)=\t^{-1}(m)$. Then~$A$
  and~$B$ are addable nodes of~$\bmu$.  If $k=l$ then it is well-known
  and easy to check that $c-r\ne d-s$. Consequently, $(k,c-r)\ne(l,d-s)$ and, hence,
  $\bigl(\bc_m(\s), \br_m(\s)\bigr)\ne\bigl(\bc_m(\t), \br_m(\t)\bigr)$
  by \autoref{D:ContentSystem}(c).
  Therefore, $(\bc(\s),\br(s))\ne(\bc(\t),\br(\t))$, giving~(a).

  Now consider (b). By assumption, $\bc(\sigma_m\s)=\bc(\t)$ and
  $\br(\sigma_m \s)=\br(\t)$, so $\sigma_m\s=\t$ by (a). Hence,
  $\s=\sigma_m\t$ as claimed.
  \end{proof}

  Part~(b) implies that if $\sigma_m\t\notin\Std(\Parts)$
  then no standard tableau has content sequence $\sigma_m\bc(\t)$ and
  residue sequence $\sigma_m\br(\t)$.

  Given $1\le m<n$ and $\t\in\Std(\bi)$, for $\bi\in I^n$, define
  scalars in~$\Kxx$ by
  \begin{equation}\label{E:Qmt}
       Q_m (\t)=
   Q^\ux_{\br_m(\t), \br_{m+1}(\t)}\bigl(\bc_m (\t), \bc_{m+1} (\t)\bigr)
   -\frac{\delta_{\br_m(\t),\br_{m+1}(\t)}}{\bigl(\bc_{m+1}(\t)-\bc_m(\t)\bigr)^2}.
  \end{equation}
  Note that $Q^\ux_{\br_m(\t), \br_{m+1}(\t)} (\bc_m (\t), \bc_{m+1}
  (\t))\in\kx$, so $Q_m(\t)\in\kx$ unless $\br_m(\t)=\br_{m+1}(\t)$.
  Further, if $\br_m(\t)=\br_{m+1}(\t)$ then $Q_m(\t)$ is well-defined
  because $\bc_m(\t)\ne\bc_{m+1}(\t)$ by \autoref{D:ContentSystem}(c) and
  \autoref{D:Contents}.
  \notation{$Q_m(\t)$}{$Q^\ux_{\br_m(\t), \br_{m+1}(\t)}(\bc_m (\t), \bc_{m+1} (\t))
        -\delta_{\br_m(\t),\br_{m+1}(\t)}/(\bc_{m+1}(\t)-\bc_m(\t))^2$}[E:Qmt]

  The following result looks innocuous but it is the key to constructing
  the seminormal representations of $\Rx(\Kxx)$.

  \begin{Lemma}\label{L:Qnonzero}
  Suppose that $\t \in \Std(\blam)$ and let $\s=\sigma_m \t$, where $1\le m<n$.
  Then $Q_m(\t) \ne 0$ if and only if~$\s\in\Std(\blam)$. Consequently,
  if $(\bc,\br)$ is a graded content system and $\s\in\Std(\blam)$
  then $Q_m(\t)$ is a nonzero homogeneous element of $\Kxx$.
  \end{Lemma}

  \begin{proof}
    For the duration of the proof set $(k,a,b)=\t^{-1}(m)$ and
    $(l,c,d) = \t^{-1} (m+1)$, so that $\bc_m(\t)=\bc(k,b-a)$,
    $\br_m(\t)=\br(k,b-a)$, $c_{m+1}(\t)=\bc(l,d-c)$ and
    $\br_{m+1}(\t)=\br(l,d-c)$

    Suppose first that $\s=\sigma_m\t\in\Std(\blam)$. If
    $\br_m(\t)=\br_{m+1}(\t)$ then $\bc_m(\t)\ne\bc_{m+1}(\t)$ by
    \autoref{L:Separation}, so that
    $Q_m(\t)=-1/(\bc_{m+1}(\t)-\bc_m(\t))^2\ne0$. Now suppose that
    $\br_m(\t)\ne \br_{m+1}(\t)$. By \autoref{E:Qmt}, $Q_m(\t)=0$ only
    if $\bc(l,d-c)$ is a root of
    $\Qx_{\br(k,a),\br(l,d-c)}(\bc(k,b-a),v)$. By axioms (b) and (c) of
    \autoref{D:ContentSystem}, $\bc(l,d-c)$ is not a root of
    $\Qx_{\br(k,b-a),\br(l,d-c)}(\bc(k,b-a),v)$ if $(k,a)\noarrow(l,c)$,
    so we can assume that $k=l$ and $d-c=b-a\pm1$ since otherwise
    $(k,a)\noarrow(l,c)$. However, if $d-c=b-a\pm1$ then $m$ and $m+1$
    are on adjacent diagonals in $\blam$, which is not possible since
    $\t$ and $\s=\sigma_m\t$ are both standard. Hence, $Q_m(\t)\ne0$
    when $\s$ is standard.

    Now, suppose that $\s\notin\Std(\blam)$. This happens if and
    only if~$m$ and~$m+1$ are in the same row or same column of the same
    component of~$\t$. That is, $k=l$ and either $a=c$ and $d=b+1$, or
    $b=d$ and $c=a+1$. That is, either $\br_{m+1}(\t)=\br(k, b-a+1)$ and
    $\bc_{m+1}(\t)=\bc(k,b-a+1)$, or $\br_{m+1}(\t)=\br(k, b-a-1)$ and
    $\bc_{m+1}(\t)=\bc(k,b-a-1)$. Hence, in both cases,
    $Q_m(\t)=Q^\ux_{\br_m(\t), \br_{m+1}(\t)}
        \bigl(\bc_m (\t)$, $\bc_{m+1}(\t)\bigr)=0$
    by \autoref{D:ContentSystem}(b).

    Finally, if $(\bc,\br)$ is a graded content system and
    $\s\in\Std(\blam)$ then $Q_m(\t)\ne0$, so it is homogeneous and
    nonzero in view of the remarks before the lemma.  Moreover,
    $Q_m(\t)$ has the expected degree by
    \autoref{E:Qpolynomials} since $\bc(k,a)$ is homogeneous of
    degree $\kill{\alpha_i}{\alpha_i}$ by \autoref{D:ContentSystem},
    where $i=\br(k,a)$.
  \end{proof}

\subsection{Seminormal forms}
  We continue to assume that $(\bc,\br)$ is a (graded) content system
  that takes values in $\kx$. Even though $(\bc,\br)$ takes values in
  $\kx$ the representations that we construct are modules for the
  $\Kxx$-algebra $\Rx(\Kxx)$ because the action of the KLR algebra on
  these modules involves the scalars $Q_m(\t)$ from \autoref{E:Qmt}, and
  these scalars typically belong to $\Kxx$, not $\kx$. To prove
  irreducibility we also use the following elements, which are not
  defined over $\kx$.

  \begin{Definition}\label{D:Ft}
    Let $\bi\in I^n$. If $\t \in \Std(\bi)$, define
    \[
       F_\t = \prod_{k=1}^n
              \prod_{\substack{\s\in\Std(\bi)\\\bc_k(\s)\ne \bc_k(\t)}}
                 \frac{y_k-\bc_k(\s)}{\bc_k(\t)-\bc_k(\s)}\cdot \ei
                                \in \Rx(\Kxx).
    \]
  \end{Definition}
  \notation{$F_\t$}
           {Semisimple idempotent in $\Rx(\Kxx)$, for $\t\in\Std(\Parts)$}
           [D:Ft]

  If $(\bc,\br)$ is a graded content system then $F_\t$ is homogeneous
  element of $\Rx(\Kxx)$ of degree~$0$ since~$\bc_k(\s)$ appears in the
  product only if $\br_k(\t)=\br_k(\s)$. Note that $\ei=\ei[\br(\t)]$, for
  $\t\in\Std(\bi)$.

  The next result gives a generalisation of Young's classical seminormal
  forms to KLR algebras with content systems.  As noted in
  \autoref{S:Rings}, $\Kxx$ is a graded field, which explains the claim
  that the module $V_\blam$ is an irreducible graded
  $\Rx(\Kxx)$-module. Recall that $\K$ is the field of fractions
  of~$\k$.

  \begin{Proposition}\label{P:SeminormalForm}
  Let $\blam \in \Parts$. Suppose that there exist scalars
  \[\set{\beta_k(\t)\in\Kxx|1\le k<n\text{ and }\t, \sigma_k\t\in\Std(\blam)}\]
  satisfying the following conditions:
  \begin{enumerate}
  \item $\beta_k (\sigma_k\t) \beta_k (\t) = Q_{k}(\t)$ if $1\le k<n$
        and $\sigma_k\t \in \Std(\blam)$;
  \item $\beta_k (\t) \beta_t (\sigma_k \t)
              = \beta_t (\t) \beta_k (\sigma_l \t)$
        if $1\le k,l<n$, $|k-l|\ne 1$ and
        $\sigma_k\t, \sigma_l \t \in \Std(\blam)$;
  \item $\beta_k (\sigma_{k+1} \sigma_k \t) \beta_{k+1} (\sigma_k \t) \beta_k (\t) =
  \beta_{k+1} (\sigma_k \sigma_{k+1} \t) \beta_k (\sigma_{k+1} \t)
    \beta_{k+1} (\t)$
  if $1\le k<n-1$ and all the tableaux appearing in this equation are
    standard.
  \end{enumerate}
  Then there exists a graded $\Rx(\Kxx)$-module $V_\blam$
  that is free as an $\Kxx$-module with homogeneous basis
  $\set{ v_\t | \t\in \Std(\blam)}$ and where $\Rx(\Kxx)$-action is determined by
  \[
    \ei v_\t = \delta_{\bi\,\br(\t)} v_\t,
        \quad y_k v_\t = \bc_k (\t) v_\t, \quad
    \psi_k v_\t  = \beta_k (\t) v_{\sigma_k \t}
    + \frac{\delta_{\br_k(\t),\br_{k+1}(\t)}}{\bc_{k+1}(\t)-\bc_k (\t)}v_\t
  \]
  for all admissible $k$, $\bi\in I^n$ and $\t\in\Std(\blam)$ and where
  $v_\s=0$ if $\s\notin\Std(\blam)$. Moreover, if $\Kxx$ is a graded
  field then $V_\blam$ is irreducible.
  \end{Proposition}

\begin{proof}
To prove that $V_\blam$ is an $\Rx(\Kxx)$-module it is enough to check that the
action of the generators of $\Rx(\Kxx)$ on $V_\blam$ respects the
relations of \autoref{D:KLR}. The action respects the cyclotomic relation
\[
            W_{i_1}(y_1) \ei =0,\qquad\text{for all } \bi\in I^n,
\]
by \autoref{D:ContentSystem}(a).  The
relations~\autoref{R:Ids}--\autoref{R:psicomm} and \autoref{R:psiy} are easily
checked by direct calculation, with condition~(b) of the proposition
used for~\autoref{R:psicomm} and relation \autoref{R:psi}
following by \autoref{L:Separation}(b).

To check relation~\autoref{R:quad}, for each $\t\in\Std(\blam)$ it is
enough to prove that
\begin{equation}\label{E:QuadRelation}
  \psi_k^2 \ei v_\t = Q^\ux_{i_k,i_{k+1}} (y_k,y_{k+1}) \ei v_\t,
  \qquad 1\le k<n\text{ and }\bi\in I^n.
\end{equation}
If $\sigma_k \t$ is not standard, then $\br_k(\t) \ne \br_{k+1}(\t)$ by
\autoref{L:Separation}(b) and $Q_k(\t)=0$ by
\autoref{L:Qnonzero}. So,
\[
  \psi_k^2 \ei v_\t = 0
    = \delta_{\bi\,\br(\t)} Q^\ux_{\br_k(\t),\br_{k+1}(\t)} (\bc_k (\t), \bc_{k+1} (\t)) v_\t
    = Q^\ux_{i_k,i_{k+1}} (y_k,y_{k+1}) \ei v_\t.
\]
On the other hand, if $\sigma_k\t$ is standard then
\begin{equation}\label{E:SeminormalQuad}
    \psi_k^2 \ei v_\t
         =  \Bigl(\beta_k (\sigma_k \t) \beta_k (\t)
              + \frac{\delta_{\br_k(\t),\br_{k+1}(\t)}}{(\bc_{k+1}(\t)-\bc_k(\t))^2} \Bigr) v_\t
         = Q^\ux_{\br_k(\t),\br_{k+1}(\t)}(y_k,y_{k+1}) \ei v_\t,
\end{equation}
where the second equality follows using condition (a) of the proposition
and the definition of $Q_k(\t)$. Hence,~\autoref{E:QuadRelation} holds in all
cases.

We now verify relation~\autoref{R:braid}. Let $\t\in \Std(\blam)$,
$1\le k<n-1$ and $\bi\in I^n$. To simplify notation, set $i=i_k$,
$i'=i_{k+1}$ and $i''=i_{k+2}$ and define
     $\t_1=\sigma_k\t$,
     $\t_2=\sigma_{k+1}\t$,
     $\t_{21}=\sigma_{k+1}\t_1$,
     $\t_{12}=\sigma_k\t_2$ and
     $\t_{121}=\sigma_k\t_{21}=\sigma_{k+1}\t_{12}$.
Note that if $\t_1\notin \Std(\blam)$, then $\t_{21}\notin \Std(\blam)$.
Similarly, $\t_{12}\notin \Std(\blam)$ if $\t_{2}\notin \Std(\blam)$ and
$\t_{121}\notin\Std(\blam)$ if either $\t_{12}\notin\Std(\blam)$ or
$\t_{21}\notin \Std(\blam)$. Using these facts and some routine,
although slightly lengthy calculations for the first equality
(cf.~\cite[Lemma~3.8]{HuMathas:SeminormalQuiver}),
shows that
\begin{equation*}
  \begin{multlined}
  (\psi_k \psi_{k+1} \psi_k - \psi_{k+1} \psi_k \psi_{k+1}) \ei v_\t\\
\quad\begin{alignedat}{2}
    &= \biggl( \delta_{ii'} \delta_{i'i''} \frac{  \bc_{k} (\t)
        + \bc_{k+2} (\t) - 2c_{k+1} (\t)}{(\bc_{k+1}(\t)
        -\bc_k(\t))^2 (\bc_{k+2}(\t)-\bc_{k+1}(\t))^2}
        +\delta_{ii''} \frac{\beta_k (\t) \beta_k (\t_1) - \beta_{k+1} (\t) \beta_{k+1} (\t_2)}
                          {\bc_{k+2} (\t) - \bc_k (\t)} \biggr) v_\t \\
    &\hspace{1cm} + \Bigl(\beta_k (\t_{21}) \beta_{k+1} (\t_1) \beta_k (\t) -
    \beta_{k+1} (\t_{12}) \beta_k (\t_2) \beta_{k+1} (\t) \Bigr) v_{\t_{121}} \\
    &=\delta_{ii''} \frac{Q_k (\t) - Q_{k+1} (\t)}{\bc_{k+2} (\t) - \bc_k (\t)} v_\t
     = \delta_{i'i''}\frac{Q^\ux_{ij} (y_{k+2},y_{k+1}) - Q^\ux_{ij} (y_k, y_{k+1})}
                          {y_k -y_{k+2}} \ei v_\t\\
    &=Q_{ii'i''}(y_k,y_{k+1},y_{k+2})\ei\v_\t
  \end{alignedat}
  \end{multlined}
\end{equation*}
where we have used conditions (a) and (c) of the proposition, and
\autoref{E:Qmt}, for the second equality. Hence,
relation~\autoref{R:braid} is satisfied. We have now shown that all of
the relations in \autoref{D:KLR} are satisfied, so $V_\blam$
is an $\Rx(\Kxx)$-module.

We next prove that $V_\blam$ is an irreducible graded
$\Rx(\Kxx)$-module when $\Kxx=\K[x^\pm]$ is a graded field. First note that
\begin{equation}\label{E:Fsvt}
    F_{\t} v_\s = \delta_{\t\s} v_\s,\qquad
    \text{ for all $\t,\s \in \Std(\Parts)$},
\end{equation}
by \autoref{D:Ft} and \autoref{L:Separation} since $v_\s$ is a
eigenvector for the $y_k$'s. Now suppose that $v\in
V_\blam$ belongs to a graded $\Rx(\Kxx)$-submodule~$M$ of~$V_\blam$ and
write $v=\sum_\s r_\s v_\s$, for $r_\s\in\Kxx$. If $r_\t\ne0$ then
$r_\t v_\t=F_\t v\in M$. Hence, $v_\t\in M$ since $M$ is a graded
submodule and~$\Kxx$ is a graded field. To show that $M=V_\blam$ it is
enough to show that $v_{\sigma_k\t} \in \Rx v_\t$ whenever $\t\in
\Std(\blam)$ and $\sigma_k \t\in \Std(\blam)$, for $1\le k<n$.  Under
these assumptions, $F_{\sigma_k\t} \psi_k v_\t = \beta_k (\t)
v_{\sigma_k\t}$. So it is enough to prove that $\beta_k (\t)\ne0$, which
follows from assumption~(a) since
$\beta_k (\t) \beta_{k} (\sigma_k\t) = Q_k (\t)$ and $Q_k(\t)\ne0$
by \autoref{L:Qnonzero}.

Finally, it remains to determine the grading on $V_\blam$. Since we have
already shown that the action of~$\Rx(\Kxx)$ on~$V_\blam$ respects the
relations and that $V_\blam$ is irreducible, and $\set{v_\s}$ is a
homogeneous basis, we can fix a grading on $V_\blam$ by fixing the
degree of one of these basis elements.  The degrees of the other basis
elements are now uniquely determined by the $\Rx(\Kxx)$-action since
$V_\blam$ is cyclic.
\end{proof}

\begin{Remark}
  Suppose that the content system $(\bc,\br)$ is not graded and takes
  values in $\k$. Then the argument of \autoref{P:SeminormalForm} shows
  that $V_\blam$ is an irreducible $\Rx(\K)$-module.
\end{Remark}

\autoref{P:SeminormalForm} constructs the module $V_\blam$ subject to
the existence of suitable scalars $\beta_k(\t)$, for $1\le k<n$ and
$\t\in\Std(\blam)$. There are two natural choices (see
\autoref{E:LGbetas}), but for now we define:
  \begin{equation}\label{E:betas}
      \beta_k(\t) = \begin{cases*}
                1&if $\sigma_k\t\Gdom \t$,\\
                Q_k (\sigma_k \t) &if  $\t\Gdom \sigma_k\t$.
            \end{cases*}
   \end{equation}

 \begin{Lemma}\label{L:betas}
   The coefficients $\beta_k(\t)$ defined by
   {\upshape\autoref{E:betas}\/} satisfy the conditions
   of \autoref{P:SeminormalForm}.
 \end{Lemma}

 \begin{proof}
   The only condition that is not obvious is that the
   $\beta$-coefficients satisfy the ``$\beta$-braid relation''
   \[
   \beta_k(\sigma_{k+1}\sigma_k\t)\beta_{k+1}(\sigma_k\t)\beta_k(\t)
              = \beta_{k+1}(\sigma_k\sigma_{k+1}\t)\beta_k(\sigma_{k+1}\t)
                  \beta_{k+1}(\t),
   \]
   for $\t\in\Std(\Parts)$ and $1\le r<d$ such that all the tableaux in
   this identity are standard. In fact, since $\beta_k(\t)$ depends only
   on the nodes $\t^{-1}(k)$ and $\t^{-1}(k+1)$, we have
   $\beta_k(\t)=\beta_{k+1}(\sigma_k\sigma_{k+1}\t)$,
   $\beta_{k+1}(\sigma_k\t)=\beta_k(\sigma_{k+1}\t)$ and
   $\beta_k(\sigma_{k+1}\sigma_k\t)=\beta_{k+1}(\t)$. These equalities
   imply the $\beta$-braid relation above.
 \end{proof}

 For each $\blam\in\Parts$ \autoref{P:SeminormalForm} constructs an
 irreducible $\Rx(\Kxx)$-module $V_\blam$.  We now fix the choice of
 $\beta$-coefficients given by \autoref{E:betas} and define $V_\blam$ to
 be the $\Rx(\Kxx)$-module defined by \autoref{P:SeminormalForm}.

 If $\t$ is a standard tableau then it is not clear from \autoref{D:Ft}
 that the element $F_\t$ is nonzero. This now follows by virtue of
 \autoref{E:Fsvt} and \autoref{L:betas}.

  \begin{Corollary}\label{C:Ftnonzero}
     Let $\t\in\Std(\blam)$, for $\blam\in\Parts$. Then $F_\t\ne 0$
     in~$\Rx(\Kxx)$.
  \end{Corollary}

  The next result shows that the representations constructed in
  \autoref{P:SeminormalForm} are pairwise non-isomorphic and, up to
  isomorphism, independent of the choice of $\beta$-coefficients in
  \autoref{P:SeminormalForm}.

  \begin{Corollary}\label{C:IsomorphicSpechts}
    Suppose that $\blam,\bmu\in\Parts$. Then  $V_\blam\cong V_\bmu$ as
    $\Rx(\Kxx)$-modules if and only if~$\blam=\bmu$.  Moreover, up to
    isomorphism, $V_\blam$ is  independent of the choice of homogeneous scalars
    $\set{\beta_k(\t)|\t\in\Std(\blam)}$ satisfying
    conditions {\upshape(a)--(c)} of \autoref{P:SeminormalForm}.
  \end{Corollary}

  \begin{proof}
    Suppose first that $\blam\ne\bmu$. By
    \autoref{L:Separation} and \autoref{E:Fsvt}, if $\t\in \Std(\blam)$
    then $F_\t V_\blam \ne 0$ and $F_{\t} V_\bmu=0$. Hence,
    $V_\blam\not\cong V_\bmu$.

    To prove the second statement suppose that $V_\blam\cong V_\bmu$ and
    that $V_\blam=\<v_\t|\t\in\Std(\blam)\>$ and
    $V_\blam'=\<v_\t'|\t\in\Std(\blam)\>$ are two $\Rx(\Kxx)$-modules with homogeneous
    structure constants $\set{\beta_r(\t)}$ and $\set{\beta_r'(\t)}$,
    respectively, satisfying the conditions of
    \autoref{P:SeminormalForm}. In particular, note that if
    $\sigma_r\t\in\Std(\blam)$ then $\beta_r(\t)$ and $\beta_r'(\t)$
    are both nonzero by \autoref{P:SeminormalForm}(a) and
    \autoref{L:Qnonzero}. Define a $\Kxx$-linear map $\theta\colon V_\blam\to
    V_\blam'$ inductively as follows. First, fix any tableau $\t_1\in\Std(\blam)$
    and set $\theta(v_{\t_1})=v_{\t_1}'$. By way of induction, suppose
    that $\theta(v_{\t_1}),\dots,\theta(v_{\t_{m-1}})$ have been defined
    and that $\t_{m}\in\Std(\blam)\setminus\set{\t_1,\dots,\t_{m-1}}$ is
    a standard tableau such that $\t_{m}=\sigma_k\t_l$, where
    $1\le k<n$ and $1\le l<m$. Set
    \[
        \theta(v_{\t_{m}})=\frac1{\beta_{k}(\t_l)}\Bigl(
           \psi_{k}-\tfrac{\delta_{\br_{k}(\t_m),\br_{k}(\t_l)}}
                            {c_{k}(\t_m)-c_{k}(\t_l)}
           \Bigr)\theta(v_{\t_l}).
    \]
    By \autoref{P:SeminormalForm}, if $\theta(v_{\t_l})\ne0$ then
    $\theta(v_{\t_{m}})\ne0$. By induction, $\theta(v_\t)$ is defined and
    nonzero for all $\t\in\Std(\blam)$. In particular,  $\theta$ is a $\Kxx$-module
    isomorphism. Moreover, $\theta (v_{\t})\in F_\t V_\blam'=\Kxx v_\t'$ by~\autoref{E:Fsvt},
    so~$\theta(v_\t)=\xi_\t v_\t'$, for some scalar
    $\xi_\t\in\Kxx$. Since $V_\blam$ and $V_\blam'$ are both $\Rx(\Kxx)$-modules,
    the construction of \autoref{P:SeminormalForm} guarantees that $\theta$ is an
    $\Rx(\Kxx)$-module homomorphism and that $V_\blam\cong V_\blam'$, as
    claimed.
  \end{proof}

  Motivated by the seminormal forms of \autoref{P:SeminormalForm}, we
  now use (graded) content systems to study the algebras $\Rx(\Kxx)$. Our
  next goal is to prove a semisimplicity result for~$\Rx(\Kxx)$, which we
  will use to study the algebras $\Rx(\kx)$ and $\Rn(\k)$.

  \subsection{Weight modules}
  This section looks at $\Rxaff(\Kxx)$=modules that are spanned by simultaneous
  eigenvectors of $y_1,\dots,y_n$. This is a first
  step towards finding a basis for $\Rx(\Kxx)$.


  Suppose that $V$ is an $\Rx(\Kxx$-module. Let
  $\bc=(c_1,\dots,c_n)\in\Kxx^n$ and $\bi\in I^n$, where $c_k$ is
  homogeneous of degree $(\alpha_{i_k}|\alpha_{i_k})$, for $1\le k\le n$. The
  \emph{$(\bc,\bi)$-weight space} of~$V$ is the $\Kxx$-module
  \[
      V_{\bc,\bi} = \set{v\in V|y_k\ei v= c_k v\text{ for }1\le k\le n}.
  \]
  A \emph{weight module} is an $\Rxaff(\Kxx)$-module that is a direct sum
  of $(\bc,\bi)$-weight spaces and is of finite rank as a
  $\Kxx$-module. For example, the module $V_\blam$ of
  \autoref{P:SeminormalForm} is an $\Rxaff(\Kxx)$-weight module.

  The next result is similar to the classification of the irreducible
  representations of the affine Hecke algebras of rank~$2$. The
  connection with the seminormal forms of \autoref{P:SeminormalForm} is
  evident in part~(b).

  \begin{Proposition}\label{P:R2affSimples}
    Let $V$ be a weight module for $\Rxaff[2](\Kxx)$ and
    suppose that $0\ne v\in V$ is a homogeneous vector such that $y_1v=c_1v$,
    $y_2v=c_2v$ and $\ei[ij]v=v$, where $c_1,c_2\in\Kxx$ and $i,j\in I$
    with~$c_1$ and~$c_2$ homogeneous of the appropriate degree.
    Then one of the following of the following mutually exclusive cases
    occurs:
    \begin{enumerate}
      \item If $\Qx_{ij}(c_1,c_2)\ne0$ then $\<v,w\>$ is an
      $\Rxaff[2](\Kxx)$-weight module of
      rank~$2$ such that $w=\psi_1 v$, $y_1w=c_2w$, $y_2w=c_1w$ and $\ei[ji]w=w$.
      \item If $i=j$ then $c_1\ne c_2$ and $V=\<v,w\>$ is an
       $\Rxaff[2](\Kxx)$-weight module of rank~$2$ such that
      $w=\bigl(\psi_1-\tfrac1{c_2-c_1}\bigr)v$, $y_1w=c_2w$, $y_2w=c_1w$
      and $\ei[ii]w=w$.
      \item If $i\ne j$ and $\Qx_{ij}(c_1,c_2)=0$ then either $V=\<v\>$
      is an $\Rxaff[2](\Kxx)$-weight module of
      rank~$1$ with $\psi_1v=0$, or $\<v,w\>$ is an
      $\Rxaff[2](\Kxx)$-weight module of rank~$2$ with
      $w=\psi_1v$ and $\psi_1w=0$.
    \end{enumerate}
  \end{Proposition}

  \begin{proof}

    As in the statement of the proposition, suppose that $v\in V$ and
    $\ei v=v$, $y_1v=c_1v$ and $y_2v=c_2v$. As in part~(a), we first
    assume that $\Qx_{ij}(c_1,c_1)\ne0$. Then $i\ne j$ since
    $\Qx_{ii}(u,v)=0$. Let $w=\psi_1 v$. Then
    $\psi_1 w=\Qx_{ij}(c_1,c_2)v\ne0$, so $w\ne0$. The remaining claims
    in~(a) now follow easily from the relations.

    Next, suppose that (b) holds, so that $i=j$. If $\psi_1v=0$ then
    $0=y_2\psi_1v=(\psi_1y_1+1)v =v$, which is a contradiction, so
    $\psi_1v\ne0$. By assumption, $V=\<v,\psi_1 v\>$ and $v$ is a weight
    vector, so $\psi_1v+av$ must be a weight vector for some $0\ne
    a\in\Kxx$. Applying the relations, $y_2\bigl(\psi_1v+av\bigr) =
    c_1\psi_1v + (ac_2+1)v$.  Since this is a weight vector, comparing
    coefficients, $ac_1=ac_2+1$.  Hence, $c_1\ne c_2$ and
    $w=\psi_1v-\tfrac1{c_2-c_1}v$ is a weight vector. The remaining
    claims in part~(b) now follow easily.

    Finally, it remains to consider (c), when $i\ne j$ and
    $\Qx_{ij}(c_1,c_2)=0$. If $w=\psi_1v\ne0$ then $\psi_1w=\psi_1^2v=0$
    since $\Qx_{ij}(c_1,c_2)=0$. In this case $\ei[ij]v=v$ and
    $\ei[ji]w=w$, so $\<v,w\>$ is $\Kxx$-free of rank~$2$. On the other
    hand, if $w=0$ then $\Kxx v$ is a $\Rxaff[2](\Kxx)$-module that is
    free of rank~$1$ as claimed.
  \end{proof}

  The symmetric group $\Sym_n$ acts on $\Kxx^n$ and $I^n$ by place
  permutations. Recall the definition of the elements
  $\phi_r\in\Rx(\Kxx)$ from \autoref{E:phi}.

  \begin{Corollary}\label{C:psiWeightSpace}
    Let $V$ be a weight module for $\Rxaff(\Kxx)$ and let
    $0\ne v\in V_{\bc,\bi}$ be homogeneous, for $\bi\in I^n$ and $\bc\in\Kxx^n$.
    Suppose that $1\le r<n$ and that $(c_r,i_r)\ne(c_{r+1},i_{r+1})$. Then
    $0\ne\phi_rv\in  V_{s_r\bc, s_r\bi}$.
  \end{Corollary}

  \begin{proof}
    By \autoref{R:psiy}, $\psi_rv\in  V_{s_r\bc, s_r\bi}+
    \delta_{i_ri_{r+1}}V_{\bc,\bi}$. In particular, $\psi_rv\in V_{s_r\bc, s_r\bi}$ if $i_r\ne i_{r+1}$.
    If $i_r=i_{r+1}$ then
    $\phi_rv\in V_{s_r\bc,s_{r+1}\bi}$ in view of \autoref{P:R2affSimples}(b)
    since $\psi_r\ei=\bigl(\psi_r(y_r-y_{r+1})+1\bigr)\ei$ in this case.
    Finally, $\phi_r$ is invertible in~$\Rx(\Kxx)$ by
    \autoref{L:KaKaKi}(a), so $\phi_rv\ne0$.
  \end{proof}

  \subsection{Content reduction}
    One of the main results of this section is \autoref{C:FsFt}, which
    shows that $\set{F_\t|\t\in\Std(\Parts)}$ is a family of pairwise
    orthogonal idempotents in $\Rx(\Kxx)$. To prove this we argue by
    induction on~$n$ to classify all weight modules for $\Rn(\Kxx)$ by
    showing that the eigenvalues of $y_1,\dots,y_n$ are given by the
    content functions on the standard tableaux.

    If $\bi\in I^n$ and $1\le m\le n$ define $\bi_{\downarrow
    m}=(i_1,\dots,i_m)\in I^m$. If $\bi\in I^m$ and $j\in I$ let
    $\bi j=(i_1,\dots,i_m,j)\in I^{m+1}$. Let
    $\StI=\set{\br(\s)| \s\in\Std(\Parts[m])}$ be the set of residue
    sequences of the standard tableaux of size~$m$. If $\bj\in I^m$ set
    \[
       \ei[\bj,n]=\sum_{\substack{\bi\in I^n\\\bi_{\downarrow m}=\bj}}
                  \ei\quad\in\Rx(\Kxx).
    \]
    By \autoref{R:Ids}, if $\bi,\bj\in I^m$ then
    $\ei[\bi,n]\ei[\bj,n]=\delta_{\bi\bj}\ei[\bi,n]$ and, moreover,
    $1_{\Rx}=\sum_{\bj\in I^m}\ei[\bj,n]$.

  Let $V$ be an $\Rx(\Kxx)$-module and suppose that $1\le m\le n$. For $\s\in
  \Std(\Parts[m])$ define $V_\s$ to be the simultaneous
  $\bc_k(\s)$-eigenspace of $y_k$ acting on~$\ei[\br(\s)]V$, for $1\le k\le m$.
  That is, $V_\s$ is the $\Kxx$-module
  \[
   V_{\s}=\set[\big]{v\in \ei[\br(\s),n]  V| y_k v = \bc_k(\s)v  \text{ for } 1\le k\le m}.
  \]
  An $\Rx(\Kxx)$-module $V$ is \textbf{$m$-content reduced} if $V$ is
  free as a $\Kxx$-module and $V=\sum_{\s\in \Std(\Parts[m])} V_\s$ as a
  $\Kxx$-module. The module $V$ is \emph{content reduced} if it is
  $n$-content reduced. If~$V$ is $m$-content reduced then the sum
  $V=\sum_{\s\in \Std(\Parts[m])} V_\s$ is necessarily direct because
  $V_\s\cap V_\t=0$, for $\s\ne\t\in\Std(\Parts[m])$. In particular,
  every content reduced module is a weight module for $\Rxaff(\Kxx)$.

  Suppose that $V$ is an $\Rx(\Kxx)$-module. We can consider $V$ as an
  $\Rxaff(\Kxx)$-module using the canonical surjection
  $\Rxaff(\Kxx)\to\Rx(\Kxx)$. By \autoref{T:KLRBasis} and
  \autoref{D:KLR}, over any ring there is an algebra embedding
  of $\Rxaff[m]$ into $\Rxaff$ that sends $\ei[\bj]$ to $\ei[\bj,n]$,
  for $\bj\in I^ m$. Therefore,~$V$ is an~$\Rxaff[m](\Kxx)$-module by
  restriction. Since $V$ is an $\Rx(\Kxx)$-module, it is killed by the
  weight polynomials $\Wbx$, so the $\Rxaff[m](\Kxx)$-action on $V$
  makes $V$ into an $\Rx[m](\Kxx)$-module. Let $\Res_{\Rx[m]}(V)$ and
  $\Res_{\Rxaff[m]}(V)$ be the restrictions of $V$ to an
  $\Rx[m](\Kxx)$-module and $\Rxaff[m](\Kxx)$-module, respectively.

  The irreducible modules $V_\blam$ of \autoref{P:SeminormalForm} are
  content reduced. Conversely, we have:

  \begin{Lemma}\label{L:SpechtDecomposition}
     Let $V$ be an $m$-content reduced $\Rx(\Kxx)$-module,
     where $1\le m\le n$. Then
     \[
        \Res_{\Rx[m]}(V)\cong\bigoplus_{\blam\in\Parts[m]}V_\blam^{\oplus a_\blam},
        \qquad\text{for some }a_\blam\ge0,
     \]
     as an $\Rx[m](\Kxx)$-module.
   \end{Lemma}

  \begin{proof}
     Since $V$ is $m$-content reduced, by definition, it is free as a
     $\Kxx$-module and has a homogeneous basis of weight vectors.
     Let~$v_\s'\in V_\s$ be such a basis vector, where $\s\in\Std(\blam)$
     and $\blam\in\Parts[m]$. To prove the lemma it is enough to show
     that $\Rx[m](\Kxx)v_\s\cong V_\blam$. Let $d_\s=\Ldt[\s]\in\Sym_n$ be the
     permutation such that $\s=d_\s\Ltlam$ and set
     $v_{\Ltlam}=\phi_{w^{-1}}v_\s'$ and $v_\t=\psi_{d_\t}v_{\Ltlam}$,
     where $\t=d_\t\Ltlam$ for $\t\in\Std(\blam)$.
     Then $v_\t$ is a nonzero element of~$V_{\t}$ by
     \autoref{C:psiWeightSpace}.  Moreover,
     $\set{v_\t|\t\in\Std(\blam)}$ is linearly independent since these
     weight spaces are disjoint. Let~$W$ be the submodule of~$V$ spanned
     by the $\set{v_\t|\t\in\Std(\blam)}$. By \autoref{P:R2affSimples}
     and \autoref{L:Separation}(b), if $\t\in\Std(\blam)$ and $1\le k<n$
     then there exist scalars~$\beta_k(\t)$ such that
     \[
        \psi_kv_\t=\beta_k(\t)v_{\sigma_k\t}
        +\tfrac{\delta_{\br_k(\t),\br_{k+1}(\t)}}{\bc_{k+1}(\t)-\bc_k(\t)}v_\t.
     \]
     In particular,~$W$ is an $\Rx[m](\Kxx)$-submodule of~$V$. Further,
     since~$W$ is an $\Rx[m](\Kxx)$-module, relations \autoref{R:quad},
     \autoref{R:psicomm} and \autoref{R:braid} imply that these
     coefficients satisfy conditions (a)--(c), respectively, of
     \autoref{P:SeminormalForm}. (In fact, the reader can check that
     $\beta_k(\t)\in\kx$ is given by \autoref{E:betas}.) Therefore,
     $W\cong V_\blam$ by \autoref{C:IsomorphicSpechts}, completing the
     proof.
  \end{proof}

  \begin{Remark}
    Using \autoref{D:KLR}, it is easy to see that if $1\le m\le
    n$ then there is a surjective algebra map from $\Rx[m](\Kxx)$ onto the
    subalgebra of $\Rx(\Kxx)$ generated by $\psi_1,\dots,\psi_{m-1}$,
    $y_1,\dots,y_m$ and $\ei[\bj,n]$, for $\bj\in I^m$. It follows from
    \autoref{C:fBases} below that this map is an isomorphism, but we
    cannot prove this yet. For now it is enough to work with $m$-content
    reduced modules, which are combinatorial shadows of these isomorphisms.
  \end{Remark}

  The next lemma can be viewed as the module theoretic origin of
  \autoref{D:Ft}. In the lemma we assume that $c_1,\dots,c_N\in\kx$ only
  because $(\bc,\br)$ takes values in $\kx$.

  \begin{Lemma}\label{L:preidempotents}
    Let $V$ be an $\Rx(\Kxx)$-module. Suppose that $
    \prod_{k=1}^N(y_r-c_k)\ei V=0 $, where $1\le r\le n$ and
    $c_1,\dots,c_N\in\kx$ are pairwise distinct and $\bi\in I^n$. Then
    \[
      \ei V=\bigoplus_{k=1}^N V_{\bi,k},
          \qquad\text{where }\quad V_{\bi,k}=\set{v\in\ei V|y_rv=c_k v},
          \text{ for }1\le k\le N.
    \]
  \end{Lemma}

  \begin{proof}
    This follows by applying the easy (polynomial) identity
    $\displaystyle\sum_{k=1}^N\prod_{l\ne k}\frac{(y_r-c_l)}{(c_k-c_l)}=1$.
  \end{proof}

  We now show that every $\Rx(\Kxx)$-module is content reduced, which is
  the linchpin of this section.

  \begin{Theorem}\label{T:ContentReduced}
    Let $V$ be a $\Kxx$-free $\Rx(\Kxx)$-module. Then $V$ is content reduced.
  \end{Theorem}

  \begin{proof}
  We argue induction on~$m$ to show that $V$ is $m$-content reduced, for
  $1\le m\le n$.

  Suppose $m=1$. Fix $\bi=(i)\in I$. By \autoref{D:ContentSystem}(a),
  \[
    \prod_{\substack{1\le l\le \ell\\\rho_l =i}}\bigl(y_1-\bc(l,0)\bigr) \ei=0
    \qquad\Longrightarrow\qquad
    \prod_{\substack{1\le l\le \ell\\\rho_l =i}}\bigl(y_1-\bc(l,0)\bigr) \ei V=0.
  \]
  In view of \autoref{D:ContentSystem}(c) and \autoref{L:Separation}(a),
  there is a self-evident bijection between the sets of standard tableaux
  $\Std(\Parts[1])$ and contents $\set{\bc(l,0)|1\le l\le\ell}$. Hence,
  the module $V$ is $1$-content reduced by \autoref{L:preidempotents}. This
  establishes the base case of our induction.

  Let $1\le m<n$. By induction, we assume that $V$ is $m$-content
  reduced. For the inductive step we show that
  $V=\bigoplus_{\t\in\Std(\Parts[m+1])}V_\t$. Fix $\s\in\Std(\Parts[m])$
  and $j\in I$ and set $V_{\s,j}= 1_{\br(\s) j,n}V_\s$.  To show that $V$ is
  $(m+1)$-content reduced it is enough to prove that
  \begin{equation}\label{E:VsjDecomp}
    V_{\s,j} = \sum_{\substack{\t\in \Std(\Parts[m+1])\\
          \t_{\downarrow m}=\s \text{ and } \br_{m+1}(\t)=j}} V_\t,
          \qquad\text{ for all $\s\in\Std(\Parts[m])$ and $j\in I$}.
  \end{equation}
  Let $\Add_{j}(\s)=\set[\big]{ \t^{-1} (m+1) | \t \in \Std(\Parts),
              \t_{\downarrow m} = \s\text{ and } \br_{m+1}(\t) = j}$
  be the set of addable $j$-nodes for $\s$. By \autoref{L:preidempotents},
  to prove \autoref{E:VsjDecomp} it suffices to show that
  \begin{equation}\label{E:VsKilling}
    \prod_{(l,r,c)\in\Add_j(\s)} \bigl(\bc(l,c-r)-y_{m+1}\bigr)V_{\s,j}=0,
  \end{equation}
  since the contents $\bc(l,c-r)$ in this product are distinct by
  \autoref{L:Separation}. By convention, empty products are~$1$, so the
  last displayed equation includes the claim that $V_{\s,j}=0$ if there
  are no standard tableaux with residue sequence $\bi=\br(\s)j$.

  Let $(k,a,b)=\s^{-1}(m)$ and set $\u=\s_{\downarrow(m-1)}\in\Std(\Parts[m-1])$.
  Define $\Add_j(\u)$ as above.

  We consider two cases.

  \Case{$j=\br_{m}(\s)$}
  By assumption, $\Add_{j}(\u) = \Add_{j}(\s)\sqcup \set{ (k,a,b)}$.
  Hence, in view of \autoref{L:KKDividedDifference}
  and~\autoref{L:DividedDiff},  it follows by induction that
  \[
      \prod_{(l,r,c) \in \Add_{j}(\u)\setminus \set{ (k,a,b)} }
               \bigl(\bc(l,c-r)-y_{m+1}\bigr) V_{\s,j} =0.
  \]
  Hence, \autoref{E:VsKilling} holds when $j=\br_m(\s)$.

  \Case{$j\ne \br_m(\s)$}
  Set $A = \set[\big]{(k,r,c)\in\Nodes|
            \br(k,r,c)=j\text{ and $(r,c)=(a+1,b)$ or $(r,c)=(a,b+1)$}}$.
  Then $|A|=-\CPair{\alpha_{\br_m(\s)}}{\alpha_j}$ and
  $\Add_j(\s)\subseteq\Add_j(\u)\sqcup A$ (disjoint union).
  By \autoref{D:ContentSystem}(b),
  \[
    Q^\ux_{\br_{m}(\s),j} (\bc_{m} (\s), v)=\epsilon\prod_{(k,r,c)\in A}
            \bigl(\bc(k,c-r)-v\bigr),
            \qquad\text{for some }\epsilon\in\k^\times.
  \]
  Hence, by induction, if $v\in V_{\s,j}$ then
  $ \psi_{m}^2 v =\epsilon \prod_{(k,r,c)\in A}\bigl(\bc(k,c-r)-y_{m+1}\bigr)v$.
  Therefore,
  \begin{align*}
   \prod_{(l,r,c) \in \Add_{j}(\u)\sqcup A} \bigl(\bc(l,c-r)-y_{m+1}\bigr) V_{\s,j}
     &= \prod_{(l,r,c) \in\Add_j(\u)}
           \bigl(\bc(l,c-r)-y_{m+1}\bigr)\cdot\psi_{m}^2V_{\s,j} \\
     &= \psi_{m}\prod_{(l,r,c) \in\Add_j(\u)}
          \bigl(\bc(l,c-r)-y_{m}\bigr) \cdot\psi_{m} V_{\s,j} \\
     &\subseteq\psi_{m} \prod_{(l,r,c) \in\Add_j(\u)}
          \bigl(\bc(l,c-r)-y_{m}\bigr)\cdot\ei[\br(\u)j\br_{m}(\s),n] V_{\u} \\
     & = 0,
  \end{align*}
  where the second equality uses \autoref{R:psiy} and the
  last equality follows by induction. In particular, \autoref{E:VsKilling} holds
  by \autoref{L:preidempotents} whenever $\Add_j(\s) = \Add_j(\u)\sqcup A$.
  We need to consider the cases when $\Add_j(\s)$ is properly contained
  in $\Add_j(\u)\sqcup A$, where \autoref{L:preidempotents} potentially
  gives weight spaces of $V_\s$ that are not indexed by standard tableaux.

  Suppose first that $(k,a,b+1)\in A$ and $(k,a,b+1)\notin\Add_j(\s)$.
  Define $c_l=\bc_l(\s)$ and $i_l=\br_l(\s)$, for $1\le l\le m$ and set
  $c_{m+1}=\bc(k,b+1-a)$ and $i_{m+1}=\br(k,b+1-a)$. Let
  $\bc=(c_1,\dots,c_{m+1})$ and~$\bi=(i_1,\dots,i_{m+1})$. By
  \autoref{L:preidempotents}, $V_{\bc,\bi}$ is a (possibly zero) summand of~$V_\s$.  By
  way of contradiction, suppose that $V_{\bc,\bi}\ne0$ and fix a nonzero
  homogeneous vector $v\in V_{\bc,\bi}$. Let $\blam=\Shape(\s)$. Then $(k,a,b+1)$ is
  not an addable node of $\blam$, so $(k,a-1,b)\in\blam$. By induction,
  $V$ is $m$-content reduced, so $V_\blam\cong\Rx[m](\Kxx)v$ as
  an $\Rx[m](\Kxx)$-module by (the proof of) \autoref{L:SpechtDecomposition}.
  Therefore, without loss of generality, we can assume that
  $\s(k,a-1,b)=m-1$. In particular, $c_{m+1}=c_{m-1}$ and
  $i_{m+1}=i_{m-1}$.  Moreover, $\psi_{m-1}v=0$ by
  \autoref{P:SeminormalForm}, since $\sigma_{m-1}\s\notin\Std(\blam)$ by
  \autoref{L:Separation}(b).  Similarly, $\psi_{m}v=0$ because~$V$ is
  $m$-content reduced and no tableau in~$\Std(\Parts[m])$ has content
  sequence $(c_1,\dots,c_{m-1},c_{m-1})$ and residue sequence
  $(i_1,\dots,i_{m-1},i_{m-1})$. Consequently,
  $\bigl(\psi_m\psi_{m-1}\psi_m-\psi_{m-1}\psi_m\psi_{m-1}\bigr)v=0$.
  Therefore, $\Qx_{i_{m-1},i_m,i_{m+1}}(y_{m-1},y_m,y_{m+1})v=0$ by
  \autoref{R:braid}. However, $Q_{i_{m-1},i_m}(c_{m-1},c_m)=0$, so
  \begin{align*}
        \Qx_{i_{m-1},i_m,i_{m+1}}(c_{m-1},c_m,y_{m+1})
           &=\frac{\Qx_{i_{m-1},i_m}(y_{m+1},c_m)}{y_{m+1}-c_{m-1}}\\
           &=\begin{cases*}
              \epsilon\bigl(\bc(k,b-1-a)-y_{m+1}\bigr) & if $\br(k,b-1-a)=i_{m-1}$,\\
              \epsilon& otherwise,
            \end{cases*}
  \end{align*}
  where $\epsilon\in\k^\times $ and the last equality follows by
  \autoref{D:ContentSystem}(b). By \autoref{D:ContentSystem}(c),
  $\bc(k,b-1-a)\ne c_{m+1}$,
  so~$Q_{i_{m-1},i_m,i_{m+1}}(y_{m-1},y_m,y_{m+1})v\ne0$, giving a
  contradiction! Hence, $V_{\bc,\bi}=0$.

  Similarly, if $(k,a,b+1)\notin A$ and $(k,a,b+1)\in\Add_j(\s)$ then
  let $\bc'=(c_1,\dots,c_m,c_{m+1}')$ and
  $\bi'=(i_1,\dots,i_m,i_{m+1}')$, where $c_{m+1}'=\bc(k,b-a-1)$ and
  $i_{m+1}'=\br(k,b-a-1)$. Then $(k,a,b-1)\in\blam$ and~$V_{\bc',\bi'}$
  is a summand of $V_\s$ by \autoref{L:preidempotents}.  Arguing as in
  the last paragraph, we deduce that $V_{\bc',\bi'}=0$.

  Consequently, if $j\ne\br_m(\s)$ then the last displayed equation,
  combined with \autoref{L:preidempotents}, shows that
  \autoref{E:VsKilling} holds.

  We have now established \autoref{E:VsjDecomp} in all cases, so $V$ is
  $(m+1)$-content reduced. This completes the proof of the inductive step
  and, hence, the proof of the proposition.
  \end{proof}

  Applying \autoref{T:ContentReduced} to the regular representation, and
  using \autoref{L:SpechtDecomposition}, shows that the algebra
  $\Rx(\Kxx)$ is completely reducible.  \autoref{P:GenericIsomorphism}
  makes this more explicit.

  \begin{Corollary}\label{C:FtProjection}
    Let $V$ be a $\Kxx$-free $\Rx(\Kxx)$-module. Then $V=\bigoplus_\t F_\t V$ as a
    $\Kxx$-module, where the sum is over $\t\in\Std(\Parts)$ such that $F_\t V\ne0$.
  \end{Corollary}

  \begin{proof}
    By \autoref{D:Ft}, if $\t\in\Std(\Parts)$ then $V_\t\subseteq\set{v\in
    V|v=F_\t v}$. On the other hand, $V=\bigoplus_\t V_\t$ by
    \autoref{T:ContentReduced}.  Therefore, $V_\t=\set{v\in V|v=F_\t v}$
    since $F_\s V\cap F_\t V=\delta_{\s\t}F_\t V$ by \autoref{L:Separation}.
  \end{proof}

  \begin{Corollary}\label{C:yFt}
    Suppose that $\t\in\Std(\Parts)$ and $1\le m\le n$. Then
    $y_mF_\t = \bc_m(\t)F_\t$ in $\Rx(\Kxx)$.
  \end{Corollary}

  \begin{proof}
    Take $V=\Rx(\Kxx)$ to be the regular representation, which is free
    as a $\Kxx$-module by base change from \autoref{P:KLRfree} since
    $\Rx(\Kxx)\cong\Kxx\otimes_{\Kx}\Rx(\Kx)$. First note that $F_\t\ne0$
    by \autoref{E:Fsvt}. By \autoref{C:FtProjection},
    $V_\t=F_\t\Rx(\Kxx)$. As~$F_\t=F_\t\cdot1\in F_\t\Rx(\Kxx)=V_\t$,
    this implies the result.
  \end{proof}

  Hence, using \autoref{L:Separation} and \autoref{D:Ft}, we obtain:

  \begin{Corollary}\label{C:FsFt}
    Let $\s,\t\in\Std(\Parts)$. Then $F_\s F_\t=\delta_{\s\t}F_\t$ in
    $\Rx(\Kxx)$.
  \end{Corollary}

  \begin{Corollary}\label{C:KLRIdempotents}
      Suppose that $\bi\in I^n$. Then, in $\Rx(\Kxx)$,
      \[ \ei= \sum_{\t\in \Std(\bi)} F_{\t}.\]
      In particular, $\ei=0$ if and only if $\bi\notin\StI$.
  \end{Corollary}

  \begin{proof}
    Take $V=\Rx(\Kxx)$ to be the regular representation of~$\Rx(\Kxx)$. By \autoref{C:FtProjection},
    \[
        \ei\Rx(\Kxx)=\bigoplus_{\t\in\Std(\bi)}F_\t\Rx(\Kxx).
    \]
    Hence, the element $\ei-\sum_{\t\in\Std(\bi)}F_\t$ acts on $\ei\Rx(\Kxx)$ as
    multiplication by zero by \autoref{C:FsFt}. Therefore, by \autoref{R:Ids},
    this element acts on $\Rx(\Kxx)$ as zero. Hence,
    $\ei=\sum_{\t\in\Std(\bi)}F_\t$ by the faithfulness of the regular
    representation.  Finally, these arguments show that if
    $\Std(\bi)=\emptyset$, then $\ei=0$. That is,~$\ei=0$ if
    and only if $\bi\notin\StI$.
  \end{proof}

  \begin{Remark}
    The last two corollaries are the main results of this section.
    Rather than the approach we have taken, these results can also be
    deduced from \autoref{P:SeminormalForm} by first showing that
    $V=\bigoplus_\blam V_\blam$ is a faithful $\Rx(\Kxx)$-module, which
    can be proved after computing the (graded) dimension of~$\Rx(\Kxx)$
    using ideas from \cite{BK:GradedDecomp,ArikiParkSpeyer:C}.  That the
    representation~$V$ is faithful now follows from
    \autoref{C:KLRIdempotents}. The next section gives a different take
    on this description of~$\Rx(\Kxx)$ as the endomorphism algebra
    of~$V$.
  \end{Remark}

  \begin{Corollary}\label{C:psiF}
  Suppose that $\br_k(\t)\ne \br_{k+1}(\t)$ for $\t \in \Std(\Parts)$ and $1\le r<n$.
  Then $y_m \psi_k F_\t = \bc_{\sigma_k (m)} (\t) \psi_k F_\t$ whenever $1\le m\le n$. In
  particular, $\psi_k F_\t=0$ if $\sigma_k\t$ is not standard.
  \end{Corollary}

  \begin{proof}
    Suppose that $\br_k(\t) \ne \br_{k+1}(\t)$. The claim that
    $y_m \psi_k F_\t = \bc_{\sigma_k (m)} (\t) \psi_k F_\t$ follows immediately
    from \autoref{R:psiy} and~\autoref{C:yFt}. For the second statement,
    if $\sigma_k\t \notin \Std(\Parts)$ then the node $\t^{-1}(k+1)$ is either
    directly to the right of, or directly below, $\t^{-1}(k)$. Therefore,
    $\br_k(\t) \ne \br_{k+1}(\t)$ by \autoref{L:AdjacentCS}.
    Consequently, by \autoref{L:Separation}(b), there is no element in
    $\Std(\Parts)$ with residue sequence $\sigma_k \br(\t)$ and content
    sequence $\sigma_k c(\t)$. Hence, $\psi_k
    F_\t=F_{\sigma_k\t}\psi_k=0$ by \autoref{C:KLRIdempotents}.
  \end{proof}

  \subsection{The algebra \texorpdfstring{$\Snl$}{Sln}}\label{S:Snell}

  This section introduces the algebra $\Snl$, which is the ``universal''
  \textit{semisimple} cyclotomic KLR algebra of level~$\ell$. In the next section
  we show if $\Rx(\Kxx)$ has a content system then it is isomorphic to
  $\Snl$. We maintain the notation of the previous sections except we
  work over the field~$\K$.

  Recall from \autoref{SS:Contents} that $\Gamma_\ell$ is the quiver of
  type $A_\infty^{\times\ell}$, with vertex set
  $\Jell=\set{1,\dots,\ell}\times\Z$. Let~$\Snl(\K)$ be the standard
  cyclotomic KLR algebra defined using the (standard) $Q$-polynomials and weight
  polynomials of \autoref{Ex:ContentSystem}(a). Let $(\bc^\Jell,\br^\Jell)$
  be the content system for $\Snl(\K)$ given in
  \autoref{Ex:ContentSystem}(a), so that $\bc^\Jell$ is identically zero
  and $\br^\Jell$ is the identity map on $\Jell$. By assumption, $\ux$ is
  the empty sequence for $\Snl$ so, by convention, $\Kxx=\K$.
  \notation{$\Snl$}{Universal level $\ell$ semisimple algebra for content system}

  To avoid confusion, if $\t\in \Std(\Parts)$ let $\br^\Jell(\t)$ be the residue
  sequence of~$\t$ with respect to the content system
  $(\bc^\Jell, \br^\Jell)$.  Explicitly, $\br^\Jell(\t) = (\br^\Jell_1(\t),\dots,\br^\Jell_n(\t))\in J_\ell^n$ where
  $\br^\Jell_m(\t) = \rho_k+b-a$ if $\t^{-1} (m) = (k,a,b)$. For convenience, set
  $\Jstd= \set{ \br^\Jell(\t) | \t\in \Std(\Parts) }$.
  By \autoref{L:Separation}, if $\bj\in\Jstd$ then there
  exists a unique standard tableau $\t\in \Std(\Parts)$ such that $\br^\Jell(\t) = \bj$
  since $\bc^\Jell$ is identically zero.

  \begin{Lemma}\label{L:SnlIdentities}
    Suppose that $1\le k< n$ and $\bj\in J_\ell^n$. Then $y_1=\dots=y_n=0$
    and $\ei[\bj]\ne0$ if and only if~$\ei[\bj]=F_\t$ for some $\t\in\Std(\Parts)$. Consequently,
    $\psi_k \ei[\bj]=0$ if $j_k \arrow j_{k+1}$ and $\ei[\bj]=0$ if $j_k=j_{k+1}$
    or $j_k=j_{k+2}$ for $1\le k<n-1$).
  \end{Lemma}

  \begin{proof}
  Let $V$ be the left regular representation of~$\Snl(\K)$.  Then
  $V=\bigoplus_{\t\in \Std(\Parts)}V_\t$ by \autoref{T:ContentReduced}.
  Since $\bc^\Jell$ is identically zero, $y_m$ acts as multiplication by
  zero on $V_\t$, for $1\le m\le n$ and  $\t\in\Std(\Parts)$. Hence,
  $y_1=\dots=y_n=0$ proving the first claim.

  Next, we show that  $\ei[\bj]\ne0$ if and only if $\ei[\bi]=F_\t$, for some $\t\in\Std(\Parts)$.
  Observe that if $\s,\t\in\Std(\Parts)$ then $\s=\t$ if and only if
  $\br^\Jell(\s)=\br^\Jell(\t)$ by \autoref{L:Separation} since $\bc_{\Jell}$
  is identically zero. Hence, $\ei[\bj]=F_\t$ for some $\t\in\Std(\Parts)$
  by \autoref{C:KLRIdempotents}. The remaining statements now follow by
  \autoref{C:KLRIdempotents} and \autoref{C:psiF}.
  \end{proof}

  \begin{Definition}
    Let $\blam\in\Parts$. For $\s,\t\in\Std(\blam)$ set
    $\Psi_{\s\t}=\psi_{w}\ei[\br^\Jell(\t)]$, where $w\in\Sym_n$ is
    the unique permutation such that $\s=w\t$.
  \end{Definition}

  \notation{$\Psi_{\s\t}$}{Basis elements of $\Snl(\K)$}

  \begin{Corollary}\label{C:SnlSpanning}
  The algebra $\Snl(\k)$ is spanned by
  $\set{ \Psi_{\s\t} | (\s,\t)\in \Std^2(\Parts)}$.
  \end{Corollary}

  \begin{proof}
    By \autoref{T:KLRBasis} and \autoref{L:SnlIdentities}, $\Snl$ is
    spanned by the set
    \[
      \set{ \psi_w \ei[\br^\Jell(\t)] | w\in\Sym_n\text{ and }\t\in \Std(\Parts)}.
    \]
    Hence, it is enough to show that if $\psi_w\ei[\br^\Jell(\t)]\ne 0$,
    for $\t\in \Std(\Parts)$ and $w\in \Sym_n$, then $w\t \in
    \Std(\Parts)$. Since~$w$ is a product of simple reflections, it is
    enough to consider the case when~$w=\sigma_k=(k,k+1)$, for $1\le k<n$.
    If~$\t$ is standard then $\sigma_k\t$ is standard unless $k$ and
    $k+1$ are in the same row, or the same column of~$\t$, in which case
    $\psi_k \ei[\bj]=0$ by \autoref{L:SnlIdentities}. Hence, if
    $\psi_k\ei\ne0$ then $\sigma_k\t\in\Std(\Parts)$ as we needed to
    show.
  \end{proof}

  Arguing by induction on $n$, it is easy to see that if $\s,\t\in
  \Std(\Parts)$ and $\br^\Jell(\s)=w\br^\Jell(\t)$, for some $w\in\Sym_n$,
  then $\Shape(\s)=\Shape(\t)$.

  Given $u,w\in \Sym_n$, write $u\preceq w$ if there is a reduced
  expression $w=\sigma_{a_1}\dots \sigma_{a_l}$ such that
  $u=\sigma_{a_1} \dots \sigma_{a_k}$, for some $0\le k<l$. (This is the
  right weak Bruhat order on $\Sym_n$.)

  \begin{Lemma}\label{L:PrefixOrder} Let $\t\in\Std(\Parts)$
    and suppose that $w\t$ is standard, for some $w\in\Sym_n$. Then $u\t$ is standard
    whenever $u\preceq w$.
  \end{Lemma}

  \begin{proof}
    If $1\le r<t\le n$ and $u(r)>u(t)$ then $w(r)>w(t)$ since $u\preceq w$.
    The result follows easily from this observation.
  \end{proof}

  \begin{Lemma}\label{L:SnlSpechtModule}
  Let $\blam\in \Std(\Parts)$. Then there exists an irreducible left
  $\Snl(\K)$-module $W_\blam$ with basis
  $\set{w_\t|\t\in\Std(\blam)}$ and where the $\Snl(\K)$-action is determined by
  \[
        \ei[\bj] w_\t = \delta_{\bj,\br^\Jell(\t)} w_\t, \quad
        y_m w_\t =0, \quad
        \psi_k w_\t = \begin{cases}
                     w_{\sigma_k\t}& \text{if } \sigma_k\t \in \Std(\blam), \\
                     0& \text{otherwise},
                 \end{cases}
  \]
  for all $\bj \in J_\ell^n$ and all admissible $k$ and $m$.
  \end{Lemma}

  \begin{proof}
    By \autoref{L:SnlIdentities},  the map $\t\mapsto \br^\Jell(\t)$ gives
    a bijection $\Std(\Parts)\bijection\Jstd$ such that $F_\t=\ei$, where
    $\bi=\br^\Jell(\t)$. Moreover, by \autoref{E:Qmt} and
    \autoref{L:SnlIdentities},
    \[  Q_k(\t)=\begin{cases*}
               1& if $\sigma_k\t\in\Std(\Parts)$,\\
               0& otherwise.
        \end{cases*}
    \]
    Therefore, in view of \autoref{E:betas}, the lemma is a special case
    of \autoref{P:SeminormalForm}.
  \end{proof}

  \begin{Remark}
    The $\Rx(\K[x^\pm])$-module $V_\blam$ is irreducible only over ~$\K[x^\pm]$.
    In contrast, it is easy to see that the module $W_\blam$ is irreducible
    over any field.
  \end{Remark}

  \begin{Remark}\label{R:SnlGrading}
    \autoref{L:SnlSpechtModule} is also a consequence
    of~\cite[Theorem 3.4]{KleshRam:PureIrredRepsKLR}. By
    \autoref{L:SnlIdentities}, the natural grading on $W_\blam$
    concentrates everything in degree~$0$.
  \end{Remark}

  We now prove that $\Snl(\K)$ is a split semisimple algebra.

  \begin{Theorem}\label{T:SSemisimple}
  The algebra $\Snl(\K)$ is a split semisimple algebra and
  $\set{W_\blam|\blam\in\Parts}$ is a complete set of pairwise
  non-isomorphic irreducible $\Snl$-modules, up to shift.
  \end{Theorem}

  \begin{proof}
    Recall from \autoref{C:SnlSpanning} that the elements
    $\set{\Psi_{\s\t}|(s,\t)\in\Std^2(\Parts)}$ span $\Snl(\k)$.  By
    \autoref{L:PrefixOrder} and \autoref{L:SnlSpechtModule}, if
    $\s,\t\in\Std(\bmu)$ then the action of $\Psi_{\s\t}$ on the module
    $W_\blam$ is given by $\Psi_{\s\t} w_\u = \delta_{\t\u} w_\s$, for
    $\u\in \Std(\blam)$. In particular, if $\bmu\ne\blam$ then~$\Psi_{\s\t}$ acts as
    zero on $W_\blam$. Moreover, this implies that the set
    $\set{\Psi_{\s\t}|(\s,\t)\in\Std^2(\Parts)}$ is linearly independent,
    and so is a basis of $\Snl(\k)$ by \autoref{C:SnlSpanning}. Extending
    scalars to~$\K$, there is a well-defined algebra isomorphism
    \[
    \mathcal{E}\map{\Snl(\K)}{\bigoplus_{\blam\in\Parts}\End_{\K}(W_\blam)};
             \Psi_{\s\t}\mapsto e_{\s\t},
    \]
    where $e_{\s\t}$ is the matrix unit given by
    $e_{\s\t}(w_\u)=\delta_{\t\u}w_\s$. It follows that $\mathcal{E}$ is
    an algebra isomorphism since $\set{\Psi_{\s\t}}$ is a basis
    of~$\Snl(\K)=\K\otimes_\k\Snl(\k)$, completing the proof.
  \end{proof}

  \begin{Remark}
    As in \autoref{R:SnlGrading}, the grading on $\Snl(\K)$ puts everything in
    degree zero. The complete set of irreducible graded $\Snl(\K)$-modules
    is $\set{q^dW_\blam|\blam\in\Parts\text{ and } d\in\Z}$. In contrast,
    if $x$ is an indeterminate, in
    degree~$1$, then the complete set of irreducible graded
    $\Snl(\K[x^\pm])$-modules is $\set[\big]{\K[x^\pm]\otimes_\K W_\blam|\blam\in\Parts}$,
    since $\K[x^\pm]$ is the unique irreducible graded $\K[x^\pm]$-module.
  \end{Remark}

  The proof of \autoref{T:SSemisimple} and \autoref{C:SnlSpanning} gives
  a basis of $\Snl(\k)$.

  \begin{Corollary}\label{C:SnlBasis}
  The algebra $\Snl(\k)$ is free as a $\k$-module with basis
  $\set[\big]{ \Psi_{\s\t}| (\s,\t)\in \Std^2(\Parts)}$.
  \end{Corollary}

  \subsection{Semisimplicity of deformed cyclotomic KLR algebras}
  This section returns to the framework of \autoref{SS:Contents}. In
  particular, we assume that $(\Qbx,\Wbx)$ is a $\kx$-deformation
  $(\Qb,\Wb)$ and that $(\bc,\br)$ is a content system for $\Rx$ with values
  in~$\kx$. This section proves that the algebras $\Rx(\Kxx)$ and
  $\Snl(\Kxx)$ are isomorphic as ungraded algebras, where $\K$ is the
  field of fractions of~$\k$.

  Recall the elements $\phi_1,\dots,\phi_{n-1}\in \Rx(K)$ defined
  in~\autoref{E:phi}.

  \begin{Lemma}\label{L:phiF}
  Suppose that $\t\in\Std(\Parts)$ and $1\le k<n$. Then, in $\Rx(\Kxx)$,
  \[
      \phi_k F_\t = \begin{cases*}
          F_{\sigma_k \t} \phi_k & if  $\sigma_k\t$ is standard, \\
          0                 & otherwise,
      \end{cases*}
  \]
  \end{Lemma}

  \begin{proof}
    By \autoref{L:KaKaKi}(d), if $1\le m<n$ then
    $\phi_k(y_m-c)=(y_{\sigma_k(m)}-c)\phi_k$. Hence, the result follows
    by \autoref{D:Ft} (and \autoref{L:Separation}).
  \end{proof}

  Let $\t\in \Std(\blam)$ and $1\le m<n$. Note that if $\bj=\br^\Jell(\t)$
  then $\br^\Jell_m(\t)\ne \br^\Jell_{m+1}(\t)$ by \autoref{L:SnlIdentities}. Recall the
  scalar $Q_m(\t)$ for $\Rx(\Kxx)$ from \autoref{E:Qmt}. Set
  \begin{equation}\label{E:qm}
  q_m(\t)= \begin{cases}
        Q_m(\t)^{-1} & \text{if } \br^\Jell_m(\t) \noarrow \br^\Jell_{m+1}(\t), \;
        \br^\Jell_m(\t) \ne \br^\Jell_{m+1}(\t) \text{ and } \sigma_m \t \vartriangleright \t, \\
        1 & \text{otherwise.}
        \end{cases}
  \end{equation}
  Note that $q_m(\t)$ is well-defined because $Q_m(\t)\ne0$ by \autoref{L:Qnonzero}.
  Moreover,
  \begin{equation}\label{Egprod}
    \text{if $\br^\Jell_m(\t) \noarrow \br^\Jell_{m+1}(\t)$ and $\br^\Jell_m(\t) \ne \br^\Jell_{m+1}(\t)$, then }
    q_m(\t) q_{m}(\sigma_m \t) = Q_m(\t)^{-1}.
  \end{equation}

  Let $\Snl(\Kxx)=\Kxx\otimes_\K\Snl(\K)$. Recall that if $A$ is graded
  then $\underline{A}$ forgets the grading on~$A$.

  \begin{Proposition}\label{P:GenericIsomorphism}
  There is an (ungraded) algebra isomorphism
  $\Theta\colon\underline{\Snl(\Kxx)}\iso \underline{\Rx(\Kxx)}$ such that
  $\Theta\bigl(y_m\bigr)=0$,
  \[
    \Theta\bigl(\ei[\bj]\bigr) = \begin{cases}
        F_{\t}&\text{if } \bj=\br^\Jell(\t) \in \Jstd, \\
        0& \text{if } \bj\notin \Jstd,
    \end{cases},\quad
    \Theta\bigl(\psi_k \ei[\bj]\bigr) =\begin{cases}
        q_k(\t) \phi_k F_{\t}& \text{if } \bj=\br^\Jell(\t) \in \Jstd, \\
        0& \text{if } \bj \notin \Jstd.
    \end{cases}
  \]
  for all $\bj\in J_\ell^n$ and all admissible $m$ and~$r$.
  \end{Proposition}

  \begin{proof}
  First, note that $\Theta(\psi_k)= \sum_\bj\Theta(\psi_k\ei[\bj])$, so the
  images of the generators of $\Snl$ under~$\Theta$ are uniquely
  determined. Hence, once we show that $\Theta$ is a homomorphism it is
  necessarily unique.  If $1\le m\le n$ then $y_m=0$, by
  \autoref{L:SnlIdentities}, so the assumption that
  $y_m\in\ker\Theta$ does not prevent $\Theta$ from being an
  isomorphism. Similarly, by \autoref{L:SnlIdentities}, if $\bj\in
  \Jell^n$ then $\ei[\bj]\ne0$ if and only if $\bj\in\Jstd$.

  To show that $\Theta$ is an algebra homomorphism it is enough to check
  that it respects the KLR relations~\autoref{R:Ids}--\autoref{R:braid}
  and the cyclotomic relation \autoref{R:cyc}. The cyclotomic
  relation~\autoref{R:cyc} is trivially satisfied and checking
  relations~\autoref{R:Ids}--\autoref{R:psicomm} and \autoref{R:psiy} is
  easy, so these are left to the reader. Relation~\autoref{R:psi} is
  routine using \autoref{L:phiF}. For relation~\autoref{R:quad} it is
  enough to show that if $\bj\in J_\ell^n$ and $1\le k<n$ then
  \[
    \Theta(\psi_k^2 \ei[\bj])
       = \Theta\bigl(Q_{j_k,j_{k+1}} (y_k,y_{k+1}) \ei[\bj]\bigr)
  \]
  By definition, the right-hand side is equal to
  \[
  \Theta\bigl(Q_{j_k,j_{k+1}} (y_k,y_{k+1}) \ei[\bj]\bigr) =
        \begin{cases*}
          F_{\t}& if  $\bj=\br^\Jell(\t) \in \Jstd$ and  $j_k \noarrow j_{k+1}$, \\
          0& otherwise.
        \end{cases*}
  \]
  If $\bj\notin\Jstd$ then $\Theta(\ei[\bj])=0$, so we may assume that
  $\bj=\br^\Jell(\t)$, for some $\t\in \Std(\Parts)$.  If $j_k\noarrow j_{k+1}$, then
  \begin{align*}
  \Theta(\psi_k^2 \ei[\bj])&= q_k(\t) q_{k}(\sigma_k \t) \phi_k^2 F_{\t}\\
     &=q_k(\t)q_{k}(\sigma_k\t)\bigl(Q_{\br^\Jell_k(\t), \br^\Jell_{k+1}(\t)}(y_k,y_{k+1})
           +\delta_{\br^\Jell_k(\t),\br^\Jell_{k+1}(\t)}\bigr)F_\t\\
     &=q_k(\t)q_{k}(\sigma_k\t)\,Q_k(\t) F_\t\\
     &= F_{\t},
  \end{align*}
  where we have used \autoref{L:KaKaKi}(f) for the second equality
  and~\autoref{Egprod} for the last equality.  On the other hand, if $j_k\arrow j_{k+1}$
  then
  $ \Theta(\psi_k^2 \ei[\bj]) = \Theta(\psi_k \ei[\sigma_k \bj]) \Theta(\psi_k \ei[\bj]) = 0$
  since $\sigma_k \bj \notin \Jstd$ (compare with \autoref{L:SnlIdentities}). Hence,
  $\Theta$ respects the quadratic relation~\autoref{R:quad}.

  Now consider the deformed braid relation~\autoref{R:braid}. Since $y_m=0$ for
  $1\le m\le n$, we need to verify that if $1\le k<n$ and $\t\in\Std(\Parts)$ and
  \[
      \Theta (\psi_k \psi_{k+1} \psi_k\ei[\br^\Jell(\t)])
        = \Theta(\psi_{k+1} \psi_k \psi_{k+1}\ei[\br^\Jell(\t)])
  \]
  If $\sigma_k\sigma_{k+1}\sigma_{k}\t=\sigma_{k+1}\sigma_k\sigma_{k+1}\t$ is
  not standard then both sides are zero, so we can  assume
  that this tableau is standard. By \autoref{L:KaKaKi}(b)
  and~\autoref{L:PrefixOrder}, it is enough to show that
  \[
  q_{k}(\sigma_{k+1} \sigma_{k} \t) q_{k+1}(\sigma_k \t) q_k(\t) =
  q_{k+1}(\sigma_k \sigma_{k+1} \t) q_{k}(\sigma_{k+1} \t) q_{k+1}(\t).
  \]
  It follows from~\autoref{E:qm} that
  $q_{k}(\sigma_{k+1} \sigma_k \t) = q_{k+1}(\t)$, $q_{k+1}(\sigma_k\t) = q_{k}(\sigma_{k+1}\t)$
  and $q_k(\t) = q_{k+1}(\sigma_k \sigma_{k+1} \t)$, so~\autoref{R:braid} is satisfied.

  We have now proved that $\Theta$ is an algebra homomorphism. By
  \autoref{C:KLRIdempotents}, to show that $\Theta$ is surjective it is
  enough to check that $\ei F_\t$, $y_k F_\t$ and $\psi_k F_\t$ belong to
  the image of~$\Theta$, for all $\bi\in I^n$, $\t\in \Std(\Parts)$ and
  all admissible~$k$.  Certainly,
  $\ei F_\t = \delta_{\bi\,\br^\Jell(\t)} F_{\t}= \delta_{\bi\,\br^\Jell(\t)}
             \Theta(\ei[\br^\Jell(\t)])\in\im\Theta$.
  Hence, $y_k F_\t\in\im\Theta$ by \autoref{C:yFt}.  Finally,
  consider $\psi_k F_\t$. If $\sigma_k \t$ is not standard, then $\psi_k
  F_\t=0$ by \autoref{C:psiF}.  Otherwise, by \autoref{E:phi} we have
  \[
      q_k(\t)^{-1} \Theta(\psi_k \ei[\br^\Jell(\t)]) = \phi_k F_\t
          = \begin{cases}
              (\bc_k (\t)-\bc_{k+1}(\t))\psi_k  F_\t  +F_\t& \text{if } \br^\Jell_k(\t)= \br^\Jell_{k+1}(\t), \\
              \psi_k F_\t& \text{if } \br^\Jell_k(\t) \ne \br^\Jell_{k+1}(\t).
            \end{cases}
  \]
  In both cases it follows that $\psi_k F_\t\in\im\Theta$, where we use
  \autoref{D:ContentSystem}(c) when $\br^\Jell_k(\t)=\br^\Jell_{k+1}(\t)$.
  Hence,~$\Theta$ is surjective.

  We have now shown that $\Theta$ is a surjective algebra homomorphism
  from $\Snl(\Kxx)$ to $\Rx(\Kxx)$. Let~$K$ be any field containing $\Kxx$.
  Extending scalars to~$K$ and using \autoref{P:SeminormalForm},
  \autoref{C:IsomorphicSpechts} and \autoref{T:SSemisimple}, the algebra
  $\Rx(K)$ has at least as many isomorphism classes of (ungraded) simple
  modules as $\Snl(K)$. Hence, by a dimension count, the induced
  map~$\Theta_K$ from~$\Snl(K)$ to~$\Rn(K)$ is an isomorphism. Therefore,
  $\Theta_K$, and hence~$\Theta$, is injective. It follows that
  $\Theta\map{\Snl(\Kxx)}\Rx(\Kxx)$ is an isomorphism of ungraded algebras,
  so the proof is complete.
  \end{proof}

  \begin{Remark}
    The isomorphism $\Theta$ of \autoref{P:GenericIsomorphism} is not
    homogeneous because, in general, the elements $\psi_k\ei[\bj]$ and
    $\Theta(\psi_k\ei[\bj])$ have different degrees.
  \end{Remark}

  Recall the irreducible graded $\Rx(\K[x^\pm])$-module $V_\blam$, for
  $\blam\in\Parts$, defined before \autoref{C:IsomorphicSpechts}.
  Combining \autoref{T:SSemisimple} and \autoref{P:GenericIsomorphism}
  shows that $\Rx(\K[x^\pm])$ is isomorphic to a direct sum of matrix algebras
  over~$\K[x^\pm]$. Hence, we have:

  \begin{Corollary}\label{C:RxIrreducibles}
    The algebra $\Rx(\K[x^\pm])$ is a split semisimple algebra over $\K[x^\pm]$
    and $\set[\big]{V_\blam|\blam\in\Parts}$ is a complete
    set of pairwise non-isomorphic irreducible graded $\Rx(\K[x^\pm])$-modules.
  \end{Corollary}

  In particular, up to isomorphism, the irreducible module $V_\blam$
  does not depend on the choice of content system~$(\bc,\br)$, for
  $\blam\in\Parts$. We already knew from \autoref{C:IsomorphicSpechts}
  that $V_\blam$ is independent of the choice of $\beta$-coefficients in
  \autoref{P:SeminormalForm}.

  \section{Cellular bases of
  \texorpdfstring{$\Rx(\Kxx)$}{Rn(k[x])}}\label{S:Cellularity}
  The main results of this paper follow from the
  construction of cellular bases for the algebra $\Rx(\kx)$, which is
  the focus of this chapter. The cellular bases that we construct are
  analogues of the $\psi$-bases of \cite{HuMathas:GradedCellular}. Using
  the results of \autoref{S:Contents} it is easy to see that the
  $\psi$-bases are linearly independent. The main difficulty is
  showing that the $\psi$-bases span the algebra $\Rx(\kx)$.

  Throughout the chapter, we continue to assume that
  $(\Gamma,\Qbx,\Wbx)$ is a $\kx$-deformation of a standard cyclotomic
  KLR datum $(\Gamma,\Qb,\Wb)$ and $(\bc,\br)$ is a (graded) content
  system with values in $\kx$ and we let $\K$ be the field of fractions
  of~$\k$. \autoref{S:Contents} studied the semisimple representation
  theory of the algebra $\Rx(\Kxx)$.

  \subsection{Integral and seminormal bases}\label{S:Bases}
  Partly inspired by \cite{HuMathas:GradedCellular, Mathas:Seminormal},
  this section defines the two new bases of $\Rx(\kx)$ that will
  ultimately allow us to prove our main results. Defining these bases is
  easy, but it will take some time to prove that they are both
  (cellular) bases over $\kx$.

  Recall from \autoref{S:Tableaux} that $\Gedom$ is the dominance order
  on~$\Parts$. If $\s\in\Std(\Parts)$ is a standard tableau and $1\le
  m\le n$ then $\s_{\downarrow m}$ is the subtableau of $\s$ that
  contains the numbers in $\set{1,\dots,m}$. Extend the dominance order
  to $\Std(\Parts)$ by defining $\s\Gedom\t$ if $\Shape(\s_{\downarrow
  m})\Gedom\Shape(\t_{\downarrow m})$, for $1\le m\le n$. Write
  $\s\Gdom\t$ if $\s\Gedom\t$ and $\s\ne\t$. Similarly, given
  $(\s,\t),(\u,\v)\in\Std^2(\Parts)$ write $(\s,\t)\Gedom(\u,\v)$ if
  $\s\Gedom\u$ and $\t\Gedom\v$. As before, write $(\s,\t)\Gdom(\u,\v)$
  if $(\s,\t)\Gedom(\u,\v)$ and $(\s,\t)\ne(\u,\v)$.

  \notation{$\s_{\downarrow m}$}{Restriction of the tableau $\s$ to $\set{1,\dots,m}$}
  \notation{$\s\Ledom\u$}{dominance on standard tableaux}
  \notation{$(\s,\t)\Ledom(\u,\v)$}{Dominance on pairs of tableaux: $\s\Ledom\u$ and $\t\Ldom\v$}

  \begin{Definition}[Residue dominance]
     Let $\s$ and $\t$ be two standard tableaux. Write $\s\GeDom\t$ if
     $\br(\s)=\br(\t)$ and $\s\Gedom\t$. If $\blam,\bmu\in\Parts$, write
     $\blam\GeDom\bmu$ if there exist $\s\in\Std(\blam)$ and
     $\t\in\Std(\bmu)$ such that $\s\GeDom\t$.
  \end{Definition}

  In what follows we could replace the posets $(\Parts,\Ldom)$ and
  $\Parts,\Gdom)$ with $(\Parts,\LDom)$ and $(\Parts,\GDom)$,
  respectively. However, doing this does not give very much additional
  information because all of our definitions are compatible with the
  block decompositions $\Rx=\bigoplus_\alpha\Rx[\alpha]$ and the residue
  dominance orderings are just the dominance ordering restricted to
  these subalgebras. We remark that in type~$\Aone$ the algebras $\Rx[\alpha]$ are
  indecomposable by \cite[(1.4)]{BK:GradedDecomp} (and
  \cite{LyleMathas:AKBlocks}). In type~$\Cone$ it is not known if
  $\Rx[\alpha]$ is indecomposable, although we expect this to be the
  case.

  Let $\blam\in\Parts$. The \emph{conjugate} of $\blam$ is the
  $\ell$-partition $\blam'=\set{(\ell-k+1,c,r)|(k,r,c)\in\blam}$.  That
  is,~$\blam'$ is the $\ell$-partition obtained from $\blam$ by
  reversing the order of the components and then swapping the rows and
  columns in each component.  As is well-known, if
  $\blam, \bmu\in\Parts$ then $\blam\Gedom\bmu$ if and only if
  $\blam'\Ledom\bmu'$.
  Similarly, the \emph{conjugate tableau} to
  $\t\in\Std(\blam)$ is the standard $\blam'$-tableau
  $\t'$ with $\t'(k,r,c)=\t(\ell-k+1,c,r)$, for $(k,r,c)\in\blam'$.
  \notation{$\blam'$}{Conjugate $\ell$-partition
    $\blam'=(\lambda^{(\ell)\prime}|\dots|\lambda^{(1)\prime})$}
  \notation{$\t'$}{Conjugate tableau: $\t'(k,r,c)=\t(\ell-k+1,c,r)$}

  It is well-known that there exist unique tableaux
  $\Gtlam$ and $\Ltlam$ such that $\Ltlam\Ledom\s\Ledom\Gtlam$, for all
  $\s\in\Std(\blam)$. Explicitly,
  $\Gtlam=(\t^{\Gdom\lambda^{(1)}}|\dots|\t^{\Gdom\lambda^{(\ell)}})$ is
  the standard $\blam$-tableau with the numbers $1,2,\dots,n$ entered in
  order from left to right along the rows of~$\t^{\Gdom\lambda^{(1)}}$,
  and then the rows of~$\t^{\Gdom\lambda^{(2)}}$ and so on. Similarly,
  $\Ltlam=(\t^{\Ldom\lambda^{(1)}}|\dots|\t^{\Ldom\lambda^{(\ell)}})$ is
  the standard $\blam$-tableau with numbers $1,2,\dots,n$ entered in
  order down the columns of the
  tableaux~$\t^{\Ldom\lambda^{(\ell)}},\dots,\t^{\Ldom\lambda^{(1)}}$.
  By construction, $\Ltlam=(\Gtlam[\blam'])'$.
  \notation{$\Ltlam,\Gtlam$}{Initial tableau with respect to $\Ldom$ and $\Gdom$}

  \begin{Definition}
    For each standard tableau $\t\in\Std(\Parts)$, let
    $\Gdt,\Ldt\in\Sym_n$ be the unique permutations such that
    $\Ldt\Ltlam=\t=\Gdt\Gtlam$. As important special cases, set
    $\Lwlam=\Ldt[\Gtlam]$ and $\Gwlam=\Gdt[\Ltlam]$.
  \end{Definition}
  \notation{$\Ldt,\Gdt$}{Permutations: $\Ldt\Ltlam=\t=\Gdt\Gtlam$, for $\t\in\Std(\Parts)$}

  Recall from \autoref{S:Quiver} that $L\map{\Sym_n}\N$ is the length
  function on~$\Sym_n$. Although normally stated using slightly
  different language, the following lemma is well-known and easy to
  prove.  See, for example, \cite[Lemma~2.18]{KMR:UniversalSpecht}.

  \begin{Lemma}\label{L:wlam}
    Suppose $\blam\in\Parts$. Then
    $\Gwlam =(\Lwlam)^{-1}$. Moreover, if $\t\in\Std(\blam)$ then
    \[
        \Lwlam = (\Gdt)^{-1}\Ldt,\qquad
        \Gwlam = (\Ldt)^{-1}\Gdt,\And
        \Ldt = \Gdt[\t'],
    \]
    with $L(\Lwlam)=L(\Ldt)+L(\Gdt)=L(\Gwlam)$.
  \end{Lemma}

  In \autoref{S:KLRproperties}, we fixed a preferred reduced expression
  $w=\sigma_{a_1}\dots \sigma_{a_l}$, for each $w\in\Sym_n$, and we defined
  $\psi_w=\psi_{a_1}\dots\psi_{a_l}$. In particular, we have preferred
  reduced expressions for the permutations
  $\Ldt$, $\Lwlam$, $\Gdt$ and $\Gwlam$ that define elements
  $\psi_{\Ldt}, \psi_{\Lwlam}, \psi_{\Gdt}, \psi_{\Gwlam}\in\Rx(\kx)$.

  Recall from \autoref{S:Tableaux} that
  $\Nodes=\set{(k,r,c)|1\le k\le\ell, r,c\ge1}$ is the
  set of nodes, which we consider as a totally ordered set under the
  lexicographic order, and that we identify an $\ell$-partition with its
  diagram $\set{(k,r,c)\in\Nodes|1\le c\le\lambda^{(k)}_r}$.

  Fix $\blam\in\Parts$. An \emph{addable} node of $\blam$ is a node
  $A=(k,r,c)\in\Nodes\setminus\blam$ such that
  $\blam\cup\set{A}\in\Parts[n+1]$. Similarly, a \emph{removable} node
  of~$\blam$ is a node $A\in\blam$ such that
  $\blam\setminus\set{A}\in\Parts[n-1]$.  If $\t\in\Std(\blam)$ let
  $\Add(\t)=\Add(\blam)$ and $\Rem(\t)=\Rem(\blam)$ be the sets of
  addable and removable nodes of~$\blam$.

  Let $\t\in\Std(\blam)$ and $1\le m\le n$ and define:
  \begin{align}\label{E:Sddable}
    \begin{split}
      \LAdd_m(\t)&=\set[\big]{A\in\Add(\t_{\downarrow m})\,|\br(A)=\br_m(\t)
                          \text{ and } A<\t^{-1}(m)}\\
      \LRem_m(\t)&=\set[\big]{A\in\Rem(\t_{\downarrow m})|\br(A)=\br_m(\t)
                          \text{ and } A<\t^{-1}(m)}\\
      \GAdd_m(\t)&=\set[\big]{A\in\Add(\t_{\downarrow m})\,|\br(A)=\br_m(\t)
                          \text{ and } A>\t^{-1}(m)}\\
      \GRem_m(\t)&=\set[\big]{A\in\Rem(\t_{\downarrow m})|\br(A)=\br_m(\t)
                          \text{ and } A>\t^{-1}(m)}.
    \end{split}
  \end{align}

  Recall from \autoref{S:KLR} that $*$ is the unique anti-isomorphism of
  $\Rx$ that fixes the generators of \autoref{D:KLR}.

  \notation{$\Lilam,\Gilam$}
     {Residue sequences: $\Lilam=\br(\Ltlam)$ and $\Gilam=\br(\Gtlam)$}[D:psist]
  \notation{$\Lylam,\Gylam$}{Polynomials $\Lylam,\Gylam\in\k[y_1,\dots,y_n]$}[D:psist]
  \notation{$\Lpsist,\Gpsist$}{The basis elements
    $\psi_{\Ldt[\s]}\Lylam \ei[\Lilam]\psi_{\Ldt}^*$ and
    $\psi_{\Gdt[\s]}\Gylam \ei[\Gilam]\psi_{\Gdt}^*$}[D:psist]

  \begin{Definition}[Integral bases]\label{D:psist}
    Let $\s,\t\in\Std(\blam)$, for $\blam\in\Parts$. Define
    \[
      \Lpsist = \psi_{\Ldt[\s]}\Lylam \ei[\Lilam]\psi_{\Ldt}^*
      \And
      \Gpsist = \psi_{\Gdt[\s]}\Gylam \ei[\Gilam]\psi_{\Gdt}^*,
    \]
    where $\Lilam=\br(\Ltlam)$, $\Gilam=\br(\Gtlam)$
    and
     \[
       \Lylam = \prod_{m=1}^n\prod_{A\in\Add_m^\Ldom(\Ltlam)}
       \bigl(y_m-\bc(A)\bigr)
       \And
       \Gylam = \prod_{m=1}^n\prod_{A\in\Add_m^\Gdom(\Gtlam)}
                  \bigl(y_m-\bc(A)\bigr).
     \]
  \end{Definition}

  By definition, if $(\s,\t)\in\Std^2(\Parts)$ then $\Lpsist$ and
  $\Gpsist$ are elements of~$\Rx(\kx)$, which depend on the choices of
  reduced expressions for $\Ldt[\s]$, $\Ldt$, $\Gdt[\s]$ and $\Gdt$.
  We will abuse notation and consider $\Lpsist$ and $\Gpsist$ as
  elements of $\Rx(\kx)$, $\Rx(\Kxx)$ and of $\Rn(\k)$. It is not yet
  clear that the elements $\Lpsist$ and $\Gpsist$ are nonzero but, if
  they are, they are homogeneous.

  To prove that $\set{\Lpsist}$ and $\set{\Gpsist}$ are bases of
  $\Rx(\kx)$ we will use some closely related \emph{seminormal bases}
  of~$\Rx(\Kxx)$.  As we will see, the seminormal bases give other
  realisations of the graded $\Rx(\Kxx)$-modules $V_\blam$ from
  \autoref{P:SeminormalForm}. In fact, this is the key to proving that
  the $\psi$-bases are linearly independent.

  \notation{$\Lfst,\Gfst$}{The basis elements $\Lfst=F_\s\Lpsist F_\t$ and
    $\Gfst=F_\s\Gpsist F_\t$, for $\s,\t\in\Std(\blam)$}[D:fst]

  \begin{Definition}[Seminormal bases]\label{D:fst}
    Let $\s,\t\in\Std(\blam)$, for $\blam\in\Parts$. Set
    \[
    \Lfst= F_\s\Lpsist F_\t \qquad\text{and}\qquad  \Gfst= F_\s\Gpsist F_\t.
    \]
  \end{Definition}

  By definition, $\Lfst, \Gfst\in\Rx(\Kxx)$ and these elements do not
  typically belong to~$\Rx(\kx)$.  We will show that $\set{\Lfst}$ and
  $\set{\Gfst}$ are cellular bases of~$\Rx(\Kxx)$. Since $\Lpsist$ and
  $\Gpsist$ are both homogeneous so are $\Lfst$ and $\Gfst$.

  Below we prove many parallel results for the elements $\set{\Lpsist}$
  and $\set{\Lfst}$, and for the elements $\set{\Gpsist}$
  and~$\set{\Gfst}$. In almost every case, the proofs are identical
  except that the $\psi^\Ldom$-basis and $f^\Ldom$-basis use the poset
  $(\Parts,\Ledom)$ whereas the $\psi^\Gdom$-basis and $f^\Gdom$-basis use
  the poset $(\Parts,\Gedom)$. For this reason, we work with a generic
  symbol $\Domin$ and write $\Dpsist$, $\Dfst$, $\Dtlam$, $\Ddt$,
  $\dots$ in place of $\Lpsist$, $\Lfst$, $\Ltlam$, $\Ldt$, $\dots$  and
  $\Gpsist$, $\Gfst$, $\Gtlam$, $\Gdt$, $\dots$, respectively.

  \begin{Lemma}\label{L:FufstFv}
    Let $\s,\t,\u,\v\in\Std(\Parts)$. Then
    $\delta_{\s\u}\delta_{\t\v}\Lfst = F_\u\Lfst F_\v$ and
    $\delta_{\s\u}\delta_{\t\v}\Gfst = F_\u\Gfst F_\v$.
  \end{Lemma}

  \begin{proof}
    This is immediate from \autoref{C:FsFt} and \autoref{D:fst}.
  \end{proof}

  In contrast, it is rarely true that
  $F_\u\Dpsist F_\v=\delta_{\s\u}\delta_{\t\v}\Dpsist$, for
  $(\s,\t),(\u,\v)\in\Std^2(\Parts)$.

  We want to show that the sets $\set{\Dpsist}$ and $\set{\Dfst}$ are bases of
  $\Rx(\Kxx)$ and that the transition matrices between the $\psi$-bases
  and the corresponding $f$-bases are unitriangular. Before we can prove
  this we need a better understanding of how $\Rx(\Kxx)$ acts on the $f$-bases
  and to do this we connect these bases to the seminormal
  representations of \autoref{S:Contents}.  Motivated
  by~\autoref{E:betas}, for $\s\in\Std(\Parts)$ and $1\le k<n$ define
  scalars $\Lbeta_k(\s), \Gbeta_k(\s)\in\kx$ by
  \begin{equation}\label{E:LGbetas}
    \Lbeta_k(\s) = \begin{cases*}
        1       & if $\s\Ldom \sigma_k\s$,\\
        Q_k(\s) & if $\sigma_k\s\Ldom \s$,
    \end{cases*}
    \And
    \Gbeta_k(\s) = \begin{cases*}
        1       & if $\s\Gdom \sigma_k\s$,\\
        Q_k(\s) & if $\sigma_k\s\Gdom \s$.
    \end{cases*}
  \end{equation}
  Repeating the argument of \autoref{L:betas} shows that:

  \begin{Lemma}
    The coefficients $\set{\Lbeta_r(\s)}$ and
    $\set{\Gbeta_r(\s)}$ satisfy conditions \upshape{(a)--(c)} of
    \autoref{P:SeminormalForm}.
  \end{Lemma}


  Hence, the coefficients $\set{\Lbeta_r(\s)}$ and
  $\set{\Gbeta_r(\s)}$ both determine irreducible graded
  $\Rx(\Kxx)$-modules $V_\blam^\Ldom$ and $V_\blam^\Gdom$, respectively.
  By \autoref{C:IsomorphicSpechts}, $V_\blam^\Ldom\cong V_\blam^\Gdom$.
  Let $\set{v^{\Dom}_\t|\t\in\Std(\blam)}$ be the basis
  of~$V^{\Dom}_\blam$ from \autoref{P:SeminormalForm}. More explicitly,
  fix a nonzero vector $v_{\Dtlam}\in F_{\Dtlam}V^{\Dom}_\blam$ and define
  $v^{\Dom}_\t$ by induction on $L(\Ddt)$ by setting
  \[
     v^\Dom_\t = \Bigl(\psi_k-\frac{\delta_{\br_k(\s)}\delta_{k+1}(\s)}{\rho_{k}(\s)}\Bigr)
                    v^\Dom_\s
  \]
  where $\Ddt=s_k\Ddt[\s]$ with $L(\Ddt)=L(\Ddt[\s])+1$, and we set
  $\rho_k(\s)=\bc_{k+1}(\s)-\bc_k(\s)\in\kx$.
  \notation{$\rho_k(\t)$}{The difference $\bc_{k+1}(\s)-\bc_k(\s)\in\kx$ }

  The next result should be compared with \autoref{P:SeminormalForm}.

  \begin{Proposition}\label{P:fstaction}
     Let $(\s,\t)\in\Std^2(\Parts)$ and suppose that $1\le k<n$, $1\le
     m\le n$ and $\bi\in I^n$. Then the elements $\Lfst$ and $\Gfst$ are
     nonzero and
     \begin{align*}
       \ei\Lfst &=\delta_{\bi\,\br(\s)}\Lfst &
       y_m\Lfst&=c_m(\s)\Lfst &
       \psi_k\Lfst &=\frac{\delta_{\br_k(\s),\br_{k+1}(\s)}}{\rho_k(\s)}\Lfst
                        +\Lbeta_k(\s)\Lfst[\u\t],\\
                        \ei\Gfst &=\delta_{\bi\,\br(\s)}\Gfst &
       y_m\Gfst&=c_m(\s)\Gfst &
       \psi_k\Gfst &=\frac{\delta_{\br_k(\s),\br_{k+1}(\s)}}{\rho_k(\s)}\Gfst
                       +\Gbeta_k(\s)\Gfst[\u\t],
     \end{align*}
     where $\u=\sigma_k\s$.
  \end{Proposition}

  \begin{proof}
    Let $\Domin$. Since $\Dfst=F_\s\Dpsist F_\t$, the formulas for
    $\ei\Dfst$ and $y_m\Dfst$ follow from \autoref{C:KLRIdempotents}
    and \autoref{C:yFt}, respectively. We use these formulas below
    without mention.


    To prove the remaining claims, fix $\t\in\Std(\blam)$ and let
    $W^{\Dom}_\t$ be the $\Kxx$-submodule of $\Rx(\Kxx)$ spanned by
    $\set{\Dfst|\s\in\Std(\blam)}$. Let $\Theta_\t\colon W^{\Dom}_\t\to
    V^{\Dom}$ be the map given by $\Theta_\t(w)=wv^{\Dom}_{\t}$, for
    $w\in W^{\Dom}_\t$.  We prove by induction on dominance order
    for~$\t$ that there exists a nonzero scalar $a_\t$, which depends
    only on $\t$, such that $\Theta_\t(\Dfst)=a_\t v^{\Dom}_\s$, for
    $\s\in\Std(\blam)$. To prove this, first consider the special case
    when $\t=\Dtlam$. By \autoref{P:SeminormalForm},
    \[
       \Dpsist[\Dtlam\Dtlam]v^{\Dom}_{\Dtlam} =\Dylam \ei[\Dilam]v^{\Dom}_{\Dtlam}
               =\prod_{m=1}^n\prod_{A\in\DAdd_m(\Dtlam)}
                \bigl(\bc_m(\Dtlam)-\bc(A)\bigr)\cdot v^{\Dom}_{\Dtlam}
               =a_{\Dtlam} v^{\Dom}_{\Dtlam},
    \]
    where $a_{\Dtlam}=\prod_m\prod_A(\bc_m(\Dtlam)-\bc(A))\in\kx$. If
    $A\in\DAdd_m(\Dtlam)$ then $\br(A)=\br_m(\Dtlam)$, so each factor
    of~$a_{\Dtlam}$ is nonzero by \autoref{D:ContentSystem}(c).  Consequently,
    $a_{\Dtlam}\ne0$. Moreover,
    $\Dfst[\Dtlam\Dtlam]v^{\Dom}_{\Dtlam}=a_{\Dtlam} v^{\Dom}_{\Dtlam\Dtlam}$ since
    $F_{\s}v^{\Dom}_{\s}=v^{\Dom}_\s$, for all $\s\in\Std(\blam)$.
    In view of \autoref{E:LGbetas} and \autoref{P:SeminormalForm}, if
    $\y\in\Std(\blam)$ then
    \[
        \Dfst[\y\Dtlam]v^{\Dom}_{\Dtlam}
        = F_\y\psi_{\Ddt[\y]}\Dfst[\Dtlam\Dtlam]v^{\Dom}_{\Dtlam}
        = a_{\Dtlam} F_\y\psi_{\Ddt[\y]}v^{\Dom}_{\Dtlam}
        = a_{\Dtlam} v^{\Dom}_{\y},
    \]
    where the last equality uses \autoref{L:FufstFv}.  It follows that
    $\Theta_{\Dtlam}$ is multiplication by $a_{\Dtlam}$. In particular,
    the map~$\Theta_{\Dtlam}$ is an $\Rx(\Kxx)$-module isomorphism
    and $W^\Dom_{\Dtlam}\cong V^\Dom_\blam$, which implies the desired
    formulas for $\psi_k\Dfst[\s\Dtlam]$ by \autoref{P:SeminormalForm}
    and \autoref{E:LGbetas}.

    Finally, suppose that $\t\ne\Dtlam$ and let
    $\Ddt=\sigma_{a_1}\dots \sigma_{a_k}$ be the preferred reduced
    expression that we fixed for the permutation~$\Ddt\in\Sym_n$ in
    \autoref{S:KLRproperties}.  Recalling the definition of~$Q_m(\t)$
    from \autoref{E:Qmt}, define
    \[
      Q(\t) = Q_{a_1}(\sigma_{a_1}\t)Q_{a_2}(\sigma_{a_2}\sigma_{a_1}\t)
          \dots Q_{a_k}(\sigma_{a_k}\dots \sigma_{a_1}\t).
    \]
    Then $Q(\t)\ne0$ by \autoref{L:Qnonzero}.  Applying
    \autoref{P:SeminormalForm}(b) $k$~times,
    \[
       \Dpsist v^{\Dom}_\t
            = \psi_{\Ddt[\s]}\Dpsist[\Dtlam\Dtlam]\psi_{\Ddt}^* v^{\Dom}_\t
            = Q(\t)\psi_{\Ddt[\s]}\Dpsist[\Dtlam\Dtlam] v^{\Dom}_{\Dtlam}
            = a_{\Dtlam} Q(\t)\psi_{\Ddt[\s]} v^{\Dom}_{\Dtlam}
            = a_{\Dtlam} Q(\t) v^{\Dom}_{\s}.
    \]
    Therefore, $\Theta_\t$ is multiplication by the scalar $a_\t=a_{\Dtlam}
    Q(\t)$, so $\Theta_\t\colon W^{\Dom}_\t\xrightarrow{\sim}V^{\Dom}$ is an
    isomorphism.  Hence, the formula for $\psi_k\Dfst$ follows from
    \autoref{P:SeminormalForm}. The proof is complete.
  \end{proof}

  Since $\Dfst=F_\s\Dpsist F_\t$, this also shows that $\Lpsist$ and
  $\Gpsist$ are nonzero, for $(\s,\t)\in\Std^2(\Parts)$.  Although we do
  not state them explicitly, applying the automorphism~$*$ to
  \autoref{P:fstaction} gives similar formulas for the right actions of
  the generators of $\Rx(\Kxx)$ on the $f$-bases.

  The first corollary of \autoref{P:fstaction} was established in its proof.

  \begin{Corollary}\label{C:Vblam}
    Let $\blam\in\Parts$ and suppose $\t\in\Std(\blam)$. Then, as
    $\Rx(\Kxx)$-modules.
    \[
        V^{\Ldom}_\blam\cong\bigoplus_\y\Kxx\Lfst[\y\t]
        \qquad\text{and}\qquad
        V^{\Gdom}_\blam\cong\bigoplus_\y\Kxx\Gfst[\y\t].
    \]
  \end{Corollary}

  \begin{Corollary}\label{C:fBases}
    The sets
    $\set{\Lfst|(\s,\t)\in\Std^2(\Parts)}$ and
    $\set{\Gfst|(\s,\t)\in\Std^2(\Parts)}$
    are bases of~$\Rx(\Kxx)$.
  \end{Corollary}

  \begin{proof}
    Let $\bi\in I^n$. By \autoref{C:KLRIdempotents}, $\ei\ne0$ if
    and only if $\bi\in\StI[n]=\set{\br(\u)|\u\in\Std(\Parts)}$.
    Moreover, if $\bi\in\StI[n]$ then $\ei=\sum_{\u\in\Std(\bi)} F_\u$.
    Hence, as $\Kxx$-modules,
    \[
      \Rx(\Kxx)=\bigoplus_{\bi,\bj\in\StI} \ei\Rx(\Kxx)\ei[\bj]
           =\sum_{\s,\t\in\Std(\Parts)} \Rx[\s\t],
           \qquad\text{where }\Rx[\s\t]=F_\s\Rx(\Kxx)F_\t.
    \]
    If $(\s,t)\in\Std^2(\Parts)$ then $\Dfst\ne0$, by
    \autoref{P:fstaction}, and $\Dfst\in \Rx[\s\t]$, by
    \autoref{C:FsFt}.  $(\s,\t)\in\Std^2(\Parts)$. Hence, $\set{\Dfst}$
    is a basis of $\Rx(\Kxx)$ and the last displayed equation becomes
    $\Rx(\Kxx)=\bigoplus_{(\s,\t)\in\Std^2(\Parts)}\Rx[\s\t]$.
  \end{proof}

  The next result shows that the idempotents $F_\t$ are scalar multiples
  of the basis elements $\Lfst[\t\t]$ and $\Gfst[\t\t]$. These scalars,
  $\Lgt$ and $\Ggt$, play an important role in what follows.

  \notation{$\Lgt,\Ggt$}{Important monomials in $\Kxx$,
      for $\t\in\Std(\Parts)$}[C:gamma]

  \begin{Corollary}\label{C:gamma}
    Suppose that $\t\in\Std(\Parts)$. Then there exist nonzero homogeneous scalars
    $\Lgt,\Ggt\in\Kxx$ such that
    \[
            \frac1{\Lgt}\Lfst[\t\t] = F_\t = \frac1{\Ggt}\Gfst[\t\t].
    \]
  \end{Corollary}

  \begin{proof}
    Let $\Domin$. By \autoref{C:fBases},
    $F_\t=\sum_{\u,v}r_{\u\v}\Dfst[\u\v]$, for some $r_{\u\v}\in\Kxx$.
    Multiplying on the left and right by $F_\t$ and applying
    \autoref{L:FufstFv} and \autoref{C:Ftnonzero} shows that
    $F_\t=r_{\t\t}\Dfst[\t\t]$. By \autoref{C:Ftnonzero}, $r_{\t\t}\ne0$.
    Therefore, setting $\Dgt=\tfrac1{r_{\t\t}}$ gives the result.
  \end{proof}

  \begin{Lemma}\label{L:fproduct}
    Suppose that $(\s,\t), (\u,\v)\in\Std^2(\Parts)$. Then
    \[
          \Lfst \Lfst[\u\v] = \delta_{\t\u}\Lgt\Lfst[\s\v]
          \And
          \Gfst \Gfst[\u\v] = \delta_{\t\u}\Ggt\Gfst[\s\v].
    \]
  \end{Lemma}

  \begin{proof}
    Let $\Domin$. If $\u\ne\t$ then
    $\Dfst\Dfst[\u\v]=\bigl(\Dfst F_\t\bigr)\Dfst[\u\v]
                     =\Dfst\bigl(F_\t \Dfst[\u\v]\bigr)=0$,
    where we have used \autoref{L:FufstFv} twice. Hence, it remains to
    consider the products $\Dfst[\s\t]\Dfst[\t\v]$. In particular, $\s$,
    $\t$ and $\v$ all have the same shape.

    By \autoref{P:fstaction}, for $\u\in\Std(\blam)$ there exist homogeneous elements
    $p_\u,q_\u\in\Rx(\Kxx)$, which are independent of~$\t$, such that
    $\Dfst[\u\t]=p_\u\Dfst[\Dtlam\t]$ and
    $\Dfst[\Dtlam\t]=q_\u\Dfst[\u\t]$. Therefore, using
    \autoref{C:gamma} and \autoref{L:FufstFv},
    \[
       \Dfst[\s\t]\Dfst[\t\v]
         = p_\s \Dfst[\Dtlam\t]\Dfst[\t\v]
         = p_\s q_\t\Dfst[\t\t]\Dfst[\t\v]
         = \Dgt p_\s q_\t F_\t\Dfst[\t\v]
         = \Dgt p_\s q_\t\Dfst[\t\v]
         = \Dgt\Dfst[\s\v],
    \]
    as required.
  \end{proof}

  We need to determine the $\gamma$-coefficients explicitly, which is
  possible because they satisfy the following recurrence relation
  involving the scalars $Q_k(\s)$ from \autoref{E:Qmt}. Note that
  $Q_k(\s)\ne0$ whenever $\sigma_k\s$ is standard by
  \autoref{L:Qnonzero}.

  \begin{Lemma}\label{L:gammaRecurrence}
    Let $\Domin$ and suppose that $\s,\t\in\Std(\Parts)$ with
    $\s\Dom\t=\sigma_k\s$, where $1\le k <n$. Then
    $\Dgt=Q_k(\s)\Dgt[\s]$.
  \end{Lemma}

  \begin{proof}
    By \autoref{E:LGbetas}, $\Dbeta_k(\s)=1$. Therefore, using
    \autoref{L:fproduct} and \autoref{P:fstaction} several times,
      \begin{align*}
        \Dgt \Dfst[\s\s]
          &= \Dfst\Dfst[\t\s]
           = \Dfst[\s\s]\Bigl(\psi_k-\frac{\delta_{\br_k(\s)\br_{k+1}(\s)}}
                                 {\rho_k(\s)}\Bigr)^2\Dfst[\s\s]\\
          &= \Dfst[\s\s]\Bigl(\psi_k^2-\frac{2\psi_k\delta_{\br_k(\s)\br_{k+1}(\s)}}
                                             {\rho_k(\s)}
               + \frac{\delta_{\br_k(\s)\br_{k+1}(\s)}}
                     {\rho_k(\s)^2}\Bigr)\Dfst[\s\s]\\
           &= \Dfst[\s\s]\Bigl(Q^\ux_{\br_k(\s)\br_{k+1}(\s)}(c_k(\s), c_{k+1}(\s))-
                   \frac{\delta_{\br_k(\s)\br_{k+1}(\s)}}
                          {\rho_k(\s)^2}\Bigr)\Dfst[\s\s]\\
           &= Q_k(\s)\Dgt[\s]\Dfst[\s\s].
      \end{align*}
      For the third equality, notice that $\psi_k\Dfst[\s\s]$ introduces
      a term involving $\Dfst[\t\s]$ but this term does not survive because
      $\Dfst[\s\s]\Dfst[\t\s]=0$ by \autoref{L:fproduct}.  The result now
      follows by \autoref{C:fBases}.
  \end{proof}

  \begin{Lemma}\label{L:gammaClosed}
  Suppose that $\t\in\Std(\blam)$, for $\blam\in\Parts$. Then
  \[
        \Lgt = \prod_{m=1}^n
           \dfrac{\prod_{A\in\LAdd_m(\t)}(\bc_m(\t)-\bc(A))}
                 {\prod_{B\in\LRem_m(\t)}(\bc_m(\t)-\bc(B))}
        \And
        \Ggt = \prod_{m=1}^n
           \dfrac{\prod_{A\in\GAdd_m(\t)}(\bc_m(\t)-\bc(A))}
                 {\prod_{B\in\GRem_m(\t)}(\bc_m(\t)-\bc(B))}.
  \]
  \end{Lemma}

  \begin{proof}
    We consider only the result for $\Lgt$ and leave the symmetric case
    of $\Ggt$ to the reader.  We argue by induction on dominance. If
    $\t=\Ltlam$ then $\Lfst[\Ltlam\Ltlam]=\Lylam \ei[\Lilam]$.
    Therefore, by
    \autoref{L:fproduct} and \autoref{P:fstaction},
    \[
        \Lgt[\Ltlam]\Lfst[\Ltlam\Ltlam]
           =\Lfst[\Ltlam\Ltlam]\Lfst[\Ltlam\Ltlam]
           =\Lylam\Lfst[\Ltlam\Ltlam]
           =\prod_{m=1}^n\prod_{A\in\LAdd_m(\Ltlam)}
              \bigl(\bc_m(\Ltlam)-\bc(A)\bigr)\Lfst[\Ltlam\Ltlam].
    \]
    As $\LRem_m(\Ltlam)=\emptyset$, for $1\le m\le n$, this gives the
    result when $\t=\Ltlam$. If $\t\Gdom\Ltlam$ then, by
    \autoref{L:gammaRecurrence}, there exists a tableau~$\s$ and an
    integer $a$, with $1\le a<n$, such that $\s\Ldom\t=\sigma_a\s$ and
    $\Lgt=Q_a(\s)\Lgt[\s]$. To complete the proof, write
    $(k,r,c)=\t^{-1}(a)$ and observe that
    $\LAdd_m(\t)=\LAdd_m(\s)$ and $\LRem_m(\t)=\LRem_m(\s)$ if $m\ne a, a+1$.
    Moreover, $\LAdd_{a}(\t)=\LAdd_{a+1}(\s)$ and $\LRem_{a}(\t)=\LRem_{a+1}(\s)$
    and
    \begin{align*}
      \LAdd_{a+1}(\t) &= \begin{cases*}
         \LAdd_{a}(\s)\setminus\set{(k,r,c)}, & if $\br_a(\s)=\br(k,r,c)$,\\
         \LAdd_{a}(\s)\cup A, & otherwise,
       \end{cases*}  \\
    \intertext{where $A$ is the set of addable $\br_a(\s)$-nodes in
    $\set{(k,r+1,c),(k,r,c-1)}$. Similarly,}
       \LRem_{a+1}(\t) &= \begin{cases*}
         \LRem_{a}(\s)\cup\set{(k,r,c)}, & if $\br_a(\s)=\br(k,r,c)$,\\
         \LRem_{a}(\s)\setminus R, & othewise,
       \end{cases*}
    \end{align*}
    where $R$ is the set of removable $\br_a(\s)$-nodes in
    $\set{(k,r+1,c),(k,r,c-1)}$.  By induction, the lemma holds for
    $\Lgt[\s]$.  Hence, recalling the definition of~$Q_a(\s)$ from
    \autoref{E:Qmt}, the lemma holds for $\Lgt$ since
    $\Lgt=Q_a(\s)\Lgt[\s]$. This completes the proof.
  \end{proof}

  We can now compute the transition matrices between the $\psi$-bases
  and the corresponding $f$-bases.

  \begin{Proposition}\label{P:psiTriangular}
    Suppose that $\s,\t\in\Std^2(\blam)$, for $\blam\in\Parts$. In $\Rx(\Kxx)$,
    \[
      \Lpsist=\Lfst +\sum_{\substack{\bmu\Ledom\blam\\(\u,\v)\in\Std^2(\bmu)}}
                  a_{\u\v}\Lfst[\u\v]
      \And
      \Gpsist=\Gfst +\sum_{\substack{\bmu\Gedom\blam\\(\u,\v)\in\Std^2(\bmu)}}
                  b_{\u\v}\Gfst[\u\v]
    \]
    for homogeneous coefficients in $\Kxx$ such that
    \begin{itemize}
      \item $a_{\u\v}\ne0$ only if
        $\br(\u)=\br(\s)$, $\br(\v)=\br(\t)$ and either
        $\bmu\Ldom\blam$, or $\bmu=\blam$, $\u\Ldom\s$ and $\v=\t$,
      \item $b_{\u\v}\ne0$ only if
        $\br(\u)=\br(\s)$, $\br(\v)=\br(\t)$ and either
        $\bmu\Gdom\blam$, or $\bmu=\blam$, $\u\Gdom\s$ and $\v=\t$.
    \end{itemize}
  \end{Proposition}

  \begin{proof}
    Let $\Domin$.  By \autoref{T:ContentReduced} and \autoref{C:gamma},
    $\ei[\Dilam]=\sum_\u F_\u =\sum_\u \frac1{\Dgt[\u]}\Dfst[\u\u]$,
    where both sums are over $\u\in\Std(\Dilam)$.
    Using \autoref{P:fstaction},
    \[
        \Dpsist[\Dtlam\Dtlam]=\Dylam \ei[\Dilam]
          =\sum_{\u\in\Std(\Dilam)}\frac1{\Dgt}\Dylam\Dfst[\u\u]
          =\sum_{\u\in\Std(\Dilam)}\frac1{\Dgt}
                \prod_{m=1}^n\prod_{A\in\DAdd_m(\Dtlam)}
                \bigl(\bc_m(\u)-\bc(A)\bigr)\Dfst[\u\u]
    \]
    for some $a_\u\in\Kxx$. If $\u=\Dtlam$ then the coefficient of
    $\Dfst[\u\u]$ in the displayed equation is~$1$ by \autoref{L:gammaClosed}.
    Now suppose that $\u\in\Std(\Dilam)$ and $\u\Domneq\Dtlam$. Let
    $m$ be minimal such that $\t_{\downarrow m}\ne(\Dtlam)_{\downarrow m}$.
    Then $A=\u^{-1}(m)\in\DAdd_m(\Dtlam)$, so~$\Dfst[\u\u]$
    appears in $\Dylam \ei[\Dilam]$ with coefficient zero. Hence,
    $\Dfst[\u\u]$ appears in $\Dpsist[\Dtlam\Dtlam]$ with nonzero
    coefficient only if $\u\Domeq\Dtlam$, so $\blam\Dom\Shape(\u)$ if
    $\u\ne\Dtlam$. This proves the base case of our induction.
    If $\s,\t\in\Std(\blam)$ then
    \[
      \Dpsist = \psi_{\Ddt[\s]}\Dpsist[\Dtlam\Dtlam]\psi^*_{\Ddt[\t]}
        = \psi_{\Ddt[\s]}\Bigl(
              \Dfst[\Dtlam\Dtlam] +\sum_{\u\Dom\Dtlam}a_\u\Dfst[\u\u]
            \Bigl)\psi^*_{\Ddt[\t]}.
    \]
    Hence, the result follows by \autoref{P:fstaction} and induction on
    $\blam$.
  \end{proof}

  By \autoref{C:fBases}, this implies that
  $\set{\Lpsist}$ and $\set{\Gpsist}$ are both bases of $\Rx(\Kxx)$.

  \subsection{Cellular algebras}\label{S:CellularAlgebras} K\"onig and
  Xi~\cite{KoenigXi:AffineCellular} introduced affine cellular algebras,
  generalising results of Graham and Lehrer~\cite{GL}.
  Following~\cite{HuMathas:GradedCellular}, this section incorporates a
  grading into this framework and at the same time allows the ground
  ring~$K$ to have a non-trivial grading.  The next section shows that
  the~$\Dfst[]$ and~$\Dpsist[]$-bases induce $K$-cellular structures on the
  algebras~$\Rx(\Kxx)$ and~$\Rx(\kx)$.

  \begin{Definition}[cf. Graham and Lehrer, K\"onig and
       Xi~\cite{GL,HuMathas:GradedCellular,KoenigXi:AffineCellular}]
    \label{graded cellular def}\label{D:cellular}
    Let $K$ be a graded commutative domain with~$1$ and suppose that $A$ is a
    graded $K$-algebra that is $K$-free and of finite rank as a $K$-module. A
    \textbf{graded $K$-cell datum} for $A$ is an ordered tuple $(P,T,a,\deg)$,
    where $(P,>)$ is the \textbf{weight poset},
    $T=\coprod_{\lambda\in P}T_\lambda$ is a finite set,
    $$a\map{\coprod_{\lambda\in P}T_\lambda\times T_\lambda}A;
              (s,t)\mapsto a_{st},
    $$
    is an injective map and $\deg\map{T}\Z$
    is a \textbf{degree function} such that:
    \begin{enumerate}[label=\upshape(C$_{\arabic*}$), ref=C$_{\arabic*}$, start=0]
      \item\label{CA:deg} If $s,t\in T_\lambda$ then
      $a_{st}$ is homogeneous of \emph{degree} $\deg(a_{st})=\deg(s)+\deg(\t)$.
      \item\label{CA:basis} The set
      $\set{a_{st}|s,t\in T_\lambda \text{ for } \lambda\in P}$
    is a $K$-basis of $A$.
      \item\label{CA:act} Let $h\in A$ be homogeneous and fix
      $s,t\in T_\lambda$, for $\lambda\in P$. There exist (homogeneous)
      scalars $r_{su}(h)\in K$, which do not depend on $t$, such that
      \[
           ha_{st} =\sum_{u\in T_\lambda}
                         r_{us}(h)a_{ut}\pmod {A^{>\lambda}},
      \]
      where $A^{>\lambda}$ is the $K$-submodule of $A$ spanned by
      $\set{a_{vw}|\mu>\lambda\text{ and }v,w\in T(\mu)}$.
      \item\label{CA:star} The $K$-linear map $*\map AA$ determined by
        $(a_{st})^*=a_{ts}$, for all $\lambda\in P$ and
       $s,t\in T_\lambda$, is an anti-isomorphism of $A$.
    \end{enumerate}
    A \emph{graded $K$-cellular algebra} is an algebra that has a
    graded $K$-cell datum. A \emph{$K$-cellular algebra} is an
    algebra that has a graded $K$-cell datum such that $\deg(t)=0$ for
    all $t\in T$.  A \emph{(graded) cellular algebra} is an algebra
    that has a (graded) $K$-cell datum when $K=K_0$ is concentrated
    in degree~$0$.
  \end{Definition}

  \begin{Remark}
    If $K=K_0$ is concentrated in degree~$0$ then a graded
    $K$-cellular algebra is a graded cellular algebra in the sense of
    \cite{HuMathas:GradedCellular}. If $K=K_0$ and $\deg(t)=0$ for all $t\in T$
    we recover the cellular algebras of Graham and Lehrer~\cite{GL}. A
    $K$-cellular algebra is a graded analogue of the
    affine cellular algebras of K\"onig and
    Xi~\cite{KoenigXi:AffineCellular} in the special case where their
    affine commutative algebra $B$ is $K$ considered as a $K_0$-algebra.
\end{Remark}

  If $L$ is a $K$-algebra, define $A(L)=L\otimes_KA$. Then $A(L)$ is
  a (graded) $L$-cellular algebra.

  Let $A=A(K)$ be a graded $K$-cellular algebra with graded $K$-cell datum
  $(P,T,c,\deg)$. As in \autoref{CA:act}, for $\lambda\in P$ let
  $A^{\ge\lambda}(K)$ be the $K$-submodule of~$A$ spanned by
  $\set{a_{st}|s,t\in T(\mu)\text{ for }\mu\ge\lambda}$. By
  \autoref{CA:act} and \autoref{CA:star}, $A^{\ge\lambda}(K)$ and
  $A^{>\lambda}(K)=\bigoplus_{\mu>\lambda}A^{\ge\mu}(K)$ are two-sided ideals
  of~$A$. Set $A_\lambda(K)=A^{\ge\lambda}(K)/A^{>\lambda}(K)$.

  For $\lambda\in P$, the
  \emph{cell module} $S_\lambda(K)$ is the free $K$-module with basis
  $\set{a_s|s\in T(\lambda)}$, where $a_s$ is homogeneous of degree
  $\deg(s)$, and where the $A$-action on $S_\lambda(K)$ is given by
  \[
     ha_{s} =\sum_{u\in T(\lambda)}r_{us}(h)a_{u},\qquad
     \text{for $h\in A$ and $s\in T(\lambda)$,}
  \]
  where $r_{us}(h)\in K$ is the scalar from \autoref{CA:act}.  If $t\in
  T(\lambda)$ then $q^{\deg t}S_\lambda(K)$ is isomorphic to the
  $A$-submodule of $A_\lambda(K)$ with basis
  $\set{a_{st}+A^{>\lambda}(K)| s\in T(\lambda)}$.

  If $L$ is a (graded) $K$-module set $S_\lambda(L)=S_\lambda(K)\otimes_K L$. For
  example, if $K=\K[x]$ and $L=q^d\K$, which is the $\K[x]$-module
  concentrated in degree $d$ on which $x$ acts as~$0$, then
  $S_\lambda(L)\cong q^dS_\lambda(\K)$.

  By \autoref{CA:act} and \autoref{CA:star}, there is a unique symmetric bilinear form
  $\<\ ,\ \>_\lambda\map{S_\lambda(L)\times S_\lambda(L)}L$ such that
  \begin{equation}\label{E:CellularForm}
    \<a_s,a_t\>_\lambda a_u=a_{us}a_t 
    \qquad\text{for $s,t,u\in T(\lambda)$.}
  \end{equation}
  Moreover, $\<\ ,\ \>_\lambda$ is homogeneous and
  $\<ax,y\>_\lambda=\<x,a^*y\>_\lambda$, for all $a\in A$ and $x,y\in
  S_\lambda(L)$. In particular, if $L$ is concentrated in degree zero
  then $\<\ ,\ \>_\lambda$ is homogeneous of degree zero. Furthermore,
  \[
    \rad S_\lambda(L)=\set[\big]{x\in S_\lambda(L)|\<x,y\>=0\text{ for all }y\in S_\lambda(L)}
  \]
  is a graded $A$-module of $S_\lambda(L)$, so that
  $D_\lambda(L)=S_\lambda(L)/\rad S_\lambda(L)$ is a graded $A$-module.

  Suppose that $K=\bigoplus_dK_d$ is a graded commutative ring such that
  $K_0$ is a field. Then $K_d$ is a finite dimensional $K_0$-vector
  space.  Let $\Irr(K)$ be a complete set of irreducible graded
  $K$-modules, up to isomorphism.  Recall from \autoref{S:Rings} that
  $q$ is the grading shift functor.

  \begin{Lemma}\label{L:Kirreducibles}
    Suppose that $K=\Kx$. Then $\Irr(\Kx)=\set{q^d\K|d\in\Z}$.
  \end{Lemma}

  \begin{proof}
    Any irreducible graded $\Kx$-module is a $\K$-vector space on which
    each $x\in\ux$ acts as multiplication by~$0$. (Compare
    \autoref{R:GradedSimples}.)
  \end{proof}

  \begin{Example}\label{Ex:kxSimples}
    Suppose that $\K$ is a field and $x$ is an indeterminate over~$\K$.
    Then $\K[x]$ is a graded field and $q^d\K[x^\pm]=\K[x^\pm]$, for
    $d\in\Z$, since $x$ has degree~$1$. Hence, $\K[x^\pm]$ is the unique
    irreducible graded $\K[x^\pm]$-module. In contrast, if $\deg y=2$
    the $\Irr(\K[y^\pm]=\set{\K[y^\pm], q\K[y^\pm]}$.  (This is why we
    define each indeterminate $x\in\ux$ to have degree~$1$.)

    Now consider $\K[x^\pm,y^\pm]$, where $y$ be a second indeterminate
    over~$\K$. Then $L=\K[z^\pm]$ becomes an irreducible graded
    $\K[x^\pm,y^\pm]$-module by letting~$x$ act as multiplication by $c_1z$
    and $y$ act as multiplication by $c_2z$, for nonzero scalars
    $c_1,c_2\in\K^\times$. Equivalently,the
    module~$L\cong\K[\x^\pm,y^\pm]/(c_2x-c_1y)$ is uniquely determined
    by the fact that $x-\tfrac{c_1}{c_2}y$ acts on~$L$ as multiplication
    by~$0$. Hence, this makes $L$ into an irreducible graded $\K[x^\pm,y^\pm]$-module
    for each $c\in\K^\times$.
  \end{Example}

  Assume that $K_0$ is a field. If $L\in\Irr(K)$ set
  $P_0(L)=\set{\lambda\in P|D_\lambda(L)\ne0}$.

  \begin{Theorem} \label{T:CellularSimples}
    Let $K$ be a graded commutative domain such that $K_0$ is a field.
    Suppose that $A$ be a graded $K$-cellular algebra. Then
    \[
       \set[\big]{D_\lambda(L)|\lambda\in P_0(L)\text{ and }L\in\Irr(K)}
    \]
    is a complete set of pairwise non-isomorphic irreducible graded
    $A$-modules.  Moreover, $D_\lambda(L)$ is self-dual as an $A$-module
    if and only if $L\in\Irr(K)$ is self-dual as a $K$-module.
  \end{Theorem}

  \begin{proof}
    By \autoref{L:Kirreducibles}, up to shift the irreducible graded
    $A$-modules are irreducible $A(L)$-modules. The result now follows
    by repeating the standard arguments for classifying the simple
    modules of cellular algebras; see
    \cite[Theorem~3.12]{KoenigXi:AffineCellular},
    \cite[Theorem~3.4]{GL}, or \cite[Theorem~2.16]{Mathas:ULect}.
  \end{proof}

  \begin{Example}\label{Ex:GradedSimples}
     Suppose that $A$ is a graded $\Kx$-cellular algebra, where $\K$ is
     a field. Define $P_0$ as above. By \autoref{L:Kirreducibles},
     $\Irr(\K[x])=\set{q^d\K|q\in\Z}$. So
     \[
       \set[\big]{D_\lambda(L)|\lambda\in P_0\text{ and }L\in\Irr(\Kx)}
          = \set[\big]{q^dD_\lambda(\K)|\lambda\in P_0\text{ and }d\in\Z}
     \]
     is a complete set of pairwise non-isomorphic irreducible graded
     $A$-modules.  Let $A(\K[x^\pm])=\K[x^\pm]\otimes_{\K[x]}A$ be the
     corresponding graded $\K[x^\pm]$-cellular algebra over $\K[x^\pm]$.  Then
     $\set{D_\lambda(\K[x^\pm])|\lambda\in P_0}$ is a complete set of pairwise
     non-isomorphic irreducible graded $A(\K[x^\pm])$-modules.
  \end{Example}

  \begin{Definition}
    Suppose that $K=\Kx$ and let $A$ be a graded $\Kx$-cellular algebra.
    Let $\lambda\in P$ and $\mu\in P_0=P_0(\K)$ and set
    $S_\lambda=S_\lambda(\K)$ and $D_\mu=D_\mu(\K)$. Then
    $\mathsf{D}=\bigl([S_\lambda:D_\mu]_q\bigr)$ is the
    \emph{graded decomposition matrix} of~$A$, where
    \[
       [S_\lambda:D_\mu]_q = \sum_{k\in\Z}\,[S_\lambda: q^kD_\mu]\,q^k
       \quad\in\N[q,q^{-1}],
    \]
    and $[S_\lambda: q^kD_\mu]$ is the multiplicity of
    $q^kD_\mu$ as a composition factor of~$S_\lambda$.
  \end{Definition}

  Standard arguments from the theory of cellular algebras now prove the
  following:

  \begin{Corollary}\label{C:TriangularDet}
    Suppose that $A$ is a graded $\Kx$-cellular algebra. Then
    \begin{enumerate}
      \item If $\lambda\in P$ and $\mu\in P_0$ then
      $[S_\mu:D_\mu]_q=1$ and $[S_\lambda:D_\mu]_q\ne0$ only
      if $\lambda\ge\mu$.
      \item The Cartan matrix of $A$ is $\mathsf{D}^T\mathsf{D}$.
    \end{enumerate}
  \end{Corollary}

  \subsection{Cellular bases for \texorpdfstring{$\Rx(\Kxx)$}{R^Λ(K[x,x^-1])}}\label{SS:RxCellular}
  This section applies the results of the last two sections to show that
  $\Rx(\Kxx)$ is a $\Kxx$-cellular algebra. We have to wait until
  \autoref{SS:kcellular} to prove that $\Rx(\kx)$ is a $\kx$-cellular
  algebra.

  We have most of the data we need to define graded
  $\Kxx$-cell data for $\Rx(\Kxx)$: we have posets $(\Parts,\Ldom)$ and
  $(\Parts,\Gdom)$ and sets of standard tableaux
  $\Std(\Parts) = \coprod_{\blam\in\Parts}\Std(\blam)$. Moreover, by the results of
  \autoref{S:Bases}, we have bases $\set{\Lfst}$, $\set{\Lpsist}$, $\set{\Gfst}$ and
  $\set{\Gpsist}$, which we view as being given by injective maps
  \[
    \Lfst[]\to\Rx(\Kxx),\quad
    \Lpsist[]\to\Rx(\Kxx),\quad
    \Gfst[]\to\Rx(\Kxx)
       \quad\text{and}\quad
    \Lpsist[]\to\Rx(\Kxx),
  \]
  which send $(\s,\t)$ to $\Lfst$, $\Lpsist$, $\Gfst$ and $\Gpsist$,
  respectively. We still need to define corresponding degree functions
  on $\Std(\Parts)$.

  For $\t\in\Std(\Parts)$, recall the homogeneous scalars
  $\Lgt,\Ggt\in\Kxx$ from \autoref{C:gamma}. As $\Kxx$ is a graded ring,
  both of these scalars have a degree in~$\Z$. Recall that
  $\deg\map\Kxx\Z$ is the degree function on~$\Kxx$ and that
  $\deg(x)=1$, for all $x\in\ux$. By \autoref{L:gammaClosed}, the scalars $\Lgt$ and
  $\Ggt$ depend on the content function~$\bc$ and are
  polynomials in~$\kx2$ and, in particular, have even degree.

  \begin{Definition}\label{D:tableauDegree}
    Let $\t\in\Std(\Parts)$. Define \textbf{degree functions},
    \[
          \Ldeg\map{\Std(\Parts)}\Z
          \qquad\text{and}\qquad
          \Gdeg\map{\Std(\Parts)}\Z,
    \]
    with respect to the posets $(\Parts,\Ldom)$ and $(\Parts,\Gdom)$,
    respectively, by
    \[
            \Ldeg(\t) = \frac12\deg\Lgt\And
            \Gdeg(\t) = \frac12\deg\Ggt,
            \qquad\text{for }\t\in\Std(\Parts).
    \]
  \end{Definition}
  \notation{$\Ldeg,\Gdeg$}{Degree functions for the $\Lpsist[]$ and
            $\Gpsist[]$ bases}[D:tableauDegree]

  When $(\bc, \br)$ is a graded content system both of these degree
  functions already exist in the literature. In type~$\Aone$, Brundan,
  Kleshchev and Wang~\cite{BKW:GradedSpecht} call $\Gdeg$ the degree of
  a tableau and $\Ldeg$ its codegree. In type $\Cone$ Ariki, Park and
  Speyer~\cite{ArikiParkSpeyer:C} use $\Gdeg$ to define the degrees of
  the basis elements of their candidates for homogeneous Specht modules. Using
  \autoref{D:tableauDegree} it is not clear that these degree functions
  coincide with those given in
  \cite{BKW:GradedSpecht,ArikiParkSpeyer:C}, however, this is immediate
  from the next result.

  Recall from \autoref{S:Quiver} that $D=\diag(\di|i\in I)$ is the
  symmetriser of the Cartan matrix of~$\Gamma$.

  \begin{Lemma}\label{L:recursiveDeg}
    Suppose that $\t\in\Std(\Parts)$. Then
    \[
        \Ldeg(\t) = \sum_{m=1}^n
             \di[\br_m(\t)]\bigl(\#\LAdd_m(\t)-\#\LRem_m(\t)\bigr)
          \And*
        \Gdeg(\t) = \sum_{m=1}^n
             \di[\br_m(\t)]\bigl(\#\GAdd_m(\t)-\#\GRem_m(\t)\bigr).
    \]
  \end{Lemma}

  \begin{proof}
    Apply \autoref{L:gammaClosed}, using the fact that $\Dgt\ne0$ and
    $\deg\bc_m(\t)=2\di[\br_m(\t)]$, which follows from \autoref{D:ContentSystem}(c)
    because $(\bc,\br)$ is a graded content system.
  \end{proof}

  We can now show that $\Rx(\Kxx)$ is a (graded) $\Kxx$-cellular algebra.

  \begin{Theorem}\label{T:Kcellular}
    Suppose that $(\bc,\br)$ is a graded content system for $\Rx(\kx)$.
    Then the algebra $\Rx(\Kxx)$ is a $\Kxx$-cellular algebra with cellular bases:
    \begin{enumerate}
      \item $\set{\Lfst|(\s,\t)\in\Std^2(\Parts)}$ with weight poset
      $(\Parts,\Ledom)$ and degree function $\Ldeg$.
      \item $\set{\Gfst|(\s,\t)\in\Std^2(\Parts)}$ with weight poset
      $(\Parts,\Gedom)$ and degree function $\Gdeg$.
      \item $\set{\Lpsist| (\s,\t)\in\Std^2(\Parts)}$ with weight poset
      $(\Parts,\Ledom)$ and degree function $\Ldeg$.
      \item $\set{\Gpsist| (\s,\t)\in\Std^2(\Parts)}$ with weight poset
      $(\Parts,\Gedom)$ and degree function $\Gdeg$.
    \end{enumerate}
  \end{Theorem}

  \begin{proof}
    Let $\Domin$. By \autoref{C:fBases}, $\set{\Dfst}$ is a $\Kxx$-basis
    of $\Rx(\Kxx)$ and by \autoref{P:fstaction} the $\Dfst[]$-basis
    satisfies \autoref{CA:act}. Recall that~$*$ is unique
    anti-isomorphism of~$\Rx(\Kxx)$ that fixes each of its generators.
    By construction, $(\Dpsist)^*=\Dpsist[\t\s]$ and $F_\s^*=F_\s$, so
    $(\Dfst)^*=\Dfst[\t\s]$ for $(\s,\t)\in\Std^2(\Parts)$. Hence,
    $\set{\Dfst}$ is a $\Kxx$-cellular basis of~$\Rx(\Kxx)$.

    Next, consider $\set{\Dpsist}$. This is a basis of $\Rx(\Kxx)$ by
    \autoref{P:psiTriangular}, so \autoref{CA:basis} is satisfied. We
    have already seen that $(\Dpsist)^*=\Dpsist[\t\s]$, verifying
    \autoref{CA:star}, so it remains to check \autoref{CA:act}. By
    \autoref{P:psiTriangular},
    \[
           \Dpsist[\s\Dtlam] \equiv \Dfst[\s\Dtlam]
                +\sum_{\u\Dom\s}r_\u\Dfst[\u\Dtlam]
                   \pmod\DRlam
    \]
    for some $r_\u\in\Kxx$ and where $\DRlam$ is the two-sided ideal
    of~$\Rx(\Kxx)$ spanned by $\set{\Dfst[\u\v]}$ where
    $\Shape(\u)=\Shape(\v)\Dom\blam$. By \autoref{P:psiTriangular},
    $\DRlam$ is also spanned by $\set{\Dpsist[\u\v]}$. Multiplying the
    last displayed equation on the left by~$a\in\Rx(\Kxx)$, and using
    \autoref{P:fstaction} and \autoref{P:psiTriangular},
    \[
          a\Dpsist[\s\Dtlam] \equiv
          a\Bigl(\Dfst[\s\Dtlam]
                +\sum_{\u\Dom\s}a_\u\Dfst[\u\Dtlam]\Bigr)
            =\sum_{\x\in\Std(\blam)}b_\x\Dfst[\x\Dtlam]
            \equiv \sum_{\x\in\Std(\blam)}c_\x\Dpsist[\x\Dtlam]
                   \mod\DRlam,
    \]
    for some homogeneous scalars $a_\u,b_\x,t_\x\in\Kxx$.
    Multiplying on the right by $\psi_{\Ddt}^*$  shows that $\Dpsist$
    satisfies \autoref{CA:act}. Hence, $\set{\Dpsist}$ is a $\Kxx$-cellular
    basis of~$\Rx(\Kxx)$.

    It remains to show that each of these bases is a graded
    $\Kxx$-cellular basis of $\Rx(\Kxx)$ when $(\bc,\br)$ is a graded
    content system. By \autoref{D:psist}, $\Dpsist$ is homogeneous, for
    $(\s,\t)\in\Std^2(\Parts)$. By \autoref{D:Ft}, $F_\t$ is homogeneous
    of degree~$0$ and $\Dfst=F_\s\Dpsist F_\t$. Hence, $\Dfst$
    is homogeneous and $\deg\Dfst=\deg\Dpsist$. Therefore, it is enough
    to show that $\deg\Dfst=\Ddeg(\s)+\Ddeg(\t)$.  Further, since~$*$ is
    homogeneous, $\deg\Dfst=\deg\Dfst[\t\s]$. So, using
    \autoref{L:fproduct},
    \[
        \deg\Dfst=\tfrac12\deg\bigl(\Dfst\Dfst[\t\s]\bigr)
                 =\tfrac12\deg\bigl(\Dgt\Dfst[\s\s]\bigr)
                 =\tfrac12\deg\bigl(\Dgt[\s]\Dgt[\t]F_\s\bigr)
                 =\deg(\s)+\deg(\t),
    \]
    as we wanted to show. This completes the proof.
  \end{proof}

  Proving that $\Rx(\Kxx)$ is a $\Kxx$-cellular algebra is nice but it
  does not directly help us in constructing a cellular basis for the KLR
  algebras~$\Rn(\k)$ and $\Rx(\kx)$. We prove that $\Rx(\kx)$ is
  $\kx$-cellular in the next section. As a prelude to doing this, for
  $\blam\in\Parts$ define $\LSlam(\Kxx)$ and $\GSlam(\Kxx)$ to be the
  graded cell modules for $\Rx(\Kxx)$-determined by the seminormal bases
  $\set{\Lfst}$ and $\set{\Gfst}$, respectively. Let $\Domin$.  By
  \autoref{P:fstaction}, $\DSlam(\Kxx)$ has basis $\set{\Dfs}$ and there is an
  isomorphism
  \[ q^{\Ddeg\Dtlam}\DSlam(\Kxx)
     \cong\Bigl(\Rx(\Kxx)\Dfst[\Dtlam\Dtlam]+\DRlam\Bigr){\bigm/}\DRlam;
         \Dfs\mapsto \Dfst[\s\Dtlam]+\DRlam.
  \]
  For $\s\in\Std(\blam)$, let $\Dpsis=\psi_{\Ddt[\s]}\Dfs[\Dtlam]$ be
  the element of $\DSlam(\Kxx)$ that is sent to
  $\Dpsist[\s\Dtlam]+\DRlam$ under this isomorphism. In view of
  \autoref{C:IsomorphicSpechts} and \autoref{P:psiTriangular}, we have:
  \notation{$\LSlam,\GSlam$}{Graded Specht modules for the $\Lpsist[]$ and $\Gpsist[]$ bases}

  \begin{Lemma}\label{L:KSpecht}
    Let $\blam\in\Parts$. As $\Kxx$-modules,
    \[
      \LSlam(\Kxx)=\bigoplus_{\s\in\Std(\blam)}\Kxx\Lpsis
      \quad\text{and}\quad
      \GSlam(\Kxx)=\bigoplus_{\s\in\Std(\blam)}\Kxx\Gpsis.
    \]
  \end{Lemma}

  By \autoref{L:FufstFv}, if
  $\theta\map{\LSlam(\Kxx)}\GSlam(\Kxx)$ is an isomorphism then
  $\theta(\Lfs)=a\Gfs$, for some $a\in\Kxx$. Comparing degrees, $a$ is
  homogeneous of degree $\Ldeg(\s)-\Gdeg(\s)$. In particular, such an
  isomorphism and its inverse are defined over $\kx$ if and only if
  $\Ldeg(\s)=\Gdeg(\s)$ for all $\s\in\Std(\blam)$.

  Let $\LSlam(\kx)=\bigoplus_\s\kx\Lpsis$ and
  $\GSlam(\kx)=\bigoplus_\s\kx\Gpsis$, where in both sums
  $\s\in\Std(\blam)$. By definition, $\LSlam(\kx)$ and~$\GSlam(\kx)$ are
  free $\kx$-modules and, by base-change, $\DSlam(\Kxx)\cong\Kxx\otimes_{\kx}\DSlam(\kx)$ by
  \autoref{L:KSpecht}. In fact, $\LSlam(\kx)$ and $\GSlam(\kx)$ are
  both $\Rx(\kx)$-modules.

  \begin{Proposition}\label{P:integral}
    Suppose that $\s\in\Std(\blam)$, for $\blam\in\Parts$. Then:
    \begin{enumerate}
      \item If $1\le k<n$ then
      $\psi_k\Lpsis\in\LSlam(\kx)$ and
      $\psi_k\Gpsis\in\GSlam(\kx)$.
      \item If $1\le m\le n$ then
      $y_m\Lpsis\in\LSlam(\kx)$ and
      $y_m\Gpsis\in\GSlam(\kx)$.
      \item If $\sigma_{b_1}\dots \sigma_{b_l}$ is a reduced expression
      for $\Ddt[\s]$ then
      \[
        \Lpsis-\psi_{b_1}\dots\psi_{b_l}\Lpsis[\Ltlam]
        \in\bigoplus_{\u\Ldom\s}\kx\Lpsis[\u]
        \And
        \Gpsis-\psi_{b_1}\dots\psi_{b_l}\Gpsis[\Gtlam]
        \in\bigoplus_{\u\Gdom\s}\kx\Gpsis[\u].
      \]
    \end{enumerate}
  \end{Proposition}

  \begin{proof}
    Let $\Domin$. To prove the proposition we argue by induction
    on the length $L(\Ddt[\s])$ of the permutation $\Ddt[\s]$.
    To start the induction, suppose that $\s=\Dtlam$, so that
    $\Ddt[\Dtlam]=1$. Then (c) is vacuously true and,
    $y_m\Dpsis[\Dtlam]=y_m\Dfs[\Dtlam]=\bc_m(\s)\Dfs[\Dtlam]$ by
    \autoref{P:fstaction}, so (b) holds. To prove~(a), applying
    \autoref{P:fstaction} shows that
    \[
       \psi_k\Dpsis[\Dtlam]=\psi_k\Dfs[\Dtlam]
       = \begin{cases*}
            \Dfs[u]=\Dpsis[\u] & if $\u=\sigma_k\Dtlam\in\Std(\blam)$,\\
            0& if $\sigma_k\Dtlam\notin\Std(\blam)$.
       \end{cases*}
    \]
    Hence, the proposition is true when $\s=\Dtlam$.

    Now suppose that $\Dtlam\Dom\s$. First, consider~(c). Let
    $\Ddt[\s]=\sigma_{a_1}\dots \sigma_{a_l}$ be the preferred reduced expression
    for $\Ddt[\s]$ that we fixed after \autoref{L:wlam}. If
    $\sigma_{b_1}\dots \sigma_{b_l}$ is a second reduced expression for~$\Ddt[\s]$
    then, by Matsumoto's theorem (see, for example,
    \cite[Theorem~1.8]{Mathas:ULect}), we can convert the reduced
    expression $\sigma_{a_1}\dots \sigma_{a_l}$ into our preferred reduced
    expression $\sigma_{b_1}\dots \sigma_{b_l}$ using only the braid relations
    of~$\Sym_n$. The~$\psi_k$ satisfy the commuting braid relations and
    by \autoref{R:braid} they satisfy the braid relations of
    length~$3$ plus or minus an ``error term'' of the form
    $\delta_{i_ki_{k+2}}Q^\ux_{i_ki_{k+1}i_{k+1}}(y_k,y_{k+1},y_{k+1})\psi_u$,
    where $u$ is smaller than $\Ddt[\s]$ in the Bruhat order and so, in
    particular, $L(u)<L(\Ddt[\s])$. Hence, by induction, part~(c) holds
    for~$\Dpsis$.

    Now consider $\psi_k\Dpsis$ as in~(a). If
    $L(\sigma_k\Ddt[\s])<L(\Ddt[\s])$ then $\Ddt[\s]$ has a reduced
    expression of the form $\sigma_k\sigma_{a_2}\dots \sigma_{a_l}$. Therefore,
    using~(c), which we have already proved,
    \begin{align*}
      \psi_k\Dpsis &=
         \psi_k\Bigl(\psi_k\psi_{a_2}\dots\psi_{a_l}\Dpsis[\Dtlam]
                   + \sum_{\u\Dom\s}r_\u\Dpsis[\u]\Bigr),
                \qquad\text{for some $r_\u\in\kx$,}\\
                   &=\psi_k^2\psi_{a_2}\dots\psi_{a_l}\Dpsis[\Dtlam]
                   + \sum_{\u\Dom\s}r_\u\psi_k\Dpsis[\u]\\
            &=Q^\ux_{\br_{k+1}(\s),\br_k(\s)}(y_k,y_{k+1})\psi_{a_2}\dots\psi_{a_l}
               \Dpsis[\Dtlam] + \sum_{\u\Dom\s}r_\u\psi_k\Dpsis[\u].
    \end{align*}
    By induction, all of these terms belong to $\DSlam(\kx)$, showing
    that $\Dpsis$ satisfies~(a).

    Finally, consider $y_m\Dpsis$. Let $\v\in\Std(\blam)$ be the unique
    standard tableau such that $\s=\sigma_{r_1}\v$. Then
    $L(\Ddt[\s])=L(\Ddt[\v])+1$, so by part~(c) and induction,
    \[
            \Dpsis[\v] = \psi_{a_2}\dots\psi_{a_l}\Dpsis[\Dtlam]
                   + \sum_{\u\Dom\v}r_\u\Dpsis[\u],
    \]
    for some $r_\u\in\kx$ (these scalars are different from those in the
    last paragraph). Therefore,
    \[
        y_m\Dpsis
         =y_m\psi_{r_1}\Bigl(\Dpsis[\v]-\sum_{\u\Dom\v}r_\u\Dpsis[\u]\Bigr).
    \]
    Applying \autoref{R:psiy} and induction now completes the proof.
  \end{proof}

%

  \subsection{Defect polynomials} The algebra $\Rx(\Kxx)$ is a split
  semisimple graded algebra, so it is naturally a symmetric algebra with
  symmetrising form given by taking the matrix trace on the
  regular representation. This form does not restrict to give a trace on
  $\Rx(\kx)$, so the aim of this section is to show how to use this trace
  form to give an ``integral'' trace form'' on~$\Rx(\kx)$. In later
  sections, these results will be used to understand the duals of some
  $\Rx(\kx)$-modules.

  We continue to assume that $(\bc,\br)$ is a graded content system for
  $\Rx(\kx)$ with values in~$\kx$. The following innocuous lemma is the
  key to constructing our trace form and to understanding the defect of
  the blocks of $\Rx(\kx)$.

  \begin{Lemma}\label{L:DefectPolynomial}
    Suppose that $\blam\in\Parts$. Then
    $\Lgt[\s]\Ggt[\s]=\Lgt\Ggt$  for all $\s,\t\in\Std(\blam)$.
  \end{Lemma}

  \begin{proof}
    It is enough to consider the case when $\s\Gdom\t=\sigma_k\s$, for
    some $1\le k<n$. In this case we have that $\Ggt=Q_k(\s)\Ggt[\s]$
    and $\Lgt[\s]=Q_k(\t)\Ggt[\t]$ by \autoref{L:gammaRecurrence}.  By
    \autoref{E:Qmt} and the symmetry of Rouquier's $Q$-polynomials,
    $Q_k(\s)=Q_k(\t)$.  Hence,
    \[
        \Lgt[\s]\Ggt[\s] = \frac{Q_k(\t)}{Q_k(\s)}\cdot\Lgt[\t]\Ggt[\t]
                         = \Lgt\Ggt,
    \]
    as required.
  \end{proof}

  \begin{Definition}\label{D:DefectPolynomial}
    Let $\blam\in\Parts$. The \textbf{$\blam$-defect polynomial}
    is $\DefPol=\Lgt\Ggt$, for any $\t\in\Std(\blam)$.
  \end{Definition}
  \notation{$\DefPol$}{The defect polynomial of $\blam\in\Parts$}[D:DefectPolynomial]

  By \autoref{L:DefectPolynomial}, the defect polynomial $\DefPol$
  depends only on $\blam$, and not on the choice of
  $\t\in\Std(\blam)$.  We will show in \autoref{C:defectPoly} that the
  degree of the defect polynomial is a block invariant. That is,
  if~$\blam,\bmu\in\Parts[\alpha]$ then $\deg\DefPol=\deg\DefPol[\bmu]$. To
  prove this, and to explain why we call this the \textit{defect}
  polynomial we need some more notation.

  For $i\in I$ and $\blam\in\Parts$ let $\Add_i(\blam)$ and
  $\Rem_i(\blam)$ be the sets of addable and removable $i$-nodes of
  $\blam$, respectively.  Recall from \autoref{S:Quiver} that
  $\set{\di|i\in I}$ is the set of symmetrisers of~$\Gamma$.

  \notation{$\alpha_\blam$}{The positive root $\sum_{A\in\blam}\alpha_{\br(A)}\in Q^+$}[D:defect]
  \notation{$\Parts[\alpha]$}{The set of $\ell$-partitions
          $\set{\blam\in\Parts|\alpha_\blam=\alpha}$}[D:defect]
  \notation{$\defect(\blam)$} {The $\Lambda$-defect of $\blam$, which is
    $\defect(\alpha_\blam)=(\Lambda,\alpha_\blam)-\tfrac12(\alpha_\blam,\alpha_\blam)$}[D:defect]
  \notation{$\Ldilam,\Gdilam$}{Number of addable minus removal $i$-nodes below/above $A$}
          [D:defect]
  \notation{$d_i(\blam)$}{Number of addable minus removable $i$-nodes of $\blam$}[D:defect]

  \begin{Definition}\label{D:defect}
    Let $\alpha\in Q^+$.
    \begin{enumerate}
      \item For $\alpha\in Q^+_n$ let
      $\Parts[\alpha]=\set{\blam\in\Parts|\alpha_\blam=\alpha}$.
      \item The \emph{$\Lambda$-defect} of $\alpha\in Q^+_n$ is
      $\defect(\alpha) = (\Lambda|\alpha) - \frac12(\alpha|\alpha)$.
      \item The \emph{$\blam$-positive root} is
      $\alpha_\blam =\sum_{A\in\blam}\alpha_{\br(A)}\in Q^+_n$.
       \item The \emph{$\Lambda$-defect} of $\blam\in\Parts[\alpha]$ is
      $\defect(\blam)=\defect(\alpha_\blam)$.
      \item Motivated by \autoref{E:Sddable}, given an addable or removable $i$-node~$A$ of $\blam$
      define
      \begin{align*}
        \Ldilam &= \di\times\bigl(\#\set{B\in\Add_i(\blam)|B<A}
                        -\#\set{B\in\Rem_i(\blam)|B<A}\bigr),\\
        \Gdilam &= \di\times\bigl(\#\set{B\in\Add_i(\blam)|B>A}
                        -\#\set{B\in\Rem_i(\blam)|B>A}\bigr),\\
      d_i(\blam)&=\di\times\bigl(\#\Add_i(\blam)-\#\Rem_i(\blam)\bigr).
      \end{align*}
    \end{enumerate}
  \end{Definition}

  By definition, $\defect(\alpha)\in\Z$. We show in
  \autoref{C:PositiveDefect} that, in fact, $\defect(\alpha)\in\N$.
  Generalising \cite[Lemma~3.11]{BKW:GradedSpecht}, we give some
  standard facts about defect.

  \begin{Lemma}\label{L:defect}
    Suppose that $\blam=\bmu+A$, where $A\in\Add_i(\bmu)$ for $i\in I$.
    Then $\alpha_\blam=\alpha_\bmu+\alpha_i$ and
    \addtocounter{equation}{-1}
    \begin{subequations}
      \begin{align}
        d_i(\blam) &= d_i(\bmu)-2\di = \Ldilam+\Gdilam+\di
        \label{L:defectSum}\\
        d_i(\blam) &=(\Lambda-\alpha_\blam|\alpha_i)
        \label{L:defectWeight}\\
        \defect(\blam) &=\defect(\bmu)+d_i(\blam)+\di
                        =\defect(\bmu)+\Ldilam+\Gdilam.
        \label{L:defectRecursive}
      \end{align}
    \end{subequations}
  \end{Lemma}

  \begin{proof}
    Equation~\autoref{L:defectSum} is just a rephrasing of \autoref{D:defect}(e).

    To prove \autoref{L:defectWeight} we argue by induction on~$n$.
    If $n=0$ then $\blam=\zero$, $\alpha_\blam=0$ and $(\Lambda|\alpha_i)=d_i(\blam)$
    is the number of addable $i$-nodes of~$\zero$. If $n>0$ then
    \[
       (\Lambda-\alpha_\blam|\alpha_i)
          = (\Lambda-\alpha_\bmu|\alpha_i) - (\alpha_i|\alpha_i)
          = d_i(\bmu) - 2\di
          = d_i(\blam),
    \]
    where the second equality follows by induction and the third
    equality from~(a). This proves \autoref{L:defectWeight}.

    Now consider~\autoref{L:defectRecursive}. As $\blam$ has a removable
    $i$-node, $\alpha_\bmu=\alpha_\blam-\alpha_i\in Q^+_{n-1}$ and
    \begin{align*}
         \defect(\blam)&= \defect(\alpha_\bmu+\alpha_i)
             =(\Lambda|\alpha_\bmu)+(\Lambda|\alpha_i)
                 -\tfrac12\bigl((\alpha_\bmu|\alpha_\bmu)
                   +2(\alpha_\bmu|\alpha_i)+(\alpha_i|\alpha_i)\bigr)\\
            &=\defect(\bmu) + (\Lambda-\alpha_\bmu|\alpha_i)-\di,
                        \hspace*{0.28\textwidth}\text{by induction,}\\
            &=\defect(\bmu) + d_i(\bmu)-\di,
            \hspace*{0.345\textwidth}\text{by~\autoref{L:defectWeight}},\\
            &=\defect(\bmu) + d_i(\blam)+\di,
    \end{align*}
    where the last equality follows by~\autoref{L:defectSum}.
    The second equality in \autoref{L:defectRecursive} follows
    by a second application of~\autoref{L:defectSum}.
  \end{proof}

  \begin{Corollary}\label{C:degCodeg}
    Suppose that $\t\in\Std(\blam)$, for $\blam\in\Parts$. Then
    $\Ldeg(\t)+\Gdeg(\t)=\defect(\blam)$.
  \end{Corollary}

  \begin{proof}
     This follows by induction on~$n$. If $n=0$ then
     $\Ldeg(\t)=\Gdeg(\t)=\defect(\blam)=0$, so the result holds.
     Suppose that $n>0$ and let $A=\t^{-1}(n)$. Set $\s=\t_{\downarrow(n-1)}$,
     $\bmu=\Shape(\s)$
     and $i=\br_n(\t)=\br(A)$. Then,
     \begin{align*}
        \Ldeg(\t)+\Gdeg(\t)
           &= \Ldeg(\s)+\Ldilam+\Gdeg(\s)+\Gdilam
               &&\text{by \autoref{L:recursiveDeg}},\\
           &=\defect(\bmu)+\Ldilam+\Gdilam
              &&\text{by induction},\\
           &=\defect(\blam),
     \end{align*}
     with the last equality coming from \autoref{L:defectRecursive}.
  \end{proof}

  We can now explain the origin of the name \textit{defect} polynomial.
  In view of \autoref{C:PositiveDefect} below, this shows that
  $\DefPol\in\kx$, for $\blam\in\Parts$. It would be interesting to
  determine these polynomials explicitly; compare
  \cite{ChlouverakiJacon:AKSchurElts}.

  \begin{Corollary}\label{C:defectPoly}
    Let $\blam\in\Parts$. Then $\DefPol$ is a homogeneous polynomial of
    degree $2\defect(\blam)$.
  \end{Corollary}

  \begin{proof}
    If $\t\in\Std(\blam)$ then, by \autoref{D:tableauDegree} and \autoref{C:degCodeg},
    the defect polynomial $\DefPol$ is homogeneous of degree
    $\deg\Lgt+\deg\Ggt=2\bigl(\Ldeg(\t)+\Gdeg(\t)\bigr)=2\defect(\blam)$.
  \end{proof}

  Although we do not need this, we note that the defect polynomial, or
  more correctly \autoref{L:DefectPolynomial}, allows us to describe the
  transition matrix between the $\Lfst[]$-basis and the $\Gfst[]$-basis,
  generalising \autoref{C:gamma}.

  \begin{Proposition}\label{P:fTransition}
    Let $\s,\t\in\Std(\blam)$, for $\blam\in\Parts$. Then
    $\Lfst=\frac{\DefPol}{\Ggt[\s]\Ggt}\Gfst$
    in $\Rx(\Kxx)$.
  \end{Proposition}

  \begin{proof}
    By \autoref{L:DefectPolynomial},
    $\Lgt[\s]/\Ggt=\Lgt/\Ggt[\s]$, so
    the statement of the proposition is equivalent to the equivalent claims that
    $\frac{\Lgt[\s]}{\Ggt}\Gfst=\Lfst=\frac{\Lgt}{\Ggt[\s]}\Gfst$.
    Since $\Dfst=(\Dfst[\t\s])^*$, it is enough to show that
    $\Lfst[\Ltlam\t]=\tfrac{\Lgt[\Ltlam]}{\Ggt}\Gfst[\Gtlam\t]$ by
    \autoref{L:fproduct}. We show this by arguing by induction on
    $L(\Ldt)$, the length of the permutation $\Ldt$.  When $\t=\Ltlam$
    the result follows from \autoref{C:gamma}. If $\t\ne\Ltlam$ then we
    can write $\t=\sigma_k\v$ with $\v\Ldom\t$ and
    $L(\Ldt[\v])=L(\Ldt)-1$. Hence, by two applications of
    \autoref{P:fstaction}, and induction,
    \[
        \Lfst[\Ltlam\t]
           = \Lfst[\Ltlam\v]\Bigl(\psi_k
           -\frac{\delta_{\br_k(\v),\br_{k+1}(\v)}}{\bc_{k+1}(\v)-c_k(\v)}\Bigr)
            = \frac{\Lgt[\Ltlam]}{\Ggt[\v]}\Gfst[\Ltlam\v] \Bigl(\psi_k
               -\frac{\delta_{\br_k(\v),\br_{k+1}(\v)}}{\bc_{k+1}(\v)-c_k(\v)}\Bigr)
           = \frac{\Lgt[\Ltlam]}{\Ggt[\v]}\Gfst[\Ltlam\t]Q_k(\v).
    \]
    This completes the proof of the inductive step, and the proposition,
    since $\Ggt[\v]=Q_k(\v)\Ggt$ by \autoref{L:gammaRecurrence}.
  \end{proof}

  By the proposition,
  $\Lfst=\frac{\DefPol}{\Ggt[\s]\Ggt}\Gfst=\frac{\Lgt[\s]}{\Ggt}\Gfst=\frac{\Lgt}{\Ggt[\s]}\Gfst$.
  In particular, the four terms in this equation have the same degree,
  which is easily checked using \autoref{C:degCodeg}.

  \subsection{A symmetrising form}
  This section uses the defect polynomials to define a symmetrising form
  on the algebra $\Rx(\kx)=\bigoplus_{\alpha\in Q^+_n}\Rx[\alpha](\kx)$, and hence
  shows that it is a graded symmetric algebra. This symmetrising form
  specialises to give a non-degenerate symmetrising form on the
  cyclotomic KLR algebra $\Rn(\k)$.

  This section is partly inspired by \cite{Mathas:gendeg}, where similar
  results were obtained for the Ariki-Koike algebras. The arguments
  given here are much shorter than those in \cite{Mathas:gendeg}, which
  is surprising both because the results here are stronger and because
  we need to prove everything from the ground up.

  \begin{Definition}
    Let $\alpha\in Q^+_n$. For $\blam\in\Parts[\alpha]$ let
    $\chi_\blam$ be the character of the irreducible $\Rx[\alpha](\Kxx)$-module
    $V_\blam(\Kxx)$. The \textbf{$\alpha$-trace form} is the map
    $\tau_\alpha\map{\Rx[\alpha](\Kxx)}\Kxx$ given by
    \[
        \tau_\alpha = \sum_{\blam\in\Parts[\alpha]}\frac1{\DefPol}\chi_\blam.
    \]
  \end{Definition}

  By \autoref{C:defectPoly}, the trace form $\tau_\alpha$ is homogeneous
  of degree $-2\defect(\alpha)$ and takes values in~$\Kxx$.

  We use the characters of $V_\blam(\Kxx)$, for $\blam\in\Parts[\alpha]$,
  in this definition to emphasise that $\tau_\alpha$ does not depend on
  a choice of basis. Note that if $\blam\in\Parts$ then
  $\LSlam(\Kxx)\cong S_\blam(\Kxx)\cong\GSlam(\Kxx)$ by
  \autoref{C:IsomorphicSpechts}.

  \begin{Example}\label{Ex:tauFt}
    Let $\s,\t\in\Std(\blam)$, where $\blam\in\Parts[\alpha]$ and
    $\alpha\in Q^+_n$.  Using \autoref{L:fproduct},
    \[
        \tau_\alpha(F_\t)=\frac{1}{\DefPol},\qquad
        \tau_\alpha(\Lfst)=\frac{\delta_{\s\t}}{\Ggt}
        \And
        \tau_\alpha(\Gfst)=\frac{\delta_{\s\t}}{\Lgt}.
    \]
  \end{Example}

  To study $\Rx(\kx)$, we use $\tau_\alpha$ to define an ``integral''
  bilinear form.  If $f(\ux)\in\k[\ux^\pm]$ is a homogeneous polynomial let
  $f_0\in\k$ be the constant term of $f(\ux)$.

  \begin{Definition}
    Let $\alpha\in Q^+_n$.
    Let $\<\ ,\ \>_\alpha\map{\Rx[\alpha](\kx)\times\Rx[\alpha](\kx)}\k$ be the
    homogeneous bilinear form on~$\Rx[\alpha](\kx)$ of degree
    $-2\defect(\alpha)$ given by
    $\<a,b\>_\alpha = \tau_\alpha(ab)_0$.
  \end{Definition}

  \notation{$\<\ ,\ \>\alpha$}{Non-degenerate symmetric bilinear form on
  $\Rx[\alpha](\kx)$}

  We leave the proof of the following easy facts about $\tau_\alpha$ and
  $\<\ ,\ \>_\alpha$ to the reader.

  \begin{Lemma}\label{L:bilinear}
    Suppose that $a,b\in\Rx[\alpha](\kx)$, for $\alpha\in Q^+_n$. Then
    \begin{enumerate}
      \item Let $a,b,\in\Rx[\alpha](\kx)$. Then $\tau_\alpha(ab)=\tau_\alpha(ba)$,
         $\tau_\alpha(a)=\tau_\alpha(a^*)$, and
         $\<a,b\>_\alpha=\<b,a\>_\alpha$.
      \item If $a,b,c\in\Rx[\alpha](\kx)$ then $\<ab,c\>_\alpha=\<a,bc\>_\alpha$.
    \end{enumerate}
  \end{Lemma}

  We want to show that $\<\ ,\ \>_\alpha$ is a homogeneous
  non-degenerate bilinear form on $\Rx[\alpha](\kx)$. The next
  results pave the way to proving this.  The first result is similar in
  spirit to \cite[Lemma~4.11]{HuMathas:GradedCellular}.

  \begin{Lemma}\label{L:tlamDominance}
    Suppose that $\blam\in\Parts$. Then there exist
    $r_{\t},s_{\t}\in\Kxx$ such that
    \[
         \Lpsist[\Ltlam\Ltlam] = \Lfst[\Ltlam\Ltlam]
         +\sum_{\t\Ldom\Ltlam}r_{\t}\Lfst[\t\t]
          \And
         \Gpsist[\Gtlam\Gtlam] = \Gfst[\Gtlam\Gtlam]
         +\sum_{\t\Gdom\Gtlam}s_{\t}\Gfst[\t\t].
    \]
  \end{Lemma}

  \begin{proof}
    Let $\Domin$.  By definition and \autoref{C:KLRIdempotents},
    \[
        \Dpsist[\Dtlam\Dtlam] = \Dylam \ei[\Dilam]
             =\Dylam\sum_{\t\in\Std(\Dilam)}\frac1{\Dgt}F_\t
             =\sum_{\t\in\Std(\Dilam)}\frac1{\Dgt}
                    \prod_{A\in\DAdd(\Dtlam)}
                        \bigl(\bc_m(\t)-\bc(A)\bigr)F_\t
    \]
    Suppose that $\t\in\Std(\Dilam)$ and that $\t\Domneq\Dtlam$. Let
    $1\le k<n$ be minimal such that
    $\t_{\downarrow k}\Dom\Dtlam[\blam,\downarrow k]$  and
    $\t_{\downarrow(k+1)}\Domneq\Dtlam[\blam,\downarrow(k+1)]$.
    Let $A=\t^{-1}(k+1)$ and $B=(\Dtlam)^{-1}(k+1)$. Abusing
    notation slightly, $B\Dom A$, so $A\in\DAdd_k(\Dtlam)$.
    That is, $A\in\DAdd(\Dtlam)$ appears in the product above,
    contributing the factor $\bc_{k+1}(\t)-\bc(A)=0$. Hence,
    $\Dfst[\t\t]=\frac1{\Dgt}F_\t$ appears in $\Dpsist[\Dtlam\Dtlam]$ only
    if $\t\Dom\Dtlam$, where dominance holds because $\br(\t)=\Dilam$.
  \end{proof}

  The next result strengthens \autoref{P:psiTriangular}. Recall from
  \autoref{S:Bases} that $(\s,\t)\Ledom(\u,\v)$ if $\s\Ledom\u$ and
  $\t\Ledom\v$.

  \begin{Lemma}\label{L:fpsiExpansions}
     Let $\s,\t\in\Std(\blam)$, for $\blam\in\Parts$. Then
     \begin{align*}
       \Lpsist &= \Lfst
       + \sum_{\substack{(\u,\v)\in\Std^2(\Parts)\\(\u,\v)\Ldom(\s,\t)}}
           a_{\u\v}\Lfst[\u\v] &
       \Gpsist &= \Gfst
       + \sum_{\substack{(\u,\v)\in\Std^2(\Parts)\\(\u,\v)\Gdom(\s,\t)}}
           b_{\u\v}\Gfst[\u\v]\\
       \Lfst &= \Lpsist
       + \sum_{\substack{(\u,\v)\in\Std^2(\Parts)\\(\u,\v)\Ldom(\s,\t)}}
           c_{\u\v}\Lpsist[\u\v] &
       \Gfst &= \Gfst
       + \sum_{\substack{(\u,\v)\in\Std^2(\Parts)\\(\u,\v)\Gdom(\s,\t)}}
           d_{\u\v}\Gpsist[\u\v],
     \end{align*}
     for some scalars $a_{\u\v}, b_{\u\v},c_{\u\v},d_{\u\v}\in\Kxx$.
  \end{Lemma}

  \begin{proof}
     Let $\Domin$. We argue by induction on the dominance
     order~$\Dom$ on~$\Parts$. Let $\blam$ be maximal with respect to
     $\Dom$. Then $\blam=(0|\dots|0|1^n)$ if $\Dom=\Ldom$ and
     $\blam=(n|0|\dots|0)$ if $\Dom=\Gdom$. In this case,
     $\Dfst[\Dtlam\Dtlam]=\Dpsist[\Dtlam\Dtlam]$, so the result holds.

     Now suppose that $\blam$ is not maximal. By
     \autoref{L:tlamDominance} the proposition holds for
     $\Dfst[\Dtlam\Dtlam]$ so, by induction, the result holds for
     $\Dpsist[\Dtlam\Dtlam]$.  Now suppose that
     $(\Dtlam,\Dtlam)\Dom(\s,\t)$, for $\s,\t\in\Std(\blam)$.  We can
     assume that $\s\ne\Dtlam$ by applying~$*$, if necessary. Pick $k$
     such that $\y=\sigma_k\s\Dom\s$. By \autoref{P:integral}(c)
     and induction,
     \[
        \Dpsist = \psi_k\Dpsist[\y\t]
          = \psi_k\Bigl(\Dfst[\y\t] +
            \sum_{(\u,\v)\Dom(\y,\t)} r_{\u\v}\Dfst[\u\v]\Bigr)
          = \Dfst[\s\t] +
            \sum_{(\u,\v)\Dom(\y,\t)} r_{\u\v}\psi_k\Dfst[\u\v],
     \]
     for some $r_{\u\v}\in\Kxx$. Consider a term $\psi_k\Dfst[\u\v]$ on
     the right-hand side and let $\w=\sigma_k\u$. If
     $L(\Ddt[\w])=L(\Ddt[\u])+1$ then $\Ddt[\w]$ is a subexpression of
     $\Ddt[s]$ since $\u\Dom\y$ and $L(\Ddt[\s])=L(\Ddt[\y])+1$, so
     $\w\Dom\s$. If $L(\Ddt[\w])=L(\Ddt[\u])+1$ then
     $\w\Dom\u\Dom\y\Dom\s$.  Therefore,~$\Dpsist$ can be written in the
     required form by \autoref{P:fstaction}. Inverting this equation,
     $\Dfst$ can also be written in the required form.  This completes
     the proof of the inductive step and hence the lemma.
  \end{proof}

  \begin{Corollary}\label{C:LpsiGpsi}
    Let $(\s,\t),(\u,\v)\in\Std^2(\Parts)$. Then
    $\Lpsist\Gpsist[\u\v]\ne0$ only if~$\t\Gedom\u$, and
    $\Gpsist[\u\v]\Lpsist\ne0$ only if~$\s\Gedom\v$.  Moreover,
    $\Lpsist\Gpsist[\t\s]=\Lfst\Gfst[\t\s]$ and
    $\Gpsist[\t\s]\Lpsist=\Gfst[\t\s]\Lfst$ are homogeneous of degree $2\defect(\blam)$.
 \end{Corollary}

  \begin{proof}
    Consider the first statement. Using \autoref{L:fpsiExpansions},
    \begin{align*}
      \Lpsist\Gpsist[\u\v]
         &=\Bigl(\Lfst+
          \sum_{\substack{(\w,\x)\in\Std^2(\Parts)\\(\w,\x)\Ldom(\s,\t)}}
             a_{\w\x}\Lfst[\w\x]\Bigr)
           \Bigl(\Gfst[\u\v]+
          \sum_{\substack{(\y,\z)\in\Std^2(\Parts)\\(\y,\z)\Gdom(\u,\v)}}
             b_{\y\z}\Gfst[\y\z]\Bigr)\\
         &=\sum_{\substack{(\w,\x)\in\Std^2(\Parts)\\(\w,\x)\Ledom(\s,\t)}}
           \sum_{\substack{(\y,\z)\in\Std^2(\Parts)\\(\y,\z)\Gedom(\u,\v)}}
              a_{\w\x}b_{\y\z}\Lfst[\w\x]\Gfst[\y\z],
    \end{align*}
    where we set $a_{\s\t}=1=b_{\u\v}$. Therefore,
    $\Lpsist\Gpsist[\u\v]\ne0$ only if $\Lfst[\w\x]\Gfst[\y\z]\ne0$ for
    some $(\w,\x),(\y,\z)\in\Std^2(\Parts)$ with $\w\Ledom\s$, $\x\Ledom\t$,
    $\y\Gedom\u$ and $\z\Gedom\v$. By \autoref{L:FufstFv},
    $\Lfst[\w\x]\Gfst[\y\z]\ne0$ only if $\x=\y$, so this forces
    $\t\Gedom\x=\y\Gedom\u$, as required.  Since
    $\Gpsist[\t\s]\Lpsist=(\Lpsist\Gpsist[\t\s])^*$, this implies that
    if $\Lpsist[\t\s]\Gpsist[\v\u]\ne0$ then $\s\Gedom\v$.  When $\u=\t$
    and $\v=\s$ the last displayed equation shows that
    $\Lpsist\Gpsist[\t\s]=\Lfst\Gfst[\t\s]$. By definition,
    $\Lpsist\Gpsist[\t\s]$ is a homogeneous element of $\Rx(\kx)$ of
    degree $2\defect(\blam)$. Similarly,
    $\Gpsist[\t\s]\Lpsist=\Gfst[\t\s]\Lfst$ is homogeneous of
    defect~$2\defect(\blam)$.
  \end{proof}

  \begin{Definition}\label{D:zlam}
    For $\blam\in\Parts$ set
    $\Lzlam=\Lpsist[\Ltlam\Ltlam]\Gpsist[\Ltlam\Ltlam]$ and
    $\Gzlam=\Gpsist[\Gtlam\Gtlam]\Lpsist[\Gtlam\Gtlam]$.
  \end{Definition}
  \notation{$\Lzlam,\Gzlam$}{Distinguished generators for Specht submodules}[D:zlam]

  By \autoref{L:wlam} we can also write
  $\Lzlam=\Lpsist[\Ltlam\Gtlam]\Gpsist[\Gtlam\Ltlam]$ and
  $\Gzlam=\Gpsist[\Gtlam\Ltlam]\Lpsist[\Ltlam\Gtlam]$. We will not need
  this, but it is not difficult to show that
  $\Lzlam=\Lpsist[\Ltlam\s]\Gpsist[\s\Ltlam]$ and
  $\Gzlam=\Gpsist[\Gtlam\s]\Lpsist[\s\Gtlam]$, for any
  $\s\in\Std(\blam)$.

  In the classical representation theory of the symmetric groups,
  elements very similar to~$\Lzlam$ and $\Gzlam$ are often used as
  distinguished generators for the semisimple Specht modules. The extra
  structure provided by the grading shows that these elements are
  ``almost'' canonical.

  \begin{Proposition}\label{P:zlam}
    Let $\blam\in\Parts[\alpha]$, for $\alpha\in Q^+$. Then
    $\Lzlam=\DefPol F_{\Ltlam}$ and $\Gzlam=\DefPol F_{\Gtlam}$.
    Consequently, $\frac1{\DefPol}\Lzlam$ and~$\frac1{\DefPol}\Gzlam$
    are (nonzero) primitive idempotents in~$\Rx(\Kxx)$ and
    $\tau_\alpha(\Lzlam)_0=1=\tau_\alpha(\Gzlam)_0$.
  \end{Proposition}

  \begin{proof}
    We give the proof only for $\Lzlam$, with the result for $\Gzlam$
    following by symmetry. Since $\Lzlam=F_{\Ltlam}\Lzlam F_{\Ltlam}$ by
    \autoref{L:FufstFv}, it follows that $\Lzlam$ is a scalar multiple
    of~$F_{\Ltlam}=\frac1{\Lgt[\Ltlam]}\Lfst[\Ltlam\Ltlam]$ by
    \autoref{C:gamma}. Then, there exist scalars
    $a_{\w\x},b_{\y\z}\in\Kxx$ such that
    \begin{align*}
        (\Lzlam)^2&= \Lpsist[\Ltlam\Ltlam]\Gpsist[\Ltlam\Ltlam]
                    \Lpsist[\Ltlam\Ltlam]\Gpsist[\Ltlam\Ltlam]
                    && \text{by \autoref{D:psist}},\\
                 &= \Lpsist[\Ltlam\Ltlam]\Bigl(\Gfst[\Ltlam\Ltlam]
     +\sum_{\mathclap{(\x,\w)\Gdom(\Ltlam,\Ltlam)}}a_{\w\x}\Gfst[\w\x]\Bigr)
                  \Bigl(\Lfst[\Ltlam\Ltlam]
     +\sum_{\mathclap{(\y,\z)\Ldom(\Ltlam,\Ltlam)}}b_{\y\z}\Lfst[\y\z]\Bigr)
                  \Gpsist[\Ltlam\Ltlam],
                  &&\text{by \autoref{L:fpsiExpansions},}\\
                &=\Lpsist[\Ltlam\Ltlam]\Gfst[\Ltlam\Ltlam]\Lfst[\Ltlam\Ltlam]
                        \Gpsist[\Ltlam\Ltlam],
                       &&\text{by \autoref{L:FufstFv},}\\
      &=\Lpsist[\Ltlam\Ltlam]\cdot\Ggt[\Ltlam]F_{\Ltlam}\cdot
                  \Lgt[\Ltlam]F_{\Ltlam}\cdot \Gpsist[\Ltlam\Ltlam],
                        &&\text{by \autoref{C:gamma},}\\
                &=\DefPol\Lzlam,&& \text{by \autoref{L:FufstFv}.}
    \end{align*}
    Hence, $\frac1{\DefPol}\Lzlam=F_{\Ltlam}$ is a primitive idempotent
    in $\Rx(\Kxx)$.  Finally, $\tau_\alpha(\Lzlam)=1$ by \autoref{Ex:tauFt}.
  \end{proof}

  Although we do not need this, it is not hard to show that
  $\Lpsist[\Ltlam\Ltlam]\Rn(\kx)\Gpsist[\Ltlam\Ltlam]=\kx\Lzlam$ is a
  free $\kx$-module of rank~$1$, giving another way to prove that
  $S_\blam(\K[x^\pm])$ is an irreducible $\Rx(\K[x^\pm])$-module.

  We have reached the main results of this section.

  \begin{Theorem}\label{T:Symmetric}
    Suppose that $(\s,\t),(\u,\v)\in\Std^2(\Parts[\alpha])$, for
    $\alpha\in Q^+_n$. Then
    \[
        \<\Lpsist, \Gpsist[\u\v]\>_\alpha = \begin{cases*}
                  1& if $(\s,\t)=(\v,\u)$,\\
                  0& if $(\s,\t)\Gnedom(\v,\u)$.
        \end{cases*}
    \]
  \end{Theorem}

  \begin{proof}
    By definition and \autoref{L:bilinear},
    $\<\Lpsist,\Gpsist[\u\v]\>_\alpha
              =\tau_\alpha(\Lpsist\Gpsist[\u\v])
              =\tau_\alpha(\Gpsist[\u\v]\Lpsist)$.
    Hence, $\<\Lpsist,\Gpsist[\u\v]\>_\alpha=0$ unless $\t\Gedom\u$ and
    $\s\Gedom\v$ by \autoref{C:LpsiGpsi}.  Now suppose that $\u=\t$ and
    $\v=\s$ and consider the inner product
    $\<\Lpsist,\Gpsist[\s\t]\>_\alpha=\tau_\alpha(\Lpsist\Gpsist[\s\t])$.
    Using \autoref{L:bilinear},
    \begin{align*}
        \<\Lpsist,\Gpsist[\t\s]\>_\alpha
           &= \tau_\alpha(\Lpsist\Gpsist[\t\s])_0
            = \tau_\alpha(\psi_{\Ldt[s]}\Lpsist[\Ltlam\t]
                    \psi_{\Gdt}\Gpsist[\Gtlam\s])_0\\
           &= \tau_\alpha(\Lpsist[\Ltlam\Ltlam]\Lpsist[\Gtlam\s]\psi_{\Ldt[s]})_0
            = \tau_\alpha(\Lpsist[\Ltlam\Ltlam]\Gpsist[\Ltlam\Ltlam])_0,
              && \text{by two applications of \autoref{L:wlam}},\\
           &=\tau_\alpha(\Lzlam)_0
            = \DefPol\tau_\alpha(F_{\Ltlam})_0,&&\text{by \autoref{P:zlam}},\\
           &= 1,
   \end{align*}
   where the last equality follows from \autoref{Ex:tauFt}.
  \end{proof}

  \subsection{Cellular bases for \texorpdfstring{$\Rx(\kx)$}{R^Λ(k[x])}}\label{SS:kcellular}
  We can now prove that $\Rx(\kx)$ is a $\kx$-cellular algebra.  In
  particular, this proves a stronger form of
  \autoref{MT:ACCellular},  our first main result from the introduction.

  \begin{Theorem}\label{T:kcellular}
    Suppose that $(\bc,\br)$ is a graded content system with values in $\kx$.
    Then $\Rx(\kx)$ is a graded $\kx$-cellular algebra with $\kx$-cellular bases:
    \begin{enumerate}
      \item $\set{\Lpsist| (\s,\t)\in\Std^2(\Parts)}$ with weight poset
      $(\Parts,\Ledom)$ and degree function $\Ldeg$.
      \item $\set{\Gpsist| (\s,\t)\in\Std^2(\Parts)}$ with weight poset
      $(\Parts,\Gedom)$ and degree function $\Gdeg$.
    \end{enumerate}
  \end{Theorem}

  \begin{proof}
     By \autoref{P:KLRfree}, $\Rx(\kx)$ is free as a $\kx$-module, so
     $\Rx(\kx)$ naturally embeds into the~$\Kxx$-algebra
     $\Rx(\Kxx)\cong\Kxx\otimes_{\kx}\Rx(\kx)$. In particular,
     the $\kx$-rank of $\Rx(\kx)$ is equal to the $\Kxx$-rank of $\Rx(\Kxx)$.

     We only show that $\set{\Lpsist}$ is a $\kx$-cellular basis of
     $\Rx(\kx)$, as the $\kx$-cellularity of $\set{\Gpsist}$ follows by
     symmetry. Since $\Rx(\kx)=\bigoplus_{\alpha\in
     Q^+_n}\Rx[\alpha](\kx)$, it is enough to show that
     $\set{\Lpsist|(\s,\t)\in\Std^2(\Parts[\alpha])}$ is a
     $\kx$-cellular basis of~$\Rx[\alpha](\kx)$, for $\alpha\in Q^+_n$.
     By \autoref{T:Kcellular}, $\set{\Lpsist|(\s,\t)\in\Parts[\alpha]}$
     is a $\Kxx$-cellular basis of $\Rx[\alpha](\Kxx)$. Therefore, to
     prove the theorem it is enough to show that
     $\set{\Lpsist(\s,\t)\in\Parts[\alpha]}$ spans $\Rx[\alpha](\kx)$
     and that the structure constants for this basis belong to $\kx$.

     Let $(\s,\t)\in\Std^2(\Parts[\alpha])$. Using \autoref{T:Symmetric}
     and Gaussian elimination to argue by induction on dominance, there
     exist homogeneous elements $\eta_{\u\v}^\Gdom\in\Rx[\alpha](\kx)$
     such that
     $\<\Lpsist,\eta^{\Gdom}_{\u\v}\>_\alpha=\delta_{(\s,\t)(\v,\u)}$
     and
     $\eta^\Gdom_{\u\v}=\Gpsist[\u\v]+\sum_{(\x,\y)\Gdom(\u,\v)}e^\Gdom_{\x\y}\Gpsist[\x\y]$,
     for homogeneous scalars $e^\Gdom_{\x\y}\in\kx$. Therefore,
     if~$h\in\Rx[\alpha](\kx)$ then
      \[
              h = \sum_{(\u,\v)\in\Std^2(\Parts[\alpha])}
              \<h,\eta^\Gdom_{\u\v}\>_\alpha\Lpsist[\u\v].
      \]
      In particular, the set
      $\set{\Lpsist|(\s,\t)\in\Std^2(\Parts[\alpha])}$ spans
      $\Rx[\alpha](\kx)$ as a $\kx$-module.  Hence, $\set{\Lpsist}$ is a
      basis of~$\Rx[\alpha](\kx)$ by \autoref{T:Kcellular}. Moreover, if
      $h\in\Rx(\kx)$ then $h\Lpsist\in\Rx(\kx)$, so
      $\<h\Lpsist,\eta^\Gdom_{\u\v}\>_\alpha\in\kx$, for
      $(\u,\v)\in]\Std^2(\Parts[\alpha])$. Therefore,
      \[
          h\Lpsist = \sum_{(\s,\t)\in\Std^2(\Parts[\alpha])}
              \<h\Lpsist,\eta^\Gdom_{\u\v}\>_\alpha\Lpsist[\u\v].
      \]
      showing that the structure constants of
      $\set{\Lpsist|(\s,\t)\in\Parts[\alpha]}$ belong to $\kx$.  Hence,
      $\set{\Lpsist|(\s,\t)\in\Parts[\alpha]}$ is a $\kx$-cellular basis
      of~$\Rx[\alpha](\kx)$ by \autoref{T:Kcellular}.
  \end{proof}

  The strategy used to prove \autoref{T:kcellular} is quite general. For
  example, an easy modifications this argument gives a streamlined proof
  of the fact that the Murphy basis of \cite[Theorem~3.26]{DJM:cyc} is a
  cellular basis of the cyclotomic Hecke algebras of
  type~$A$~\cite[Theorem~3.26]{DJM:cyc}.

  \begin{Remark}
    In type $\Aone$, even in the ungraded world, pairs of dual bases for
    the algebras $\Rx(\kx)$ are not known. It seems hard to explicitly
    describe the basis $\set{\eta^\Gdom_{\s\t}}$ that is dual to
    $\set{\Lpsist}$. Similarly, it is hard to describe the basis
    $\set{\eta^\Ldom_{\s\t}}$ that is dual to $\set{\Gpsist}$. On the
    other hand, using \autoref{T:kcellular}, it is straightforward to
    check that $\set{\eta^\Gdom_{\s\t}}$ and $\set{\eta^\Ldom_{\s\t}}$
    are $\kx$-cellular bases of $\Rx(\kx)$.
  \end{Remark}

  As noted in \autoref{Ex:ContentSystem}, content systems $(\bc,\br)$
  do not always exist in positive characteristic. Nonetheless, by
  base-change, \autoref{T:kcellular} gives cellular bases over other
  rings. Indeed, since \autoref{Ex:ContentSystem} gives content systems
  with values in~$\Z[x]$ for quivers of types $\Aone$ and $\Cone$, we
  obtain cellular bases over~$\k[x]$ for arbitrary rings~$\k$.

  \begin{Corollary}\label{C:RxCellular}
    Suppose that $(\bc,\br)$ is a graded content system with values in $\kx$ and
    let $K$ be  commutative domain with~$1$ that is a $\kx$-algebra. Then
    $\Rx(K)$ is a graded $K$-cellular algebra with cellular
    bases:
    \begin{enumerate}
      \item $\set{\Lpsist| (\s,\t)\in\Std^2(\Parts)}$ with weight poset
      $(\Parts,\Ledom)$ and degree function $\Ldeg$.
      \item $\set{\Gpsist| (\s,\t)\in\Std^2(\Parts)}$ with weight poset
      $(\Parts,\Gedom)$ and degree function $\Gdeg$.
    \end{enumerate}
  \end{Corollary}

  \begin{proof}
    This is immediate from \autoref{T:kcellular} since
    $\Rx(K)\cong K\otimes_{\kx}\Rx(\kx)$.
  \end{proof}

  Essentially as an important special case, this implies that the
  (standard) cyclotomic KLR algebras $\Rn(K)$ of type $\Aone$ or $\Cone$
  are cellular over any ring $K$.

  \begin{Corollary}\label{C:RnCellular}
    Let $K$ be commutative domain with~$1$ and suppose that $\Rn(K)$ is a
    cyclotomic KLR algebra of type $\Aone$, $A_\infty$,  $\Cone$ or
    $C_\infty$.  Then $\Rn(K)$ is a graded cellular algebra with cellular
    bases:
    \begin{enumerate}
      \item $\set{\Lpsist| (\s,\t)\in\Std^2(\Parts)}$ with weight poset
      $(\Parts,\Ledom)$ and degree function $\Ldeg$.
      \item $\set{\Gpsist| (\s,\t)\in\Std^2(\Parts)}$ with weight poset
      $(\Parts,\Gedom)$ and degree function $\Gdeg$.
    \end{enumerate}
  \end{Corollary}

  \begin{proof}
    For quivers of type $\Aone$ of $\Cone$, by
    \autoref{L:IntegralContentSystems} there exist graded content system
    $(\bc,\br)$ with values in $\Z[x]$ for a deformed cyclotomic KLR
    algebra $\Rx(|Z[x])$. Therefore,
    $\Rn(K)\cong K\otimes_{\Z[x]}\Rx(\Z[x])$ as $K$-algebras, where $K$
    is considered as a $\Z[x]$-algebra by letting $x$ act as
    multiplication by~$0$, so the result follows by \autoref{T:kcellular}.
    For quivers of type $A_\infty$ of $C_\infty$, by taking $e$
    sufficiently large, this implies that the cyclotomic KLR algebras of
    type $A_\infty$ and $C_\infty$ are cellular; compare with
    \cite[Corollary~2.10]{HuMathas:SeminormalQuiver}.
  \end{proof}

  \begin{Remark}
    For the cyclotomic KLR algebras of type $\Aone$
    \autoref{C:RnCellular} recovers, with considerably less effort, the
    main theorem of Li~\cite{GeLi:IntegralKLR}, which generalises the
    main theorem of \cite{HuMathas:GradedCellular} to give an integral
    basis of~$\Rn(\Z)$. The papers
    \cite{Bowman:ManyCellular,MathasTubbenhauer:Subdivision} use
    Websters' diagrammatic KLRW algebras to construct different cellular
    bases for the cyclotomic KLR algebras of types~$\Aone$ and~$\Cone$,
    which depend on a choice of ``loading''. In type~$\Aone$,
    Bowman~\cite[Proposition~7.3]{Bowman:ManyCellular} has shown that
    the transition matrix between the $\Gpsist[]$-basis of~$\Rn(\k)$ and
    the ``asymptotic Webster diagram basis'' is unitriangular. In
    type~$\Cone$, we do not know the relationship between
    the cellular bases considered in this paper and those
    in~\cite{MathasTubbenhauer:Subdivision}, although it seems likely
    that Bowman's arguments generalise to show that the transition
    matrices between these bases is unitriangular in the ``asymptotic
    case''.
  \end{Remark}

  \begin{Remark}
     The cellular bases in \autoref{T:kcellular} give graded Specht
     modules for the cyclotomic KLR algebras $\Rn(\k)$. In type~$\Aone$
     this recovers the results of
     \cite{BKW:GradedSpecht,HuMathas:GradedCellular}. Ariki, Park
     and Speyer~\cite{ArikiParkSpeyer:C} have given a conjectural
     construction of graded Specht modules in type~$\Cone$ using
     analogues of the homogeneous Garnir relations from
     \cite{KMR:UniversalSpecht}, and they have proved these conjectures
     in type $C_\infty$. As shown in~\cite{Mathas:Intertwiners}, it is
     easy to prove the conjectures of \cite{ArikiParkSpeyer:C} using
     \autoref{T:kcellular}.
  \end{Remark}

  It is very difficult to do calculations with the cyclotomic KLR
  algebras $\Rn$. In contrast, it is very easy to calculate with the
  $\psi$-bases of $\Rx(\kx)$ because the transition matrices to the
  corresponding seminormal bases are unitriangular by
  \autoref{P:psiTriangular} and the action of $\Rx(\kx)$ on the
  seminormal bases is completely described by \autoref{P:fstaction}.
  The rest of this paper can be viewed as theoretical applications of
  this observation. In a different direction, this observation is used
  in~\cite{ChungMathasSpeyer:TypeCDecomp,Mathas:SageKLR} to implement
  the cyclotomic KLR algebras of types $\Aone$ and $\Cone$ in
  \textsc{SageMath}~\cite{sage}.

  An $R$-algebra $A$ is a \emph{graded symmetric algebra} algebra if
  there is a non-degenerate homogeneous bilinear form $\<\ ,\
  \>\map{A\times A}R$ of degree~$d$ such that $\<ab,c\>=\<a,bc\>$, for
  all $a,b,c\in A$; compare \cite[Definition~66.1]{C&R}.  Hence,
  combining \autoref{T:Symmetric} and \autoref{T:kcellular} yields:

  \begin{Corollary}
    Let $\alpha\in Q^+_n$. Then $\Rx[\alpha](\kx)$ is a graded symmetric
    algebra with homogeneous trace form of degree $-2\defect(\alpha)$.
  \end{Corollary}

  The bilinear form $\<\ ,\ \>_\alpha$ is defined over $\k$. So, in view
  of \autoref{L:IntegralContentSystems}, we obtain the corresponding results for
  the cyclotomic KLR algebras $\Rn(\Z)$.

  \begin{Corollary}\label{C:Symmetric}
    Let $\alpha\in Q^+_n$. Then $\Rn[\alpha](\Z)$ is a graded symmetric
    algebra with homogeneous trace form of degree $-2\defect(\alpha)$.
    In particular, the cyclotomic Hecke algebras of type $\Aone$ and
    $\Cone$ are graded symmetric algebras over any ring.
  \end{Corollary}

  For the cyclotomic KLR algebras of type $\Aone$, \autoref{C:Symmetric}
  was first proved as \cite[Corollary~6.18]{HuMathas:GradedCellular}.
  Later, Kashiwara~\cite{Kashiwara:KLRBiadjointness} and
  Webster~\cite[Remark~3.19]{Webster:HigherRep} used categorical and
  diagrammatic arguments, respectively, to show that cyclotomic KLR
  algebras of symmetrisable type are graded symmetric algebras.

  As our first application of the trace form on $\Rx(\kx)$ we show that
  the graded Specht modules $\LSlam(\kx)$ and $\GSlam(\kx)$ are dual to
  each other, up to shift.

  \begin{Proposition}\label{P:SpechtDual}
    Suppose that $K$ is a $\kx$-module and let $\blam\in\Parts[\alpha]$,
    for $\alpha\in Q^+_n$. Then
    \[
      \LSlam(K)\cong q^{\defect(\blam)}\GSlam(K)^\circledast
        \And
      \GSlam(K)\cong q^{\defect(\blam)}\LSlam(K)^\circledast
    \]
    as $\Rx(\kx)$-modules.
  \end{Proposition}

  \begin{proof}
    The two isomorphisms are equivalent so we  prove only the first
    one. For $\s\in\Std(\blam)$ let
    $\theta_\s\in q^{\defect(\blam)}\GSlam(K)^\circledast$ be the unique
    $K$-linear map such that
    \[
      \theta_\s(\Gpsis[\t])=\<\Lpsist[\Ltlam\s],\Gpsist[\t\Ltlam]\>_\alpha,
      \qquad\text{for }\t\in\Std(\blam).
    \]
    Define a homomorphism
    $\Theta\map{\LSlam(K)}{\GSlam(K)^\circledast}$ by
    $\Theta(\Lpsis)=\theta_\s$, for $\s\in\Std(\blam)$.
    By \autoref{C:degCodeg},
    $\Ldeg(\Ltlam)+\Gdeg(\Gtlam)=\defect(\blam)$, so $\Theta$ is a
    homogeneous map of degree zero into
    $q^{\defect(\blam)}\bigl(\GSlam(K)\bigr)^\circledast$. In view of
    \autoref{L:bilinear}, $\Theta$ is an $\Rx(K)$-module homomorphism
    and, by \autoref{T:Symmetric}, it is an isomorphism of $K$-modules.
  \end{proof}

  In particular, the specialisation of $\ux$ to~$0$, which corresponds
  to taking $K=\k$, shows that
  \[
      \LSlam(\k)\cong q^{\defect(\blam)}\GSlam(\k)^\circledast
        \And
      \GSlam(\k)\cong q^{\defect(\blam)}\LSlam(\k)^\circledast
  \]
  as $\Rn(\k)$-modules. In view of \autoref{L:IntegralContentSystems}, and base
  change, $\k$ can be an arbitrary ring. In type~$\Aone$, this recovers
  \cite[Proposition~6.19]{HuMathas:GradedCellular}.

  As the last result in this section, we note that combining
  \autoref{L:fpsiExpansions} and \autoref{T:kcellular} gives the
  following useful strengthening of \autoref{P:integral}(b).

  \begin{Corollary}\label{C:ypsi}
    Suppose that $1\le m\le n$ and $\s,\t\in\Std(\blam)$, for
    $\blam\in\Parts$. Then
    \[
         y_m\Lpsist = \bc_m(\s)\Lpsist
         + \sum_{(\u,\v)\Ldom(\s,\t)} c_{\u\v}\Lpsist[\u\v]
              \And
         y_m\Gpsist = \bc_m(\s)\Gpsist
         + \sum_{(\u,\v)\Gdom(\s,\t)} d_{\u\v}\Gpsist[\u\v]
    \]
    for some $c_{\u\v},  d_{\u\v}\in\kx$ such that
    \begin{itemize}
      \item $c_{\u\v}\ne0$ only if
    $\br(\u)=\br(\s)$, $\br(\v)=\br(\t)$ and either
    $\bmu\Ldom\blam$, or $\bmu=\blam$, $\v=\t$ and $\u\Ldom\s$,
      \item $d_{\u\v}\ne0$ only if
    $\br(\u)=\br(\s)$, $\br(\v)=\br(\t)$ and either
    $\bmu\Gdom\blam$, or $\bmu=\blam$, $\v=\t$ and $\u\Gdom\s$.
    \end{itemize}
  \end{Corollary}

  Notice, in particular, that the coefficients of the leading term
  $\Dpsist$ are zero in the standard KLR algebras~$\Rn(\k)$ since
  $\bc_m(\s)$ is a polynomial in $\kx$ with zero constant term by the
  degree requirements of \autoref{D:ContentSystem}. Hence, it follows
  that $y_r^{|\Std(\bi)|}\ei=0$ in $\Rx(\k)\cong\Rn(\k)$, generalising
  \cite[Corollary~4.31]{HuMathas:SeminormalQuiver}.

\section{Graded Specht and simple modules}
  This chapter uses the cellular bases of \autoref{T:kcellular} to
  construct complete sets of graded simple modules for $\Rx(\Kx)$. We
  prove some identities relating the decomposition matrices associated
  to the different bases and over different fields. Some of these
  results will be instrumental in the next chapter when we show than the
  algebra $\bigoplus_{n\ge0}\Rx(\Kx)$ categorifies the integral highest
  weight module $L(\Lambda)$ of the corresponding Kac-Moody algebra.

  In this chapter we slightly weaken the assumptions of the last two
  chapters and assume that $(\Gamma,\Qbx,\Wbx)$ is a $\kx$-deformation
  of a standard cyclotomic KLR datum $(\Gamma,\Qb,\Wb)$ and $(\bc,\br)$
  is a (graded) content system with values in~$\kx$. Assume that $\K$ is
  a field that is a $\k$-algebra, so that $\Rx(\Kx)$ is a graded
  $\Kx$-cellular algebra by \autoref{C:RxCellular}. As explained below,
  the results in this chapter apply to the standard cyclotomic KLR
  algebras of type $\Aone$, $A_\infty$, $\Cone$ and $C_\infty$ since the
  graded irreducible $\Rx(\Kx)$-modules and the graded irreducible
  $\Rn(\K)$-modules coincide.

  \subsection{Irreducible modules}\label{SS:GradedSimples}
  This section describes the irreducible graded $\Rx$-modules, both as
  subquotients and as submodules of $\Rx$. Recall that $\K$ is a field
  that is a $\k$-algebra.

  Let $L$ be a $\kx$-module. Fix $\blam\in\Parts$. Via
  \autoref{E:CellularForm}, the $\kx$-cellular algebra framework equips
  the Specht modules $\LSlam(L)$ and~$\GSlam(L)$ with homogeneous
  symmetric associative bilinear forms that are characterised by
  \begin{equation}\label{E:SpechtForm}
      \<\Lpsis,\Lpsis[\t]\>^\Ldom_\blam\Lpsist[\u]=\Lpsist[\u\s]\Lpsis[\t]
        \And
      \<\Gpsis,\Gpsis[\t]\>^\Gdom_\blam\Gpsist[\u]=\Gpsist[\u\s]\Gpsis[\t],
  \end{equation}%
  \notation{$\<\ ,\ \>^\Ldom_\blam$, $\<\ ,\ \>^\Gdom_\blam$}
           {Bilinear forms on $\LSlam$ and $\GSlam$}%
  for $\s,\t\u\in\Std(\blam)$. The \emph{radicals} of
  the graded Specht modules are the submodules defined by:
  \begin{align*}
    \rad\LSlam(L) &= \set{a\in\LSlam(L)|
            \<a,b\>^\Ldom_\blam=0 \text{ for all }b\in\LSlam(L)}\\
    \rad\GSlam(L) &= \set{a\in\LSlam(L)|
            \<a,b\>^\Gdom_\blam=0 \text{ for all }b\in\GSlam(L)}.
  \end{align*}
  Note that these definitions make sense for any (graded) $\kx$-module~$L$.

  \begin{Definition}\label{D:RnSimples}
    Let $\bmu\in\Parts[\alpha]$, for $\alpha\in Q^+_n$. Let $L$ be a
    $\kx$-module and define
    \[
      \LDmu(L)=\LSlam[\bmu](L)/\rad\LSlam[\bmu](L)
        \And
      \GDnu(L)=\GSlam[\bmu](L)/\rad\GSlam[\bmu](L)
    \]
    If $K=\Kx$ then $\LDmu(\K)$ and $\GDnu(\K)$ are
    $\Rx(\Kx)$-modules. Set
    \[
      \LKlesh[\alpha]=\set{\bmu\in\Parts[\alpha]|\LDmu(\K)\ne0}
      \qquad\text{and}\qquad
      \GKlesh[\alpha]=\set{\bmu\in\Parts[\alpha]|\GDnu(\K)\ne0}.
    \]
    Let $\LKlesh=\bigcup_{\alpha\in Q^+_n}\LKlesh[\alpha]$ and
    $\GKlesh=\bigcup_{\alpha\in Q^+_n}\GKlesh[\alpha]$.
  \end{Definition}
  \notation{$\LDmu,\GDnu$}{Simple $\Rx$-modules defined by the $\Lpsist$
            and $\Gpsist$ bases}[D:RnSimples]
  \notation{$\LKlesh,\GKlesh$}{Indexing sets for simple $\Rx$-modules}

  When the choice of $L$ is clear (usually, $L=\K)$, then
  we write $\LDmu$ and $\GDnu$.

  As $\K$-vector spaces, with respect to the $\ux$-grading, $\LDmu(\K)$
  is the degree zero component of $\LDmu(\Kx)$ and $\GDnu(\K)$ is the
  degree zero component of $\GDnu(\Kx)$. The modules $\LDmu(\Kx)$ and
  $\GDnu(\Kx)$ are free $\Kx$-modules, and so infinite dimensional
  $\K$-vector spaces if $\ux\ne\emptyset$, whereas $\LDmu(\K)$ and
  $\GDnu(\K)$ are finite dimensional $\K$-vector spaces upon which each
  $x\in\ux$ acts as multiplication by~$0$.

  Even though our notation does not reflect this, the sets $\LKlesh$ and
  $\GKlesh$ depend on~$\charge$ and, \textit{a priori}, on the
  field~$\K$.  In type $\Aone$ the sets $\LKlesh$ and $\GKlesh$ have
  already been determined \cite{BK:GradedDecomp,Ariki:class}.  In
  \autoref{T:KleshchevSimples} below we give a uniform characterisation
  of $\LKlesh$ and $\GKlesh$ in types $\Aone$ and $\Cone$. In
  particular, this result shows that the sets $\LKlesh$ and $\GKlesh$ do
  not depend on the choice of field~$\K$.

  Combining \autoref{T:Kcellular} and \autoref{T:CellularSimples} we obtain:

  \begin{Theorem}\label{T:GradedSimples}
    Let $\Domin$ and suppose that $K=\Kx$. Then
    $\set[\big]{q^z\DDmu(\K)|\bmu\in\DKlesh\text{ and } z\in\Z}$
    is a complete set of pairwise non-isomorphic irreducible graded
    $\Rx(\Kx)$-modules. Moreover, $\DDmu(\K)$ is a graded self-dual
    $\Rx(\Kx)$-module, for $\bmu\in\DKlesh$.
  \end{Theorem}

  By \autoref{C:RnCellular} and \autoref{Ex:GradedSimples}, the set of
  isomorphism classes of irreducible graded $\Rn(\K)$-module coincides
  with the set of isomorphism classes of irreducible $\Rx(\Kx)$-modules.
  The point is that if $L$ is a $\Kx$-module and some $x\in\ux$ does not
  act on~$L$ as multiplication by zero then $\DDmu(L)$ is not
  irreducible.

  We next show how to realise the graded simple modules of
  $\Rx(\K[\ux])$ as submodules of $\Rx(\K[\ux])$, up to shift.  To do
  this we first need a similar description of the Specht modules, for
  which we use the elements~$\Lzlam$ and~$\Gzlam$ from \autoref{D:zlam}.
  Extending the definition of $\Dzlam$, for $\s\in\Std(\blam)$ set
  \[
     \Lzlams=\psi_{\Ldt[\s]}\Lzlam=\Lpsist[\s\Ltlam]\Gpsist[\Ltlam\Ltlam]
     \qquad\text{and}\qquad
     \Gzlams=\psi_{\Gdt[\s]}\Gzlam=\Gpsist[\s\Gtlam]\Lpsist[\Gtlam\Gtlam].
  \]

  \begin{Lemma}\label{L:SpechtSubmodule}
    Let $\blam\in\Parts$. Then there are $\Rx(\kx)$-module
    isomorphisms
    \[
      \Rx(\kx)\Lzlam\cong q^{\defect(\blam)+\Gdeg\Ltlam}\LSlam
        \And
      \Rx(\kx)\Gzlam\cong q^{\defect(\blam)+\Ldeg\Gtlam}\GSlam.
    \]
    Moreover, these modules have bases $\set{\Lzlam|\s\in\Std(\blam)}$
    and $\set{\Gzlam|\s\in\Std(\blam)}$, respectively.
  \end{Lemma}

  \begin{proof}
    Let $\Doming$. By \autoref{C:LpsiGpsi}, there is a well-defined,
    homogeneous, $\Rx(\kx)$-module homomorphism
    $\pi^\Dom_\blam\map{q^{\defect(\blam)+\Ddeg\Dtlam}\DSlam}\Rx(\kx)\Dzlam$
    given by
    \[
         \pi^\Dom_\blam\bigl(\Dpsist[\s\Dtlam]+\DRlam\bigr)
             =\Dpsist[\s\Dtlam]\psi^{\doM}_{\Dtlam\Dtlam}
             =\Dzlams, \qquad\text{for }\s\in\Std(\blam).
    \]
    By \autoref{T:Kcellular}, $\pi^\Dom_\blam$ is homogeneous of degree
    zero.  The set $\set{\Dzlams|\s\in\Std(\blam)}$ is a basis for the
    image of $\pi^{\Dom}$ since multiplying by the idempotents $F_\t$,
    for $\t\in\Std(\blam)$, shows that these elements are linearly
    independent. Hence, $\Rx(\kx)\Dzlam=\im\pi^\Dom$ in view of
    \autoref{P:zlam}. The result follows.
  \end{proof}

  By \autoref{D:zlam},
  $\Gpsist[{\Gtlam\Ltlam}]\Lzlam =\Gzlam\Gpsist[{\Ltlam\Gtlam}]$
  and $\Lpsist[{\Ltlam\Gtlam}]\Gzlam =\Lzlam\Gpsist[{\Ltlam\Gtlam}]$,
  for $\blam\in\Parts$. Applying \autoref{L:wlam},
  \begin{equation}\label{E:psiZ}
      \Lpsist[\Ltlam\Gtlam]\Gzlam
        =\Lpsist[\Ltlam\Gtlam]\cdot\Gpsist[\Gtlam\Gtlam]\Lpsist[\Gtlam\Gtlam]
        =\Lpsist[\Ltlam\Ltlam]\Gpsist[\Ltlam\Ltlam]\cdot\Lpsist[\Ltlam\Gtlam]
        =\Lzlam\Lpsist[\Ltlam\Gtlam]
  \end{equation}
  and, similarly,
  $\Gpsist[\Gtlam\Ltlam]\Lzlam=\Gzlam\Gpsist[\Gtlam\Ltlam]$. The next
  result, which has its origins in the work of James~\cite[\S11]{James},
  shows that these elements generate the simple $\Rx(\kx)$-modules.

  \begin{Theorem}\label{T:SimpleSubmodule}
    Suppose $\bmu\in\LKlesh[\alpha]$ and $\bnu\in\GKlesh[\alpha]$,
    for $\alpha\in Q^+$. As $\Rx(\Kx)$-modules,
    \[
      q^{2\defect(\bmu)+\Ldeg\Gtlam}\LDmu(\K)
          \cong\Rx(\K)\Lzlam[\bmu]\Lpsist[{\Ltlam[\bmu]\Gtlam[\bmu]}]
          \quad\text{and}\quad
      q^{2\defect(\bnu)+\Gdeg\Ltlam}\GDnu(\K)
          \cong\Rx(\K)\Gzlam[\bnu]\Gpsist[{\Gtlam[\bnu]\Ltlam[\bnu]}]
    \]
    In particular, $\LDmu(\K)\ne0$ if and only
    if~$\Lzlam[\bmu]\Lpsist[{\Ltlam[\bmu]\Gtlam[\bmu]}]\ne0$
    and $\GDnu(\K)\ne0$ if and only
    if~$\Gzlam[\bnu]\Gpsist[{\Gtlam[\bnu]\Ltlam[\bnu]}]\ne0$
    in~$\Rx(\Kx)$.
  \end{Theorem}

  \begin{proof}
    We prove only the first isomorphism as the second isomorphism
    follows by symmetry. We first prove some related results
    over~$\kx$. As in the proof of \autoref{P:SpechtDual}, define
    $\theta_\t\in\LSlam[\bmu](\kx)^\circledast$ by
    $\theta_\t(\Lpsis[\u])
       =\tau_\alpha(\Gpsist[{\u\Ltlam[\bmu]}]\Lpsist[{\Ltlam[\bmu]\t}])$,
    for $\t,\u\in\Std(\bmu)$. Using \autoref{E:psiZ},
    \autoref{L:SpechtSubmodule} and \autoref{P:SpechtDual}, there are
    homogeneous $\Rx(\kx)$-module homomorphisms (the reader is welcome to
    determine the degrees of these maps),
    \[
    \LSlam[\bmu](\kx)
         \xrightarrow{\ f \ } \Rx(\kx)\Lzlam[\bmu]
         \xrightarrow{\ g \ } \Rx(\kx)\Gzlam[\bmu]
         \xrightarrow{\ h \ } \LSlam[\bmu](\kx)^\circledast,
    \]
    given by $f(\Lpsis)=\Lzlams=\psi_{\Ldt[\s]}\Lzlam[\bmu]$,
    $g(a)=a\Lpsist[{\Ltlam[\bmu]\Gtlam[\bmu]}] $ and
    $h(\Gzlams[\t])=\theta_\t$, for tableaux $\s,\t\in\Std(\bmu)$ and
    $a\in\Rx(\kx))$.  By \autoref{L:SpechtSubmodule} and the
    proof of \autoref{P:SpechtDual}, $f$ and $h$ are isomorphisms. Let
    $\theta=h\circ g\circ f$ be the composition of these three maps. To
    determine $\theta$, for $\s\in\Std(\bmu)$ write
    \[
       \Lzlams = \Lpsist[{\s\Ltlam[\bmu]}]
           \Gpsist[{\Ltlam[\bmu]\Ltlam[\bmu]}]
         = \sum_{(\u,\v)\in\Std^2(\Parts)} a_{\u\v}\Gpsist[\u\v],
           \qquad\text{ for }a_{\u\v}\in L.
    \]
    By \autoref{CA:act} and \autoref{T:kcellular}, $a_{\u\v}\ne0$ only
    if $\Shape(\v)\Ledom\bmu$, with equality only if $\v=\Ltlam[\bmu]$.
    Therefore,
    \[
         \theta(\Lpsis)=h(\Lzlams\Lpsist[{\Ltlam[\bmu]\Gtlam[\bmu]}])
         = h\Bigl(\ \sum_{\mathclap{(\u,\v)\in\Std^2(\Parts)}}\ a_{\u\v}
              \Gpsist[\u\v]\Lpsist[{\Ltlam[\bmu]\Gtlam[\bmu]}]\,\Bigr)
         = \sum_{\u\in\Std(\bmu)} a_{\u\Ltlam[\bmu]}
               h\bigl(\Gpsist[{\u\Ltlam[\bmu]}]\Lpsist[{\Ltlam[\bmu]
                     \Gtlam[\bmu]}]\bigr)
         = \sum_{\u\in\Std(\bmu)} a_{\u\Ltlam[\bmu]}\theta_{\u},
    \]
    where we have used \autoref{C:LpsiGpsi}, for the third equality, and \autoref{L:wlam}
    for the last equality together with the identity
    $\Gzlams[\u]=\Gpsist[{\u\Gtlam[\bmu]}] \Lpsist[{\Gtlam[\bmu]\Gtlam[\bmu]}]
                =\Gpsist[{\u\Ltlam[\bmu]}] \Lpsist[{\Ltlam[\bmu]\Gtlam[\bmu]}]$.
    Consequently, since $\tau_\alpha$ is a trace form,
    \begin{align*}
        \theta(\Lpsis)(\Lpsis[\t]) &= \sum_{\u\in\Std(\bmu)}
             a_{\u\Ltlam[\bmu]} \theta_\u(\Lpsis[\t])
        = \sum_{\u\in\Std(\bmu)} a_{\u\Ltlam[\bmu]}\tau_\alpha\bigl(
             \Gpsist[{\u\Ltlam[\bmu]}]\Lpsist[{\Ltlam[\bmu]\t}]\bigr)
        =\tau_\alpha\Bigr(\sum_{\u\in\Std(\bmu)}a_{\u\Ltlam[\bmu]}
             \Gpsist[{\u\Ltlam[\bmu]}]\Lpsist[{\Ltlam[\bmu]\t}]\Bigr)\\
       &=\tau_\alpha\Bigr(\sum_{(\u,\v)\in\Std^2(\Parts)}a_{\u\v}
             \Gpsist[{\u\v}]\Lpsist[{\Ltlam[\bmu]\t}]\Bigr),
             \qquad \text{by \autoref{C:LpsiGpsi}},\\
       &=\tau_\alpha\bigl(\Lzlams\Lpsist[{\Ltlam[\bmu]\t}]\bigr)
        =\tau_\alpha\bigl(\Lpsist[{\s\Ltlam[\bmu]}]
           \Gpsist[{\Ltlam[\bmu]\Ltlam[\bmu]}]\Lpsist[{\Ltlam[\bmu]\t}]\bigr)
        =\tau_\alpha\bigl( \Gpsist[{\Ltlam[\bmu]\Ltlam[\bmu]}]
            \Lpsist[{\Ltlam[\bmu]\t}]\Lpsist[{\s\Ltlam[\bmu]}]\bigr)\\
       &=\<\Lpsis[\t],\Lpsis\>_\blam^\Ldom\tau_\alpha\bigl(
            \Gpsist[{\Ltlam[\bmu]\Ltlam[\bmu]}]
            \Lpsist[{\Ltlam[\bmu]\Ltlam[\bmu]}]\bigr)
        =\<\Lpsis[\t],\Lpsis\>_\blam^\Ldom\tau_\alpha(\Lzlam)
        =\<\Lpsis[\t],\Lpsis\>_\blam^\Ldom,
    \end{align*}
    where the first equality on the last line uses \autoref{C:LpsiGpsi},
    and the definition of the inner product on~$\LSlam[\bmu](\kx))$, and
    the last equality follows by \autoref{P:zlam}. Hence, ignoring the
    degree shift, $\theta$ is the natural $\kx$-linear map from
    $\LSlam[\bmu](\kx))\to\LSlam[\bmu](\kx))^\circledast$ induced by the
    bilinear form $\<\ ,\ \>_\blam^\Ldom$ on~$\LSlam[\bmu](\kx))$.

    Finally, to identify $\LDmu(\K)$, consider $\K$ as a $\Kx$-module by
    letting each $x\in\ux$ act as zero. Tensoring with~$\K$, the
    calculations above show that, for the induced maps after base
    change, $\theta\ne0$ if and only if~$\LDmu(\K)\ne0$. By
    construction, the maps~$f$ and~$h$ are both isomorphisms, so
    $\LDmu(\K)\ne0$ if and only if $g\ne0$, which is if and only
    if~$\Lzlam[\bmu]\Lpsist[{\Ltlam[\bmu]\Gtlam[\bmu]}]\ne0$.  Further,
    if $\LDmu(\K)\ne0$ then
    $q^d\LDmu(\K)\cong\im(g\circ f)
        =\Rx(\K)\Lzlam[\bmu]\Lpsist[{\Ltlam[\bmu]\Gtlam[\bmu]}]$, for
    some $d\in\Z$. Inspection of the maps, using \autoref{L:defectSum},
    shows that $d=2\defect(\blam)+\Ldeg\Gtlam$.
  \end{proof}

  \begin{Remark}
    If $\bmu\in\LKlesh$ then the simple module
    $\Rx(\K)\Lzlam[\bmu]\Lpsist[{\Ltlam[\bmu]\Gtlam[\bmu]}]$ is the socle
    of a projective cover of~$\LDmu(\K)$, up to shift. The module
    $\Rx(\K)\Lzlam[\bmu]\Lpsist[{\Ltlam[\bmu]\Gtlam[\bmu]}]$ is spanned
    by $\set{\Lzlams\Lpsist[{\Ltlam[\bmu]\Gtlam[\bmu]}]|\s\in\Std(\bmu)}$.
  \end{Remark}

  \subsection{Graded decomposition numbers}\label{SS:GradedDecomp}
  This section introduces graded decomposition matrices together with
  the key result that these matrices are unitriangular. This will be
  used in the next chapter to construct bases in the Grothendieck groups
  of~$\Rx(\Kx$), which we use to prove \autoref{MT:Simples} from the
  introduction.

  If $M$ is an $\Rx(\Kx)$-module and $D$ is an irreducible
  $\Rx(\Kx)$-module then the \textbf{graded decomposition multiplicity} of~$D$
  in~$M$ is the Laurent polynomial
  \[
    [M:D]_q = \sum_{k\in\Z} [M:q^kD]\,q^k\quad\in\N[q,q^{-1}],
  \]
  where $[M:q^kD]\in\N$ is equal to the number of composition factors
  of~$M$ that are isomorphic to $q^kD$.

  The \textbf{graded decomposition numbers} of $\Rx(\Kx)$
  are the decomposition multiplicities
  \begin{equation}\label{E:GradedDec}
   \Ldlammu =[\DSlam(\K): \DDmu(\K)]_q
   \qquad\text{and}\qquad
   \Gdlamnu =[\GSlam(\K): \GDnu(\K)]_q
  \end{equation}
  for $\blam\in\Parts$, $\bmu\in\LKlesh$ and $\bnu\in\GKlesh$.
  The \emph{graded decomposition matrices} of $\Rx(\Kx)$ are the matrices
  \[
    \LDec=\bigl(\Ldlammu\bigr)
      \And
    \GDec=\bigl(\Gdlamnu\bigr),
  \]
  The most important result that we need about the decomposition
  matrices of~$\Rx(\Kx)$ is the following.
  \notation{$\Ldlammu,\Gdlamnu$}{Graded decomposition numbers for $\Rx$}[E:GradedDec]

  \begin{Theorem}\label{T:TriangularDecomp}
    Suppose that $\K$ is a field and that $\blam\in\Parts$.
    \begin{enumerate}
      \item If $\bmu\in\LKlesh$ then $\Ldlammu[\bmu\bmu]=1$ and
      $\Ldlammu\ne0$ only if $\blam\Ledom\bmu$ and
      $\alpha_\blam=\alpha_\bmu$.
      \item If $\bnu\in\GKlesh$ then $\Gdlamnu[\bnu\bnu]=1$ and
      $\Gdlamnu\ne0$ only if $\blam\Gedom\bnu$ and
      $\alpha_\blam=\alpha_\bnu$.
    \end{enumerate}
  \end{Theorem}

  \begin{proof}
    Let $\Domin$, $\blam\in\Parts$ and $\bmu\in\DKlesh$. The theory of
    graded cellular algebras, via \autoref{T:CellularSimples}, shows
    that the decomposition matrix $\DDec$ is unitriangular when the rows
    and columns are ordered with respect to any total order that refines
    $\Dom$-dominance.  Hence, $\Ddlammu[\bmu\bmu]=1$ and $\Ddlammu\ne0$
    only if~$\blam\Domeq\bmu$.  The remaining claim follows because the
    cellular bases of \autoref{T:kcellular} give the decomposition
    $\Rx(\K)=\bigoplus_{\alpha\in Q^+_n}\Rx[\alpha](\Kx)$ of $\Rx(\K))$ into
    a direct sum of two-sided ideals.
  \end{proof}

  For $\bmu\in\LKlesh$ let $\LYmu$ be the \textbf{projective cover} of
  $\LDmu$ as an $\Rx(\K)$-module. Similarly, let~$\GYnu$ be the
  projective cover of $\GDnu$ as an $\Rx(\K)$-module, for
  $\bnu\in\GKlesh$.  \notation{$\LYmu,\GYnu$}{Projective covers of
  $\LDmu$ and $\GDnu$, respectively}

  \begin{Proposition}\label{P:PIMs}
    Let $\K$ be a field.
    \begin{enumerate}
      \item Let $\bmu\in\LKlesh$. Then $\LYmu$ has a filtration
      $\LYmu=\LYmu[\bmu,1]\supset\LYmu[\bmu,2]\supset\dots
            \supset\LYmu[\bmu,z]$
      such that there exist $\ell$-partitions
      $\blam_1,\dots,\blam_z\in\Parts$ with
      $\LYmu[\bmu,k]/\LYmu[\bmu,k+1]\cong\Ldlammu[\blam\blam_k]
            \LSlam[\blam_k]$
      and $k>l$ whenever $\blam_k\Ldom\blam_l$.
      \item Let $\bmu\in\GKlesh$. Then $\GYnu$ has a filtration
      $\GYnu=\GYnu[\bmu,1]\supset\GYnu[\bmu,2]\supset\dots
            \supset\GYnu[\bmu,z]$
      such that there exist $\ell$-partitions
      $\blam_1,\dots,\blam_z\in\Parts$ with
      $\GYnu[\bmu,k]/\GYnu[\bmu,k+1]\cong\Gdlamnu[\blam\blam_k]
        \GSlam[\blam_k]$
      and $k>l$ whenever $\blam_k\Gdom\blam_l$.
    \end{enumerate}
  \end{Proposition}

  \begin{proof}
    This comes from the general theory of (graded) cellular
    algebras; see~\cite[Theorem~3.7]{GL} or
    \cite[Lemma~2.25]{HuMathas:GradedCellular}.
  \end{proof}

  Define \emph{graded Cartan matrices} $\LCar=\bigl(\Lclammu\bigr)$ and
  $\GCar=\bigl(\Gclamnu\bigr)$ by
  \[
    \Lclammu = [\LYmu:\LDmu[\bnu]]
    \And
    \Gclamnu = [\GYnu:\GDnu].
  \]
  If $M$ is matrix let $M^T$ be its transpose.

  Standard arguments now show that the $\Kx$-cellular algebra $\Rx(\Kx)$ enjoys
  the following graded analogue of Brauer–Humphreys reciprocity; compare
  \cite[Theorem~2.17]{HuMathas:GradedCellular}.

  \begin{Corollary}\label{C:Cartan}
    Suppose that $\K$ is a field. Then $\LCar=\bigl(\LDec\bigr)^T\LDec$
    and $\GCar=\bigl(\GDec\bigr)^T\GDec$.
  \end{Corollary}

  \subsection{Adjustment matrices}\label{SS:Adjustment}
  Following \autoref{L:IntegralContentSystems}, in this section we
  assume that $\k=\Z$, so the content system $(\bc,\br)$ is defined over $\Z[\ux]$.
  By assumption, $\K$ is a field that is a $\k$-algebra, which means
  that we are assuming that $\K$ is a field. Then $\Rx(\Kx)\cong
  \Kx\otimes_{\Z[\ux]}\Rx(\Z[\ux])$ is a graded $\Kx$-cellular algebra
  by \autoref{T:kcellular}. The main result of this section compares the
  decomposition matrices of the two algebras $\Rx(\Q[\ux])$ and
  $\Rx(\Kx)$.

  Let $\A[I^n]$ be the free
  $\A$-module generated by $I^n$. The \textbf{$q$-character} of a
  finite dimensional  $\Rx(\Kx)$-module~$M$ is
  \[\ch M=\sum_{\bi\in I^n}(\gdim M_\bi)\bi \in\A[I^n],\]
  where $M_\bi=\ei M$, for $\bi\in I^n$. For example,
  $\ch\DSlam(\Kx)=\sum_{\t\in\Std(\blam)}q^{\Ddeg(\t)}\br(\t)$.

  The \emph{bar involution} is the $\Z$-linear involution on~$\A$
  given by setting $\overline{f(q)}=f(q^{-1})$, for $f(q)\in\Z$. Extend
  the bar involution to an automorphism of $\A[I^n]$ by declaring
  that $\overline{\,\bi\,}=\bi$, for $\bi\in I^n$. It is easy to see
  that $\ch(M^\circledast)=\overline{\ch M}$.
  \notation{$\ch M$}{Formal character in $\A[I^n]$, for the $\Rx$-module $M$}
  \notation{$\overline{\phantom{m}}$}{The bar involution on
  $\A+\Z[q,q^{-1}]$ given by $\overline{f(q)}=f(q^{-1})$}

  The following result is well-known and is easily proved by induction
  of the height of $\alpha\in Q^+$. This result is stated as
  \cite[Theorem~3.17]{KhovLaud:diagI}, with the reader being invited to
  repeat the proof of \cite[Theorem~3.3.1]{Klesh:book}.

  \begin{Theorem} \label{T:InjectiveCh}
    Let $\K$ be a field. Then the character map
    $\ch\map{[\Rep\Rx(\Kx)]}\A[I^n]$ is injective.
  \end{Theorem}

  The definition of the modules $\LDmu(L)$ and $\GDnu(L)$, and the
  radicals of the Specht modules, makes sense for any $\Z[\ux]$-module
  $L$. For $\bmu\in\DKlesh$ and $\bnu\in\GKlesh$ define
  \[
        \LEmu(L)=L\otimes_{\Z[\ux]}\LDmu(\Z[\ux])
        \qquad\text{and}\qquad
        \GEnu(L)=L\otimes_{\Z[\ux]}\GEnu(\Z[\ux]).
  \]
  For $\blam\in\Parts$, let
  $G^\Dom_\blam=\Bigl(\<\Dpsis,\Dpsis[\t]\>_\blam^\Dom\Bigr)_{\s,\t\in\Std(\blam)}$
  be the Gram matrix of the bilinear form \autoref{E:SpechtForm} on the
  Specht module $S^\Dom_\blam$. By considering the Smith normal form of
  $G^\Dom_\blam$, it is straightforward to prove the following. (Compare
  with \cite[Theorem~3.7.4]{Mathas:Singapore}.)

  \begin{Lemma}\label{L:SimpleCharacters}
    Let $\bmu\in\Parts$ and $\Domin$. Then
    $\DEmu(\Z[\ux])$ is a $\Z[\ux]$-free $\Rx(\Z[\ux])$-module.
    Moreover, $\DDmu(\Q)\cong\DEmu(\Q)$.
  \end{Lemma}

  The following polynomials define a map between the Grothendieck groups
  of $\Rx(\Q[\ux])$ and $\Rx(\Kx)$.

  \begin{Definition}\label{D:Adjustment}
    Let $\K$ be a field, $\Domin$ and $\bmu,\bnu\in\DKlesh$. Define
    Laurent polynomials $\Danumu$ by
    \[
      \Danumu=\sum_{q\in\Z} [\DEmu[\bnu]:q^d\DDmu(\K)]q^d \in\N[q,q^{-1}].
    \]
    The matrix $\DAdj=\bigl(\Danumu\bigr)$ is the \emph{graded
    adjustment matrix} of $\Rx(\Kx)$.
  \end{Definition}

  \begin{Theorem}\label{T:Adjustment}
    Suppose that $\K$ is a field and let $\Domin$.
    \begin{enumerate}
      \item If $\bmu,\bnu\in\DKlesh$ then $\Danumu\ne0$ only if
      $\bmu\Domeq\bnu$ and $\alpha_\bmu=\alpha_\bnu$. Moreover,
      $\overline{\Danumu}=\Danumu$.
      \item As matrices, $\DDec=\DDec(\Q)\DAdj$. That is, if $\blam\in\Parts$
      and $\bmu\in\DKlesh$ then
      \[
         \Ddlammu(\Kx)=\sum_{\bnu\in\DKlesh}\Ddlammu[\blam\bnu](\Q)\Danumu[\bnu\bmu].
      \]
    \end{enumerate}
  \end{Theorem}

  \begin{proof}Every composition factor of $\DEmu$ is a
    composition factor of $\DSlam[\bmu](\K)$, so the first statement
    in~(a) follows from \autoref{T:TriangularDecomp}. By
    \autoref{L:SimpleCharacters}, the adjustment matrix induces a
    well-defined map of Grothendieck groups
    $\DAdj\map{[\Rep\Rn(\Q[\ux)]}[\Rep\Rn(\Kx)]$ given by
    \[
       \DAdj\bigl([\DDmu[\bnu](\Q)]\bigr)=[\DEmu[\bnu](\K)]
       =\sum_{\bmu\in\DKlesh}\Danumu{}[\DDmu(\K)].
    \]
   Taking $q$-characters, $\ch\DDmu(\Q)=\sum_\bnu\Danumu\ch\DDmu[\bmu](\K)$.
   Applying $\circledast$ to both sides, the self-duality of the
   simple modules now implies that $\overline{\Danumu}=\Danumu$, which
   completes the proof of part~(a). To prove~(b), observe that
   \begin{align*}
      \sum_{\bmu\in\DKlesh}\Ddlammu[\blam\bmu](\K)\ch\DDmu(\K)
         &= \ch\DSlam(\K)
          = \ch\DSlam(\Q)\\
        & =\sum_{\bnu\in\DKlesh}\Ddlammu[\blam\bnu](\Q)\ch\DDmu[\bnu](\Q)\\
        & =\sum_{\bnu\in\DKlesh}\Ddlammu[\blam\bnu](\Q)\ch\DEmu[\bnu](\K)\\
        & =\sum_{\bnu\in\DKlesh}\Ddlammu[\blam\bnu](\Q)
            \sum_{\bmu\in\DKlesh} \Danumu\ch\DDmu(\K).
   \end{align*}
   Comparing the coefficient of $\ch\DDmu(\K)$ on both sides using
   \autoref{T:InjectiveCh} proves part~(b).
  \end{proof}

  We prove in \autoref{T:KleshchevSimples} below that
  $\DKlesh(\K)=\DKlesh(\Q)$ for any field $\K$, which implies that~$\DAdj$
  is a square unitriangular matrix.

  \subsection{A Mullineux-like involution}\label{S:Mullineux}
  \autoref{T:GradedSimples} gives two descriptions of the
  simple $\Rx(\K)$-modules $\set{q^z\LDmu(\K)}$ and $\set{q^z\GDnu(\K)}$.
  The aim of this section is set up the machinery for comparing these
  different constructions of the simple $\Rx(\K)$-modules. We start
  with a definition.

  \begin{Definition}\label{D:Mullineux}
    Let $\mull\map{\LKlesh}\GKlesh$ be the unique bijection such that
    $\LDmu(\K)\cong\GDnu[\mull(\bmu)](\K)$, for $\bmu\in\LKlesh$.
  \end{Definition}
  \notation{$\mull(\bmu)$}{Bijection $\mull\map\LKlesh\GKlesh$ such
      that $\LDmu\cong\GDnu[\mull(\bmu)]$}[D:Mullineux]

  If $\bmu\in\LKlesh$ and $\bnu\in\GKlesh$ then, by
  \autoref{T:GradedSimples}, the modules $q^z\LDmu(\K)$ and
  $q^y\GDnu(\K)$ are self-dual if and only if $z=0$ and $y=0$,
  respectively. Hence, the map $\mull$ of \autoref{D:Mullineux} is
  well-defined.

  Like the sets $\LKlesh$ and $\GKlesh$, \textit{a priori}, the map
  $\mull$ depends on~$\Lambda$, $\charge$, and the field~$\K$. We give an
  explicit description of~$\mull$ in \autoref{C:SimpleClassification}
  below, which shows that~$\mull$ is independent of~$\K$. In the next
  section we show that~$\mull$ is closely related to the sign
  isomorphism. In particular, in the special case of the symmetric
  groups, the map $\bmu\mapsto\mull(\bmu)'$ is the Mullineux map~\cite{Mull}.

  Recall from \autoref{SS:GradedDecomp} that $\DYmu$ is the projective
  cover of $\DDmu$, for $\bmu\in\DKlesh$. Hence, we have:

  \begin{Lemma}\label{L:mullY}
     Let $\bmu\in\LKlesh$. Then $\LYmu\cong\GYnu[\mull(\bmu)]$.
  \end{Lemma}

  Using $\mull$ we can give the precise relationship between the graded
  decomposition numbers $\Ldlammu$ and $\Gdlamnu$. In particular, this
  shows that the graded decomposition matrices $\LDec$ and $\GDec$
  encode equivalent information.

  Recall from the last section that the \emph{bar involution} is the
  $\Z$-linear automorphism of $\A$ given by
  $\overline{f(q)} = f(q^{-1})$.

  \begin{Proposition}\label{P:DecompComparision}
    Suppose that $\K$ is a field.
    \begin{enumerate}
      \item If $\blam\in\Parts$ and $\bmu\in\LKlesh$ then
      $\Ldlammu = q^{\defect{\blam}}\overline{\Gdlamnu[\blam\mull(\bmu)]}$.
      \item If $\blam\in\Parts$ and $\bmu\in\LKlesh$ then
      $\Ldlammu\ne0$ only if $\mull(\bmu)\Ledom\blam\Ledom\bmu$.
      \item If $\blam\in\Parts$ and $\bnu\in\GKlesh$ then
      $\Gdlamnu\ne0$ only if $\mull^{-1}(\bnu)\Gedom\blam\Gedom\bnu$.
    \end{enumerate}
  \end{Proposition}

  \begin{proof}
    Using formal characters and \autoref{P:SpechtDual}, we have
    \begin{align*}
      \sum_{\bmu\in\LKlesh}\Ldlammu{}\ch\LDmu(\K)
            &= \ch\LSlam(\K)
             = q^{\defect(\blam)}\ch\GSlam(\K)^\circledast
             = q^{\defect(\blam)}\overline{\ch\GSlam(\K)}\\
            &= q^{\defect(\blam)}
               \overline{\sum_{\bnu\in\GKlesh}\Gdlamnu{\ch\GDnu(\K)}}\\
            &= q^{\defect(\blam)}\sum_{\bnu\in\GKlesh}
               \overline{\Gdlamnu}\ch\GDnu(\K)\\
            &= q^{\defect(\blam)}\sum_{\bmu\in\LKlesh}
               \overline{\Gdlamnu[\blam\mull(\bmu)]}\ch\GDnu[\mull(\bmu)](\K)
    \end{align*}
    where the second last equality follows because $\GDnu(\K)$ is
    self-dual by \autoref{T:GradedSimples}. Part~(a) follows by
    comparing the coefficient of $\ch\LDmu(\K)$ on both sides using
    \autoref{T:InfinityCrystalGraph}.

    For~(b), if $\Ldlammu\ne0$ then $\blam\Ledom\bmu$ by
    \autoref{T:TriangularDecomp}. Moreover,
    $\Gdlamnu[\blam\mull(\bmu)]\ne0$ by~(a), so $\blam\Gedom\mull(\bmu)$
    by \autoref{T:TriangularDecomp}. The proof of~(c) is similar.
  \end{proof}

  Recalling the adjustment matrices of \autoref{SS:Adjustment}, we  obtain:

  \begin{Corollary}
    Let $\K$ be a field and $\bmu,\bnu\in\LKlesh$. Then
      $\Lanumu=\overline{\Ganumu[\mull(\bnu)\mull(\bmu)]}$.
  \end{Corollary}

  \begin{proof}
    Using \autoref{T:Adjustment}(b), twice, and
   \autoref{P:SpechtDual},
   \begin{align*}
     \sum_{\bnu,\bmu\in\LKlesh}\Ldlammu[\blam\bnu](\Q)\Lanumu\ch\LDmu(\K)
         &= \ch\LSlam(\K) = q^{\defect(\blam)}\overline{\ch\GSlam(\K)}\\
         &= q^{\defect(\blam)}\sum_{\bsig,\btau\in\GKlesh}
             \overline{\Gdlamnu[\blam\bsig](\Q)}\, \overline{\Ganumu[\bsig\btau]}
                 \ch\GDnu[\btau](\K)\\
         &= \sum_{\bmu\in\LKlesh}
         \sum_{\bnu\in\LKlesh}\Ldlammu[\blam\bnu](\Q)\,
           \overline{\Ganumu[\mull(\bnu)\mull(\bmu)]}\ch\LDmu(\K),
   \end{align*}
   where the last equality uses \autoref{P:DecompComparision}(a), where we set
   $\bsig=\mull(\bnu)$ and $\btau=\mull(\bmu)$.  The result
   follows by \autoref{T:InjectiveCh}.
  \end{proof}

  Part~(a) and \autoref{T:TriangularDecomp} imply that if
  $\bmu\in\LKlesh$ then
  $\Ldlammu[\mull(\bmu)\bmu]=q^{\defect(\bmu)}=\Gdlamnu[\bmu\mull(\bmu)]$.

  \begin{Example}
    Suppose that $\Gamma$ is a quiver of type $\Cone[2]$,
    $\Lambda=\Lambda_0$ and $n=6$. Direct
    calculation shows that the graded decomposition numbers
    of~$\Rx<\Lambda_0>[6](\Kx)$ are:
    \begin{center}
      \hspace*{12mm}
        \begin{DecompositionMatrix}
                    &    (6)&  (5,1)&  (4,2)&(4,1^2)&(3,2,1)\\
          (6)       &    1&     &     &     &     \\
          (5,1)     &    q&    1&     &     &     \\
          (4,2)     &    q&q^{2}&    1&     &     \\
          (4,1^2)   &    .&    .&    .&    1&     \\
          (3^2)     &q^{2}&    .&    q&    .&     \\
          (3,2,1)   &    .&    .&    .&    .&    1\\
          (3,1^3)   &    .&    .&    .&    q&    .\\
          (2^3)     &    q&    .&q^{2}&    .&    .\\
          (2^2,1^2) &q^{2}&    q&q^{3}&    .&    .\\
          (2,1^4)   &q^{2}&q^{3}&    .&    .&    .\\
          (1^6)     &q^{3}&    .&    .&    .&    .\\
        \end{DecompositionMatrix}
      \hspace*{20mm}
        \begin{DecompositionMatrix}
                    &    (1^6)&  (2,1^4)&(2^2,1^2)&  (3,1^3)&  (3,2,1)\\
          (1^6)     &      1&       &       &       &       \\
          (2,1^4)   &      q&      1&       &       &       \\
          (2^2,1^2) &      q&  q^{2}&      1&       &       \\
          (2^3)     &  q^{2}&      .&      q&       &       \\
          (3,1^3)   &      .&      .&      .&      1&       \\
          (3,2,1)   &      .&      .&      .&      .&      1\\
          (3^2)     &      q&      .&  q^{2}&      .&      .\\
          (4,1^2)   &      .&      .&      .&      q&      .\\
          (4,2)     &  q^{2}&      q&  q^{3}&      .&      .\\
          (5,1)     &  q^{2}&  q^{3}&      .&      .&      .\\
          (6)       &  q^{3}&      .&      .&      .&      .\\
        \end{DecompositionMatrix}
      \hspace*{12mm}
        \\[2mm]
      \hspace*{12mm}
        Graded decomposition matrix $\LDec[6](\Kx)$
      \hspace*{18mm}
        Graded decomposition matrix $\GDec[6](\Kx)$
    \end{center}
    In particular,  these decomposition matrices are independent of the
    characteristic and, in this example, the map $\mull$ sends a
    partition to its conjugate, as defined in \autoref{S:Bases}.
  \end{Example}

  \begin{Remark}\label{R:}
    If $\K$ is a field of characteristic zero, and if $\Rx(\Kx)$ is an
    algebra of type $\Aone$, then \autoref{P:DecompComparision} implies
    that if $\blam\ne\bmu$ then $0<\deg\Ddlammu\le \defect(\bmu)$, with
    equality if and only if $\blam=\mull(\bmu)$; see
    \cite[Corollary~3.6.7]{Mathas:Singapore}.  This result follows
    because in this case $\Ddlammu\in\delta_{\blam\bmu}+q\N[q]$ by
    \autoref{C:Ariki} below.  In positive characteristic, and in
    type~$\Cone$, this is no longer true. Even in type~$\Aone$,
    combining \autoref{P:DecompComparision} and
    \cite[Corollary~5]{Evseev:BadAdjustment} (and
    \cite[Example~3.7.13]{Mathas:Singapore}), shows that the degrees of
    the graded decomposition numbers are not bounded by the defect in
    positive characteristic.
  \end{Remark}

  \subsection{The sign isomorphism}\label{S:SignAutomorphism}
  A \textit{sign isomorphism} of the KLR algebras of type~$\Aone$ was
  introduced in \cite[(3.14)]{KMR:UniversalSpecht}. This section
  generalises this map to include the quivers of type $\Cone$ and it
  describes its effect on the Specht modules and simple modules of
  $\Rx$. In type~$\Aone$, many of the results in this section are graded
  analogues of results in \cite[\S3]{HuMathas:FayersConjecture}.

  \begin{Definition}\label{D:sign}
    The \emph{sign automorphism} of~$\Gamma$ is quiver automorphism
    $\sgn\map\Gamma\Gamma$ given by
    \[
        \sgn(i) = \begin{cases*}
          e-i\pmod e& for type $\Aone$,\\
          e-1-i     & for type $\Cone$,
        \end{cases*}
    \]
    for $i\in I$.
    If $\bi=(i_1,\dots,i_n)\in I^n$ let
    $\bi^\sgn=\bigl(\sgn(i_1),\dots,\sgn(i_n)\bigr)\in I^n$.
  \end{Definition}
  \notation{$\eps$}{Sign automorphism of $\Gamma$ and associated maps
     on $\Rx$, $\Uq$, $\ldots$}[D:sign]

  It is straightforward to check that $\cij=\cij[\sgn(i)\sgn(j)]$, for
  all $i,j\in I$, showing that $\sgn$ is a quiver automorphism
  of~$\Gamma$.  The sign automorphism of~$\Gamma$ induces automorphisms
  of the lattices~$P^+$ and~$Q^+$, given by
  $\Lambda\mapsto \Lambda^{\sgn}$ and $\alpha\mapsto \alpha^{\sgn}$,
  that are uniquely determined by
  \[
  \kill{\alpha^\vee_i}{\Lambda^{\sgn}}=\kill{\alpha^\vee_{\sgn(i)}}{\Lambda}
     \quad\text{and}\quad
     \kill{\alpha^\vee_j}{\alpha^{\sgn}}=\kill{\alpha^\vee_{\sgn(j)}}{\alpha},
     \qquad\text{for }i,j\in I,
  \]
  respectively.

  By definition, the algebra $\Rx[\alpha](\kx)$ depends on the families
  polynomials $\Wbx$ and $\Qbx$ from \autoref{N:WbxQbx}. Define
  polynomials $\Wbxs=\bigl(W^{\ux,\sgn}_i(u)\bigr)_{i\in I}$ and
  $\Qbxs=\bigl(Q^{\ux,\sgn}_{ij}(u,v)\bigr)_{i,j\in I}$ by
  \begin{equation}\label{E:SignedPolys}
    W^{\ux,^\sgn}_i(u)=W^\ux_{\eps(i)}(-u)\And*
    Q^{\ux,\sgn}_{ij}(u,v)=Q^\ux_{\eps(i)\eps(j)}(-u,-v),
    \qquad\text{for $i,j\in I$}.
  \end{equation}
  Set $\sRx=\Rx<\Lambda^\sgn>[\alpha^\sgn](\Qb^\sgn,\Wb^\sgn)$. If
  $(\bc,\br)$ is a (graded) content system for $\Rx$ then
  $(-\bc,\sgn\circ\br)$ is a graded content system with values in~$\kx$ for~$\sRx$.

  If $\charge=(\kappa_1,\dots,\kappa_\ell)$ is an $\ell$-charge for
  $\Lambda$ then $\charge^{\sgn}=(-\kappa_\ell,\dots,-\kappa_1)$ is the
  corresponding \emph{signed charge}.

  \begin{Proposition}\label{P:SignedKLR}
    Let $\Lambda\in P^+$ and $\alpha\in Q^+$. Then there is a unique
    graded algebra isomorphism $\sgn\map{\Rx[\alpha](\kx)}{\sRx(\kx)}$
    such that
    \[
       \sgn\bigl(\ei) = \ei[\bi^\sgn],\qquad
       \sgn(\psi_k) = -\psi_k\And
       \sgn(y_m) = -y_m,
    \]
    for $\bi\in I^n$, $1\le k<n$ and $1\le m\le n$.
  \end{Proposition}

  \begin{proof}
    Checking the relations in \autoref{D:KLR} shows that there
    is a well-defined surjective homomorphism isomorphism
    $\sgn\map{\Rx[\alpha](\kx)}{\sRx(\kx)}$ of graded algebras.
    By symmetry, there is also a well-defined surjective graded algebra
    homomorphism $\sgn'\map{\sRx}{\Rx[\alpha]}$. By definition,
    $\sgn\circ\sgn'$ and $\sgn'\circ\sgn$ are identity maps, so the
    result follows. (Hereafter, we abuse notation and use $\sgn$ for
    both of these isomorphisms.)
  \end{proof}

  The isomorphism $\sgn\map{\Rx[\alpha](\kx)}{\sRx(\kx)}$ of
  \autoref{P:SignedKLR} is the \emph{sign isomorphism}. This generalises
  the sign automorphism of the group algebra of the symmetric group,
  which corresponds to the special case when
  $\Lambda=\Lambda_0$ in type~$\Aone$ for $\Rn(\K)$, when~$\K$ is a field.
  By base change, \autoref{P:SignedKLR} induces isomorphisms
  $\Rx[\alpha](L)\bijection\sRx(L)$ for any $\kx$-algebra $L$.
  Setting $\ux=0$ we obtain an analogous isomorphism
  $\sgn\map{\Rn[\alpha](\k)}{\sRn(\k)}$.

  If $M$ is an $\sRx$-module let $M^\sgn$ be the $\sgn$-\emph{twisted}
  $\Rx[\alpha](\kx)$-module that is equal to $M$ as a $\kx$-module and
  where the $\Rx[\alpha]$-action is twisted by $\sgn$, so that $a\cdot
  m=\sgn(a)m$, for $a\in\Rx[\alpha](\kx)$ and $m\in M$. By
  \autoref{P:SignedKLR}, this induces an equivalence of categories
  $\Rep\sRx(\kx)\to\Rep\Rx[\alpha](\kx)$ given by $M\mapsto M^\sgn$. In
  the special case of the symmetric groups, this is the equivalence of
  categories induced by tensoring with the sign representation.  This
  follows because if $\K$ is a field then there is an isomorphism
  $\Rn<\Lambda_0>(\K)\cong\K\Sym_n$ by the main result of
  \cite{BK:GradedKL} and in this case $\sgn$ induces an auto-equivalence
  of $\Rep\Rn<\Lambda_0>(\K)$. More generally, $\sgn$ induces an
  auto-equivalence of $\Rep\Rx(\Kx)$ whenever $\Lambda=\Lambda^\sgn$.

  Most of our notation so far implicitly
  depends on $\Lambda$ and sometimes $\alpha$ and $\charge$. To avoid
  ambiguity, we decorate our notation with $\sgn$ whenever it is applied to
  objects associated with the algebra $\sRx(\kx)$, and we
  continue to use our existing notation for the algebras
  $\Rx[\alpha](\kx)$. In particular, $\DSlam*$ and $\DDmu*$ are the
  graded Specht and simple $\sRx(\kx)$--modules .
  The main results of this section explore the twisted modules $(\DSlam*)^\sgn$
  and $(\DDmu*)^\sgn$, for $\blam\in\Parts[\alpha]$ and
  $\bmu\in\DKlesh[\alpha]$.

  We need ``sign adapted'' combinatorics for the KLR algebras. As
  suggested by the terminology, in the representation theory of the
  symmetric groups this is given by conjugate partitions and tableaux,
  as defined in \autoref{S:Bases}.

  Extending the definition of the conjugate of an $L$-partition from
  \autoref{S:Bases}, the \emph{conjugate} of the node $A=(m,r,c)$ is the
  node $A'=(\ell-m+1,c,r)$. In particular, if $\blam\in\Parts$ then its
  conjugate is $\blam'=\set{A'|A\in\blam}$ and the conjugate of
  $\t\in\Std(\blam)$ is the tableau $\t'\in\Std(\blam')$ given by
  $\t'(A)=\t(A')$, for $A\in\blam'$. If $A$ is a node then $(A')'=A$, so
  conjugation is an involution on the sets of $\ell$-partitions and
  standard tableaux.

  A straightforward walk through the definitions reveals that the
  following identities hold.

  \begin{Lemma}\label{L:SignedCombinatorics}
    Let $\blam\in\Parts[\alpha]$, for $\alpha\in Q^+$. If
    $A\in\blam'$ then
    \[
      \Ldilam*(\blam')=\Gdilam[A'],\qquad
      \Gdilam*(\blam')=\Ldilam[A'],\qquad
      d^\sgn_i(\blam')=d_{\sgn(i)}(\blam) \And
      \defect^\sgn(\blam')=\defect(\blam).
    \]
    Moreover, if $\s\in\Std(\blam)$ then
    $\br(\s')=\br(\s)^\sgn$,
    $\deg^{\Ldom}_\sgn(\s')=\Gdeg(\s)$ and
    $\deg^{\Gdom}_\sgn(\s')=\Ldeg(\s)$.
  \end{Lemma}

  \begin{Proposition}\label{P:SignedBasis}
    Suppose that $\s,\t\in\Std(\blam)$, for $\blam\in\Parts[\alpha]$.
    Then
    \[
        \sgn\bigl(\Lpsist\bigr)=\pm\psi^{\Gdom\sgn}_{\s'\t'}
        \And
        \sgn\bigl(\Gpsist\bigr)=\pm\psi^{\Ldom\sgn}_{\s'\t'}.
    \]
  \end{Proposition}

  \begin{proof}
    This is a straightforward exercise in the definitions. Observe that
    $\Ltlam=\t^{\Gdom \sgn}_{\blam'}$ and
    $\Gtlam=\t^{\Ldom \sgn}_{\blam'}$. Consequently,
    if $\u\in\Std(\blam)$ then $d_\u^{\Ldom\sgn}=\Gdt[\u']$ and
    $d_\u^{\Gdom\sgn}=\Ldt[\u']$. By \autoref{D:psist} and
    \autoref{E:SignedPolys}, $y^{\Ldom\sgn}_{\blam'}=\pm\Gylam$
    and $y^{\Gdom\sgn}_{\blam'}=\pm\Lylam$, implying the result.
  \end{proof}

  For the Specht modules of the symmetric groups,
  James~\cite[Theorem~8.15]{James} proved the famous result that
  $S^{\lambda'}\cong\text{sgn}\otimes S^\lambda$, where $S^\lambda$ is a
  Specht module for the symmetric group $\Sym_n$ and $\text{sgn}$ is its
  sign representation. This next result generalises James' theorem.

  \begin{Corollary}\label{C:SignedSpecht}
    Suppose that $\blam\in\Parts[\alpha]$, for $\alpha\in Q^+$. Then
    $\LSlam\cong\bigl(\GSlam*[\blam']\bigr)^\sgn$ and
    $\GSlam\cong\bigl(\LSlam*[\blam']\bigr)^\sgn$.
  \end{Corollary}

  \begin{proof}
    By \autoref{P:SignedBasis},
    $\LeRlam<\alpha>\cong\bigl(\GeRlam*<\alpha>[\blam']\bigr)^\sgn$ and
    $\GeRlam<\alpha>\cong\bigl(\LeRlam*<\alpha>[\blam']\bigr)^\sgn$,
    implying the result.
  \end{proof}

  This allows us to identify the twisted simple $\sRx$-modules as
  $\Rx[\alpha]$-modules. The result says that these modules are
  isomorphic once you conjugate the $\ell$-partitions and interchange
  the $\Ldom$-simple modules and the $\Gdom$-simple modules. The simple
  modules are defined over the field $\K$.

  \begin{Corollary}\label{C:SignedSimples}
    Let $\bmu\in\LKlesh[\alpha]$ and $\bnu\in\GKlesh[\alpha]$. Then
    $\LDmu\cong\bigl(\GDnu*[\bmu']\bigr)^\sgn$
    and $\GDnu\cong\bigl(\LDmu*[\bnu']\bigr)^\sgn$.
  \end{Corollary}

  \begin{proof}
     Let $\head(M)$ be the head of $M$, which is its
     maximal semisimple quotient. Then, using \autoref{C:SignedSpecht},
     $\LDmu\cong \head(\LSlam[\bmu])
             \cong \bigl(\head\GSlam*[\bmu']\bigr)^\sgn
             \cong \bigl(\GDnu*[\bmu']\bigr)^\sgn$.
     The second isomorphism is proved in exactly the same way.
  \end{proof}

  Recall from \autoref{D:Mullineux} that $\mull\map\LKlesh\GKlesh$ is
  the map given by $\LDmu\cong\GDnu[\mull(\bmu)]$, for $\bmu\in\LKlesh$.
  In the special case of the symmetric groups the next result says that
  the map $\bmu\mapsto\mull(\bmu)'$ is the Mullineux map.

  \begin{Corollary}\label{C:MullinuxMap}
    Let $\bmu\in\LKlesh[\alpha]$. Then
    \[
      \LDmu\cong\bigl(\LDmu*[\mull(\bmu)']\bigr)^\sgn,\quad
      \GDnu[\mull(\bmu)]\cong\bigl(\GDnu*[\bmu']\bigr)^\sgn,\quad
      \LYmu\cong\bigl(Y^{\Ldom\sgn}_{\mull(\bmu)'}\bigr)^\sgn
      \quad\text{and}\quad
      \GYnu[\mull(\bmu)]\cong\bigl(Y^{\Gdom\sgn}_{\bmu'}\bigr)^\sgn,\quad
    \]
    In particular,
    $\set{D^{\Ldom\sgn}_{\bmu}|\bmu'\in\GKlesh}$ and
    $\set{D^{\Gdom\sgn}_{\bnu}|\bnu'\in\LKlesh}$ are both complete sets
    of pairwise non-isomorphic self-dual irreducible graded
    $\sRx$-modules.
  \end{Corollary}

  \begin{proof}
    Using \autoref{C:SignedSimples},
    $\LDmu\cong\GDnu[\mull(\bmu)]\cong\bigl(\LDmu*[\mull(\bmu)']\bigr)^\sgn$. The
    proof of the second isomorphism is similar and the remaining
    isomorphisms follow by the uniqueness of projective covers.
  \end{proof}

  If $M$ is an $\Rx$-module then its \emph{socle}, $\soc M$, is its
  maximal semisimple submodule. Dually, the \emph{head} of $M$, $\head
  M$, is the maximal semisimple subquotient of $M$.
  \notation{$\soc M$}{The socle of $M$}
  \notation{$\head M$}{The head of $M$}

  \begin{Corollary}\label{C:SpechtSocle}
    Let $\bmu\in\LKlesh[\alpha]$ and $\bnu\in\GKlesh[\alpha]$. Then
    \[
      \soc\LSlam[\bmu]\cong q^{\defect(\bmu)}\bigl(\LDmu*[\mull(\bmu)']\bigr)^\sgn
      \And
      \soc\GSlam[\bnu]\cong q^{\defect(\bnu)}\bigl(\GDnu*[\mull(\bnu)']\bigr)^\sgn.
    \]
  \end{Corollary}

  \begin{proof}
    Using \autoref{P:SpechtDual},
    \[
      \soc\LSlam[\bmu]
      \cong\soc\bigl(q^{\defect(\bmu)}{\GSlam[\bmu]}^\circledast\bigr)
         \cong q^{\defect(\bmu)}\head\bigl(\GSlam[\bmu]\bigr)^\circledast
         \cong q^{\defect(\bmu)}\GDnu[\bmu]
         \cong q^{\defect(\bmu)}\bigl(\LDmu*[\mull(\bmu)']\bigr)^\sgn.
    \]
    where the last isomorphism follows from \autoref{C:MullinuxMap}.
    The second isomorphism is similar.
  \end{proof}

  The last result in this section can be viewed as a generalisation of
  \cite[Theorem~7.2]{LLT}.

  \begin{Corollary}\label{C:SignedDecomposition}
    Let $\blam\in\Parts$ and $\bmu\in\LKlesh$ and $\bnu\in\GKlesh$. Then
    \[
      [S^{\Ldom\sgn}_\blam: D^{\Ldom\sgn}_{\mull(\bmu)'}]_q
           =q^{\defect{\blam'}}\overline{[\LSlam[\blam']: \LDmu]_q}
           \qquad\text{and}\qquad
      [S^{\Gdom\sgn}_\blam: D^{\Gdom\sgn}_{\bmu'}]_q
           =q^{\defect{\blam'}}\overline{[\GSlam[\blam']: \GDnu[\mull(\bmu)]]_q}.
    \]
  \end{Corollary}

  \begin{proof}
    We prove only the second identity. Using \autoref{C:SignedSpecht}
    and \autoref{C:SignedSimples},
    \[
      [S^{\Gdom\sgn}_\blam: D^{\Gdom\sgn}_{\bmu'}]_q
         = [\bigl(S^{\Gdom\sgn}_\blam\bigr)^\sgn: \bigl(D^{\Gdom\sgn}_{\bmu'}\bigr)^\sgn]_q
         = [\LSlam[\blam']:\LDmu]_q
         = q^{\defect\blam'}\overline{[\GSlam[\blam']: \GDnu[\mull(\bmu)]]_q}
    \]
    where the last equality follows from
    \autoref{P:DecompComparision}(a) and
    \autoref{L:SignedCombinatorics}.
  \end{proof}

\section{Categorification}\label{S:Categorification}

  This chapter brings together all of our previous work to prove that
  the algebras $\Rx(\K[x])$ categorify the integrable highest weight
  modules of the corresponding Kac-Moody algebras, which is
  \autoref{MT:Categorification} from the introduction. As applications,
  we classify the simple $\Rx(\K[x])$-modules (\autoref{MT:Simples}), and
  prove their modular branching rules (\autoref{MT:ModularBranching}).
  To do this we first use the algebras $\Rx(\k[x^\pm])$ to prove
  the branching rules for the graded Specht modules of~$\Rx(\k[x])$, which
  leads almost directly to our categorification theorem. We then use the
  representation theory of $\Rx(\K[x])$ to describe the canonical bases of
  the highest weight modules, which gives us a way of studying the
  simple modules of~$\Rx(\K[x])$.

  Throughout this chapter we continue to assume that $(\bc,\br)$ is a
  (graded) content system with values in $\k[x]$ for a cyclotomic KLR
  algebra $\Rx(\k[x])$, and $\K$ is a field that is a $\k$-algebra so
  that $\Rx(\K[x])$ is a graded $\K[x]$-cellular algebra by
  \autoref{C:RxCellular}. In particular, as discussed in the last
  chapter, \autoref{C:RnCellular} implies that the results in this
  chapter apply to the (standard) cyclotomic KLR algebras of types $\Aone$,
  $A_\infty$, $\Cone$ and $C_\infty$.

  \subsection{Branching rules}\label{SS:Branching}
  This section proves analogues of the classical branching rules of the
  symmetric groups for the $\Rx$-Specht modules. That is, we describe
  the modules obtained by inducing and restricting the graded Specht
  modules. The strategy is to first prove the branching rules for the
  semisimple algebras $\Rx(\K[x^\pm])$ and then to use this result to prove the
  branching rules for $\Rx(\k[x])$, after which the branching rules for
  $\Rx$ and $\Rn$ follow by specialisation. In the next section we use these
  results to show that~$\Rx$ categorifies the integral highest weight
  modules of~$\Uq$.

  \notation{$\Rep_\K\Rx(\K[x])$}{Category of graded $\Rx$-modules, which are finite dimensional over~$\K$}
  \notation{$\Proj_\K\Rx(\K[x])$}{Full subcategory of $\Rep_\K\Rx(\K[x])$ of projective modules}
  \notation{$\iRes$}{The $i$-restriction functor $\Rep_\K\Rx[\alpha+\alpha_i]to\Rep_\K\Rx[\alpha]$}
  \notation{$\iInd$}{The $i$-induction functor $\Rep_\K\Rx[\alpha]\to\Rep_\K\Rx[\alpha+\alpha_i]$}

  Before we can begin, we need to define the categories that we are
  going to work in.  Fix $\alpha\in Q^+_n$.  Let $\Rep\Rx[\alpha](\k[x])$
  be the category of finitely generated graded
  $\Rx[\alpha](\k[x])$-modules, and similarly define
  $\Rep\Rx[\alpha](\K[x])$. Let $\Rep_\K\Rx[\alpha](\K[x])$ be the full
  subcategory of $\Rep\Rx[\alpha](\K[x])$ consisting of graded
  $\Rx[\alpha](\K[x])$-modules that are \textit{finite dimensional} as
  $\K$-vector spaces.  Let $\Proj\Rx[\alpha](\k[x])$ and
  $\Proj_\K\Rx[\alpha](\K[x])$ be the additive subcategories of graded
  projective modules in $\Rep\Rx[\alpha](\k[x])$ and
  $\Rep_\K\Rx[\alpha](\K[x])$, respectively.  Similarly, let
  $\Rep\Rn[\alpha](\k)$ and $\Proj\Rn[\alpha](\k)$ $\Rep\Rx[\alpha](\K)$
  and $\Proj\Rx[\alpha](\K)$ be the corresponding subcategories of
  graded $\Rn[\alpha](\k)$-module.  and graded
  $\Rx[\alpha](\K)$-modules, respectively.

  Set $\Rep\Rx(\k[x])=\bigoplus_{\alpha\in Q^+_n}\Rep\Rx[\alpha](\k[x])$, and
  similarly for the other categories defined above.

  Ultimately, we are most interested in the category $\Rep_\K\Rx(\K[x])$,
  which is quite different to $\Rep\Rx(\K[x])$. For example, the graded
  Specht module $\DSlam(\K[x])$ does not belong to $\Rep_\K\Rx(\K[x])$ but
  it does belong to $\Rep\Rx(\K[x])$. The categories $\Rep_\K\Rx(\K[x])$ and
  $\Rep\Rn(\K)$ are also not equivalent but they have isomorphic
  Grothendieck groups by the remarks after \autoref{T:GradedSimples}.

  Let $i\in I$ and $\alpha\in Q^+_n$. Set
  $\ei[\alpha,i] = \sum_{\bj\in I^\alpha}\ei[\bj i]$.
  Define \emph{$i$-restriction} and \emph{$i$-induction} functors:
  \begin{align*}
    \iRes \map{\Rep{\Rx[\alpha+\alpha_i]\k[x])}}{\Rep{\Rx[\alpha](\k[x])}};
         M\mapsto \ei[\alpha,i]\Rx[\alpha+\alpha_i](\k[x])
              \otimes_{\Rx[\alpha+\alpha_i]}M,\\
    \iInd \map{\Rep{\Rx[\alpha](\k[x])}}{\Rep{\Rx[\alpha+\alpha_i](\k[x])}};
     M\mapsto \Rx[\alpha+\alpha_i](\k[x])\ei[\alpha,i]\otimes_{\Rx[\alpha](\k[x])}M.
  \end{align*}
  Abusing notation, we also write $\iRes\map{\Rep{\Rx[n+1]}}{\Rep{\Rx}}$
  and $\iInd\map{\Rep{\Rx}}{\Rep{\Rx[n+1]}}$ for the corresponding
  induced functors on these module categories.  These functors can be
  defined as the direct sum of the functors defined above or they can be
  defined directly by replacing each occurrence of $\ei[\alpha,i]$ in
  the definitions above with $\ei[n,i]=\sum_{\alpha\in Q^+_n}\ei[\alpha,i]$.
  We further abuse notation and use $\iRes$ and $\iInd$ for the
  induced functors on all of the categories defined above.

  \begin{Proposition}\label{P:exactness}
    Let $i\in I$. There is a (non-unital) embedding of graded algebras
    $\iota_{n,i}\colon\Rx\hookrightarrow\Rx[n+1]$ such that
    \[
            \ei[\bj]\mapsto \ei[\bj i],\quad
            \psi_r\ei[\bj]\mapsto \psi_r\ei[\bj i]\And*
            y_m\ei[\bj]\mapsto y_m\ei[\bj i],
    \]
    for $\bj\in I^n$, $1\le r<n$ and $1\le m\le n$.
    Moreover, if $M\in\Rep{\Rx[n+1]}$ then $\iRes(M)=\ei[n,i]M$ and
    if~$N\in\Rep{\Rx}$ then $\iInd(N)=\Rx[n+1]\ei[n,i]N$, so
    $\iRes$ and $\iInd$ are exact functors.
  \end{Proposition}

  \begin{proof}
    The relations \autoref{D:KLR}, together with
    \autoref{T:kcellular}, imply that there is a unique non-unital
    algebra embedding $\iota_{\alpha,\alpha_i}\colon\Rx[\alpha]
    \hookrightarrow\Rx[\alpha+\alpha]$ such that
     \[
            \ei[\bj]\mapsto \ei[\bj i],\quad
            \psi_r\ei[\bj]\mapsto \psi_r\ei[\bj i]\And*
            y_m\ei[\bj]\mapsto y_m\ei[\bj i],
     \]
     for $\bj\in I^\alpha$, $1\le r<n$ and $1\le m\le n$. In particular,
     $\iRes$ is an exact functor.
     Kashiwara \cite[Corollary~3.3]{Kashiwara:KLRBiadjointness} proves
     that $\iInd$ is exact.
  \end{proof}

  The aim of this section is to describe the modules $\iRes\LSlam$ and
  $\iInd\GSlam$, for $\blam\in\Parts$.  We start with the easier case of
  restriction, following \cite{Mathas:SpechtRestriction}. If $\Domin$
  then \autoref{P:psiTriangular}, $\DSlam(\K[x^\pm])$ has an $\Dfst[]$-basis
  and a~$\Dpsist[]$-basis, for which the transition matrices are
  unitriangular.  Note that $\LSlam(\K[x^\pm])\cong\GSlam(\K[x^\pm])$ in view of
  \autoref{C:IsomorphicSpechts} and \autoref{P:SeminormalForm}.

  If $\t\in\Std(\blam)$ let $\t_{\downarrow}=\t_{\downarrow(n-1)}$. Let
  $\K'$ be the field of fractions of~$\k$.

  \begin{Lemma}\label{L:SpechtRrestriction}
    Suppose that $\blam\in\Parts[\alpha+\alpha_i]$. Then, as
    $\Rx[\alpha](\K'[x^\pm])$-modules,
    \[
       \iRes\bigl(\LSlam(\K'[x^\pm])\bigr)\cong
           \bigoplus_{\mathclap{B\in\Rem_i(\blam)}}\LSlam[\blam{-}B](\K'[x^\pm]).
           \quad\text{and}\quad
       \iRes\bigl(\GSlam(\K'[x^\pm])\bigr)\cong
           \bigoplus_{\mathclap{B\in\Rem_i(\blam)}}\GSlam[\blam{-}B](\K'[x^\pm]).
    \]
  \end{Lemma}

  \begin{proof}
    This follows from \autoref{L:SpechtDecomposition} but to understand
    how the Specht modules restrict over $\k[x]$ we need to describe the
    isomorphism explicitly.  Let $\Domin$. By \autoref{T:Kcellular},
    $\iRes\bigl(\DSlam(\K'[x^\pm])$ has basis
    $\set{\Dfs|\s\in\Std(\blam)\text{ and }\br_n(\t)=i}$, which is in
    bijection with the set of tableaux $\bigcup_B\Std(\blam{-}B)$ where
    $B\in\Rem_i(\blam)$.  Define a $\K'[x^\pm]$-linear map
    \begin{equation}\label{E:Restriction}
      \theta\map{\iRes\bigl(\DSlam(\K'[x^\pm])\bigr)}\bigoplus_{B\in\Rem_i(\blam)}
          \DSlam[\blam{-}B](\K'[x^\pm]);\quad \Dfs\mapsto\Dfs[\s_{\downarrow}],
          \qquad\text{ for }\s\in\Std(\blam).
    \end{equation}
    By \autoref{P:fstaction} this is an isomorphism of $\Rx(\K'[x^\pm])$-modules.
  \end{proof}

  There are no grading shifts in \autoref{L:SpechtRrestriction} because
  $\K'[x^\pm]\cong q^d\K'[x^\pm]$ as a $\Z$-graded ring, for $d\in \Z$. The analogue
  of this result over $\k[x]$ requires grading shifts that are given by
  the integers~$\Ldilam$ and~$\Gdilam$ from \autoref{D:defect}.

  \begin{Proposition}\label{P:iRestriction}
    Suppose that $\blam\in\Parts[\alpha+\alpha_i]$ and let
    $A_1>\dots>A_z$ be the removable $i$-nodes of $\blam$. Then
    there exist $\Rx[\alpha](\k[x])$-module filtrations
    \begin{align*}
        \iRes\bigl(\LSlam(\k[x])\bigr)=\LSlam[\blam,z](\k[x])
           \supset\LSlam[\blam,z-1](\k[x])
           \supset\dots\supset\LSlam[\blam,2](\k[x])
           \supset\LSlam[\blam,1](\k[x])\supset0\\
        \iRes\bigl(\GSlam(\k[x])\bigr)=\GSlam[\blam,1](\k[x])
           \supset\LSlam[\blam,2](\k[x])
           \supset\dots\supset \GSlam[\blam,z-1](\k[x])
           \supset\GSlam[\blam,z](\k[x])\supset0
    \end{align*}
    with
    $\LSlam[\blam,k](\k[x])/\LSlam[\blam,k-1](\k[x])
         \cong q^{\Ldilam[A_k]}\LSlam[\blam-A_k](\k[x])$
    and
    $\GSlam[\blam,k](\k[x])/\GSlam[\blam,k+1](\k[x])
         \cong q^{\Gdilam[A_k]}\GSlam[\blam-A_k](\k[x])$,
    for $1\le k\le z$.
  \end{Proposition}

  \begin{proof}
    Consider $\iRes\bigl(\LSlam\bigr)$.  As in
    \autoref{L:SpechtRrestriction}, the module
    $\iRes\bigl(\DSlam(\k[x])\bigr)$ has basis
    \[
        \set{\Lpsis|\s\in\Std(\blam)\text{ and }\br_n(\s)=i}
           =\bigcup_{k=1}^z\set{\Lpsis|\s_\downarrow\in\Std(\blam-A_k)}.
    \]
    For $1\le k\le z$, define
    $\LSlam[\blam,k](\k[x]) =
       \<\Lpsis\mid\s_{\downarrow}\in\Std(\blam-A_s)
         \text{ for }1\le s\le k\>$.
    Then
    $\iRes\bigl(\LSlam(\k[x])\bigr)=\LSlam[\blam,z](\k[x])\supset\dots
        \supset\LSlam[\blam,1](\k[x])\supset0$
    is an $\Rx[\alpha](\k[x])$-module filtration of
    $\iRes\bigl(\LSlam(\k[x])\bigr)$ by \autoref{P:integral} and
    \autoref{C:ypsi}. In view of
    \autoref{P:psiTriangular}, it follows easily by induction on
    dominance that the $\Rx(\K[x^\pm])$-module isomorphism
    $\theta$ defined in \autoref{E:Restriction} induces
    $\Rx(\k[x])$-module isomorphisms
    \[
        \theta_k\map{\LSlam[\blam,k](\k[x])/\LSlam[\blam,k-1](\k[x])}
                   q^{\Ldilam}\LSlam[\blam-A_k](\k[x]);
                \Lpsis\mapsto\Lpsis[\s_{\downarrow}].
    \]
    This completes the proof for $\iRes\bigl(\LSlam(\k[x])\bigr)$. The
    filtration of $\iRes\bigl(\GSlam(\k[x])\bigr)$ can be constructed in
    exactly the same way.  Alternatively, it can be deduced from the
    filtration of $\iRes\bigl(\LSlam(\k[x])\bigr)$ using
    \autoref{P:SpechtDual} and \autoref{L:defectSum}.
  \end{proof}

  By base change, we obtain the corresponding result over any ring~$L$
  that is a $\k[x]$-module.

  \begin{Corollary}\label{C:SpechtRestriction}
    Suppose that $L$ is a $\k[x]$-module,
    $\blam\in\Parts[\alpha+\alpha_i]$ and let $A_1>\dots>A_z$ be the
    removable $i$-nodes of $\blam$. Then there exist
    $\Rx[\alpha](L)$-module filtrations
    \begin{align*}
        \iRes\bigl(\LSlam(L)\bigr)=\LSlam[\blam,z](L)\supset\LSlam[\blam,x-1](L)
           \supset\dots\supset\LSlam[\blam,2](L)
           \supset\LSlam[\blam,1](L)\supset0\\
        \iRes\bigl(\GSlam(L)\bigr)=\GSlam[\blam,1](L)\supset\LSlam[\blam,2](L)
           \supset\dots\supset \GSlam[\blam,z-1](L)
           \supset\GSlam[\blam,z](L)\supset0
    \end{align*}
    with
    $\LSlam[\blam,k](L)/\LSlam[\blam,k-1](L)
         \cong q^{\Ldilam[A_k]}\LSlam[\blam-A_k](L)$
    and
    $\GSlam[\blam,k](L)/\GSlam[\blam,k+1](L)
         \cong q^{\Gdilam[A_k]}\GSlam[\blam-A_k](L)$,
    for $1\le k\le z$.
  \end{Corollary}

  In view of \autoref{P:Specialisation}, a special case of
  \autoref{C:SpechtRestriction} gives Specht filtrations of the Specht
  modules~$\DSlam(L)$ for the standard cyclotomic algebras $\Rn(L)$, for
  $\Domin$.  In type~$\Aone$ this recovers
  \cite[Theorem~4.11]{BKW:GradedSpecht} when $L$ is a field and
  \cite[\S5]{Mathas:SpechtRestriction} for general~$L$.

  Next we consider the induced modules $\iInd\bigl(\LSlam\bigr)$ and
  $\iInd\bigl(\GSlam\bigr)$ using ideas that go back to
  Ryom-Hansen~\cite{RyomHansen:GradedTranslation}.  First, some
  notation. Let $\Domin$ and suppose $A\in\Add_i(\blam)$. Let
  $\Dtlam[\blam,A]\in\Std(\blam{+}A)$ be the unique standard tableau
  such that $(\Dtlam[\blam,A])_{\downarrow}=\Dtlam$.  Note that this
  forces $\Dtlam[\blam,A](A)=n+1$.

  The following example is suggestive of how the graded induction
  formulas are proved for the Specht modules are proved over~$\k[x]$.

  \begin{Example}
    Let $\blam=(3^2,2)$ and consider the quivers $\Aone[2]$
    and $\Cone[2]$. The residues in $\blam$ are:
    \begin{center}
      \begin{tikzpicture}
        \pic (A) at (0,0) {atableaux tableau={{0,1,2},{2,0,1},{1,2}}};
        \node at ([xshift=-30]A-2-1) {$\Aone[2]$};
        \node[blue] at (2.25,-0.25){$A_3$};
        \node[blue] at (1.75,-1.25){$A_2$};
        \node[blue] at (0.75,-1.75){$A_1$};
      \end{tikzpicture}
      \qquad
      \begin{tikzpicture}
        \pic (C) at (0,0) {atableaux tableau={{0,1,2},{1,0,1},{2,1}}};
        \node at ([xshift=-30]C-2-1) {$\Cone[2]$};
        \node[blue] at (2.25,-0.25){$A_3$};
        \node[blue] at (1.75,-1.25){$A_2$};
        \node[blue] at (0.75,-1.75){$A_1$};
      \end{tikzpicture}
    \end{center}
    In type $\Aone[2]$, take $i=0$ so that
    $\Add_i(\blam)=\set{A_1,A_2,A_3}$ where, as above, $A_1=(4,1)$, $A_2=(3,2)$
    and $A_3=(1,4)$.  The standard tableaux~$\Ltlam[\blam,A_r]$
    and~$\Gtlam[\blam,A_r]$ are:
    \begin{align*}
      \Gtlam[\blam,A_1] &= \Tableau{{1,2,3},{4,5,6},{7,8},{9}}&
      \Gtlam[\blam,A_2] &= \Tableau{{1,2,3},{4,5,6},{7,8,9}}  &
      \Gtlam[\blam,A_3] &= \Tableau{{1,2,3,9},{4,5,6},{7,8}}  \\
      \Ltlam[\blam,A_1] &= \Tableau{{1,4,7},{2,5,8},{3,6},{9}}&
      \Ltlam[\blam,A_2] &= \Tableau{{1,4,7},{2,5,8},{3,6,9}}  &
      \Ltlam[\blam,A_3] &= \Tableau{{1,4,7,9},{2,5,8},{3,6}}
    \end{align*}
    In type $\Cone[2]$, take $i=1$ so that
    $\Add_i(\blam)=\set{A_1,A_3}$.
%
  \end{Example}

  \begin{Lemma}\label{L:SemisimpleInduction}
    Suppose that $\blam\in\Parts[\alpha]$, for $\alpha\in Q^+_n$.
    Then, as $\Rx[\alpha+\alpha_i](\K'[x^\pm])$-modules,
    \[
      \iInd\bigl(\LSlam(\K'[x^\pm])\bigr)\cong
            \bigoplus_{A\in\Add_i(\blam)}\LSlam[\blam{+}A](\K'[x^\pm]).
            \quad\text{and}\quad
      \iInd\bigl(\GSlam(\K'[x^\pm])\bigr)\cong
            \bigoplus_{A\in\Add_i(\blam)}\GSlam[\blam{+}A](\K'[x^\pm]).
    \]
  \end{Lemma}

  \begin{proof}
    Let $\Domin$.  By \autoref{L:SpechtSubmodule},
    $\DSlam\cong\Rx(\K'[x^\pm])\Dzlam$. Hence, it is enough to describe
    \[
      \iInd\bigl(\Rx(\K'[x^\pm])\Dzlam\bigr)=\Rx[\alpha+\alpha_i](\K'[x^\pm])\Dzlam.
    \]
    Let $\iota_{\alpha,i}\map{\Rx[\alpha](\K'[x^\pm])}\Rx[\alpha+\alpha_i](\K'[x^\pm])$
    be the embedding of \autoref{P:exactness}.  Now
    $\Dzlam=\DefPol F_{\Dtlam}$ by \autoref{P:zlam}, so
    \begin{align}
      \begin{split}
      \iota_{\alpha,i}(\Dzlam) & = \DefPol F_{\Dtlam} \ei[\Dilam i]
        = \DefPol F_{\Dtlam}\sum_{\t\in\Std(\Dilam i)}\frac1{\Dgt}F_\t\\
       &= \DefPol \sum_{\substack{\t\in\Std(\Dilam i)\\\t_{\downarrow}=\Dtlam}}
                                  \frac1{\Dgt}F_\t
        = \sum_{A\in\Add_i(\Dtlam)}
          \frac{\DefPol}{\Dgt[{\Dtlam[\blam,A]}]}F_{\Dtlam[\blam,A]},
        \end{split}\label{E:zlamup}
    \end{align}
    where the second equality follows from \autoref{L:Separation} and
    \autoref{P:fstaction}. Note that the coefficients in the last
    equation are homogeneous and, hence, invertible in~$\K'[x^\pm]$. Therefore,
    by \autoref{L:FufstFv}, the induced module
    $\iInd\bigl(\DSlam(\K'[x^\pm])\bigr)$ is spanned by the elements
    $\set{\Dfst[\s{\Dtlam[\blam,A]}]|\s\in\Std(\blam{+}A)
                   \text{ and }A\in\Add_i(\blam)}$.
    \autoref{C:Vblam} now implies the result.
  \end{proof}

  The second last line of the proof of \autoref{L:SemisimpleInduction}
  is the reason why we are working over the polynomial rings $\k[x]$ and
  $\K'[x^\pm]$ in this section rather than over the multivariate
  polynomial rings $\kx$ and $\K'[\ux^\pm]$.

  \begin{Proposition}\label{P:iInduction}
    Suppose that $\blam\in\Parts[\alpha]$ and let
    $A_1>\dots>A_z$ be the addable $i$-nodes of $\blam$. Then
    there exist $\Rx[\alpha+\alpha_i](\k[x])$-module filtrations
    \begin{align*}
        \iInd\bigl(\LSlam(\k[x])\bigr)=\LSlam[\blam,1](\k[x])
           \supset\LSlam[\blam,2](\k[x])
           \supset\dots\supset\LSlam[\blam,z-1](\k[x])
           \supset\LSlam[\blam,z](\k[x])\supset0\\
        \iInd\bigl(\GSlam(\k[x])\bigr)=\GSlam[\blam,z](\k[x])
           \supset\LSlam[\blam,z-1](\k[x])
           \supset\dots\supset \GSlam[\blam,2](\k[x])
           \supset\GSlam[\blam,1](\k[x])\supset0
    \end{align*}
    such that
    $\LSlam[\blam,k](\k[x])/\LSlam[\blam,k+1](\k[x])
         \cong q^{\Ldilam[A_k]}\LSlam[\blam+A_k](\k[x])$
    and
    $\GSlam[\blam,k](\k[x])/\GSlam[\blam,k-1](\k[x])
         \cong q^{\Gdilam[A_k]}\GSlam[\blam+A_k](\k[x])$,
    for $1\le k\le z$.
  \end{Proposition}

  \begin{proof}
    If $\Add_i(\blam)=\emptyset$ then $\iInd(\LSlam(\k[x]))=0$ by
    \autoref{L:SemisimpleInduction}, so we can assume
    $\Add_i(\blam)\ne\emptyset$.  We only consider
    $\iInd\bigl(\LSlam(\k[x])\bigr)$. Set
    $\LZlamu
    =q^{-\defect(\blam)-\Ldeg(\Ltlam)}\Rx[\alpha+\alpha_i]\iota_{n,i}(\Lzlam)$.
    Then $\iInd\bigl(\LSlam(\k[x])\bigr)\cong\LZlamu$, by
    \autoref{L:SpechtSubmodule}, so, it is enough to show that $\LZlamu$
    has the required filtration. To do this we first construct a basis
    for~$\LZlamu$.

    By \autoref{T:kcellular},
    $\iota_{n,i}(\Gpsist[\Ltlam\Ltlam])
    =\sum_{(\s,\t)\in\Std^2(\Parts[n+1])}a_{\s\t}\Gpsist$,
    for $a_{\s\t}\in\k[x]$. Therefore, if $h\in\Rx[n+1](\k[x])$ then
    \[
        \iota_{n,i}(h\Lzlam) = \sum_{(\s,\t)\in\Std^2(\Parts[n+1])}
               a_{\s\t}h\Lylam\ei[\Lilam i]\Gpsist
    \]
    By \autoref{E:zlamup}, we may assume that $a_{\s\t}\ne0$ only if
    $\t=\Ltlam[\blam,A_k]$, for $1\le k\le z$. Further, by
    \autoref{C:ypsi}, if $\s\ne\Ltlam[\blam,A_k]$ then
    $\Lylam\ei[\Lilam i]\Gpsist$ can be written as a linear combination of
    more dominant terms, so we can assume that $\s=\t$. That is,
    \[
        \iota_{n,i}(h\Lzlam) = \sum_{k=1}^z
        a_kh\Lylam\ei[\Lilam i]\Gpsist[{\Ltlam[\blam,A_k]\Ltlam[\blam,A_k]}],
        \qquad\text{for }a_k\in\k[x].
    \]
    By \autoref{C:LpsiGpsi}, the product
    $\Lpsist[\u\v]\Gpsist[{\Ltlam[\blam,A_k]\Ltlam[\blam,A_k]}]\ne0$ only if
    $\Ltlam[\blam,A_k]\Gedom\v$. Since we also need
    $\br(\v)=\br(\Ltlam[\blam,A_k])$, the term
    $\Lpsist[\u\v]\Gpsist[{\Ltlam[\blam,A_k]\Ltlam[\blam,A_k]}]$
    is nonzero only if
    $\v=\Ltlam[\blam,A_l]$ for $1\le l\le k$.

    For $1\le k\le z$ let $n_k=\Ltlam[\blam{+}A_k](A_k)\in\set{1,\dots,n}$,
    $\psi_{n..n_k}=\psi_{n}\dots\psi_{n_k}$ if $n_k<n+1$ and
    set $\psi_{n..n_k}=1$ if $n_k=n+1$. Observe that
    $\Ltlam[\blam,A_k]=\psi_{n..n_k}\Ltlam[\blam+A_k]$. Therefore, in
    $\Rx[n+1](\k[x])$,
    \[
      y_{n+1}^{\Ldilam[A_k]}\psi_{n..n_k}\iota_{n,i}(\Lpsist[\Ltlam\Ltlam])
        =y_{n+1}^{\Ldilam[A_k]}\psi_{n..n_k}\Lylam\ei[\Lilam i]
        =\Lylam[\blam+A_k]\ei[\Lilam i]\psi_{n..n_k}
        =\Lpsist[{\Ltlam[\blam+A_k]\Ltlam[\blam,A_k]}].
    \]
    For $\s\in\Std(\Parts[\blam+A_k])$ set
    $\Lzlamus=\Lpsist[{\s\Ltlam[\blam,A_l]}]
           \Gpsist[{\Ltlam[\blam,A_k]\Ltlam[\blam,A_k]}]$.
    Then we have shown that
    \[
       y_{n+1}^{\Ldilam[A_k]}\psi_{n..n_k}\iota_{n,i}(\Lpsist[\Ltlam\Ltlam])
           \Gpsist[{\Ltlam[\blam,A_k]\Ltlam[\blam,A_l]}]
        = \sum_{l=k}^za_l\Lpsist[{\s\Ltlam[\blam,A_l]}]
           \Gpsist[{\Ltlam[\blam,A_k]\Ltlam[\blam,A_l]}]
           =a_k\Lzlamus,
    \]
    where the equality follows from \autoref{C:LpsiGpsi}. In particular,
    $a_k\Lzlamus\in\iInd\LSlam$, whenever $\s\in\Std(\blam+A_k)$ and
    $1\le k\le z$.

    Let $M$ be the free $\k[x]$-module spanned by
    $\set{\Lzlamus|\s\in\Std(\blam+A_k)\text{ and } 1\le k\le z}$. We
    claim that $M=\LZlamu=\iInd\LSlam(\k)$, which is equivalent claiming
    that $a_k\in\k^\times$, for $1\le k\le z$. If $x$ divides
    some $a_k$ then the $\K'$-dimension of $\LZlamu\otimes_{\k[x]}\K'$ is
    strictly smaller than the $\K'[x^\pm]$-rank of $\iInd\LSlam(\K'[x^\pm])$ by
    \autoref{L:SemisimpleInduction}, which is a contradiction.
    Therefore, $a_k\in\k$ for $1\le k\le z$. An easy argument using
    Nakayama's lemma (cf.~\cite[Proposition~4.6]{HuMathas:GradedInduction}),
    now shows that $M=\LZlamu$. In particular, this shows that
    $\set{\Lzlamus|\s\in\Std(\blam+A_k)\text{ and } 1\le k\le z}$ is a
    basis of $\LZlamu$.

    We can construct the promised filtration of $\LZlamu$. Define
    \[
         \LSlam[\blam,k](\k[x]) = \<\Lzlamus\mid\s\in\Std(\blam{+}A_m)
            \text{ for }1\le m\le k\>, \qquad\text{ for }0\le k\le z.
    \]
    Then
    $\LZlamu=\LSlam[\blam,1](\k[x]) \supset\LSlam[\blam,2](\k[x])
           \supset\dots\supset\LSlam[\blam,z-1](\k[x])
           \supset\LSlam[\blam,z](\k[x])\supset0$
    and each $\LSlam[\blam,k](\k[x]) $ is an $\Rx(\k[x])$-submodule of~$\LZlamu$
    by \autoref{T:kcellular}. By \autoref{C:LpsiGpsi}, for $1\le k\le z$
    define homogeneous $\Rx(\k[x])$-module homomorphisms
    $\pi_k\map{q^{\Ldilam[A_k]}\LSlam[\blam{+}A_k](\k[x])}
          \LSlam[\blam,k](\k[x])/\LSlam[\blam,k-1](\k[x])$  by
    \[
    \pi_k\Bigl(\Lpsist[{\s\Ltlam[\blam,A_k]}]+\LRlam[(\blam{+}A_k)]\Bigr)
            =\Lpsist[\s{\Ltlam[\blam,A_k]}]
              \Gpsist[{\Gtlam[\blam{+}A_z]\Ltlam[\blam,A_z]}]
              +\LSlam[\blam,k-1](\k[x])
            =\Lzlamus+\LSlam[\blam,k-1](\k[x]),
    \]
    for $\s\in\Std(\blam{+}A_k)$. By construction, these maps are
    surjective and hence bijective in view of
    \autoref{L:SemisimpleInduction}. To complete the proof we need to
    check that the map $\pi_k$ is homogeneous of degree~$0$.  Now,
    $\Ldeg(\Gtlam[\blam,A])=\Ldeg(\Gtlam)+\Gdilam$ and
    $\Gdeg(\Ltlam[\blam,A])=\Ldeg(\Ltlam)+\Ldilam$. Recalling the degree
    shifts in the definition of~$\LZlamu$,
    \begin{align*}
        \deg\pi_k &=
           \deg\bigl(\Gpsist[{\Gtlam[\blam{+}A_z]\Ltlam[\blam,A_z]}]\bigr)
               +\Ldeg(\Ltlam[\blam,A_k])
               -\bigl(\defect(\blam)+\Gdeg(\Gtlam)\bigr)-\Ldilam[A_k]=0,
    \end{align*}
    where we have once again used \autoref{C:degCodeg}.
  \end{proof}

  \begin{Corollary}\label{C:SpechtInduction}
    Suppose that $L$ is a $\k[x]$-module, $\blam\in\Parts[\alpha]$ and let
    $A_1>\dots>A_z$ be the addable $i$-nodes of $\blam$. Then
    there exist $\Rx[\alpha+\alpha_i](L)$-module filtrations
    \begin{align*}
        \iInd\bigl(\LSlam(L)\bigr)=\LSlam[\blam,1](L)
           \supset\LSlam[\blam,2](L)
           \supset\dots\supset\LSlam[\blam,z-1](L)
           \supset\LSlam[\blam,z](L)\supset0\\
        \iInd\bigl(\GSlam(L)\bigr)=\GSlam[\blam,z](L)
           \supset\LSlam[\blam,z-1](L)
           \supset\dots\supset \GSlam[\blam,2](L)
           \supset\GSlam[\blam,1](L)\supset0
    \end{align*}
    such that
    $\LSlam[\blam,k](L)/\LSlam[\blam,k+1](L)
         \cong q^{\Ldilam[A_k]}\LSlam[\blam+A_k](L)$
    and
    $\GSlam[\blam,k](L)/\GSlam[\blam,k-1](L)
         \cong q^{\Gdilam[A_k]}\GSlam[\blam+A_k](L)$,
    for $1\le k\le z$.
  \end{Corollary}

  In particular, this result includes filtrations of the induced Specht
  modules for the cyclotomic KLR algebras $\Rn(\k)$.
  In type~$\Aone$, this includes the main theorem of
  \cite[Theorem~4.11]{HuMathas:GradedInduction}, which describes Specht
  filtrations of the $\Rn<\alpha+\alpha_i>(L)$-modules
  $\iInd\bigl(\DSlam(L)\bigr)$ for $\Domin$.

  Finally, we note that we obtain the graded branching rules for the
  Specht modules of $\Rx(\K[x])$ by taking $L=\K$, or $L=\K[x]$, in
  \autoref{C:SpechtRestriction} and \autoref{C:SpechtInduction}.

  \subsection{Two dualities}
  As in \autoref{SS:Branching}, we continue to assume that $(\bc,\br)$ is
  a content system with values in $\k[x]$ and let $\K$ be a
  field that is a $\k$-algebra. In this section we work in the categories
  $\Rep_\K\Rx(\K[x])$ and $\Proj_\K\Rx(\K[x])$ of graded $\Rx(\K[x])$-modules
  that are finite dimensional as $\K$-vector spaces.

  Recall from \autoref{E:Dual} that $\circledast$ defines a graded
  duality on $\Rx(\K[x])$-modules. Similarly,  define
  $\#$ to be the graded functor given by
  \begin{equation}\label{E:HashDual}
    M^\#=\HOM_{\Rx(\K[x])}(M,\Rx(\K[x])),\qquad\text{ for $M\in\Rep_\K\Rx(\K[x])$},
  \end{equation}
  with the natural action of $\Rx(\K[x])$ on $M^\#$.  Consider
  $\circledast$ and $\#$ as endofunctors of $\Rep_\K\Rx(\K[x])$ and
  $\Proj_\K\Rx(\K[x])$.  As noted in \cite[Remark~4.7]{BK:GradedDecomp},
  \autoref{T:Symmetric} implies that these two functors agree up to
  shift.
  \notation{$M^\#$}{The projective dual: $M^\#=\Hom_{\Rx(\K[x])}(M,\Rx(\K[x]))$}[E:HashDual]

  \begin{Lemma}\label{L:SameDuals}
    Let $\alpha\in Q^+$. Then
    $\#\cong q^{2\defect(\alpha)}\circ\circledast$
    as endofunctors of $\Rep_\K\Rx[\alpha](\K[x])$.
  \end{Lemma}

  \begin{proof}
    By \autoref{T:Symmetric},
    $\Rx[\alpha](\K[x])\cong q^{2\defect(\alpha)}(\Rx[\alpha](\K[x]))^\circledast$.
    If $M\in\Rep_\K\Rx[\alpha](\K[x])$ then
    \begin{align*}
       M^\# & =\HOM_{\Rx[\alpha](\K[x])}(M,\Rx[\alpha](\K[x]))
       =\HOM_{\Rx[\alpha]}\bigl(M,
            q^{2\defect(\alpha)}(\Rx[\alpha](\K[x]))^\circledast\bigr)\\
       &\cong\HOM_{\Rx[\alpha]}\bigl(M,
            q^{2\defect(\alpha)}\HOM_{\K[x]}(\Rx[\alpha](\K[x]),\K[x])\bigr)\\
       &\cong q^{2\defect(\alpha)}
            \HOM_{\K[x]}\bigl(M\otimes_{\Rx[\alpha](\K[x])}\Rx[\alpha](\K[x]),
                            \K[x]\bigr)\\
       &\cong q^{2\defect(\alpha)}M^\circledast,
    \end{align*}
    where the third isomorphism is the standard hom-tensor adjointness.
    All of these isomorphisms are functorial, so the lemma follows.
  \end{proof}

  As $M$ is a finite dimensional $\K$-vector space,
  $(M^\circledast)^\circledast\cong M$ for all $M\in\Rep_\K\Rx(\K[x])$.
  Hence, $(M^\#)^\#\cong M$ by \autoref{L:SameDuals}.  Therefore,
  $\circledast$ and $\#$ define self-dual equivalences on the module
  categories $\Rep_\K\Rx(\K[x])$ and $\Proj_\K\Rx(\K[x])$.

  \begin{Proposition}\label{P:CommutingDuals}
    Suppose that $i\in I$. Then there are functorial
    isomorphisms
    \begin{align*}
    \circledast\circ\iRes&\cong\iRes\circ\circledast
          \map{\Rep_\K\Rx[n+1](\K[x])}\Rep_\K\Rx(\K[x]),\\
      \#\circ\iInd&\cong\iInd\circ\#\map{\Proj_\K\Rx(\K[x])}\Proj_\K\Rx[n+1](\K[x]).
    \end{align*}
  \end{Proposition}

  \begin{proof}
    The isomorphism $\circledast\circ\iRes\cong\iRes\circ\circledast$ is
    immediate from the definitions. For the second isomorphism, recall
    that if $P\in\Proj_\K\Rx(\K[x])$ then
    $\HOM_{\Rx(\K[x])}(P,M)\cong\HOM_{\Rx(\K[x])}(M,\Rx(\K[x]))\otimes_{\Rx(\K[x])}M$,
    for any $\Rx(\K[x])$-module~$M$. Now,
    \[
       (\Rx[n+1](\K[x])\ei[n,i])^\#
        =\HOM_{\Rx[n+1]\K[x])}(\Rx[n+1](\K[x])\ei[n,1],
                       \Rx[n+1](\K[x]))\cong \Rx[n+1](\K[x])\ei[n,i],
    \]
    where the last isomorphism follows because $\ei[n,i]^*=\ei[n,i]$. Therefore,
    \begin{align*}
      \iInd(P^\#)
        &=\HOM_{\Rx(\K[x])}\Bigl(P,\Rx(\K[x])\Bigr)
                  \otimes_{\Rx(\K[x])}\Rx[n+1](\K[x])\ei[n,i]\\
        &\cong\HOM_{\Rx(\K[x])}\Bigl(P,\Rx[n+1](\K[x])\ei[n,i]\Bigr)\\
        &\cong\HOM_{\Rx(\K[x])}\Bigl(P,
           \HOM_{\Rx[n+1](\K[x])}\bigl(\ei[n,i]\Rx[n+1](\K[x]),
               \Rx[n+1](\K[x])\bigr)\Bigr)\\
        &\cong\HOM_{\Rx[n+1](\K[x])}\Bigl(
            P\otimes_{\Rx(\K[x])}\Rx[n+1](\K[x])\ei[n,i],
                      \Rx[n+1](\K[x])\Bigr)\\
        &\cong(\iInd P)^\#,
    \end{align*}
    where the second last isomorphism is the usual tensor-hom adjointness.
  \end{proof}

  It follows from \autoref{P:CommutingDuals} and \autoref{L:SameDuals}
  that the functors $\circledast$ and $\iInd$, and $\#$ and $\iRes$,
  commute up to shift.

\subsection{Grothendieck groups and the Cartan pairing}
\label{SS:GrothendieckGroups}
  We are now ready to prove the categorification theorems from the
  introduction, which will allow us to classify the simple
  $\Rx(\K[x])$-modules and prove our modular branching rules.  As in the
  last two sections we continue to assume that $\Rx(\K[x])$ is defined
  using a graded content system with values in~$\k[x]$, where the field
  $\K$ is a $\k$-algebra. In particular,
  this means that the graded branching rules for the Specht modules for
  $\Rx(\K[x])$ are given by the results in \autoref{SS:Branching}.

  Recall that $q$ is an indeterminate over~$\Z$ and that
  $\A=\Z[q,q^{-1}]$.  Let $[\Rep_\K\Rx(\K[x])]$, $[\Proj_\K\Rx(\K[x])]$,
  be the Grothendieck groups of the corresponding categories of graded
  $\Rx(\K[x])$-modules, which are categories of finite dimensional
  $\K$-vector spaces.  We consider each of these Grothendieck groups as
  $\A$-modules, where $q$ acts by grading shift. If $M$ is a
  module in one of these categories, let $[M]$ be its image in the
  corresponding Grothendieck group.  Since~$q$ is the grading shift
  functor, which is exact, $[qM]=q[M]$.

  \notation{$[\Rep_\K\Rx(\K[x])]$}{Grothendieck group of $\Rep_\K\Rx(\K[x])$}
  \notation{$[\Proj_\K\Rx(\K[x])]$}{Grothendieck group of $\Proj_\K\Rx(\K[x])$}
  \notation{$\RepN$}{$\bigoplus_{n\ge0}[\Rep_\K\Rx(\K[x])]$}
  \notation{$\ProjN$}{$\bigoplus_{n\ge0}[\Proj_\K\Rx(\K[x])]$}

  Rather than considering the Grothendieck groups in isolation it is
  advantageous to consider all of them together.  Define
  \[
    \RepN  = \bigoplus_{n\ge0}\RepK\And*
    \ProjN = \bigoplus_{n\ge0}\ProjK.
  \]
  These Grothendieck groups are independent of the choice of cellular
  basis in \autoref{T:kcellular}, however, we give parallel categorification
  results for the two $\psi$-bases of~$\Rx(\K[x])$.

  By \autoref{P:exactness}, the induction and restriction functors $\iInd$
  and $\iRes$ are exact and send projectives to projectives.  Therefore
  they induce $\A$-linear automorphisms of the Grothendieck groups
  $\RepN$ and $\ProjN$, which are given by
  \[
      \iInd{}[M] = [\iInd M] \And
      \iRes{}[M] = [\iRes M]
  \]
  for all modules $M$ and $i\in I$.

  Let $M$ and $N$ be free $\A$-modules.  A \emph{semilinear} map of
  $\A$-modules is a $\Z$-linear map $\theta\map MN$ such that
  $\theta(q^dm)=q^{-d}\theta(m)=\overline{q^d}\theta(m)$, for all
  $d\in\Z$ and $m\in M$.  A \emph{sesquilinear} map $f\map{M\times N}\A$
  is a function that is semilinear in the first variable and linear in
  the second.

  Let $\CPair{\ }{\ }\map{\ProjN\times\RepN}\A$ be the \emph{Cartan pairing}, which is determined by
  \begin{equation}\label{E:CartanPairing}
    \CPair[\big]{[P]}{[M]}=\delta_{mn}\gdim\HOM_{\Rx(\K[x])}\bigl(P,M\bigr),
  \end{equation}
  for $P\in\Proj_\K\Rx[m](\K[x])$ and $M\in\Rep_\K\Rx[n](\K[x])$.
  The Cartan pairing is sesquilinear because
  \[
  \HOM_{\Rx(\K[x])}(q^{-k}P,M)\cong\HOM_{\Rx(\K[x])}(P,q^{k}M)
         \cong q^{k}\HOM_{\Rx(\K[x])}(P,M),
         \qquad\text{for any }k\in\Z.
  \]
  The Cartan pairing is characterised by either of the two properties:
  \begin{equation}\label{E:CParingYmuDnu}
    \CPair[\big]{[\LYmu[\blam]]}{[\LDmu]} = \delta_{\blam\bmu}
    \qquad\text{or}\qquad
    \CPair[\big]{[\GYnu]}{[\GDnu[\bsig]} = \delta_{\bnu\bsig}
  \end{equation}
  for $\blam,\bmu\in\LKlesh$ or $\bnu,\bsig\in\GKlesh$, respectively.
  \notation{$\CPair{\ }{\ }$}
           {Cartan pairing $\ProjN\times\RepN\to\A$}
           [E:CartanPairing]

  \begin{Remark}
    By the remarks after \autoref{T:GradedSimples}, as abelian groups,
    \[
        \RepK\cong[\Rep\Rn(\K)]
        \quad\text{ and }\quad
        \ProjK\cong[\Proj\Rn(\K)].
    \]
    In what follows, we could work with the Grothendieck groups
    $[\Rep\Rn(\K)]$ and $[\Proj\Rn(\K)]$.
  \end{Remark}

  \subsection{Fock spaces}\label{SS:Fock}
  This section proves that $\ProjN$ and $\RepN$ categorify the integral
  form and its dual, respectively, of an irreducible integrable highest
  weight module of the quantised Kac-Moody algebra $\Uq$. We start by
  recalling the results and definitions that we need from the Kac-Moody
  universe.  The arguments in this section are mostly standard, and
  follow (and correct) \cite{Mathas:Singapore}. Our approach is similar
  to \cite{BK:GradedDecomp} except that we use the representation theory
  of the KLR algebras to construct the canonical bases, rather than vice
  versa. What is non-standard is that these arguments apply
  simultaneously in types $\Aone$ and~$\Cone$.

  Recall $\A=\Z[q,q^{-1}]$. Set $\AA=\Q(q)$.
  For $i\in I$ and $k\in\Z$ let $[k]_i = (q_i^k-q_i^{-k})/(q_i-q_i^{-1})$,
  where $q_i=q^{\di}$. If $k>0$ set $[k]_i!=[1]_i[2]_i\dots[k]_i$.
  For non-commuting indeterminates $u$ and $v$ and $i\in I$ set
  \[
    (\adq u)^c(v)=\sum_{d=0}^c(-1)^d\tfrac{[c]_i!}{[c-d]_i![k]_i!}u^{c-d}vu^d.
  \]
  \notation{$q_i$}{For $i\in I$, $q_i=q^{\di}$}
  \notation{$[k]_i$}{For $k\in\Z$, $[k]_i$ is the quantum integer $(q_i^k-q_i^{-k})/(q_i-q_i^{-1})\in\A$}
  \notation{$[k]_i!$}{For $k>0$, $[k]_i!$ is the quantum factorial $[1]_i\dots[k]_i\in\A$}
  \notation{$\Uq$}{Quantum group of the Kac-Moody algebra $\mathfrak{g}_\Gamma$}
  \notation{$E_i,F_i,K_i^{\pm}$}{Generators of $\Uq$}

  \begin{Definition}
    The \textbf{quantum group} $\Uq$ is the $\AA$-algebra with
    generators $E_i$, $F_i$, $K^\pm_i$, for $i\in I$, and relations:
    \bgroup
       \setlength{\abovedisplayskip}{1pt}
       \setlength{\belowdisplayskip}{2pt}
      \begin{align*}
        K_iK_j &=K_jK_i,
           &K_iK_i^{-1}&=1,
           &[E_i,F_j]&=\delta_{ij}\tfrac{K_i-K_i^{-1}}{q_i-q_i^{-1}},
      \end{align*}
      \begin{align*}
          K_iE_jK_i^{-1}&=q^{\cij}E_j,
          &K_iF_jK_i^{-1}&=q^{-\cij}F_j,
      \end{align*}
      \begin{align*}
        (\adq E_i)^{1-\cij}(E_j)=0=(\adq F_i)^{1-\cij}(F_j),\qquad
        \text{for }i\ne j.
      \end{align*}
   \egroup
    The quantum group $\Uq$ is a Hopf algebra with coproduct determined by
    \[
      \Delta(K_i)=K_i\otimes K_i,\quad
      \Delta(E_i)=E_i\otimes K_i+1\otimes E_i\And*
      \Delta(F_i)=F_i\otimes1+ K_i^{-1}\otimes F_i,
    \]
    for $i\in I$.
  \end{Definition}

  We will only need basic facts about highest weight theory and
  canonical bases for~$\Uq$.  Detailed accounts of the representation
  theory of $\mathfrak{g}_\Gamma$ and $\Uq$ can be found in
  \cite{Ariki:book,Lusztig:QuantBook,Kac}.

  \notation{$\LFock,\GFock$}{$\Uq$-Fock spaces associated to the
  $\Lpsis[]$ and $\Gpsist[]$ bases}
  \notation{$\Lfock,\Gfock$}{Basis elements of the Fock spaces $\LFock$ and $\GFock$}

  \begin{Definition}
    Let $\Lambda\in P^+$. The \emph{combinatorial Fock spaces}
    $\LFock[\A]$ and $\GFock[\A]$ are the free $\A$-modules with basis
    the sets of symbols $\set{\Lfock|\blam\in\Parts*}$ and
    $\set{\Gfock|\blam\in\Parts*}$, respectively.

    Set $\LFock[\AA]=\AA\otimes_\A\LFock[\A]$ and
    $\GFock[\AA]=\AA\otimes_\A\GFock[\A]$.
  \end{Definition}

  By definition, $\LFock[\AA]$ and $\GFock[\AA]$ are infinite dimensional
  $\AA$-vector spaces.  For $\Domin$, identify $\Dfock$ with
  $1_{\AA}\otimes_\A\Dfock$, for $\blam\in\Parts$.  Then
  $\set{\Dfock|\blam\in\Parts*}$ is an $\AA$-basis $\DFock$.

  Let $\zero=(0|\dots|0)\in\Parts*$ be the empty $\ell$-partition.
  Recall the integers $\Ldilam$, $\Gdilam$, and $d_i(\blam)$ from
  \autoref{D:defect}. Note that these definitions depend on
  $(\Lambda,\charge)$.

  \begin{Theorem}[Hayashi~\cite{Hayashi}, Misra-Miwa~\cite{MisraMiwa},
                  Premat~\cite{Premat:TypeCFock}]\label{T:FockSpace}
    Let $\Lambda\in P^+$.
    \begin{enumerate}
      \item The Fock space $\LFock$ is an integrable $\Uq$-module with $\Uq$-action
      determined by
      \[
         E_i\cdot\Lfock = \sum_{B\in\Rem_i(\blam)}
                            q^{-\Gdilam[B]}\Lfock[\blam{-}B],
         \quad
         F_i\cdot\Lfock=\sum_{A\in\Add_i(\blam)}q^{\Ldilam}\Lfock[\blam{+}A],
         \And*
         K_i\cdot\Lfock=q^{-d_i(\blam)}\Lfock,
      \]
      for $i\in I$ and $\blam\in\Parts$.

      \item The Fock space $\GFock$ is an
      integrable $\Uq$-module with $\Uq$-action determined by
      \[
         E_i\cdot\Gfock = \sum_{B\in\Rem_i(\blam)}
                            q^{-\Ldilam[B]}\Gfock[\blam{-}B],
         \quad
         F_i\cdot\Gfock=\sum_{A\in\Add_i(\blam)}q^{\Gdilam}\Gfock[\blam{+}A],
         \And*
         K_i\cdot\Gfock=q^{-d_i(\blam)}\Gfock,
      \]
      for $i\in I$ and $\blam\in\Parts$.

    \end{enumerate}
  \end{Theorem}

  \begin{proof}
    To prove (a) and (b) it is enough to verify that these actions
    respect the relations of~$\Uq$. Recall the sign automorphism of
    \autoref{S:SignAutomorphism}. In particular, by
    \autoref{L:SignedCombinatorics}, $\Ldilam=\Gdilam*[A'](\blam')$,
    where if~$A\in\Add(\blam)\cup\Rem(\blam)$ then $\Ldilam$ is computed
    with respect to $(\Lambda,\charge)$ and $\Gdilam*[A']$ is computed with
    respect to $(\Lambda^\sgn,\charge^\sgn)$.  Hence, parts~(a) and~(b)
    are equivalent and it suffices to prove~(b).

    If $\Gamma$ is a quiver of type $\Aone$ then (b) is due to
    Hayashi~\cite{Hayashi} in level~$1$, with the result in higher
    levels following by applying the coproduct, as was observed by Misra
    and Miwa~\cite{MisraMiwa}.  For quivers of type~$\Cone$, this was
    proved by Premat~\cite[Theorem~3.1]{Premat:TypeCFock} in level~$1$
    (see also Kim and Shin~\cite{KimShin:TypeCLLT}), with
    the result in higher levels again following by applying the
    coproduct, as noted already in \cite[\S1]{ArikiParkSpeyer:C}.
  \end{proof}

  \autoref{T:FockSpace} does not give the $\Uq$-actions on the Fock
  spaces that we want because this action does not commute with the bar
  involution on $L(\Lambda)$, which is introduced in
  \autoref{S:CanonicalBases} below. Let $\tau\map\Uq\Uq$ be anti-linear
  anti-automorphism given by
  \[
        \tau(K_i)=K_i^{-1},\quad
        \tau(E_i)=q^{\di}F_iK_i^{-1}\quad\text{and}\quad
        \tau(F_i)=q^{-\di}K_iE_i
        \qquad\text{for }i\in I.
  \]
  This map is not an involution but it is invertible.
  Twisting the $\Uq$-action from \autoref{T:FockSpace} by $\tau$ gives
  the $\Uq$-action on the Fock space that we need.

  \begin{Corollary}\label{C:FockSpace}
    Suppose that $\Lambda\in P^+$.
    \begin{enumerate}
      \item The Fock space $\LFock$ is an integrable $\Uq$-module with
      $\Uq$-action determined by
      \[
         E_i\Lfock = \sum_{B\in\Rem_i(\blam)} q^{\Ldilam[B]}\Lfock[\blam{-}B],
         \quad
         F_i\Lfock=\sum_{A\in\Add_i(\blam)}q^{-\Gdilam}\Lfock[\blam{+}A],
         \And*
         K_i\Lfock=q^{d_i(\blam)}\Lfock,
      \]
      for $i\in I$ and $\blam\in\Parts$.

      \item The Fock space $\GFock$ is an
      integrable $\Uq$-module with $\Uq$-action determined by
      \[
         E_i\Gfock = \sum_{B\in\Rem_i(\blam)} q^{\Gdilam[B]}\Gfock[\blam{-}B],
         \quad
         F_i\Gfock=\sum_{A\in\Add_i(\blam)}q^{-\Ldilam}\Gfock[\blam{+}A],
         \And*
         K_i\Gfock=q^{d_i(\blam)}\Gfock,
      \]
      for $i\in I$ and $\blam\in\Parts$.

    \end{enumerate}
  \end{Corollary}

  \begin{proof}
    We consider only (a) and leave part~(b) to the reader since this is
    similar.  Using \autoref{T:FockSpace}, and the fact that $\tau$ is an
    anti-isomorphism of~$\Uq$, we can define a new action of $\Uq$
    on~$\LFock$ by~$E_i\Lfock=\tau(F_i)\cdot\Lfock$,
    $F_i\Lfock=\tau(E_i)\cdot\Lfock$ and
    $K_i\Lfock=\tau(K_i)\cdot\Lfock$, for $i\in I$ and
    $\blam\in\Parts$. Therefore,
    \begin{align*}
      E_i\Lfock &= \tau(F_i)\cdot\Lfock = q^{-\di}K_iE_i\cdot\Lfock
            = \sum_{B\in\Rem_i(\blam)}
                   q^{\di+d_i(\blam)-\Gdilam[B]}\Lfock[\blam{-}B]\\
           &= \sum_{B\in\Rem_i(\blam)} q^{\Ldilam[B]}\Lfock[\blam{-}B],
    \end{align*}
    where the last equality follows from  \autoref{L:defectSum}.
    The other identities are similar.
  \end{proof}

  In what follows we always use the $\Uq$-action on the Fock spaces
  $\LFock$ and $\GFock$ from \autoref{C:FockSpace}.  We work with both
  Fock spaces because they are closely intertwined and by using both
  Fock spaces we will be able to determine the labelling of the simple
  $\Rx(\K[x])$-modules and the map $\mull$ from \autoref{D:Mullineux}.  As
  our notation suggests, the Fock spaces $\LFock$ and $\GFock$ can be
  naturally associated with the $\Lpsist[]$ and $\Gpsist[]$-bases of
  $\Rx(\K[x])$, respectively. To make this connection precise we need a
  little more notation.

  A vector $v$ in a $\Uq$-module  has \emph{weight} $\wt(v)=\theta$ if
  $K_iv=q^{(\theta|\alpha_i)}v$, for all $i\in I$.
  \autoref{C:FockSpace}, and \autoref{L:defectWeight}, imply that if
  $\blam\in\Parts[\alpha]$ then
  \begin{equation}\label{E:weight}
    \wt(\Lfock)=\Lambda-\alpha=\wt(\Gfock),\qquad
        \text{for all }\blam\in\Parts[\alpha].
  \end{equation}
  In particular,
  $\LFock$ and $\GFock$ are both integrable highest weight modules
  for~$\Uq$ and $\Lfock[\zero]$ and $\Gfock[\zero]$ are highest weight
  vectors of weight $\Lambda$.

  \notation{$\wt(v)$}{Weight of an element in a Fock space}[E:weight]
  \notation{$L(\Lambda)$}{Irreducible integrable highest weight module
  for $\Uq$ of weight $\Lambda$}

  Let $L(\Lambda)_\AA$ be the irreducible integrable highest weight
  module for $\Uq$ with highest weight~$\Lambda$. Then
  $L(\Lambda)_\AA=\Uq v_\Lambda$, where $v_\Lambda$ is a highest weight
  vector of weight~$\Lambda$.

  \begin{Corollary}
     Let $\Lambda\in P^+$. Then
     $\Uq\Lfock[\zero]\cong L(\Lambda)_\AA\cong\Uq\Gfock[\zero]$
     as $\Uq$-modules,
  \end{Corollary}

  \begin{proof}
    By \autoref{C:FockSpace} and \autoref{E:weight},
    the vectors $\Lfock[\zero]\in\LFock$ and
    $\Gfock[\zero]\in\GFock$ are both highest weight vectors of
    weight~$\Lambda$. Therefore,
    $\Uq\Lfock[\zero]\cong L(\Lambda)_\AA\cong\Uq\Gfock[\zero]$
    required.
  \end{proof}

  To make use of this result, recall from
  \autoref{SS:GrothendieckGroups} that $\RepN$ and $\ProjN$ are the
  direct sums of Grothendieck groups of graded $\Rx(\K[x])$-modules and
  graded projective $\Rx(\K[x])$-modules, respectively, for $n\ge0$.
  In particular, $\RepN$ and $\ProjN$ are free an $\A$-modules.
  \notation{$\Parts*$}{The set $\bigcup_{n\ge0}\Parts$}
  \notation{$\LKlesh*,\GKlesh*$}
      {The sets $\bigcup_{n\ge0}\LKlesh$ and $\bigcup_{n\ge0}\GKlesh$}

  Let $\Parts*=\bigcup_{n\ge0}\Parts$, $\LKlesh*=\bigcup_{n\ge0}\LKlesh$
  and $\GKlesh*=\bigcup_{n\ge0}\GKlesh$. By \autoref{T:GradedSimples} and
  \autoref{T:TriangularDecomp}, $\RepN$ comes equipped
  with four distinguished bases:
  \begin{equation}\label{E:RepBasis}
     \set[\big]{[\LDmu]|\bmu\in\LKlesh*},\quad
     \set[\big]{[\LSlam[\bmu]]|\bmu\in\LKlesh*}\quad
     \set[\big]{[\GDnu]|\bnu\in\GKlesh*},\And*
     \set[\big]{[\GSlam[\bnu]]|\bmu\in\GKlesh*}
  \end{equation}
  Here, $\LDmu=\LDmu(\K)$, $\LSlam=\LSlam(\K)$, $\GDnu=\GDnu(\K)$ and
  $\GSlam[\bnu]=\GSlam[\bnu](\K)$ are finite dimensional $\K$-modules.
  In contrast, the projective Grothendieck group $\ProjN$ has only two
  natural bases:
  \begin{equation}\label{E:ProjBases}
     \set[\big]{[\LYmu]|\bmu\in\LKlesh*}\And
     \set[\big]{[\GYnu]|\bnu\in\GKlesh*},
  \end{equation}
  where, as in \autoref{SS:GradedDecomp}, $\LYmu=\LYmu(\K)$ and
  $\GYnu=\GYnu(\K)$ are the projective covers of $\LDmu$ and $\GDnu$,
  respectively.  Define elements $\set{\Lymu|\bmu\in\LKlesh}$ and
  $\set{\Gynu[\bmu]|\bnu\in\GKlesh}$ of $\LFock$ and $\GFock$,
  respectively, by setting
  \begin{equation}\label{E:ymu}
        \Lymu = \sum_{\blam\in\Parts} \Ldlammu\Lfock
        \And
        \Gynu = \sum_{\blam\in\Parts}\Gdlamnu[\blam\bnu]\Gfock.
  \end{equation}
  \notation{$\Lymu,\Gynu$}{Images of $[\LYmu]$ and $[\GYnu]$ in $\LFock$
  and $\GFock$}[E:ymu]

  Set $\RepN_\AA=\AA\otimes_\A\RepN$ and $\ProjN_\AA=\AA\otimes_\A\RepN$.

  \begin{Proposition}\label{P:Fock}
    Suppose that $\Lambda\in P^+$. Identify $E_i$ and $\iRes$, and $F_i$
    and $\iInd\circ q^{\di} K_i^{-1}$, for $i\in I$. Then there
    are $\Uq$-module embeddings
    \begin{align*}
      \Ldec_T&\map{\ProjN_\AA}\LFock[\AA]; [\LYmu]\mapsto\Lymu &
      \Gdec_T&\map{\ProjN_\AA}\GFock[\AA]; [\GYnu]\mapsto\Gynu
    \intertext{and $\Uq$-module surjections}
        \Ldec&\map{\LFock[\AA]}\RepN_\AA; \Lfock\mapsto[\LSlam] &
        \Gdec&\map{\GFock[\AA]}\RepN_\AA; \Gfock\mapsto[\GSlam]
    \end{align*}
    Consequently, $\ProjN_\AA\cong L(\Lambda)\cong \RepN_\AA$ as
    $\Uq$-modules.
  \end{Proposition}
  \notation{$\Ldec,\Gdec$}{Surjective decomposition maps $\Ddec\map{\DFock}{\RepN}$}[P:Fock]
  \notation{$\Ldec_T,\Gdec_T$}{Injective decomposition maps $\Ddec_T\map{\ProjN}{\DFock}$}[P:Fock]

  \begin{proof}
    Let $\Doming$. By \autoref{T:TriangularDecomp} and
    \autoref{P:PIMs}, there are well-defined $\AA$-linear
    maps~$\Ddec_T$ and~$\Ddec$, with $\Ddec_T$ injective and~$\Ddec$
    surjective. It remains to check that these maps are
    homomorphisms of $\Uq$-modules.

    Let $i\in I$. By \autoref{P:exactness},  the functors $\iRes$ and
    $\iInd$ are exact, and send projective modules to projective
    modules, so they both induce $\AA$-linear endomorphisms of the
    Grothendieck groups $\ProjN$ and $\RepN$. Taking $L=\K$ in
    \autoref{C:SpechtRestriction} and \autoref{C:SpechtInduction},
    \begin{align*}
      E_i[\DSlam] &= [\iRes \DSlam]= \sum_{B\in\Rem_i(\blam)}
          q^{\Ddilam[B]}[\DSlam[\blam-B]],\\
          F_i[\DSlam] &=[\iInd\circ q^{\di}K_i^{-1}\DSlam]
                   =\sum_{A\in\Add_i(\blam)}q^{\Ddilam+\di-d_i(\blam)}
                        [\DSlam[\blam+A]]
                   =\sum_{A\in\Add_i(\blam)}q^{-\Dodilam}[\DSlam[\blam+A]],
    \end{align*}
    where the last equality uses \autoref{L:defectSum}. Therefore, by
    identifying $E_i$ with the functor $\iRes$, and~$F_i$ with the
    functor $\iInd \circ q^{\di} K_i^{-1}$, the linear maps~$\Ldec_T$
    and~$\Ldec$ become well-defined $\Uq$-module homomorphisms by
    \autoref{C:FockSpace}. As $\Uq$-modules, $\RepN$ and $\ProjN$ are
    both cyclic because they are both generated
    by~$[\DYmu[\zero]]=[\DSlam[\zero]]=[\DDmu[\zero]]$.  By definition,
    $\Ldec_T([\LYmu[\zero]])=\Lfock[\zero]$ and
    $\Ldec(\Lfock[\zero])=[\LSlam[\zero]]$, so the proposition follows
    since $\Uq\Lfock[\zero]\cong L(\Lambda)\cong\Uq\Gfock[\zero]$ is an
    irreducible $\Uq$-module.
  \end{proof}

  Since $K_i\Dfock=q^{d_i(\blam)}\Dfock$, for $\blam\in\Parts*$, we view
  $K_i$ as a grading shift functor on $\Rep_\K\Rx(\K[x])$, for $i\in I$.
  Hereafter, for $i\in I$ we identify $E_i$ and $\iRes$, and $F_i$ and
  $\iInd\circ q^{\di}K_i^{-1}$, as functors on
  $\Rep_\K\Rx[\bullet](\K[x])$ and $\Proj_\K\Rx[\bullet](\K[x])$.

  \begin{Remark}
    Let $\Domin$. Then \autoref{P:Fock} can be interpreted as saying
    that there is a commutative diagram of $\Uq$-modules:
    \[
      \begin{tikzpicture}[>=stealth,->,shorten >=2pt,looseness=.5,auto]
        \matrix (M)[matrix of math nodes,row sep=1cm,column sep=16mm]{
          \ProjN_\AA  &  \DFock[\AA] \\ & \RepN_\AA \\
         };
         \draw[->](M-1-1)--node[above]{$\Ddec_T$}(M-1-2);
         \draw[->](M-1-2)--node[right]{$\Ddec$}(M-2-2);
         \draw[->](M-1-1)--node[below]{$\Dcar$}(M-2-2);
      \end{tikzpicture}
    \]
    The map $\Dcar\map{\ProjN_\AA}\RepN_\AA$ is given by the Cartan
    matrix, which is the natural embedding of $\ProjN_\AA$ into
    $\RepN_\AA$. Of course, $\Ddec$ is the decomposition map and
    $\Ddec_T$ is its transpose. Hence, \autoref{C:Cartan} categorifies
    \autoref{P:Fock} .
  \end{Remark}

  \begin{Remark}
    Let $\sgn$ be the sign automorphism of~$\Gamma$ from
    \autoref{D:sign}. Abusing notation slightly, the quiver automorphism
    $\sgn$ induces a unique automorphism of $\Uq$ such that
    \[
        \sgn(E_i)=E_{\sgn(i)},\qquad
        \sgn(F_i)=F_{\sgn(i)}\And
        \sgn(K_i)=K_{\sgn(i)},\qquad
        \text{for all }i\in I
    \]
    Let $\LFock*=\<\s^{\Ldom\sgn}_\blam|\blam\in\Parts*\>_\A$ and
    $\GFock*=\<\s^{\Gdom\sgn}_\blam|\blam\in\Parts*\>_\A$ be the Fock spaces
    with $\UA$-action defined using the functions
    $d^{\Ldom\sgn}_A(\blam)$ and $d^{\Gdom\sgn}_a(\blam)$ from
    \autoref{S:SignAutomorphism}. Then \autoref{L:SignedCombinatorics}
    implies that there are $\Uq$-module isomorphisms
    $\TT^\sgn_\Ldom\colon\LFock\cong\GFock*$ and $\TT^\sgn_\Gdom\colon\GFock\cong\LFock*$
    given by $\TT^\sgn_\Ldom(\Lfock)=s^{\Gdom\sgn}_{\blam'}$ and\
    $\TT^\sgn_\Gdom(\Gfock)=s^{\Ldom\sgn}_{\blam'}$, for
    $\blam\in\Parts*$. Equivalently, there are $\Uq$-module isomorphisms
    $\LFock\cong\bigl(\GFock*\bigr)^\sgn$ and
    $\GFock\cong\bigl(\LFock*\bigr)^\sgn$, where the $\Uq$ actions on
    $\LFock*$ and $\GFock*$ are twisted by~$\sgn$.
    These results should be compared with \autoref{C:SignedSpecht}.
  \end{Remark}

  We need to prove an ``integral'' version of the $\Uq$-module isomorphisms
  in \autoref{P:Fock} over~$\A$. To do this recall that Lusztig's
  $\A$-form of $\Uq$ is the $\A$-subalgebra $\UA$ of~$\Uq$ that
  generated by the quantised divided powers $E_i^{(k)}=E_i^k/\qint k!$
  and~$F_i^{(k)}=F_i^k/\qint k!$, for $i\in I$ and~$k\ge0$.
  For any $\A$-module $A$ set $\UA[A]=A\otimes_\A\UA$.

  \autoref{C:FockSpace} implies that $\UA$ acts on the
  $\A$-submodule~$\DFock$ of~$\DFock[\AA]$; compare with
  \cite[Lemma~6.15]{Mathas:ULect} and \cite[Lemma~6.2]{LLT}. Set
  \begin{equation}\label{E:IntegralFock}
        \LLa=\UA\Lfock[\zero] \And \GLa=\UA\Gfock[\zero].
  \end{equation}
  Then  \autoref{P:Fock} implies that
  $\AA\otimes_\A\LLa\cong L(\Lambda)\cong\AA\otimes_\A\GLa$, as
  $\Uq$-modules, and that:
  \notation{$\LLa,\GLa$}{Highest weight modules as submodules of
  $\LFock$ and $\GFock$}

  \begin{Corollary}\label{C:ProjK}
    Suppose that $\Lambda\in P^+$. Then $\LLa\cong\ProjN\cong\GLa$
    as $\UA$-modules.
  \end{Corollary}

  The analogue of this result for $\RepN$ requires some
  Lie theory.  Define symmetric bilinear forms
  $\LPair{\ }{\ }\map{\LFock\times\LFock}\A$ and
  $\GPair{\ }{\ }\map{\GFock\times\GFock}\A$ by
  \begin{equation}\label{E:FockPairing}
    \LPair[\big]{\Lfock}{\Lfock[\bmu]}=\delta_{\blam\bmu}q^{\defect\blam}
    \And
    \GPair[\big]{\Gfock}{\Gfock[\bmu]}=\delta_{\blam\bmu}q^{\defect\blam}
      \qquad \text{for }\blam,\bmu\in\Parts*,
  \end{equation}
  and extending linearly. By definition, both of these bilinear forms
  are non-degenerate.  By restriction, we consider $\LPair{\ }{\ }$ and
  $\GPair{\ }{\ }$ as (possibly degenerate) bilinear forms on~$\LLa$ and
  $\GLa$, respectively.
  \notation{$\LPair{\ }{\ },\GPair{\ }{\ }$}{Semilinear pairings on $\LFock$
  and $\GFock$}[E:FockPairing]

  \begin{Lemma}\label{L:ContrBiadjoint}
    Let $\Domin$.  The bilinear form $\DPair{\ }{\ }$ on $\DLa$ is
    characterised by the properties:
    \[
       \DPair[\big]{\Dfock[\zero]}{\Dfock[\zero]}=1,\quad
       \DPair[\big]{E_iu}{v}=\DPair[\big]{u}{F_iv}\And*
       \DPair[\big]{F_iu}{v}=\DPair[\big]{u}{E_iv},
    \]
    for all $i\in I$ and $u,v\in\DLa$.
  \end{Lemma}

  \begin{proof}
     By definition, $\DPair[\big]{\Dfock[\zero]}{\Dfock[\zero]}=1$. Let
     $i\in I$. To show that~$E_i$ and~$F_i$ are biadjoint with respect
     to~$\DPair{\ }{\ }$ it is enough to consider the cases when
     $u=\Dfock[\bmu]$ and $v=\Dfock$, for $\blam,\bmu\in\Parts*$. By
     \autoref{C:FockSpace},
     $\DPair[\big]{F_i\Dfock[\bmu]}{\Dfock}
         =0=\DPair[\big]{\Dfock[\bmu]}{E_i\Dfock}$
     unless $\blam=\bmu+A$ for some $A\in\Add_i(\blam)$. Moreover,
     if~$A\in\Add_i(\bmu)$ and $\blam=\bmu+A$ then using
     \autoref{C:FockSpace} and \autoref{L:defect},
     \begin{align*}
       \DPair[\big]{F_i\Dfock[\bmu]}{\Dfock}
                &=q^{\defect(\blam)-\Dodilam(\bmu)}
                 =q^{\defect(\blam)-d_i(\bmu)+\di+\Ddilam(\bmu)}
                 =q^{\defect(\bmu)+\Ddilam(\bmu)}
                 =\DPair[\big]{\Dfock[\bmu]}{E_i\Dfock}.
     \end{align*}
    Similarly, $\DPair[\big]{E_i\Dfock}{\Dfock[\bmu]}
          =\DPair[\big]{\Dfock}{F_i\Dfock[\bmu]}$,
    for all $\blam,\bmu\in\Parts$. As~$\Dfock[\zero]$ is the highest
    weight vector of weight~$\Lambda$ in the irreducible module
    $\AA\otimes_\A\DLa$, it follows by induction on weight that these
    three properties uniquely determine the bilinear form $\DPair{\ }{\ }$
    on~$\DLa$.
  \end{proof}

  As the next result shows, the pairings $\LPair{\ }{\ }$ and $\GPair{\
  }{\ }$ are closely related to the Cartan pairing defined in
  \autoref{E:CartanPairing}. Recall the functor $\#$ from
  \autoref{E:HashDual}.

  \begin{Lemma}\label{L:FormComparision}
     Suppose that $u\in\ProjN$ and $v\in\LFock$ with $\wt(v)=\beta$.
     Then
      \[
         \LPair[\big]{\Ldec_T(u^\#)}{v}
            = q^{\defect(\beta)}\CPair[\big]{u}{\Ldec(v)}
         \And
         \GPair[\big]{\Gdec_T(u^\#)}{v}
            = q^{\defect(\beta)}\CPair[\big]{u}{\Gdec(v)}
      \]
  \end{Lemma}

  \begin{proof}
    Let $\Domin$. It is enough to check this when $x=q^a[\DYmu]$
    and~$v=\Dfock$, for $a\in\Z$, $\bmu,\blam\in\DKlesh*$ and
    $\blam\in\Parts*$. As $\CPair{\ }{\ }$ is sesquilinear, and
    $\DPair{\ }{\ }$ is bilinear,
    \begin{align*}
      q^{\defect(\blam)}\CPair[\big]{q^a[\DYmu]}{\Ddec(\Dfock)}
        & = q^{\defect(\blam)-a}\sum_{\bnu\in\DKlesh*}\Ddlammu[\blam\bnu]
            \CPair[\big]{[\DYmu]}{[\DDmu[\bnu]]}\\
        &= q^{\defect(\blam)-a}\Ddlammu
         = q^{-a}\sum_{\bnu\in\DKlesh*}\Ddlammu[\bnu\bmu]
              \DPair[\big]{\Dfock[\bnu]}{\Dfock[\blam]}\\
        &=q^{-a}\DPair[\Big]{\Ddec_T\bigl([\DYmu]\bigr)}{\Dfock}
         =\DPair[\Big]{\Ddec_T\bigl([q^{a}\DYmu]^\#\bigr)}{\Dfock}.
    \end{align*}
    The last equality follows because $[q^a\DYmu]^\#=q^{-a}[\DYmu]$, by
    \autoref{E:HashDual}, since $\DYmu$ is projective.
  \end{proof}

  We can now show that the Cartan pairing is biadjoint with respect to
  $\iInd$ and $\iRes$, for $i\in I$.

  \begin{Theorem}\label{T:Biadjointness}
      Let  $u\in\ProjN$, $v\in\RepN$,and $i\in I$. Then
      \[
        \CPair[\big]{\iInd u}{v}=\CPair[\big]{u}{\iRes v} \And
        \CPair[\big]{\iRes u}{v}=\CPair[\big]{u}{\iInd v}.
      \]
  \end{Theorem}

  \begin{proof}Let $\Domin$. Since $\Ddec$ is surjective, we can write
    $v=\Ddec(\dot v)$ where $\dot v\in L_\A(\Lambda)$ and $\wt(\dot
    v)=\Lambda-\alpha$. Then $\CPair{\iRes u}{v}=0$ unless
    $\wt(u)=\Lambda-\alpha+\alpha_i$, in which case we compute
    \begin{align*}
      \CPair[\big]{\iRes u}{v} &=\CPair[\big]{E_i u}{\Ddec(\dot v)} \\
        &=q^{-\defect(\alpha)}\DPair[\Big]{\Ddec_T\bigl((E_iu)^\#\bigr)}{\dot v}
            &&\text{by \autoref{L:FormComparision}},\\
        &=q^{\defect(\alpha)}\DPair[\Big]{E_i\Ddec_T(u^\circledast)}{\dot v},
            &&\text{by \autoref{L:SameDuals} and \autoref{P:CommutingDuals}},\\
        &=q^{\defect(\alpha)}\DPair[\big]{\Ddec_T(u^\circledast)}{F_i\dot v},
            &&\text{by \autoref{L:ContrBiadjoint}},\\
        &=q^{-\defect(\alpha)-2\defect(\alpha-\alpha_i)}
             \DPair[\big]{\Ddec_T(u^\#)}{F_i\dot v},
            &&\text{by \autoref{L:SameDuals}},\\
        &=q^{-\defect(\alpha-\alpha_i)}\CPair[\big]{u}{F_iv},
            &&\text{by \autoref{L:FormComparision}},\\
        &=\CPair[\big]{u}{\iInd v},
    \end{align*}
    where the last equality uses \autoref{L:defectRecursive} and the
    identifications of~$F_i$ and~$\iInd\circ q^{-\di}K^{-1}_i$ from
    \autoref{P:Fock}.  A similar calculation
    shows that $\CPair[\big]{u}{\iRes v}=\CPair[\big]{\iInd u}{v}$.
  \end{proof}

  \begin{Remark}
    Working over a positively graded ring, Kashiwara
    \cite[Theorem~3.5]{Kashiwara:KLRBiadjointness} shows that
    $(\iRes,\iInd)$ is a biadjoint pair, which implies
    \autoref{T:Biadjointness}.  \autoref{L:FormComparision} can be
    interpreted as saying that the Cartan pairing categorifies the
    Shapovalov form; compare \cite[Lemma~3.1 and Theorem~4.18(4)]{BK:GradedDecomp}.
  \end{Remark}

  The modules $\LLa$ and $\GLa$ are \emph{standard $\A$-forms} of the
  irreducible $\Uq$-module $L(\Lambda)$. The corresponding
  \emph{costandard} $\A$-forms of $L(\Lambda)$ are the dual lattices:
  \begin{align*}
    \LLa^*&=\set[\big]{v\in \LLa|\CPair uv\in\A\text{ for all } u\in \LLa}
    \\
    \GLa^*&=\set[\big]{v\in \GLa|\CPair uv\in\A\text{ for all } v\in \GLa}
  \end{align*}
  By \autoref{L:FormComparision},
  $\DLa^*=\set{v\in \AA\otimes_\A\DLa|\DPair uv\in\A\text{ for all }
                u\in \DLa}$.
  \notation{$\LLa^*,\GLa^*$}{Dual highest weight modules as submodules of
  $\LFock$ and $\GFock$}

  We can now prove the main result of this section. Categorical analogues of
  this result have been obtained by Brundan and
  Kleshchev~\cite[Theorem~4.18]{BK:GradedDecomp} in type $\Aone$ and
  Kang and Kashiwara~\cite[Theorem~6.2]{KangKashiwara:CatCycKLR} for all
  symmetrisable Kac-Moody algebras. The following theorem
  provides an explicit bridge between the graded representation theory
  of $\Rx(\K[x])$ and the representation theory of $\UA$, which will be
  exploited in the following sections.

  \begin{Theorem}[Cyclotomic categorification]\label{T:CyclotomicCat}
    Suppose that $\Lambda\in Q^+$. Then, as $\UA$-modules,
    \[
       \LLa \cong\ProjN\cong\GLa
       \And
       \LLa^*\cong\RepN\cong\GLa^*.
    \]
  \end{Theorem}

  \begin{proof}
    The two isomorphisms for $\ProjN$ were already noted in
    \autoref{C:ProjK}. Let $\Domin$. Using the fact that $\DLa\cong\ProjN$,
    together with \autoref{E:CartanPairing} and
    \autoref{T:Biadjointness}, shows that $\DLa^*\cong\RepN$ as
    $\UA$-modules.
  \end{proof}

  In particular, note that \autoref{T:CyclotomicCat} implies that the
  sets $\LKlesh$ and $\GKlesh$ are independent of the field~$\K$. (In
  fact, this already follows from \autoref{P:Fock}.) We will soon give
  recursive descriptions of these sets.

  \subsection{Canonical bases}\label{S:CanonicalBases}
  A key feature of integrable highest weight modules is that they come
  equipped with the closely related canonical bases and crystal
  bases. This section connects the natural bases of~$\ProjN$ and~$\RepN$
  with canonical bases of $L_\A(\Lambda)$ and $L_\A(\Lambda)^*$.

  \begin{Lemma}\label{L:EiFiBarCommute}
    Let $i\in I$. Then $E_i\circ\circledast\cong\circledast\circ E_i$
    and $F_i\circ\circledast\cong\circledast\circ F_i$ as functors on
    $\Rep_\K\Rx[\bullet](\K[x])$.
  \end{Lemma}

  \begin{proof}
    By \autoref{P:CommutingDuals}, $\iRes$ commutes
    with $\circledast$ as functors on $\Rep_\K\Rx[\bullet](\K[x])$.
    Therefore, it is enough to show that
    $F_i\circ\circledast\cong\circledast\circ F_i$ as functors on
    $\Rep_\K\Rx[\alpha](\K[x])$, for $\alpha\in Q^+$.  As in
    \autoref{P:Fock}, identify~$F_i$ with the functor $\iInd\circ
    q^{\di}K_i^{-1}=q^{-\di}K_i\circ\iInd$ on $\Rep_\K\Rx[\alpha](\K[x])$.
    Then there are isomorphisms
    \begin{align*}
      F_i\circ\circledast
      &\cong q^{\di}\iInd K^{-1}_i\circ q^{-2\defect\alpha}\#
           && \text{by \autoref{L:SameDuals}},\\
        &\cong q^{\di-d_i(\alpha)-2\defect\alpha}\iInd\circ\#
           &&\text{where $d_i(\alpha)=(\Lambda-\alpha|\alpha_i)$,}\\
        &\cong q^{\di-d_i(\alpha)-2\defect\alpha}\#\circ\iInd
           && \text{by \autoref{P:CommutingDuals}},\\
        &\cong q^{-2\defect(\alpha+\alpha_i)}\#\circ q^{d_i(\alpha)-\di}
             \circ\iInd
           &&\text{by \autoref{L:defect}},\\
        &\cong\circledast\circ q^{-\di}K_i\iInd\cong\circledast\circ F_i,
           &&\text{by \autoref{L:SameDuals}}.
    \end{align*}
    So, $E_i$ and $F_i$ commute with $\circledast$
    when acting on~$\Rep_\K\Rx[\bullet](\K[x])$ (and as functors
    on~$\Proj_\K\Rx[\bullet](\K[x])$).
  \end{proof}

  In contrast, $E_i$ and $F_i$ do \textit{not} commute with $\#$ --- and
  nor do the functors $\iInd$ and $\circledast$.

  The functors $\#$ and $\circledast$ of \autoref{E:HashDual} and
  \autoref{E:Dual}, respectively, induce semilinear automorphisms
  of $\ProjN$ and~$\RepN$, which are given by:
  \[
     [P]^\#=[P^\#], \And
     [M]^\circledast=[M^\circledast]
  \]
  for $M\in\Rep_\K\Rx(\K[x])$ and $P\in\Proj_\K\Rx(\K[x])$.  \autoref{L:SameDuals} shows
  that these automorphisms are closely related.  By restriction, we
  consider $\circledast$ as a semilinear automorphism of $\ProjN$.

  The \emph{bar involution} on $\overline{\space}\map\Uq\Uq$ is the
  unique semilinear involution such that
  \[
        \overline{E_i}=E_i,\qquad
        \overline{F_i}=F_i \And
        \overline{K_i}=K_i^{-1},
        \qquad\text{for all }i\in I.
  \]
  Recall that $\Lambda\in P^+$ is a dominant weight and that
  $L(\Lambda)=\Uq v_\Lambda$ is an integrable highest weight module,
  where~$v_\Lambda$ a highest weight vector of weight $\Lambda$. The bar
  involution of~$\Uq$ induces a unique semilinear bar
  involution~$\overline{\phantom{m}}$ on~$L(\Lambda)$ such that
  $\overline{v_\Lambda}=v_\Lambda$ and
  $\overline{av}=\overline{a}\,\overline{v}$, for all $a\in\Uq$ and
  $v\in L(\Lambda)$.  \notation{$\overline{v}$}{Bar involution applied
  to an element $v$ of an integrable $\Uq$-module}

  \begin{Corollary}\label{C:BarStar}
      Let $u\in\LLa$, $v\in\GLa$ and $p\in\ProjN$. Then
      \[
        \Ldec(u)^\circledast= \Ldec(\overline{u}),\quad
        \Gdec(v)^\circledast= \Gdec(\overline{v}),\quad
        \Ldec_T(p^\#)=q^{2\defect(\alpha)}\overline{\Ldec_T(p)}\And*
        \Gdec_T(p^\#)=q^{2\defect(\alpha)}\overline{\Gdec_T(p)}.
      \]
  \end{Corollary}

  \begin{proof}
    Let $\Domin$. Since
    $\overline{\Dfock[\zero]}=\Dfock[\zero]={\Dfock[\zero]}^\circledast$
    is the highest weight vector in~$\DLa$, arguing by induction on weight using
    \autoref{L:EiFiBarCommute},  it
    follows that $\Ddec(\overline{f})=\bigl(\Ddec(f)\bigr)^\circledast$,
    for all $f\in\DLa$. As $\ProjN$ embeds into $\RepN$,
    $\Ddec_T(p^\circledast)=\overline{\Ddec_T(p)}$, for all
    $p\in\ProjN$. Hence,
    $\Ddec_T(p^\#)=q^{2\defect{\alpha}}\overline{\Ddec_T(p)}$ since
    $\#\cong q^{2\defect(\alpha)}\circ\circledast$ by
    \autoref{L:SameDuals}.
  \end{proof}

  That is, $\circledast$ categorifies the bar involution on the Fock space.

  \begin{Remark}
    The Fock spaces $\LFock$ and $\GFock$ are both integrable highest
    weight modules. Hence, both Fock spaces come equipped with bar
    involutions that are unique up to a choice of scalars, corresponding
    to the choice of highest weight vectors. Motivated by
    \autoref{P:SpechtDual}, let $\TT\map{\LFock}{\GFock}$  be the unique
    linear map such that
    $\TT(\Lfock)=q^{\defect\blam}\overline{\Gfock}$, for
    $\blam\in\Parts*$. Then \autoref{C:FockSpace}, \autoref{P:Fock} and
    \autoref{L:defect} imply that $\TT$ is a $\Uq$-module isomorphism
    and that
    $\TT\circ\overline{\phantom{m}}=\overline{\phantom{m}}\circ\TT$.
    Similarly, the map $\TT'\colon\GFock\to\LFock$, which sends $\Gfock$
    to $q^{\defect\blam}\overline{\Lfock}$ for $\blam\in\Parts*$, is a
    $\Uq$-module isomorphism and
    $\TT'\circ\overline{\phantom{m}}=\overline{\phantom{m}}\circ\TT'$.
    Moreover, $\TT\circ\TT'$ and $\TT'\circ\TT$ are both identity maps.
    We will not use these observations in what follows, except
    implicitly in the sense that, as this remark suggests, working with
    the two Fock spaces, $\LFock$ and $\GFock$, serves as a replacement
    for giving an explicit description of the bar involution on either Fock
    space.
  \end{Remark}

  \begin{Lemma}\label{L:Shapovalov}
    Suppose that $P\in\Proj\Rn(F)$ and $M\in\Rep\Rn<m>(F)$. Then
    \[
       \CPair[\big]{[P]}{[M]^\circledast}
          =\overline{\CPair[\big]{[P]^\#}{[M]}}.
    \]
  \end{Lemma}

  \begin{proof}
    This is a standard tensor-hom adjointness argument; see, for
    example, \cite[Lemma~2.5]{BK:GradedDecomp}.
  \end{proof}

  By \autoref{E:CartanPairing}, with respect to the Cartan pairing, the
  bases $\set{[\LYmu]|\bmu\in\LKlesh[\bullet]}$ and
  $\set{[\GYnu]|\bnu\in\GKlesh[\bullet]}$ of $\ProjN$ are dual to the
  bases $\set{[\LDmu]|\bmu\in\LKlesh[\bullet]}$ and
  $\set{[\GDnu]|\bnu\in\GKlesh[\bullet]}$ of $\RepN$, respectively. The
  projective Grothendieck group $\ProjN$ comes equipped with only one
  natural basis $\set{[\DYmu]|\bmu\in\DKlesh[\bullet]}$. In contrast,
  the Grothendieck group $\RepN$ has two quite different bases,
  $\set{[\DDmu]|\bmu\in\DKlesh[\bullet]}$ and
  $\set{[\DSlam[\bmu]]|\bmu\in\DKlesh[\bullet]}$, given by the simple
  modules and the Specht modules. To define  a second basis of $\ProjK$,
  which turns out to be dual to the dual Specht modules,
  define the \emph{inverse graded decomposition numbers} to be the
  Laurent polynomials $\Lelammu,\Gelamnu[\bsig\bnu]\in\A$ given by
  \begin{equation}\label{E:InverseDecomp}
    \bigl(\Lelammu\bigr) = \bigl(\Ldlammu\bigr)^{-1} \And
    \bigl(\Gelamnu[\bsig\bnu]\bigr) = \bigl(\Gdlamnu[\bsig\bnu]\bigr)^{-1},
  \end{equation}
  where $\blam,\bmu\in\LKlesh$, $\bnu,\bsig\in\GKlesh$ and the rows and
  columns of these matrices are ordered using the lexicographic orders
  $\Llex$ and $\Glex$, respectively. These polynomials are well-defined
  because these submatrices of the decomposition matrices of $\Rx(\K[x])$
  are lower unitriangular square matrices by \autoref{T:TriangularDecomp}.
  For $\bmu\in\LKlesh$ and $\bnu\in\GKlesh$ define virtual projective
  modules by
  \begin{equation}\label{E:Xmu}
      \LXmu = \sum_{\blam\Ledom\bmu}{\Lelammu}\,[\LYmu[\blam]]
      \And
      \GXnu = \sum_{\bsig\Gedom\bnu}{\Gelamnu[\bsig\bnu]}
               \,[\GYnu[\bsig]],
  \end{equation}
  where $\blam\in\LKlesh$ and $\bsig\in\GKlesh$ in the sums.  As the
  matrices in \autoref{E:InverseDecomp} are invertible,
  $\bigcup_{n\ge0}\set{\LXmu|\bmu\in\LKlesh}$ and
  $\bigcup_{n\ge0}\set{\GXnu|\bnu\in\GKlesh}$ are both $\A$-bases
  of~$\ProjN$. The definition of the $\DXmu[]$-bases suggests that these
  elements depend on~$\K$ but the next result shows that these elements
  are independent of~$\K$.
  \notation{$\Lelammu,\Gelamnu$}{Entries of the inverse graded decomposition matrices}[E:InverseDecomp]
  \notation{$\LXmu,\GXnu$}{Fake projective modules, which give bases of $\ProjN$}[E:Xmu]

  \begin{Lemma}\label{L:Xmu}
    Suppose that $\bmu,\blam\in\LKlesh$ and $\bnu,\bsig\in\GKlesh$. Then
    $\CPair[\big]{\LXmu[\bmu]}{[\LSlam[\blam]]^\circledast}
         =\delta_{\blam\bmu}$
    and
    $\CPair[\big]{\GXnu}{[\GSlam[\bsig]]^\circledast}=\delta_{\bnu\bsig}$.
  \end{Lemma}

  \begin{proof}
    It is enough to prove the first statement as the second follows by
    symmetry. By the definitions,
    \begin{align*}
        \CPair[\big]{\LXmu}{[\LSlam[\bsig]]^\circledast}
            &=\CPair[\Big]{ \sum_{\blam\Ledom\bmu}{\Lelammu}
                 \,[\LYmu[\blam]]}{[\LSlam[\bsig]]^\circledast}
             =\sum_{\blam\Ledom\bmu}\overline{\Lelammu}
                  \CPair[\big]{[\LYmu[\blam]]}{[\LSlam[\bsig]]^\circledast}\\
            &=\sum_{\blam\Ledom\bmu}\overline{\Lelammu}\, \CPair[\Big]{[\LYmu[\blam]]}
            {\sum_{\btau\Gedom\bsig}\overline{\Ldlammu[\bsig\btau]}[\LDmu[\btau]]^\circledast}\\
            &=\sum_{\substack{\btau\Gedom\bsig\\\blam\Ledom\bmu}}
                \overline{\Ldlammu[\bsig\btau]}\,\overline{\Lelammu}
                \CPair[\big]{[\LYmu[\blam]]}{[\LDmu[\btau]]}\\
            &=\sum_{\bsig\Ledom\blam\Ledom\bmu}
                \overline{\Ldlammu[\bsig\blam]}\,\overline{\Lelammu},
    \end{align*}
    where the last equality follows by \autoref{E:CParingYmuDnu}. Note
    that in these sums, $\blam,\btau\in\DKlesh$. The result now follows
    by \autoref{E:InverseDecomp}.
  \end{proof}

  Applying the two bar involutions $\#$ and $\circledast$ shows that
  if $\Domin$ then
  \begin{equation}\label{E:DYBarInvariance}
       [\DYmu]^\# = [\DYmu]  \quad \text{and}\quad
       [\DDmu]^\circledast = [\DDmu],
       \qquad\text{for }\bmu\in\DKlesh,
  \end{equation}
  with the $\#$-identities following because $\LYmu$ and $\GYnu$ are
  projective and the $\circledast$-identities coming from
  \autoref{T:GradedSimples}. It is less clear what these involutions
  do to the other bases of $\ProjN$ and $\RepN$.

  \begin{Lemma}\label{L:BarTriangular}
    Let $\bmu\in\LKlesh$ and $\bnu\in\GKlesh$. Then
    \begin{align*}
       \bigl(\LXmu\bigr)^\# & = \LXmu
        + \sum_{\blam\Ldom\bmu}x^{\Ldom}_{\blam\bmu}(q)\LXmu[\blam],&
       [\LSlam[\bmu]]^\circledast & = [\LSlam[\bmu]]
        + \sum_{\blam\Ldom\bmu}s^{\Ldom}_{\bmu\blam}(q)[\LSlam],\\
        \bigl(\GXnu\bigr)^\# & = \GXnu
        + \sum_{\bsig\Gdom\bnu}x^{\Gdom}_{\bsig\bnu}(q)\GXnu[\bsig],&
       [\GSlam[\bnu]]^\circledast & = [\GSlam[\bnu]]
        + \sum_{\bsig\Gdom\bnu}s^{\Gdom}_{\bnu\bsig}(q)[\GSlam[\bsig]].
    \end{align*}
    for Laurent polynomials
    $x^{\Ldom}_{\blam\bmu}(q), s^{\Ldom}_{\blam\bmu}(q),
     x^{\Gdom}_{\bsig\bnu}(q), s^{\Gdom}_{\bsig\bmu}(q)\in\A$
    with $\blam\in\LKlesh$ and $\bsig\in\GKlesh$.
  \end{Lemma}

  \begin{proof}
    Let $\bsig\in\LKlesh$. Using
    \autoref{T:TriangularDecomp} and \autoref{E:InverseDecomp},
    \begin{align*}
      [\LSlam[\bmu]]^\circledast
        & =\Bigl(\sum_{\balp\Ledom\bmu}\Ldlammu[\bmu\balp]
              [\LDmu[\balp]]\Bigr)^\circledast
          =\sum_{\balp\Ledom\bmu}\overline{\Ldlammu[\bmu\balp]}\,
              [\LDmu[\balp]]\\
         &=\sum_{\balp\Ledom\bmu}\overline{\Ldlammu[\bmu\balp]}\,
              \sum_{\blam\Ledom\balp}\Lelammu[\balp\blam][\LSlam]\\
         &=[\LSlam[\bmu]] + \sum_{\blam\Ldom\bmu}\Bigr(
              \sum_{\substack{\balp\in\LKlesh\\\blam\Ledom\balp\Ledom\bmu}}
              \overline{\Ldlammu[\bmu\balp]}\Lelammu[\balp\blam]\Bigr)[\LSlam],
    \end{align*}
    where the last equality follows because
    $\Ldlammu[\bmu\bmu]=1=\Lelammu[\bmu\bmu]$ by
    \autoref{T:TriangularDecomp}. This proves the result
    for $[\LSlam[\bmu]]^\circledast$, which this implies that
    $\LXmu{}^\#$ has the required expansion by \autoref{L:Xmu} and
    \autoref{L:Shapovalov}. The remaining claims are similar.
  \end{proof}

  \notation{$\LGmu,\GGnu$}{$\#$-canonical basis vectors in $\ProjN$}[T:CanonicalBasis]
  \notation{$\LGlammu,\GGlamnu$}{Transition matrices between the $\set{[\DXmu]}$ and $\set{\DGmu}$ bases}[T:CanonicalBasis]
  \notation{$\LHmu,\GHnu$}{$\circledast$-canonical basis vectors in $\RepN$}[T:CanonicalBasis]
  \notation{$\LHlammu,\GHlamnu$}{Transition matrices between the $\set{[\DSlam]}$ and $\set{\DHmu}$ bases}[T:CanonicalBasis]

  \begin{Theorem}\label{T:CanonicalBasis}
    Let $\bmu\in\LKlesh$ and $\bnu\in\GKlesh$. Then there exist
    bases $\set{\LGmu|\bmu\in\LKlesh}$ and $\set{\GGnu|\bnu\in\GKlesh}$
    of $\ProjN$, and $\set{\LHmu|\bmu\in\LKlesh}$ and
    $\set{\GHnu|\bnu\in\GKlesh}$ of $\RepN$, that are uniquely
    determined by the conditions:
    \begin{align*}
      \bigl(\LGmu\bigr)^\# &= \LGmu \And*
      \LGmu=\LXmu +\sum_{\blam\Ldom\bmu}\LGlammu\LXmu[\blam]\\
      \bigl(\GGnu\bigr)^\# &= \GGnu \And*
      \GGnu[\bmu]=\GXnu[\bmu] +\sum_{\blam\Gdom\bnu}\GGlamnu\GXnu[\blam]\\
      \bigl(\LHmu\bigr)^\circledast &= \LHmu \And*
      \LHmu=[\LSlam[\bmu]] +\sum_{\blam\Ldom\bmu}\LHlammu{} [\LSlam]\\
      \bigl(\GHnu\bigr)^\circledast &=\GHnu\And*
      \GHnu=[\GSlam[\bnu]]+\sum_{\blam\Gdom\bnu}\GHlamnu{}[\GSlam].
    \end{align*}
    for polynomials $\LGlammu,\LHlammu\in\delta_{\blam\bmu}+q\Z[q]$ and
    $\GGlamnu[\blam\bnu],\GHlamnu[\blam\bnu]\in\delta_{\blam\bnu}+q\Z[q]$,
    for $\bmu\in\LKlesh$ and $\bnu\in\GKlesh$.
    In particular, the basis elements  $\LGmu$, $\GGnu$,
    $\LHmu$ and $\GHnu$, and these polynomials, are independent of
    the field~$\K$.
  \end{Theorem}

  \begin{proof}
    Given \autoref{L:BarTriangular}, this result is a consequence of
    \textit{Lusztig's Lemma}~\cite[Lemma 24.2.1]{Lusztig:QuantBook},
    which is easily proved by induction on dominance using Gaussian
    elimination and \autoref{L:BarTriangular}. See
    \cite[Proposition~3.5.6]{Mathas:Singapore} for a proof that uses
    very similar language to that used here.
  \end{proof}

  A key point in \autoref{T:CanonicalBasis} is that the coefficients
  appearing in \autoref{L:BarTriangular} belong to~$\A$.  As the
  notation suggests, the polynomials $\DGlammu$ are related to the
  decomposition matrices of $\Rx(\K[x])$ and the polynomials~$\DHlammu$
  are related to the inverse decomposition matrices. See
  \autoref{T:CanonicalToDY} below for a precise statement.

  By \autoref{T:CanonicalBasis}, $\set{\LGmu|\bmu\in\LKlesh[\bullet]}$ and
  $\set{\GGnu|\bnu\in\GKlesh[\bullet]}$ are bases of $\ProjN$ and
  $\set{\LHmu|\bmu\in\LKlesh[\bullet]}$ and $\set{\GHnu|\bnu\in\GKlesh[\bullet]}$ are
  bases of $\RepN$.

  \begin{Definition}\hfill
    \begin{enumerate}
      \item The $\circledast$-\emph{canonical bases} of $\RepN$ are the
      two bases $\set{\LHmu|\bmu\in\LKlesh[\bullet]}$ and
      $\set{\GHnu|\bnu\in\GKlesh[\bullet]}$.
      \item The $\#$-\emph{canonical bases} of $\ProjN$ are the two
      bases $\set{\LGmu|\bmu\in\LKlesh[\bullet]}$ and
      $\set{\GGnu|\bnu\in\GKlesh[\bullet]}$.
    \end{enumerate}
  \end{Definition}

  We frequently call these four bases \textit{canonical bases} of $\RepN$ and
  $\ProjN$.  In \autoref{T:KleshchevSimples} below we show that, up to scaling,
  these bases coincide with Lusztig's (dual) canonical
  bases~\cite[\S14.4]{L:can} and Kashiwara's (upper and lower) global
  bases~\cite{Kashiwara:CrystalBases} of $L(\Lambda)$.

  For now we note that \autoref{T:CanonicalBasis} and
  \autoref{L:SameDuals} imply:

  \begin{Corollary}\label{C:CanonicalBases}
    Suppose that $\bmu\in\LKlesh$ and $\bnu\in\GKlesh$. Then
    \[
      \bigl(\LGmu\bigr)^\circledast = q^{-2\defect\bmu}\LGmu,\quad
      \bigl(\GGnu\bigr)^\circledast = q^{-2\defect\bnu}\GGnu,\quad
      \bigl(\LHmu\bigr)^\# = q^{2\defect\bmu}\LHmu\quad\text{and}\quad
      \bigl(\GHnu\bigr)^\# =q^{2\defect\bnu}\GHnu.
    \]
  \end{Corollary}

  The next result shows that these bases of $\ProjK$ and $\RepK$ are
  dual with respect to the Cartan pairing.  The matrix identities in the
  next result should be compared with \autoref{E:InverseDecomp}.

  \begin{Corollary}\label{C:CanonicalDual}
    Suppose that $\blam,\bmu\in\LKlesh$ and $\bnu,\bsig\in\GKlesh$. Then
    $\CPair{\LGmu[\blam]}{\LHmu}=\delta_{\blam\bmu}$ and
    $\CPair{\GGnu}{\GHnu[\bsig]}=\delta_{\bnu\bsig}$.
    Equivalently, the two matrix identities hold
    \[
      \bigl(\LHlammu[\blam\bmu]\bigr) = \bigl(\LGlammu\bigr)^{-1}
      \And
      \bigl(\GHlamnu[\bsig\bnu]\bigr) = \bigl(\GGlamnu[\bsig\bnu]\bigr)^{-1}.
    \]
  \end{Corollary}

  \begin{proof}
    Let $\Domin$. Let $\balp,\bbet\in\DKlesh$. Direct calculation reveals
    that
    \begin{align*}
      \CPair[\big]{\DGmu[\balp]}{\DHmu[\bbet]}
         &=\CPair[\big]{[\DGmu[\balp]]}{[\DHmu[\bbet]]^\circledast}
          =\CPair[\Big]{\sum_{\bsig\in\DKlesh}\DGlammu[\bsig\balp]\DXmu[\bsig]}
               {\sum_{\btau\in\DKlesh}\overline{\DHlammu[\bbet\btau]}
                              [\DSlam[\btau]]}\\
         &=\sum_{\bsig,\btau\in\DKlesh}
         \overline{\DGlammu[\bsig\balp]}\,\overline{\DHlammu[\bbet\btau]}
          \CPair[\big]{\DXmu[\bsig]}{[\DSlam[\btau]]^\circledast}\\
          &=\sum_{\bsig\in\DKlesh}
               \overline{\DHlammu[\bbet\bsig]\DGlammu[\bsig\balp]},
    \end{align*}
    where the last equality follows by \autoref{L:Xmu}. Therefore,
    $\CPair{\DGmu[\balp]}{\DHmu[\bbet]}\in\delta_{\balp\bbet}+q^{-1}\Z[q^{-1}]$.
    However, by \autoref{L:Shapovalov},
    \[
        \CPair{\DGmu[\balp]}{\DHmu[\bbet]}
        = \overline{\CPair{\DGmu[\balp]{}^\#}{\DHmu[\bbet]{}^\circledast}}
        = \overline{\CPair{\DGmu[\balp]}{\DHmu[\bbet]}}
        \in\delta_{\balp\bbet}+qZ[q].
    \]
    Hence, $\CPair{\DGmu[\balp]}{\DHmu[\bbet]}=\delta_{\balp\bbet}$. The
    calculation in the first displayed equation shows that this is
    equivalent to the matrix identity in the statement of the corollary.
  \end{proof}

  In particular, this shows that the $\#$-canonical bases of $\ProjN$ and the
  $\circledast$-canonical bases of $\RepN$ encode equivalent information.

  \begin{Lemma}\label{L:CanonicaBasisCoeffs}
    Let $\blam, \bmu\in\LKlesh$ and $\bsig, \bnu\in\GKlesh$. Then
    \[
        \LGlammu=\CPair{\LGmu}{[\LSlam]},\quad
        \GGlamnu=\CPair{\GGnu}{[\GSlam[\bsig]]},\quad
        \LHlammu=\CPair{\LXmu[\blam]}{\LHmu}
        \quad\text{and}\quad
        \GHlamnu=\CPair{\GXnu[\bsig]}{\GHnu}.
    \]
  \end{Lemma}

  \begin{proof}
    Let $\Domin$ and $\bmu\in\DKlesh$ and $\blam\in\Parts$. Using
    \autoref{L:FormComparision} and \autoref{T:CanonicalBasis},
    \[
      \CPair{\DGmu}{[\DSlam]}
        =\overline{\CPair{\DGmu{}^\#}{[\DSlam]^\circledast}}
        =\overline{\CPair{\DGmu}{[\DSlam]^\circledast}}
        = \sum_{\btau\in\DKlesh}\DGlammu[\btau\bmu]
        \overline{\CPair[\big]{\DXmu[\btau]}{[\DSlam]^\circledast}}
        =\DGlammu[\blam\bmu],
    \]
    where the last equality comes from \autoref{L:Xmu}. The proof of the
    other identities are similar.
  \end{proof}

  For $\bmu\in\LKlesh$, $\bnu\in\GKlesh$ and $\blam,\bsig\in\Parts$
  define Laurent polynomials
  \begin{equation}\label{E:CanonicalBasisCoeffs}
    \LGlammu=\CPair{\LGmu}{[\LSlam]}
    \qquad\text{and}\qquad
    \GGlamnu[\bsig\bnu]=\CPair{\GGnu}{[\GSlam[\bsig]]}.
  \end{equation}
  By \autoref{L:CanonicaBasisCoeffs}, if $\blam,\bmu\in\DKlesh$ then
  $\DGlammu$ coincides with the polynomial defined in
  \autoref{T:CanonicalBasis}. In particular, if $\blam\in\DKlesh$ then
  $\DGlammu\in\delta_{\blam\bmu}+q\Z[q]$ by \autoref{T:CanonicalBasis}.
  We will show in \autoref{C:FullDecompDegs} below that this is still
  true when $\blam\in\Parts\setminus\DKlesh$. Moreover, we show that
  $\DGlammu\in\delta_{\blam\bmu}+q\N[q]$ in type~$\Aone$.

  \notation{$\Lglammu,\Gglamnu$}{Transition matrices between the $\set{[\DGmu]}$ and $\set{\DYmu}$ bases}[T:CanonicalToDY]
  \notation{$\Lhlammu,\Ghlamnu$}{Transition matrices between the $\set{[\DHmu]}$ and $\set{\DDmu}$ bases}[T:CanonicalToDY]

  \begin{Theorem}\label{T:CanonicalToDY}
     For $\bmu,\blam\in\LKlesh$ and $\bnu,\bsig\in\GKlesh$, there exist
     bar invariant polynomials
     $\Lglammu,\Gglamnu[\bsig\bnu],\Lhlammu,\Ghlamnu[\bnu\bsig]\in\A$
     such that
     \begin{align*}
        [\LYmu] &= \LGmu + \sum_{\blam\Ldom\bmu}\Lglammu\LGmu[\blam],&
        [\GYnu] &= \GGnu
           + \sum_{\bsig\Gdom\bnu}\Gglamnu[\bsig\bnu]\GGnu[\bsig],\\
        [\LDmu] &= \LHmu + \sum_{\blam\Gdom\bmu}\Lhlammu[\bmu\blam]\LHmu[\blam],&
        [\GDnu] &= \GHnu
           + \sum_{\bsig\Ldom\bnu}\Ghlamnu[\bnu\bsig]\GHnu[\bsig].
     \end{align*}
     Moreover,
     for $\bsig,\blam\in\Parts$, the following matrix identities hold:
     \begin{align*}
        \bigl(\Lhlammu\bigr) &= \bigl(\Lglammu\bigr)^{-1}, &
        \bigl(\Ghlamnu\bigr) &= \bigl(\Gglamnu\bigr)^{-1}, \\
        \bigl(\Ldlammu\bigr) &= \bigl(\LGlammu\bigr)\bigl(\Lglammu\bigr),&
        \bigl(\Gdlamnu[\bsig\bnu]\bigr)
           &= \bigl(\GGlamnu[\bsig\bnu]\bigr)\bigl(\Lglammu\bigr).
     \end{align*}

  \end{Theorem}

  \begin{proof}Let $\Domin$.
    By \autoref{E:DYBarInvariance}, $[\DYmu]$ is a $\#$-invariant
    element of $\ProjN$ and $[\DDmu]$ is a $\circledast$-invariant
    element of $\RepN$. Hence, the first four identities follow by
    \autoref{E:InverseDecomp} and \autoref{L:BarTriangular}. (These
    four identities describe the transition matrices between the
    $\set{[\DYmu]}$ and $\set{[\DGmu]}$ bases and between the $\set{[\DDmu]}$
    and $\set{[\DHmu]}$ bases.) Since
    $\CPair{[\DYmu]}{[\DDmu[\bnu]]}=\delta_{\bmu\bnu}$, by
    \autoref{E:CartanPairing}, these transition matrices are inverse to
    each other by \autoref{C:CanonicalDual}.  Finally, if
    $\blam\in\Parts$ and $\bmu\in\DKlesh$ then
    \begin{align*}
      \Ddlammu &= \CPair[\big]{[\DYmu]}{[\DSlam]}
        =\CPair[\Big]{\sum_{\bnu\in\DKlesh}\Dglammu[\bnu\bmu]\DGmu[\bnu]}
           {[\DSlam]}
        =\sum_{\bnu\in\DKlesh}\Dglammu[\bnu\bmu]\CPair[\big]{\DGmu[\bnu]}
           {[\DSlam]}\\
        &=\sum_{\bnu\in\DKlesh}\DGlammu[\blam\bnu]\Dglammu[\bnu\bmu],
    \end{align*}
    where the third equality follows because
    $\overline{\Dglammu[\bnu\bmu]}=\Dglammu[\bnu\bmu]$ is bar invariant.
    This gives the required factorisation of the decomposition matrices
    $\Ddec$.
  \end{proof}

  As a consequence, we recover the Ariki-Brundan-Kleshchev
  categorification theorem.

  \begin{Corollary}[{%
    Brundan and Kleshchev~\cite[Theorem~5.3 and Corollary~5.15]{BK:GradedDecomp}}]
     \label{C:Ariki}
     Let $\Gamma$ be a quiver of type $\Aone$ and suppose that $\K$ is a
     field of characteristic~$0$. Then
     \[
        [\LYmu] = \LGmu,\qquad
        [\GYnu] = \GGnu,\qquad
        [\LDmu] = \LHmu \And
        [\GDnu] = \GHnu.
     \]
     for all $\bmu\in\LKlesh$ and all $\bnu\in\GKlesh$. Consequently,
     if $\blam\in\Parts$, $\bmu\in\LKlesh$ and $\bnu\in\GKlesh$ then
     \[
       \Ldlammu = \CPair{\LGmu}{[\LSlam]} \And
       \Gdlamnu = \CPair{\GGnu}{[\GSlam]}.
     \]
     In particular,
     $\Ldlammu=\LGlammu\in\delta_{\blam\bmu}+q\N[q]$ if $\blam\in\LKlesh$ and
     $\Gdlamnu=\GGlamnu\in\delta_{\blam\bnu}+q\N[q]$ if $\blam\in\GKlesh$.
   \end{Corollary}

  \begin{proof}
    Let $\Domin$. The algebras $\Rn(\K)\cong\Rx(\K)$ are cellular by
    \autoref{C:RnCellular}, so every field is a splitting field
    for~$\Rx(\K)$, so we can assume that $\K=\CC$. In type~$\Aone$, Brundan
    and Kleshchev~\cite{BK:GradedKL} proved that the cyclotomic KLR
    algebra $\Rn(\CC)$ is isomorphic to a (degenerate) Ariki-Koike
    algebra $\Hn(\CC)$.  Ariki~\cite[Theorem~4.4(2)]{Ariki:can}, and
    Brundan and Kleshchev~\cite[Theorem~3.10]{BK:DegenAriki} in the
    degenerate case, proved that the dual canonical basis of $\RepK[\CC]$
    at $q=1$ coincides with the basis of
    $[\Rep\Hn<\bullet>]=\bigoplus_{n\ge0}[\Rep\Hn(\CC)]$ given by the
    images of the irreducible $\Hn$-modules. Therefore, $\DHmu=[\DDmu]$,
    for $\bmu\in\DKlesh$, since the simple module $\DDmu$ is self-dual
    by \autoref{T:CellularSimples}. The remaining claims now follow in
    view of \autoref{T:CanonicalBasis} and \autoref{L:Xmu}.
  \end{proof}

  \begin{Example}\label{Ex:DecompMatrices}
    Given \autoref{C:Ariki}, in type $\Cone$ it is natural to ask if the
    $\circledast$-canonical bases of~$L_\A(\Lambda)^*$ coincide with the
    bases of simple modules, and the $\#$-canonical bases with the bases
    of principal indecomposable $\Rx(\K)$-modules when $\K$ is a field of
    characteristic zero. It is shown in
    \cite{ChungMathasSpeyer:TypeCDecomp} that this first fails for the
    principal block of~$\Rx<\Lambda_0>[8](\CC)$ when~$\Gamma$ is a
    quiver of type~$\Cone[2]$. Several other examples are given where
    the canonical bases do not coincide with the natural bases of these
    Grothendieck groups in type~$C$, including an example when $n=13$
    that shows that the graded decomposition numbers of $\Rx(\K[x])$ are
    not necessarily polynomials, even in characteristic zero.
  \end{Example}

  The transition matrices $\bigl(\Lglammu\bigr)$,
  $\bigl(\Gglamnu\bigr)$, $\bigl(\Lhlammu\bigr)$ and
  $\bigl(\Ghlamnu\bigr)$ in \autoref{T:CanonicalToDY} are analogues of
  the adjustment matrices of \autoref{D:Adjustment}. These matrices
  express the decomposition matrices of $\Rx(\K[x])$ in terms of the
  canonical bases and dual canonical bases. By taking inverses, similar
  ``adjustment matrix'' identities hold for the inverse decomposition
  matrices.

  Recall the Mullineux involution $\mull\map\LKlesh\GKlesh$ from
  \autoref{D:Mullineux}. The next result should be compared with
  \autoref{P:DecompComparision}.

  \begin{Proposition}\label{P:CanonicalMullineux}
    Let $\bmu\in\LKlesh$. Then $\LGmu=\GGnu[\mull(\bmu)]$ and
    $\LHmu=\GHnu[\mull(\bmu)]$. Moreover, if $\blam\in\Parts$ then
    \(
       \LGlammu[\blam\bmu]
         = q^{\defect\blam}\overline{\GGlamnu[\blam\mull(\bmu)]}.
    \)
  \end{Proposition}

  \begin{proof}
  By \autoref{D:Mullineux}, $[\LDmu]=[\GDnu[\mull(\bmu)]]$ and
  $[\LYmu]=[\GYnu[\mull(\bmu)]]$. Hence, $\LGmu=\GGnu[\mull(\bmu)]$ and
  $\LHmu=\GHnu[\mull(\bmu)]$ by \autoref{T:CanonicalToDY} and the
  uniqueness of the canonical basis elements established in
  \autoref{T:CanonicalBasis}. To prove the remaining claim, if
  $\bmu\in\LKlesh$ and $\blam\in\Parts$ then
  \begin{align*}
      \LGlammu
       &= \CPair[\big]{\LGmu[\bmu]}{[\LSlam]}
       = q^{\defect\blam}
         \CPair[\big]{\GGnu[\mull(\bmu)]}{[\GSlam]^\circledast}
       =q^{\defect\blam}\overline{%
         \CPair[\big]{\GGnu[\mull(\bmu)]}{[\GSlam]}}
         =q^{\defect\blam}\overline{\GGlamnu[{\mull(\bmu)\blam}]},
  \end{align*}
  where we have used \autoref{P:SpechtDual} and \autoref{L:Shapovalov}.
  \end{proof}

  Combining \autoref{T:CanonicalBasis} and
  \autoref{P:CanonicalMullineux}, we obtain.

  \begin{Corollary}\label{C:DecompDegs}
    Let $\bmu\in\LKlesh$, $\bnu\in\LKlesh$ and
    $\blam,\bsig\in\LKlesh\cup\GKlesh$.
    \begin{enumerate}
      \item If $\LGlammu\ne0$ then $\bmu\Ledom\blam\Ledom\mull(\bmu)$ and
      $\alpha_\blam=\alpha_\bmu$.
      Moreover, $\LGlammu[\bmu\bmu]=1$,
      $\LGlammu[\mull(\bmu)\bmu]=q^{\defect\bmu}$ and
      if~$\mull(\bmu)\Ldom\blam\Ldom\bmu$ then
      $0<\deg\LGlammu<\defect\bmu$.
      \item If $\GGlamnu\ne0$ then $\bmu\Gedom\blam\Gedom\mull(\bmu)$ and
      $\alpha_\blam=\alpha_\bmu$.
      Moreover, $\GGlamnu[\bmu\bmu]=1$,
      $\GGlamnu[\mull^{-1}(\bmu)\bmu]=q^{\defect\bmu}$ and
      if~$\mull^{-1}(\bmu)\Gdom\blam\Gdom\bmu$ then
      $0<\deg\GGlamnu<\defect\bmu$.
    \end{enumerate}
  \end{Corollary}

  \begin{proof}
    If $\blam,\bmu\in\LKlesh$ then
    $\LGlammu\in\delta_{\blam\bmu}+q\Z[q]$ by
    \autoref{T:CanonicalBasis}. Hence, the only claim in~(a) that is not
    immediate from \autoref{P:CanonicalMullineux} is that
    $0<\deg\LGlammu<\defect\bmu$ when $\blam\in\GKlesh$
    and~$\blam\notin\set{\bmu,\mull(\bmu)}$. In this case,
    $\GGlamnu[\blam\mull(\bmu)]\in\delta_{\blam\mull(\bmu)}+q\Z[q]$, so
    $0<\deg\LGlammu<\defect\bmu$ by \autoref{P:CanonicalMullineux}. This
    proves~(a). The proof of~(b) is similar.
  \end{proof}

  Later, we will show that this result is true for
  $\blam,\bsig\in\Parts$. There are similar identities for the
  polynomials $\LHlammu$ and $\GHlamnu$, which we leave for the reader.

  \begin{Corollary}\label{C:PositiveDefect}
    Let $\blam\in\Parts[\alpha]$, for $\alpha\in Q^+_n$. Then $\defect\alpha=\defect\blam\ge0$.
  \end{Corollary}

  \begin{proof}
    This is implicit in \autoref{C:DecompDegs} since $\LGlammu$ and
    $\GGlamnu$ are polynomials.
  \end{proof}

  \subsection{Crystal bases of Fock spaces}
  The categorification results of the last few sections imply that
  the number of self-dual graded simple modules is independent of the
  characteristic, but we have not yet determined the sets $\LKlesh$ and
  $\GKlesh$ that index the simple $\Rx(\K[x])$-modules. To do this we now
  describe the crystal graphs of~$\LLa$ and~$\GLa$. We start by
  recalling Kashiwara's theory of global and crystal bases and Lusztig's
  theory of canonical bases.


  Suppose that $V$ be an integrable highest weight module for $\Uq$.
  If $i\in I$ then~$E_i$ and~$F_i$ act on $V$ as locally nilpotent
  linear operators.  Therefore, by \cite[16.1.4]{Lusztig:QuantBook},
  each weight vector $v\in V$ can be written uniquely in the form
  \[
        v = \sum_{r\ge0} F_i^{(r)}v_r
  \]
  such that $E_iv_r=0$ and $K_iv_r=q^{\CPair{\wt(v_r)}{\alpha_i}+r\di}v_r$,
  for $r\ge0$. For $i\in I$, the \emph{Kashiwara operators} $e_i$ and
  $f_i$ are the linear endomorphisms of~$V$ defined by
  \begin{equation}\label{E:KashiwaraOps}
           e_i v = \sum_{r\ge1}F_i^{(r-1)}v_r\And
           f_i v = \sum_{r\ge0}F_i^{(r+1)}v_r.
  \end{equation}
  \notation{$e_i,f_i$}{Kashiwara's crystal operators, for $i\in I$}[E:KashiwaraOps]
  For $\bi\in I^n$ set $e_\bi=e_{i_n}\dots e_{i_2}e_{i_1}$ and
  $f_\bi=f_{i_n}\dots f_{i_2}f_{i_1}$.

  \notation{$\AA_0$}{Ring of rational functions regular at $0$}
  \notation{$\AA_\infty$}{Ring of rational functions regular at $\infty$}
  Let $\AA_0$ be the subring of rational functions $\AA=\Q(q)$ that are
  regular at zero and let $\AA_\infty$ be the rational function that are
  regular at infinity. To allow us to work with these two rings
  simultaneously, if $\omega\in\set{0,\infty}$ set
  \[
        q_\omega = \begin{cases*}
            q      & if $\omega=0$,\\
            q^{-1} & if $\omega=\infty$.
        \end{cases*}
  \]
  \notation{$q_\omega$}{Shorthand notation with $q_0=q$ and $q_\infty=q^{-1}$}

  \begin{Definition}[{Kashiwara~\cite[Definition 2.3.1]{Kashiwara:CrystalBases}}]
    \label{D:CrystalBasis}
    Let $V$ be an integrable $\Uq$-module. Fix $\omega\in\set{0,\infty}$.
    A \emph{$\omega$-crystal base}
    of~$V$ is a pair $(\L[\omega], \B[\omega])$ such that:
    \begin{enumerate}
      \item The module $\L[\omega]$ is a free $\AA_\omega$-submodule of~$V$ such that
      $V\cong\AA\otimes_{\AA_\omega}\L[\omega]$ and~$\L[\omega]$ is a direct sum of
      $\Uq$-weight spaces and it is invariant under the actions of $e_i$
      and $f_i$, for $i\in I$.
      \item The set $\B[\omega]$ is a basis of the $\Q$-vector space
      $\L[\omega]/q_\omega\L[\omega]=\<\B[\omega]\>_\Q$.
      \item The elements of $\B[\omega]$ are images of weight
      vectors under the map $\L[\omega]\to\L[\omega]/q_\omega\L[\omega]$.
      \item If $i\in I$ then $e_i\B[\omega]\subset\B[\omega]\cup\set{0}$ and
      $f_i\B[\omega]\subset\B[\omega]\cup\set{0}$.
      \item If $b,b'\in\B[\omega]$ and $i\in I$ then $e_ib=b'$ if and only if
      $f_ib'=b$.
    \end{enumerate}
  \end{Definition}

  This section describes the $0$-crystal base $(\L,\B)$ and the
  $\infty$-crystal base $(\L[\infty], \B[\infty])$ of~$L(\Lambda)$.

  If $V=\Uq v_\Lambda$ is an integrable highest weight module with
  highest weight vector~$v_\Lambda$ then, as in
  \autoref{S:CanonicalBases}, the bar involution on~$V$ is defined to be
  the unique semilinear automorphism such that
  $\overline{v_\Lambda}=v_\Lambda$ and
  $\overline{av}=\overline{a}\,\overline{v}$, for all $v\in V$ and $a\in
  \Uq$.

  \begin{Theorem}[{Lusztig \cite[\S14.4]{Lusztig:QuantBook},
    Kashiwara \cite{Kashiwara:CrystalBases}}]
    \label{T:KashiwaraLusztig}
    Let $V$ be an integrable $\Uq$-module. Fix $\omega\in\set{0,\infty}$
    and suppose that $(\L[\omega],\B[\omega])$ is an $\omega$-crystal
    basis for~$V$.  Then there exists a unique $\A$-basis
    $ \B[\omega](\Lambda)=\set{G_{\omega,b}|b\in\B[\omega](\Lambda)} $
    of $V_\A(\Lambda)$ such that $\overline{G_{\omega,b}}=G_{\omega,b}$
    and $G_{\omega,b}\equiv b\pmod{q_\omega\L[\omega](\Lambda)}$, for
    $b\in\B[\omega](\Lambda)$.
  \end{Theorem}

  The basis $\B(\Lambda)$  of~$V(\Lambda)$ is Lusztig's
  \emph{dual canonical basis}, or Kashiwara's \emph{lower global basis}
  and the basis $\B[\infty](\Lambda)$ is Lusztig's
  \emph{canonical basis}, or Kashiwara's \emph{upper global basis}.

  To apply these results to the combinatorial Fock spaces
  $\LLa$ and $\GLa$, and the Grothendieck groups $\ProjN$ and $\RepN$, we
  first generalise the integers $\Ldilam$ and $\Gdilam$ from
  \autoref{D:defect}.  If $\blam,\bmu\in\Parts*$ and $i\in I$ write
  $\blam\iarrow[i^r]\bmu$ if $|\bmu|=|\blam|+r$ and
  $\bmu=\blam\cup\set{A_1,\dots,A_r}$, where
  $\set{A_1,\dots,A_r}\subseteq\Add_i(\blam)$, and define
  \begin{align*}
    \Ldilam[\bmu] &= \di\sum_{s=1}^r\Bigl(\#\set[\big]{B\in\Add_i(\bmu)|B<A_s}
                    -\#\set[\big]{B\in\Rem_i(\blam)|B<A_s}\Bigr),\\
   \Gdilam[\bmu] &= \di\sum_{s=1}^r\Bigl(\#\set[\big]{B\in\Add_i(\blam)|B>A_s}
                    -\#\set[\big]{B\in\Rem_i(\blam)|B>A_s}\Bigr).
  \end{align*}
  By definition, if $\bmu=\blam\cup\set{A}$, for $A\in\Add_i(\blam)$,
  then $\Ldilam[\bmu]=\Ldilam$ and $\Gdilam[\bmu]=\Gdilam$.

  \begin{Lemma}\label{L:DividedPowers}
    Let $\blam\in\Parts$ and $i\in I$. Then, for $r\ge0$,
    \[
        F_i^{(r)}\Lfock = \sum_{\blam\iarrow[i^k]\bmu}
                        q^{-\Gdilam[\bmu]}\Lfock[\bmu]
        \And
        F_i^{(r)}\Gfock = \sum_{\blam\iarrow[i^k]\bmu}
                        q^{-\Ldilam[\bmu]}\Gfock[\bmu]
    \]
  \end{Lemma}

  \begin{proof}
    This follows easily by induction on $r$ using the fact that
    $F^{(r+1)}_i=\qint{r+1}F^{(r)}_i$; see
    \cite[Lemma~6.15]{Mathas:ULect} for a similar argument. The base
    case for the induction is given by \autoref{C:FockSpace}.
  \end{proof}

  \notation{$\zero\Dgood\bmu$}{A $\Dom$-good node sequence from $\zero$ to $\bmu$}[D:NormalGood]

  \begin{Definition}[Normal and good nodes]\label{D:NormalGood} Let $\blam\in\Parts$ and
    $i\in I$.
    \begin{enumerate}
      \item A removable $i$-node $A\in\Rem_i(\blam)$ is
      \emph{$\Ldom$-normal} if $\Ldilam\le0$ and
      $\Ldilam<\Ldilam[B]$ if~$B<A$, for~$B\in\Rem_i(\blam)$.
      \item A normal $i$-node $A$ is \emph{$\Ldom$-good} if $A\le B$
      whenever $B$ is a $\Ldom$-normal $i$-node. Equivalently, $A$ is a
      $\Ldom$-good $i$-node if $\Ldilam\le\Ldilam[B]$ for all
      $B\in\Rem_i(\blam)$ with equality only if $A\le B$.
      \item A removable $j$-node $A\in\Rem_j(\blam)$ is
      \emph{$\Gdom$-normal} if $\Gdilam\le0$ and
      $\Gdilam<\Gdilam[B]$ if~$B>A$, for~$B\in\Rem_j(\blam)$.
      \item A normal $j$-node $A$ is \emph{$\Gdom$-good} if $A\ge B$
      whenever $B$ is a $\Gdom$-normal $j$-node.  Equivalently, $A$ is a
      good $i$-node if $\Gdilam\le\Gdilam[B]$ for all
      $B\in\Rem_i(\blam)$ with equality only if $A\ge B$.
    \end{enumerate}
    If $\bmu=\blam+A$ write $\blam\Lgood[i]\bmu$ if $A$ is a
    $\Ldom$-good $i$-node of~$\bmu$ and write $\blam\Ggood[j]\bnu$ if $A$ is an
    $\Gdom$-good $j$-node of~$\bnu$. More generally, if $\bmu,\bnu\in\Parts$ and
    $\bi,\bj\in I^n$, write $\zero\Lgood[\bi]\bmu$ and
    $\zero\Ggood[\bj]\bnu$ if there exist $\ell$-partitions
    $\bmu_1,\dots,\bmu_n=\bmu$ and $\bnu_1,\dots,\bnu_n=\bnu$ such that
     \[
            \zero\Lgood[i_1]\bmu_1\Lgood[i_2]\dots\Lgood[i_n]\bmu_n=\bmu
            \And
            \zero\Ggood[j_1]\bnu_1\Ggood[j_2]\dots\Ggood[j_n]\bnu_n=\bnu,
     \]
     respectively.
  \end{Definition}

  There is a dual definition for conormal and cogood nodes.

  \begin{Definition}[Conormal and cogood nodes]\label{D:CoNormalGood}
    Let $\blam\in\Parts$ and $i\in I$.
    \begin{enumerate}
      \item An addable $i$-node $A\in\Add_i(\blam)$ is
      \emph{$\Ldom$-conormal} if $\Ldilam\ge0$ and $d_A(\blam)>d_B(\blam)$
      if $A<B$, for $B\in\Add_i(\blam)$.
      \item A normal $i$-node $A$ is \emph{$\Ldom$-cogood} if $A\ge B$
      whenever $B$ is a $\Ldom$-normal $i$-node.
      \item An addable $j$-node $A\in\Add_j(\blam)$ is
      \emph{$\Gdom$-conormal} if $\Gdilam[j]\ge0$ and
      $d_A(\blam)>d_B(\blam)$ if $A>B$, for $B\in\Add_j(\blam)$.
      \item A normal $j$-node $A$ is \emph{$\Gdom$-cogood} if $A\le B$
      whenever $B$ is a $\Gdom$-normal $j$-node.
    \end{enumerate}
  \end{Definition}

  In particular, if $\bmu=\blam\cup A$ then $A$ is a good $i$-node of
  $\bmu$ if and only if~$A$ is a cogood $i$-node of~$\blam$.

  Normal and conormal nodes are often defined by listing the addable and
  removable $i$-nodes for~$\blam$ lexicographically and then recursively
  deleting all adjacent addable-removable pairs for $\Ldom$-normal
  nodes, and removable-addable pairs for $\Gdom$-normal nodes. After all
  such pairs have been removed, the normal nodes are the removable nodes
  that remain and the conormal nodes are the addable nodes. It is
  slightly tedious, but straightforward, to check that these
  descriptions of normal and conormal nodes are equivalent to the two
  definitions above; compare with \cite[Lemma~11.2]{Ariki:book}.

  \begin{Example}
    Consider the partition $\lambda=(4,3,1)$ for the algebra
    $\Rx<\Lambda_0>[6](\K[x])$ of type $\Cone[2]$. The type~$\Cone[2]$
    residues in $\lambda$ are given by the diagram:
    \[
    \Tableau{{0,1,2,1},{1,0,1},{2}}
    \]
    Then
    $\zero\Lgood[0](1)\Lgood[1](2)\Lgood[1](2,1)\Lgood[0](2^2)
        \Lgood[2](3,2)\Lgood[2](3,2,1)\Lgood[1](4,2,1)\Lgood[1](4,3,1)$. It
        follows from \autoref{T:KleshchevSimples} below that
        $\LDmu[(4,3,1)]\ne0$. In contrast,
    $(3)\Ggood[1](3,1)\Ggood[0](3,2)\Ggood[1](3^2)\Ggood[1](4,3)
        \Ggood[2](4,3,1)$.
    The partition $(3)$ does not have any $\Gdom$-normal nodes, so
    $\GDnu[(4,3,1)]=0$ by \autoref{T:KleshchevSimples}.
  \end{Example}

  Analogues of the next result are well-known. Given its importance to
  the main results of this paper we give the proof, following
  \cite[Theorem~6.17]{Mathas:ULect}. Perhaps unexpectedly, the result
  mixes up the dominance and reverse dominance partial orders.

  \begin{Theorem}\label{T:InfinityCrystalGraph}
    Let $\blam,\bmu\in\Parts$ and $i\in I$.
    \begin{enumerate}
      \item If $\blam$ does not have a $\Gdom$-good $j$-node then
      $e_j\Lfock\in q^{-1}\LFock[\AA_\infty]$.
      \item If $\blam\Ggood[j]\bmu$ then
      $e_j\Lfock[\bmu]=\Lfock\pmod{q^{-1}\LFock[\AA_\infty]}$
      and $f_j\Lfock=\Lfock[\bmu]\pmod{q^{-1}\LFock[\AA_\infty]}$.
      \item If $\blam$ does not have a $\Gdom$-good $i$-node then
      $e_i\Gfock\in q^{-1}\GFock[\AA_\infty]$.
      \item If $\blam\Lgood[i]\bmu$ then
      $e_i\Gfock[\bmu]=\Gfock\pmod{q^{-1}\GFock[\AA_\infty]}$
      and $f_i\Gfock=\Gfock[\bmu]\pmod{q^{-1}\GFock[\AA_\infty]}$.
    \end{enumerate}
  \end{Theorem}

  \begin{proof}
    We prove only parts (a) and (b) as the proofs of (c) and (d) follow
    by symmetry. First suppose that $\blam$ does not have a $\Gdom$-good
    $i$-node. If $A\in\Rem_i(\blam)$ then $\Gdilam(\blam)>0$, so there
    are at least as many addable $i$-nodes below $A$ as there are
    removable $i$-nodes. Let $\check A$ be the highest addable
    $i$-node of $\blam$ such that $A<\check A$ and
    $\Ldilam[\check A](\blam)=\Ldilam[A](\blam)+1$. As
    $\Ldilam[A](\blam)>0$ the node $\check A$ always exists and if
    $A,B\in\Rem_i(\blam)$ then $\check A=\check B$ if and only if $A=B$.
    If $M\subseteq\Rem_i(\blam)$ let $\check\blam_M=\blam{-}M{+}\check
    M$, where $\check M=\set{\check A|A\in M}$. That is,
    $\check\blam_M$ is the $\ell$-partition obtained from~$\blam$ by
    removing the $i$-nodes in $M$ from~$\blam$ and then adding on the
    nodes in $\check M$. In particular, $|\check\blam_M|=|\blam|$. Now
    set
    \[
      \check\Omega_i(\Lfock)=\sum_{M\subseteq\Rem_i(\blam)} (-q)^{-|M|}
            \Lfock[\check\blam_M]\in\LFock[\AA_\infty].
    \]
    By \autoref{C:FockSpace}, $\Lfock[\bnu]$ appears in
    $E_i\check\Omega_i(\Lfock)$ only if
    $\Rem_i(\bnu)=\check M\cup N$ where $\Rem_i(\blam)=M\sqcup
    N\sqcup\set{A}$ (disjoint union). Now, $\Lfock[\bnu]$ appears
    in $E_i\Lfock[\blam_{M}]$ and in $E_i\Lfock[\blam_{M\cup\set{A}}]$,
    and its coefficient in~$E_i\Omega_i(\Lfock)$ is
    \[
        (-q)^{-|M|+d^\Ldom_A(\blam_M)}
           +(-q)^{-|M|-1+d^\Ldom_{\check A}(\blam_{M\cup\set{A}})}=0,
    \]
    where the last equality follows because $d^\Ldom_{\check
    A}(\blam_M)=\Ldilam[\check A](\blam)=\Ldilam+1=d^\Ldom_A(\blam_M)$,
    which is the key identity underpinning this theorem.  Hence,
    $E_i\check\Omega_i(\Lfock)=0$ and, consequently,
    $e_i\check\Omega_i(\Lfock)=0$ by \autoref{E:KashiwaraOps}.
    Therefore,
    \[
            e_i\Lfock\equiv e_i\check\Omega_i(\Lfock)
            =0 \pmod{q^{-1}\LFock[\AA_\infty]},
    \]
    proving (a).

    To prove (b) we continue to assume that $\blam$ has no
    $\Gdom$-normal $i$-nodes and compute $f_i^r\Lfock$, for $r\ge0$.
    Using the notation above, set
    \[
      \Ncal_i(\blam)
        = \set[\big]{A\in\Add_i(\blam)|A\ne\check B
                       \text{ for any }B\in\Add_i(\blam)}
        = \set{A_1>\dots>A_z}.
    \]
    Observe that $z=\#\Ncal_i(\blam)=d_i(\blam)$ and that
    $s=\Gdilam[A_s]$, for $1\le s\le z$. So, $\Ncal_i(\blam)$ is the set of
    $\Gdom$-conormal $i$-nodes of~$\blam$.

    For $K\subseteq\Add_i(\bnu)$ let $\bnu{+}K$ be the $\ell$-partition
    $\bnu\cup K$. Using \autoref{E:KashiwaraOps} for the first
    congruence, and \autoref{L:DividedPowers} for the following equality,
    \begin{align*}
      f_i^r\Lfock &\equiv F_i^{(r)}\check\Omega_i(\Lfock)
      \quad\pmod{q^{-1}\LFock[\AA_\infty]}\\
      &=\sum_{M\subseteq\Rem_i(\blam)}(-q)^{-|M|}
        \sum_{\substack{K\subseteq\Add_i(\check\blam_M)\\|K|=r}}
        q^{-\Gdilam[\check\blam_M{+}K](\check\blam_M)}\Lfock[\check\blam_M{+}K]\\
      &=\sum_{M\subseteq\Rem_i(\blam)}(-q)^{-|M|}
        \sum_{\substack{K\subseteq\Add_i(\blam)\setminus\check M\\|K|=r}}
        q^{-\Gdilam[\check\blam_M{+}K](\check\blam_M)}\Lfock[\check\blam_M{+}K]\\
      &=\sum_{\substack{K\subseteq\Add_i(\blam)\\|K|=r}}
        \sum_{\substack{M\subseteq\Rem_i(\blam)\\\check M\cap K=\emptyset}}
          (-q)^{-|M|-\Gdilam[\blam{+}K](\blam)} \Lfock[\check\blam_M{+}K]\\
     &\equiv\strut\begin{cases*}
       \Lfock[\blam{+}\set{A_1,\dots,A_r}]& if $1\le r\le z$,\\
       0& otherwise,
      \end{cases*}
    \end{align*}
    where the last equation, which is modulo $q^{-1}\LFock[\AA_\infty]$, follows because
    if $K\ne\set{A_1,\dots,A_r}$ or $M\ne\emptyset$ then
    $|M|-\Gdilam[\blam{+}K](\blam)>0$.  To complete the proof of~(b) it
    remains to observe that $A_r$ is the $\Gdom$-good $i$-node
    of~$\blam{+}\set{A_1,\dots,A_{r-1}}$.
  \end{proof}

  \begin{Definition}\label{D:CrystalGraph}
    Suppose that $\Lambda\in P^+$. Define
    \begin{align*}
        \LCry &=\set{\bmu|\bmu\in\Parts\text{ and }\zero\Lgood\bmu
          \text{ for some }\bi\in I^n \text{ and }n\ge0 }
    \intertext{and}
        \GCry &=\set{\bnu|\bnu\in\Parts\text{ and }\zero\Ggood\bnu
          \text{ for some }\bj\in I^n \text{ and }n\ge0 }
    \end{align*}
    and set
    $\LCrystal*=\set{\Lfock[\bnu]+q^{-1}\LLa[\AA_\infty]|\bnu\in\GCry}$
    and
    $\GCrystal*=\set{\Gfock[\bmu]+q^{-1}\GLa[\AA_\infty]|\bmu\in\LCry}$.
  \end{Definition}

  By definition, $\LCrystal*$ is contained in
  $\LLa[\AA_\infty]/q^{-1}\LLa[\AA_\infty]$ and, similarly, $\GCrystal*$ is
  contained in $\GLa[\AA_\infty]/q^{-1}\GLa[\AA_\infty]$.

  \begin{Corollary}\label{C:InfinityCrystalGraph}
     Let $\Lambda\in P^+$. Then $\bigl(\LLa[\AA_\infty], \GCrystal*\bigr)$
     and $\bigl(\GLa[\AA_\infty], \LCrystal*\bigr)$ are $\infty$-crystal
     bases of~$L(\Lambda)$.
  \end{Corollary}
  \notation{$\LCry,\GCry$}{The sets $\set{\bmu\in\Parts*|\zero\Dgood\bmu}$}[D:CrystalGraph]

  \begin{proof}
    We only prove the result for $\bigl(\LLa[\AA_\infty],
    \GCrystal*\bigr)$.  The only condition in \autoref{D:CrystalBasis}
    that is not clear from \autoref{T:InfinityCrystalGraph} is
    that~$\GCrystal*$ is a $\Q$-basis of~$\LLa[\AA
    _\infty]/q^{-1}\LLa[\AA _\infty]$. Since $\LLa[\AA _\infty]$ is a
    highest weight module,
    \[
        \LLa[\AA _\infty]/q^{-1}\LLa[\AA _\infty]
           = \bigl\<\, f_\bi\Lfock[\zero]+q^{-1}\LLa[\AA _\infty]\bigm
                   |\bi\in I^n \text{ for }n\ge0\,\bigr\>_{\AA _\infty}.
    \]
    Hence, it is enough to show that
    $\set{f_\bi\Lfock[\zero]+q^{-1}\LLa[\AA _\infty]|\bi\in I^n}$ is spanned by
    \[
    \set{\Lfock[\bmu]+q^{-1}\LLa[\AA _\infty]|\bmu\in\GCry\cap\Parts\text{ for }n\ge0}.
    \]
    We argue by induction on~$n$. If $n=0$ there is nothing
    to prove since $\Lfock[\zero]$ is a highest weight vector
    in~$\LLa[\AA _\infty]$. By way of induction, suppose that the claim
    is true for~$n$ and consider the statement for~$n+1$. Fix
    $\bmu\in\GCry$ and $\bi\in I^n$ such that $\zero\Ggood[\bi]\bmu$. By
    \autoref{T:InfinityCrystalGraph}, $f_i\Lfock\in q^{-1}\LLa[\AA
    _\infty]$ if and only if $\bmu$ has no $\Gdom$-conormal $i$-nodes
    and, moreover, if $A$ is the $\Gdom$-cogood $i$-node then
    $f_i\Lfock[\bmu]\equiv\Lfock[\bmu{+}A]\pmod{q^{-1}\LLa[\AA _\infty]}$.
    This completes the proof of the inductive step and hence proves the corollary.
  \end{proof}

  For $i,j\in I$ and $\blam\in\Parts$ define functions
  $\Leps_i,\Lphi_i\map{\LCry}\Z$ and
  $\Geps_j,\Gphi_j\map{\GCry}\Z$ by
  \begin{align}\label{E:CrystalData}
    \begin{aligned}
        \Leps_i(\bmu)&=\#\set{A\in\Add_i(\bmu)|A\text{ is $\Ldom$-normal}}&
        \Geps_i(\bnu)&=\#\set{A\in\Add_i(\bnu)|A\text{ is $\Gdom$-normal}}\\
        \Lphi_i(\bmu)&=\#\set{A\in\Rem_i(\bmu)|A\text{ is $\Ldom$-conormal}}&
        \Gphi_i(\bnu)&=\#\set{A\in\Rem_i(\bnu)|A\text{ is $\Gdom$-conormal}}
    \end{aligned}
  \end{align}
  for $\bmu\in\LCry*$ and $\bnu\in\GCry*$. Let $i,j\in I$. These
  definitions readily imply that if $i\in I$ then
  \begin{equation}\label{E:EpsPhi}
        d_i(\bmu)=\Lphi_i(\bmu)-\Leps_i(\bmu)
        \quad\text{and}\quad
        d_i(\bnu)=\Gphi_i(\bnu)-\Geps_i(\bnu),
        \qquad\text{for $\bmu\in\LCry$ and $\bnu\in\GCry$.}
  \end{equation}
  Abusing notation, if $\blam, \bmu\in\LCry$ and $\blam\Lgood[i]\bmu$
  we write $e_i\bmu=\blam$ and $f_i\blam=\bmu$. Similarly, if $\bsig,
  \bnu\in\GCry$, write $e_j\bnu=\bsig$ and $f_j\bsig=\bnu$ if
  $\bsig\Ggood[j]\bnu$. If $\Deps_i(\blam)=0$ set $
  e_i\blam=0$ and if $\Dphi_i(\blam)=0$ set $f_i\blam=0$.
  \notation{$\Leps_i(\bmu),\Geps_i(\bmu)$}{The number of $\Dom$-normal $i$-nodes, for $i\in I$}
  \notation{$\Lphi_i(\bmu),\Gphi_i(\bmu)$}{The number of $\Dom$-conormal $i$-nodes, for $i\in I$}

  By \autoref{C:InfinityCrystalGraph}, if $m$ is a non-negative integer and
  $\blam\in\DCry$ then $e_i^m\blam\ne0$ if and only
  if~$m\le\Deps_i(\blam)$ and $f_i^m\blam\ne0$ if and only if
  $m\le\Dphi_i(\blam)$. Therefore, following
  \cite[\S7.2]{Kashiwara:CrystalSurvey}, the datum
  $(\LCry,e_i,f_i,\Leps,\Lphi,\wt)$ uniquely determines Kashiwara's
  upper crystal graph of~$\GLa$, where $\wt$ is the weight function
  of~\autoref{E:weight}. Similarly, the datum
  $(\GCry,e_i,f_i,\Geps,\Gphi,\wt)$ determines the upper crystal graph
  of~$\LLa$.

  Using \autoref{T:KashiwaraLusztig}, the crystal bases
  $\LCrystal*$ and $\GCrystal*$ lift to canonical bases
  \[
    \set[\big]{\LCan*[\bnu]|\bnu\in\GCry}
    \qquad\text{and}\qquad
    \set[\big]{\GCan*[\bmu]|\bmu\in\LCry}
  \]
  of $\LLa^*$ and $\GLa^*$, respectively, that
  are uniquely determined by the properties:
  \begin{align}\label{E:InfinityFockCanonical}
    \begin{split}
      \overline{\LCan*[\bnu]}&=\LCan*[\bnu]\And
      \LCan*[\bnu]\equiv\Lfock[\bnu] \pmod{q^{-1}\LLa[\AA_\infty]}\\
      \overline{\GCan*[\bmu]}&=\GCan*[\bmu]\And
      \GCan*[\bmu]\equiv\Gfock[\bmu] \pmod{q^{-1}\GLa[\AA_\infty]}.
    \end{split}
  \end{align}
  for $\bnu\in\GCry$ and $\bmu\in\LCry$.

  When combined with \autoref{T:GradedSimples}, the next result proves
  \autoref{MT:Simples} from the introduction. As remarked at the start
  of \autoref{S:Categorification}, \autoref{C:RnCellular},
  this result applies to all (standard) cyclotomic KLR algebras of types
  $\Aone$, $A_\infty$, $\Cone$ and $C_\infty$.

  \begin{Theorem}\label{T:KleshchevSimples}
    Let $\Lambda\in P^+$. Then
    $\LKlesh=\LCry$ and $\GKlesh=\GCry$. Moreover,
    if $\bmu\in\LKlesh$ then
    \[
        \Ldec_T(q^{-\defect\bmu}\LGmu[\bmu])=\LCan*[\mull(\bmu)]
        \quad\text{and}\quad
        \Gdec_T(q^{-\defect\bmu}\GGnu[\mull(\bmu)])=\GCan*[\bmu].
    \]
  \end{Theorem}

  \begin{proof}
    By working with $\LLa$ we prove that $\GCry=\GKlesh$ and that
    $\Ldec_T(q^{-\defect\bmu}\LGmu[\bmu])=\LCan*[\mull(\bmu)]$ for
    $\bmu\in\LKlesh$. The remaining results are proved in exactly the
    same way and are left as an exercise for the reader.  By
    \autoref{C:BarStar} and \autoref{L:EiFiBarCommute}, the functor
    $\circledast$ categorifies the bar involution on $\LLa$, so
    $\set{q^{-\defect\bmu}\LGmu|\bmu\in\LKlesh[\bullet]}$ is the
    $\infty$-canonical basis of $\ProjN$. By
    \autoref{T:KashiwaraLusztig}, the $\infty$-canonical basis is
    uniquely determined by the choice of highest weight vector, and $\Ldec_T$
    sends $\LGmu[\zero]$ to $\Lfock[\zero]$. Hence, if $\bmu\in\LKlesh$
    then $\Ldec_T(q^{-\defect\bmu}\LGmu)=\LCan*[\bnu]$, for some
    $\bnu\in\GCry$. To determine the $\ell$-partition $\bnu$, we compute
    in~$\LLa[\AA_\infty]$:
    \begin{align*}
        \Ldec_T(q^{-\defect\bmu}\LGmu)
        &= q^{-2\defect\bmu}\sum_{\blam\in\Parts}
            \LPair[\big]{\Ldec_T(\LGmu)}{\Lfock}\Lfock
            & \text{by \autoref{E:FockPairing} and \autoref{P:Fock},}\\
        &= q^{-\defect\bmu}\sum_{\blam\in\Parts}
            \CPair[\big]{\LGmu}{[\LSlam]}\Lfock
               & \text{by \autoref{L:FormComparision},}\\
        &= q^{-\defect\bmu}\sum_{\blam\in\Parts}\LGlammu\Lfock
               & \text{by \autoref{L:CanonicaBasisCoeffs},}\\
        &= \Lfock[\mull(\bmu)]+ \sum_{\blam\in\Parts}
        \overline{\LGlammu[\blam\mull(\bmu)]}\Lfock\pmod{q^{-1}\LFock[\AA_\infty]}
          & \text{by \autoref{P:CanonicalMullineux}}\\
        &\equiv \Lfock[\mull(\bmu)]+
        \sum_{\mathclap{\blam\in\Parts\setminus(\LKlesh\cup\GKlesh)}}
        \overline{\LGlammu[\blam\mull(\bmu)]}\Lfock\pmod{q^{-1}\LFock[\AA_\infty]},
    \end{align*}
    where the last equality comes from \autoref{C:DecompDegs}.
    Therefore, \autoref{T:InfinityCrystalGraph} and
    \autoref{E:InfinityFockCanonical} force
    $\bnu=\mull(\bmu)$ and $\overline{\LGlammu[\blam\mull(\bmu)]}
        =q^{-\defect\bmu}\LGlammu\in\delta_{\blam\mull(\bmu)}+q^{-1}Z[q^{-1}]$,
    for $\blam\in\Parts$. That is,
    \[
        \Ldec_T(q^{-\defect\bmu}\LGmu)=\LCan*[\mull(\bmu)]
        \quad\text{and}\quad
        \bnu=\mull(\bmu)\in\GKlesh.
    \]
    In particular, this shows that
    $\GCry=\set{\mull(\bmu)|\bmu\in\LKlesh[\bullet]}=\GKlesh[\bullet]$,
    where the last equality is \autoref{D:Mullineux}. This completes the
    proof.
  \end{proof}

  \autoref{T:KleshchevSimples} completes the
  classification of the simple $\Rx(\K[x])$-modules from
  \autoref{T:GradedSimples} by giving a description of the
  sets~$\LKlesh$ and $\GKlesh$. The crystal graphs of $L(\Lambda)$ allow
  us to strengthen this characterisation of $\LKlesh$ and $\GKlesh$.

  \begin{Corollary}\label{C:SimpleClassification}
    Let $\K$ be a field and suppose that $\bmu\in\Parts$.
    \begin{enumerate}
      \item The $\Rx(\K[x])$-module $\LDmu(F)\ne0$ if and only if
      $\bmu\in\LKlesh$.
      \item The $\Rx(\K[x])$-module $\GDnu[\bmu](F)\ne0$ if and only if
      $\bmu\in\GKlesh$.
      \item The $\ell$-partition $\bmu\in\LKlesh$ if and only if
      $\zero\Lgood\bmu$ for some $\bi\in I^n$.
      \item The $\ell$-partition $\bmu\in\GKlesh$ if and only if
      $\zero\Ggood[\bi]\bmu$ for some $\bi\in I^n$.
      \item If $\bmu\in\LKlesh$ and $\bi\in I^n$ then
      $\zero\Lgood\bmu$ if and only if $\zero\Ggood[\bi]\mull(\bmu)$.
    \end{enumerate}
  \end{Corollary}

  \begin{proof}
    By invoking \autoref{T:KleshchevSimples} and \autoref{T:GradedSimples},
    parts (a)--(d) are restatements of the identities $\LKlesh=\LCry$ and
    $\GKlesh=\GCry$. For part~(e), if $\bmu\in\LKlesh$ then
    $\zero\Lgood\bmu$ if and only if the sequence~$\bi$ labels a path in
    the crystal graph of $\LLa$ from~$\zero$ to~$\bmu$. By
    \autoref{T:CyclotomicCat}, the
    $\Uq$-modules $\LLa^*$ and $\GLa^*$ have isomorphic crystal graphs.
    Any crystal isomorphism preserves the labels on the paths, so
    $\zero\Lgood\bmu$ is a path in the crystal graph of~$\LLa$ if and
    only if $\zero\Ggood[\bi]\bnu$ is a path in the crystal graph
    of~$\GLa$, for some $\bnu\in\GKlesh$. Applying
    \autoref{T:KleshchevSimples} twice,
    \[
        \LCan*[\mull(\bmu)]
           =\Ldec_T(q^{-\defect\bnu}\LGmu)
           \quad\text{and}\quad
        \GCan*[\bnu]
           =\Gdec_T(q^{-\defect\bnu}\GGnu[\mull(\bnu)])
    \]
    By \autoref{P:CanonicalMullineux},
    $\LGmu=\GGnu[\mull(\bmu)]$, so the map $\Gdec_T\circ(\Ldec_T)^{-1}$
    induces a crystal isomorphism
    \[
        \bigl(\LLa[\AA_\infty], \GCrystal*\bigr)
            \to \bigl(\GLa[\AA_\infty], \LCrystal*\bigr),
    \]
    which sends $\LCan*[\mull(\bmu)]+q^{-1}\LLa[\AA_\infty]$ to
    $\GCan*[\bmu]+q^{-1}\GLa[\AA_\infty]$. Hence,
    part~(e) follows
  \end{proof}

  We have now proved a strong form of \autoref{MT:Simples} from the
  introduction.

  Notice that \autoref{C:SimpleClassification} gives a description of
  the map $\bmu\mapsto\mull(\bmu)$, for $\mull\map\LKlesh\GKlesh$.
  Explicitly, if $\bmu\in\LKlesh$ then we can find $\bi\in I^n$ such
  that $\zero\Lgood\bmu$ is a path in the crystal graph of $\LLa^*$
  from~$\Lfock[\zero]$ to~$\Lfock[\bmu]$.  Then $\mull(\bmu)\in\GKlesh$
  is the unique $\ell$-partition such that $\zero\Ggood\mull(\bmu)$ in
  the crystal graph of $\GLa^*$. In view of \autoref{C:SignedSimples},
  if $\Gamma$ is a quiver of type~$\Aone$ and $\Lambda=\Lambda_0$, this
  gives a variation on Kleshchev's description of the Mullineux map of
  the symmetric group, which is the function $\bmu\mapsto\mull(\bmu)'$,
  for $\bmu\in\LKlesh$.

  The proof of \autoref{T:KleshchevSimples} gives the following
  strengthening of \autoref{C:DecompDegs}.

  \begin{Corollary}\label{C:FullDecompDegs}
    Let $\bmu\in\LKlesh$, $\bnu\in\LKlesh$ and
    $\blam,\bsig\in\Parts$.
    \begin{enumerate}
      \item If $\LGlammu\ne0$ then $\bmu\Ledom\blam\Ledom\mull(\bmu)$ and
      $\alpha_\blam=\alpha_\bmu$.
      Moreover, $\LGlammu[\bmu\bmu]=1$,
      $\LGlammu[\mull(\bmu)\bmu]=q^{\defect\bmu}$ and
      if~$\mull(\bmu)\Ldom\blam\Ldom\bmu$ then
      $0<\deg\LGlammu<\defect\bmu$.
      \item If $\GGlamnu\ne0$ then $\bmu\Gedom\blam\Gedom\mull(\bmu)$ and
      $\alpha_\blam=\alpha_\bmu$.
      Moreover, $\GGlamnu[\bmu\bmu]=1$,
      $\GGlamnu[\mull^{-1}(\bmu)\bmu]=q^{\defect\bmu}$ and
      if~$\mull^{-1}(\bmu)\Gdom\blam\Gdom\bmu$ then
      $0<\deg\GGlamnu<\defect\bmu$.
    \end{enumerate}
  \end{Corollary}

  By \autoref{C:Ariki}, $\DGlammu=[\DSlam:\DDmu]_q$ in type $\Aone$
  when~$\K$ is a field of characteristic zero, so
  $\DGlammu\in\delta_{\blam\bmu}+q\N[q]$ in this case. In type~$\Cone$,
  we can only say that $\DGlammu\in\delta_{\blam\bmu}+q\Z[q]$, and that
  these polynomials approximate the graded decomposition numbers in the
  sense of \autoref{T:CanonicalToDY}.

  The final results in this section describe the $0$-canonical bases of
  $\LLa$ and $\GLa$. To do this we retrace our steps and prove a
  variation of \autoref{T:InfinityCrystalGraph}.

  \begin{Theorem}\label{T:ZeroCrystalGraph}
    Let $\blam,\bmu\in\Parts$ and $i\in I$.
    \begin{enumerate}
      \item If $\blam$ does not have a $\Ldom$-good $i$-node then
      $e_i\Lfock\in q\LFock[\AA_0]$.
      \item If $\blam\Lgood\bmu$ then
      $e_i\Lfock[\bmu]=\Lfock\pmod{q\LFock[\AA_0]}$
      and $f_i\Lfock=\Lfock[\bmu]\pmod{q\LFock[\AA_0]}$.
      \item If $\blam$ does not have a $\Gdom$-good $j$-node then
      $e_j\Gfock\in q\GFock[\AA_0]$.
      \item If $\blam\Ggood\bmu$ then
      $e_j\Gfock[\bmu]=\Gfock\pmod{q\GFock[\AA_0]}$
      and $f_j\Gfock=\Gfock[\bmu]\pmod{q\GFock[\AA_0]}$.
    \end{enumerate}
  \end{Theorem}

  \begin{proof}
     The proof is almost identical to the proof of
     \autoref{T:InfinityCrystalGraph}. For (a), suppose that $\blam$
     does not have a $\Ldom$-good $i$-node. For $A\in\Rem_i(\blam)$
     define $\hat A$ to be the lowest addable $i$-node of~$\blam$ such
     that $A>\hat A$ and $\Ldilam[\hat A]=\Ldilam+1$. If
     $M\subseteq\Rem_i(\blam)$ set $\hat\blam_M=\blam-M+\hat M$, where
     $\hat M=\set{\hat A|A\in M}$, and define
     \[
       \hat\Omega_i(\Lfock) =
          \sum_{M\subseteq\Rem_i(\blam)}(-q)^{|M|}\Lfock[\hat\blam_M]
     \]
     Exactly we before, it now follows that
     $e_i\Lfock\in q\LFock[\AA_0]$ proving (a) with (b) following
     similarly. We leave the details to the reader.
  \end{proof}

  As before, set
    $\LCrystal=\set{\Lfock[\bnu]+q^{-1}\LLa[\AA_\infty]|\bnu\in\GCry}$
    and
    $\GCrystal=\set{\Gfock[\bmu]+q^{-1}\GLa[\AA_\infty]|\bmu\in\LCry}$.
  The argument of \autoref{C:InfinityCrystalGraph} now yields:

  \begin{Corollary}\label{C:ZeroCrystalBasis}
     Let $\Lambda\in P^+$. Then $\bigl(\LLa[\AA_0], \LCrystal\bigr)$ and
     $\bigl(\GLa[\AA_0], \GCrystal\bigr)$ are $0$-crystal bases of~$L(\Lambda)$.
  \end{Corollary}

  By \autoref{T:KashiwaraLusztig}, the crystal bases $\LCrystal$ and
  $\GCrystal$ lift to canonical bases $\set{\LCan|\bmu\in\LCrystal}$
  of~$\LLa$, and $\set{\GCan|\bnu\in\GCrystal}$ of~$\GLa$,
  that are uniquely determined by the properties:
  \begin{align}\label{E:FockCanonical}
    \begin{split}
      \overline{\LCan}&=\LCan\And \LCan\equiv\Lfock[\bmu]
      \pmod{q\LLa[\AA_0]}\\
      \overline{\GCan}&=\GCan\And \GCan\equiv\Gfock[\bnu]
      \pmod{q\GLa[\AA_0]}.
    \end{split}
  \end{align}
  for $\bmu\in\LCrystal$ and $\bnu\in\GCrystal$.
  Now set
    $\LCrystal=\set{\Lfock[\bnu]+q\LLa[\AA_0]|\bnu\in\LKlesh}$
  and
    $\GCrystal=\set{\Gfock[\bmu]+q\GLa[\AA_0]|\bmu\in\GKlesh}$.

    \begin{Theorem}\label{T:ZeroCrystal}
      Suppose that $\bsig,\bmu\in\LKlesh$ and $\blam,\bnu\in\GKlesh$.
      Then $\Ldec(\LCan)=\LHmu$, $\Gdec(\GCan)=\GHnu$,
      \[
        \LPair{\LCan*[\blam]}{\LCan}=\delta_{\blam\mull(\bmu)}
        \qquad\text{and}\qquad
        \GPair{\GCan*[\bsig]}{\GCan}=\delta_{\mull(\bsig)\bnu}.
      \]
  \end{Theorem}

  \begin{proof}
    By \autoref{T:KleshchevSimples}, $\LCry=\LKlesh$. Therefore, by
    \autoref{L:FormComparision} and the uniqueness of canonical bases
    from \cite[Theorem~5]{Kashiwara:CrystalBases}, if $\bnu\in\LKlesh$
    then we can write
    $\Ldec(\LCan[\bnu])=\LHmu[\bmu]$, for some $\bmu\in\LKlesh$. By
    \autoref{T:CanonicalBasis}, if $\bmu\in\LKlesh$ then
    \[
          (\LHmu)^\circledast=\LHmu\qquad\text{and}\qquad
          \LHmu\equiv[\LSlam[\bmu]]\pmod{q\RepN}.
    \]
    Hence, $\Ldec(\LCan)=\LHmu$. Similarly, $\Gdec(\GCan)=\GHnu$.
    Using \autoref{T:KleshchevSimples} and \autoref{L:FormComparision},
    if $\btau\in\LKlesh$ then
    \[
        \LPair[\big]{\LCan*[\mull(\btau)]}{\LCan}
          = \LPair[\big]{\Ldec_T(q^{-\defect\btau}\LGmu[\btau])}{\LCan}
          = q^{\defect\btau}\CPair{(q^{-\defect\btau}\LGmu[\btau])^\#}{\LHmu}
          = \CPair{\LGmu[\btau]}{\LHmu}
          = \delta_{\btau\bmu},
    \]
    where the last equality follows by \autoref{T:CanonicalToDY} and
    \autoref{E:CParingYmuDnu}.  Setting $\blam=\mull(\btau)$ gives the
    first inner product in the displayed equation. The inner product
    $\GPair{\GCan*[\bsig]}{\GCan}$ can be computed in the same way.
  \end{proof}

  \subsection{Modular branching rules}\label{S:ModularBranching}
  This section uses the results of the last section, and \autoref{T:KLRBasis},
  to prove precise forms of the modular branching theorem, which is
  \autoref{MT:ModularBranching} from the introduction. That is, we prove
  that the modular branching rules for~$\Rx(\K[x])$ categorify the crystal
  graph of~$L(\Lambda)$. In principle, this result has already been
  proved by Lauda and Vazirani~\cite{LaudaVazirani:CatCrystals},
  however, their theorem does not imply our result because it is not
  clear how to relate their labelling of the irreducible
  $\Rx(\K[x])$-modules to \autoref{C:SimpleClassification}.  On the other
  hand, our results do imply those of \cite{LaudaVazirani:CatCrystals}
  for the cyclotomic KLR algebras of types $\Aone$ and $\Cone$.
  Moreover, our approach to the modular branching rules is considerably
  shorter than the other routes in the literature because we have
  already established the link between the representation theory
  of~$\Rx(\K[x])$ and the crystal graph of~$L(\Lambda)$.

  Suppose that $M$ is an $\Rx(\K[x])$-module. Recall from
  \autoref{S:SignAutomorphism} that $\head M$ and $\soc M$ are the head
  of socle of $M$, respectively. For $i\in I$ and $k\ge0$ inductively
  define $\Rx(\K[x])$-modules $\tei^kM$ and $\tfi^kM$ by setting
  $\tei^0M=M=\tfi^0M$ and if $k\ge0$ define
  \[
    \tei^{k+1}M=\soc\bigl(E_i(\tei^kM)\bigr)
    \qquad\text{and}\qquad
    \tfi^{k+1}M=\head\bigl(F_i(\tfi^kM)\bigr).
  \]
  Using these operators attach the following non-negative integers to
  $M$:
  \[
    \eps_i(M) = \max\set{k\ge0|\tei^kM\ne0}
    \qquad\text{and}\qquad
    \phi_i(M) = \max\set{k\ge0|\tfi^kM\ne0}.
  \]
  The key result that we need is the following, which lifts some of the
  easy preliminary results from Grojnowski's approach to the modular
  branching rules into our setting.

  \begin{Proposition}\label{P:Grojnowski}
    Let $\bmu\in\LKlesh$, $\bnu\in\GKlesh$ and $i,j\in I$ and assume
    that $\eps_i(\LDmu)>0$ and $\eps_j(\GDnu)>0$.
    \begin{enumerate}
      \item As $\Rx[n-1](\K[x])$-modules, $E_i(\LDmu)$ is
      self-dual and $\tei\LDmu$ is irreducible with
       $\eps_i(\tei\LDmu)=\eps_i(\LDmu)-1$.
       Moreover,  if $[E_i\LDmu:L]>0$ and $L\not\cong q^b\tei\LDmu$ as
       $\Rx[n-1]$-modules, then $\eps_i(L)<\eps_i(\tei\LDmu)$.
      \item As $\Rx[n-1](\K[x])$-modules, $E_j(\GDnu)$ is
      self-dual and $\tei[j]\GDnu$ is irreducible with
       $\eps_j(\tei[j]\GDnu)=\eps_j(\GDnu)-1$. Moreover,  if
       $[E_j\GDnu:L]>0$ and $L\not\cong q^b\tei[j]\GDnu$ as
       $\Rx[n-1]$-modules, then $\eps_j(L)<\eps_j(\tei[j]\GDnu)$.
       \item Let $M$ be an irreducible $\Rx(\K[x])$-module. Then
       $y_n$ acts nilpotently on $E_iM$ with nilpotency index
       $\eps_i(M)$.
    \end{enumerate}
  \end{Proposition}

  \begin{proof}
    The modules $E_i(\LDmu)$ and $E_j(\GDnu)$ are self-dual by
    \autoref{P:CommutingDuals}.  The remaining claims in~(a) follow from
    \cite{ChuangRouq:sl2,KhovLaud:diagI}. In more detail, by
    construction any irreducible $\Rx[m](\K[x])$-module is an irreducible
    $\Rxaff[m](\K[x])$-module. Hence, $\tei\DDmu=\soc(E_i\DDmu)$ is an
    irreducible $\Rx[n-1](\K[x])$-module by
    \cite[Corollary~3.12]{KhovLaud:diagI}, which also shows that
    $\eps_i(\tei\DDmu)=\eps_i(\DDmu)-1$. The remaining statements follow
    from \cite[Lemma~3.9]{KhovLaud:diagI}.  (The paper
    \cite{KhovLaud:diagI} assumes that the quiver~$\Gamma$ is
    simply-laced but the arguments apply without change in
    type~$\Cone$.)

    Parts (b) now follows by symmetry.

    Now consider~(c). Since $y_n$ has positive degree, it is a nilpotent
    element of~$\Rx(\K[x])$, so the real claim here is that~$y_n$ has
    nilpotency index $\eps_i(M)$ when acting on $E_iM$. This can be
    proved by repeating the argument of \cite[Theorem~3.5.1]{Klesh:book}
    using Lemma~2.1 and Lemma~3.7 of \cite{KhovLaud:diagI}.
  \end{proof}

  \begin{Corollary}\label{C:SimpleHead}
   Suppose that  $\blam,\bmu\in\LKlesh[\bullet]$ and $\bsig,\bnu\in\GKlesh[\bullet]$ and
   fix $i,j\in I$ and $a,b\in\Z$.
   \begin{enumerate}
     \item If $\soc(E_i\LDmu)\cong q^a\LDmu[\blam]$, then
       $\head(F_i\LDmu[\blam])\cong q^{\di-d_i(\blam)-a}\LDmu$
    \item If $\soc(E_j\GDnu)\cong q^b\GDnu[\bsig]$, then
   $\head(F_j\GDnu[\bsig])\cong q^{\di[j]-d_i(\bsig)-b}\GDnu$.
   \end{enumerate}
  \end{Corollary}

  \begin{proof} Let $\Domin$ and suppose that
    $\blam,\bmu\in\DKlesh[\bullet]$ and
    $i\in I$. By tensor-hom adjointness,
    \[
        \Hom_{\Rx(\K[x])}\Bigl(q^{a}\iInd \DDmu[\blam],\DDmu\Bigl)
          \cong\Hom_{\Rx[n-1](\K[x])}\Bigl(q^{a}\DDmu[\blam],\iRes \DDmu\Bigl).
    \]
    By assumption, the right-hand
    hom-space is nonzero if and only if $\soc(E_i\DDmu)\cong
    q^a\DDmu[\blam]$. On the other hand,
    $\iInd\DDmu[\blam]=q^{\di-d_i(\blam)}F_i\DDmu[\blam]$ and
    $F_i\DDmu[\blam]$ is self-dual by \autoref{P:CommutingDuals}.
    Therefore, the left-hand hom-space is nonzero if and only if
    $q^{\di-d_i(\blam)-a}\DDmu$ is a quotient of
    $F_i\DDmu[\blam]$. Moreover, since $\soc(E_i\DDmu)$ is
    irreducible by \autoref{P:Grojnowski}, it follows that
    $\head(F_i\DDmu[\blam])$ is irreducible, so this completes the proof.
  \end{proof}

  By \autoref{P:Grojnowski}, if $L$ is a composition factor of $E_i\DDmu$
  then $\eps_i(L)<\eps(\tei\LDmu)$, so we also obtain:

  \begin{Corollary}\label{C:epsBound}
      Suppose that $i,j\in I$ and let $\bmu\in\LKlesh$ and $\bnu\in\GKlesh$. Then
      \[
            \eps_i(\LDmu) = \max\set[\big]{k\ge0|E_i^k\LDmu\ne0}
            \quad\text{and}\quad
            \eps_i(\GDnu) = \max\set[\big]{k\ge0|E_i^k\GDnu\ne0}.
      \]
  \end{Corollary}

  Recall the definition of the quantum integers $[k]_i$ and quantum
  factorials $[k]_i!$ from \autoref{SS:Fock}.

  Kashiwara's theory of global crystal bases, combined with
  \autoref{C:ZeroCrystalBasis} and \autoref{T:ZeroCrystalGraph}, gives:

  \begin{Lemma}[{Kashiwara~\cite[Lemma~12.1]{Kashiwara:CrystalSurvey}}]\label{L:EiFiCan}
      Suppose that $i,j\in I$ and let $\bmu\in\LKlesh$ and $\bnu\in\GKlesh$. Then
      \begin{align*}
        E_i\LHmu &= [\Leps_i(\bmu)]_i\LHmu[e_i\bmu]
        +\sum_{\substack{\blam\in\LKlesh[n-1]\\
                         \mathclap{\Leps_i(\blam)<\Leps_i(\bmu)-\di}}}
              a^{\Ldom,i}_{\blam\bmu} \LHmu[\blam],
        &
        E_j\GHnu &= [\Geps_j(\bnu)]_j\GHnu[e_j\bnu]
          +\sum_{\substack{\bsig\in\GKlesh[n-1]\\
                          \mathclap{\Geps_j(\bsig)<\Geps_j(\bnu)-\di}}}
              a^{\Gdom,i}_{\bsig\bnu} \GHnu[\bsig],\\
        F_i\LHmu &= [\Lphi_i(\bmu)]_i\LHmu[f_i\bmu]
          +\sum_{\substack{\blam\in\LKlesh[n+1]\\
                          \mathclap{\Lphi_i(\blam)<\Lphi_i(\bmu)-\di[j]}}}
              b^{\Ldom,j}_{\blam\bmu} \LHmu[\blam],
        &
        F_j\GHnu &= [\Gphi_j(\bnu)]_j\GHnu[f_j\bnu]
          +\sum_{\substack{\bsig\in\GKlesh[n+1]\\
                          \mathclap{\Gphi_j(\bsig)<\Gphi_j(\bnu)-\di[j]}}}
              b^{\Gdom,j}_{\bsig\bnu} \GHnu[\bsig].
      \end{align*}
      for bar invariant Laurent polynomials
      $a^{\Ldom,i}_{\blam\bmu},a^{\Gdom,i}_{\blam\bmu},
      b^{\Ldom,j}_{\blam\bmu},b^{\Gdom,j}_{\blam\bmu}\in\A$.
  \end{Lemma}

  Similar to \autoref{C:epsBound},  we can use \autoref{L:EiFiCan} to argue
  by induction to determine the crystal data statistics $\Deps_i(\bmu)$
  and $\Dphi_i(\bmu)$ from \autoref{E:CrystalData}, for
  $\bmu\in\DKlesh$:
  \begin{equation}\label{E:epsBound}
    \Deps_i(\bmu) = \max\set[\big]{k\ge0|E_i^k\DHmu\ne0}
    \quad\text{and}\quad
    \Dphi_i(\bmu) = \max\set[\big]{k\ge0|F_i^k\DHmu\ne0},
  \end{equation}
  Using the last two results we can prove the ``modular restriction rules''
  for the simple $\Rx(\K[x])$-modules. By \autoref{P:Grojnowski}, we
  already know that $\tei D_\bmu$ is irreducible so the next result
  precisely identifies which irreducible it is. We remind the reader
  that this result applies to any cyclotomic KLR algebra of type $\Aone$,
  $A_\infty$, $\Cone$ or $C_\infty$ by \autoref{C:RnCellular}.

  For $\Domin$ define $\Dow$ to be the minimal element
  of~$\Parts$ with respect to the partial order~$\Dom$. That is,
  $\Low=(n|0|\dots|0)$ when $\Dom=\Ldom$, and $\Gow=(0|\dots|0|1^n)$
  when $\Dom=\Gdom$.

  \notation{$\Gow$}
      {The minimal $\ell$-partition $(0|\dots|0|1^n)$ in $(\Parts,\Gedom)$}
  \notation{$\Low$}
      {The minimal $\ell$-partition $(n|0|\dots|0)$ in $(\Parts,\Ledom)$}

  \begin{Theorem}\label{T:ModularRestriction}
      Suppose that $i,j\in I$, $\bmu\in\LKlesh$ and $\bnu\in\GKlesh$. Then
      $\eps_i(\LDmu)=\Leps_i(\bmu)$ and $\eps_j(\GDnu)=\Geps_j(\bnu)$.
      If $\eps_i(\bmu)\ne0$ and $\eps_j(\bnu)\ne0$, respectively, then
      as $\Rx[n-1](\K[x])$-modules,
      \[
        \tei\LDmu\cong q^{\di(\Leps_i(\bmu)-1)}\LDmu[e_i\bmu]
        \qquad\text{and}\qquad
        \tei[j]\GDnu\cong q^{\di[j](\Geps_j(\bnu)-1)}\GDnu[e_j\bnu].
      \]
  \end{Theorem}

  \begin{proof}
    It is enough to consider case $\tei\LDmu$, because the result for
    $\tei\GDnu$ is then implied by symmetry.  We argue, first, by induction on
    $n$ and then on the $\Ldom$-dominance order to show that
    $\eps_i(\DDmu)=\Leps_i(\bmu)$ and that, \textit{up to shift},
    $\tei\LDmu\cong\LDmu[e_i\bmu]$. First, suppose that
    $\bmu=\Gow=(n|0|\ldots|0)$, which is the maximal element of $\LKlesh$
    under dominance. Then $\LDmu$ is the one dimensional trivial module
    of $\Rx(\Kx)$ and $[\LDmu]=\LHmu$ by
    \autoref{T:CanonicalToDY}. Hence, $\eps_i(\LDmu)=\eps_i(\bmu)$ and
    $\tei\LDmu=\LDmu[e_i\bmu]$ if $\eps_i(\LDmu)\ne0$, which is if and
    only if~$i=\br_n(\bmu)$, $e_i\Gow=e_i\bmu=\Gow[n-1]$ and $\eps_i(\bmu)=1$.
    Therefore, the theorem holds when $\bmu=\Gow$.

    Now suppose that $\bmu\ne\Gow$ is not maximal with respect to
    dominance in $\LKlesh$. By induction we can assume that, up to
    shift, $\tei\LDmu[\bsig]=\LDmu[e_i\bsig]$ whenever $\bsig\in\LKlesh$
    and $\bsig\Gdom\bmu$. Set $\eps=\eps_i(\LDmu)$. By
    \autoref{C:epsBound} and \autoref{P:Grojnowski}, there exists
    $\bnu\in\LKlesh[n-\eps]$ and a polynomial $f(q)\in\N[q,q^{-1}]$ such
    that $E_i^{(\eps)}\LDmu=f(q)[\LDmu[\bnu]]$. We will show that
    $\bnu=e_i^\eps\bmu$. By \autoref{T:CanonicalToDY}, we can write
    \[
        [\LDmu] = \LHmu+\sum_{\bsig\Gdom\bmu}
          \Lglammu[\bsig\bmu]\LHmu[\bsig].
    \]
    Let $\eps'=\max\set{\Leps_i(\bsig)|\Lglammu[\bsig\bmu]\ne0}$.    If
    $\eps'>\eps$ then, by \autoref{L:EiFiCan},
    \[
      E^{(\eps')}_i[\LDmu]=
        \sum_{\substack{\bsig\Gedom\bmu\\\Leps_i(\bsig)=\eps'}}
          \Lglammu[\bsig\bmu]\LHmu[e_i^{\eps'}\bsig].
    \]
    In particular, $E^{(\eps')}_i[\LDmu]\ne0$, a contradiction.
    Similarly, if $\eps'<\eps$ then $E^{(\eps')}_i[\LDmu]=0$, giving a
    second contradiction. Hence, $\eps'=\eps$ and we have
    \[
      f(q)[\LDmu[\bnu]]= E^{(\eps)}_i[\LDmu]
         = \sum_{\substack{\bsig\Gedom\bmu\\\Leps_i(\bsig)=\eps}}
               \Lglammu[\bsig\bmu]\LHmu[e_i^\eps\bsig].
    \]
    If $\eps_i(\bmu)<\eps=\eps_i(\LDmu)$ then $\bnu=e_i^\eps\bsig$, for
    some $\bsig\Gedom\bmu$. Applying \autoref{C:SimpleHead} and
    induction, it follows that
    $\LDmu\cong\tfi^\eps\LDmu[\bnu]\cong\LDmu[\bsig]$, up to shift.
    This is a contradiction since $\bsig\Gdom\bmu$. Therefore,
    $\eps_i(\bmu)=\eps_i(\LDmu)$ and $\tei\LDmu=\LDmu[e_i\bmu]$, up to
    shift, completing the proof of the inductive step.

    We have now shown that $\eps_i(\LDmu)=\Leps_i(\bmu)$ and if
    $\eps_i(\bmu)>0$ then $\tei\LDmu\cong q^d\LDmu[e_i\bmu]$, for some
    $d\in\Z$, and it remains to show that $d=\di(\Leps_i(\bmu)-1)$. To
    complete the proof, observe that because
    $\eps_i(\LDmu)=\Leps_i(\bmu)$, Kashiwara's \autoref{L:EiFiCan}
    implies that $[E_i\LDmu:\LDmu[e_i\bmu]]_q=[\Leps_i(\bmu)]_i$.
    By
    \autoref{R:psiycomm}, $y_n$ commutes with $\Rx[n-1](\K[x])$, so
    multiplication by~$y_n$ defines an $\Rx[n-1](\K[x])$-module
    endomorphism of~$E_i\LDmu$. By \autoref{P:Grojnowski}(c), the
    nilpotency index of~$y_n$ acting on~$E_i\LDmu$ is~$\Leps_i(\bmu)$.
    Therefore,
    \begin{equation}\label{E:yn}
        [y_n^k\LDmu/y_n^{k+1}\LDmu: \LDmu[e_i\bmu]]_q\ne0,
        \qquad\text{for $0\le k<\eps_i(\bmu)$.}
    \end{equation}
    Moreover, every composition factor of $E_i\LDmu$
    isomorphic to $\LDmu[e_i\bmu]$, up to shift, arises uniquely in this
    way by the remarks above.  The module $E_i\LDmu$ is self-dual by
    \autoref{P:Grojnowski}(a).  Consequently, $\head(E_i\LDmu)\cong
    q^d\LDmu[e_i\bmu]$, for some $d\in\Z$. Moreover,
    $\tei\LDmu=\soc(E_i\LDmu)\cong
    q^{d+2\di(\Leps_i(\bmu)-1)}\LDmu[e_i\bmu]$ by \autoref{E:yn}.
    Hence, using self-duality again, $d=-\di(\Leps_i(\bmu)-1)$, so
    $\tei\LDmu=q^{\di(\Leps_i(\bmu)-1)}\LDmu[e_i\bmu]$ as claimed.
\end{proof}

\begin{Corollary}\label{C:ModularInduction}
  Let $i,j\in I$, $\bmu\in\LKlesh$ and $\bnu\in\GKlesh$. Then
  $\phi_i(\LDmu)=\Lphi_i(\bmu)$, $\phi_j(\GDnu)=\Gphi_j(\bnu)$
  and
  \[
    \tfi\LDmu\cong q^{\di(1-\Dphi_i(\bmu))}\LDmu[f_i\bmu]
    \qquad\text{and}\qquad
    \tfi[j]\GDnu\cong q^{\di[j](1-\Dphi_j(\bnu))}\GDnu[f_j\bnu]
  \]
  as $\Rx[n+1](\K[x])$-modules.
\end{Corollary}

\begin{proof}
  Let $\Domin$. By \autoref{E:EpsPhi},
  $d_i(\bmu)=\Dphi_i(\bmu)-\Deps_i(\bmu)$, so
  $\tfi\DDmu\cong q^{\di(1-\Dphi_i(\bmu))}\DDmu[f_i\bmu]$ by
  \autoref{T:ModularRestriction} and \autoref{C:SimpleHead}.
  In turn, this implies that $\phi_i(\DDmu)=\Dphi_i(\bmu)$.
\end{proof}

Since $\eps_i(\DDmu)=\Deps_i(\bmu)$ by \autoref{T:ModularRestriction},
and $\phi_i(\DDmu)=\Dphi_i(\bmu)$ by \autoref{C:ModularInduction},
\autoref{L:EiFiCan} now implies:

\begin{Corollary}\label{C:EiFiCan}
    Let $i,j\in I$, $\bmu\in\LKlesh$ and $\bnu\in\GKlesh$. Then
    \begin{align*}
      E_i[\LDmu] &= [\Leps_i(\bmu)]_i[\LDmu[e_i\bmu]]
      +\sum_{\substack{\blam\in\LKlesh[n-1]\\
                       \mathclap{\Leps_i(\blam)<\Leps_i(\bmu)-\di}}}
            c^{\Ldom,i}_{\blam\bmu} [\LDmu[\blam]],
      &
      E_j[\GDnu] &= [\Geps_j(\bnu)]_j[\GDnu[e_j\bnu]]
        +\sum_{\substack{\bsig\in\GKlesh[n-1]\\
                        \mathclap{\Geps_j(\bsig)<\Geps_j(\bnu)-\di}}}
            c^{\Gdom,i}_{\bsig\bnu}[\GDnu[\bsig]],\\
      F_i[\LDmu] &= [\Lphi_i(\bmu)]_i[\LDmu[f_i\bmu]]
        +\sum_{\substack{\blam\in\LKlesh[n+1]\\
                        \mathclap{\Lphi_i(\blam)<\Lphi_i(\bmu)-\di[i]}}}
            d^{\Ldom,i}_{\blam\bmu} [\LDmu[\blam]],
      &
      F_j[\GDnu] &= [\Gphi_j(\bnu)]_j[\GDnu[f_j\bnu]]
        +\sum_{\substack{\bsig\in\GKlesh[n+1]\\
                        \mathclap{\Gphi_j(\bsig)<\Gphi_j(\bnu)-\di[j]}}}
            d^{\Gdom,j}_{\bsig\bnu}[\GDnu[\bsig]].
    \end{align*}
    for bar invariant Laurent polynomials
    $c^{\Ldom,i}_{\blam\bmu},c^{\Gdom,i}_{\blam\bmu},
    d^{\Ldom,j}_{\blam\bmu},d^{\Gdom,j}_{\blam\bmu}\in\N[q,q^{-1}]$.
\end{Corollary}

Many people have observed that the last result implies that the
dimension of $\DDmu$ is at least the number of paths in the
$\Dom$-crystal graph from $\zero$ to $\bmu$, but we can do much better.
If $\bmu\in\DKlesh$ and $\zero\Dgood\bmu$ is a good node sequence, define
the bar invariant polynomial $[\eps_\bi]\in\N[q,q^{-1}]$ recursively by
setting
\[
        [\Deps_\bi(q)] =
          \begin{cases*}
            [\Deps_{i_n}(\bmu)]_{i_n}[\eps_{\bi'}(q)],
                & if $n>0$ and $\bi'=(i_1,\dots,i_{n-1})$,\\
             1  & if $n=0$.
          \end{cases*}
\]
Given two characters $\chi,\chi'\in\N[q,q^{-1}][I^n]$ write
$\chi\ge\chi'$ if $\chi-\chi'\in\N[q,q^{-1}][I^n]$.

\begin{Corollary}\label{C:gdim}
   Let $\bmu\in\LKlesh$ and $\bnu\in\GKlesh$. Then
   \[
     \ch\LDmu\ge\displaystyle\sum_{\zero\Lgood\bmu}[\Leps_\bi(q)]\bi
     \qquad\text{and}\qquad
     \ch\GDnu\ge\displaystyle\sum_{\zero\Ggood\bnu}[\Geps_\bj(q)]\bj.
   \]
\end{Corollary}

\begin{proof}
  This follows easily from \autoref{C:gdim} by induction
  on~$n$.
\end{proof}

This result is rarely sharp.  When $\Rn(F)$ is semisimple and
$\DSlam=\DDmu[\blam]$ is concentrated in degree zero, then the
$\Dom$-good residue sequences are in bijection with the standard
$\blam$-tableaux and $[\Deps_\bi(q)]=1$ (cf.
\cite[Proposition~2.4.6]{Mathas:Singapore}). It follows that the
right-hand side is the graded character of the Specht module, which is
concentrated in degree zero in the semisimple case, so in this case
$\DDmu=\DSlam[\bmu]$ and both bounds in corollary are sharp.

\begin{Corollary}
    Let $i,j\in I$, $\bmu\in\LKlesh$ and $\bnu\in\GKlesh$. Then
    \[
      \END_{\Rx[n-1](F)}(\iRes\LDmu) \cong F[y_n]/(y_n^{\Leps_i(\bmu)})
      \quad\text{and}\quad
      \END_{\Rx[n-1](F)}(\iRes\GDnu) \cong F[y_n]/(y_n^{\Geps_i(\bnu)}).
    \]
    as $\Z$-graded algebras.
\end{Corollary}

\begin{proof}
  Let $\Domin$.  As observed in the proof of
  \autoref{T:ModularRestriction}, multiplication by $y_n$ defines an
  $\Rx[n-1](F)$-module homomorphism of $E_i\DDmu=\iRes\DDmu$ and $y_n$
  acts on $E_i\DDmu$ as a nilpotent operator of index~$\Deps_i(\bmu)$.
  Hence, the image of~$y_n$ in the endomorphism ring
  $\END_{\Rx[n-1](F)}(E_i\LDmu)$ generates a subalgebra isomorphic
  to~$F[y_n]/(y_n^{\Leps_i(\bmu)})$. By \autoref{E:yn}, the image of the
  endomorphism given by multiplication by $y_n^k$ has head isomorphic to
  $q^{\di(2k+1-\Deps_i(\bmu))}\DDmu[e_i\bmu]$, for $0\le k<\Deps_i(\bmu)$. On
  the other hand, if $\phi$ is a (homogeneous) $\Rx[n-1](\K[x])$-module
  endomorphism of $E_i\DDmu$ then~$\phi$ then
  $\head(\im\phi)\cong q^k\DDmu[e_i\bmu]$, for some $k\in\Z$. As
  $[E_i\DDmu:\Deps_i(\bmu)]_q=[\Deps_i(\bmu)]_i$, it follows that
  $\phi(m)=y_n^km$, for some $k$.
\end{proof}

We are missing a description of the endomorphism rings
$\END_{\Rx[n+1](F)}(\iInd\LDmu)$ and
$\END_{\Rx[n+1](F)}(\iInd[j]\LDmu)$, for $\bmu\in\LKlesh$,
$\bnu\in\LKlesh$ and $i,j\in I$. Naively, we might expect that
\[
    \END_{\Rx[n+1](F)}(\iInd\LDmu)\cong F[c_{n+1}]/(c_{n+1}^{\Lphi_i(\bmu)})
    \quad\text{and}\quad
    \END_{\Rx[n+1](F)}(\iInd\GDnu)\cong F[c_{n+1}]/(c_{n+1}^{\Gphi_j(\bnu)}),
\]
where $c_{n+1}=y_1+y_2+\dots+y_{n+1}$. In type~$\Aone$, this result was
proved by Brundan and Kleshchev~\cite[Theorem~4.9]{BK:GradedDecomp}.
Unfortunately, in type~$\Cone$, the element $c_{n+1}$ is rarely
homogeneous, so this statement needs to be modified. In any case, we do
not see how to obtain a description of these endomorphism rings using the
results of this paper.

\bigskip
\IndexOfNotation

\let\u=\uaccent

\end{document}